\documentclass{amsart}
\usepackage[margin=1.3in,heightrounded]{geometry}
\usepackage{setspace}
\usepackage[T1]{fontenc}
\usepackage[utf8]{inputenc}
\usepackage{lmodern}
\usepackage{microtype}
\usepackage[dvipsnames]{xcolor}
\usepackage{graphicx}
\usepackage{mathtools}
\usepackage{amssymb}
\usepackage{amsthm}
\usepackage{xparse}
\usepackage{xspace}
\usepackage{mathrsfs}
\usepackage{xurl}
\usepackage[backref=page]{hyperref}
\hypersetup{
 colorlinks,
 linkcolor={red!50!black},
 citecolor={green!50!black},
 urlcolor={blue!80!black}
}
\usepackage[all]{hypcap}
\usepackage[capitalize, noabbrev]{cleveref}

\usepackage{subcaption}
\usepackage{tikz-cd}
\usepackage{quiver}

\DeclareMathAlphabet{\mathbbold}{U}{bbold}{m}{n}

\newcommand{\Z}{\mathbb Z}
\newcommand{\Q}{\mathbb Q}
\newcommand{\R}{\mathbb R}
\newcommand{\C}{\mathbb C}
\newcommand{\F}{\mathbb F}

\newcommand{\Sph}{\mathbb S}

\newcommand{\GL}{\mathrm{GL}}
\newcommand{\SL}{\mathrm{SL}}
\newcommand{\Spin}{\mathrm{Spin}}
\newcommand{\Pin}{\mathrm{Pin}}
\newcommand{\SO}{\mathrm{SO}}
\newcommand{\rO}{\mathrm{O}}

\newcommand{\Hom}{\mathrm{Hom}}
\newcommand{\Ext}{\mathrm{Ext}}

\newcommand{\Aut}{\mathrm{Aut}}
\newcommand{\End}{\mathrm{End}}
\newcommand{\cH}{\mathcal H}
\newcommand{\coker}{\mathrm{coker}}

\newcommand{\AtoA}{\Phi}

\newcommand{\bl}{\text{--}}
\newcommand{\cA}{\mathcal A}

\newcommand{\inj}{\hookrightarrow}

\newcommand{\id}{\mathrm{id}}

\newcommand{\MTSpin}{\mathit{MTSpin}}
\newcommand{\MTPin}{\mathit{MTPin}}

\DeclarePairedDelimiter\paren{(}{)}
\DeclarePairedDelimiter\ang{\langle}{\rangle}
\DeclarePairedDelimiter\abs{\lvert}{\rvert}

\DeclarePairedDelimiter\bkt{[}{]}
\DeclarePairedDelimiter\set{\{}{\}}

\makeatletter
\let\oldparen\paren
\def\paren{\@ifstar{\oldparen}{\oldparen*}}
\let\oldbkt\bkt
\def\bkt{\@ifstar{\oldbkt}{\oldbkt*}}
\makeatother

\DeclareDocumentCommand{\shortexact}{s O{} O{} mmmm}{
\IfBooleanTF{#1}{ 
 \begin{tikzcd}[ampersand replacement=\&]
        1 \& {#4}
        \&  {#5}
        \& {#6}
        \& {1 #7}
        \arrow[from=1-1, to=1-2]
        \arrow["#2", from=1-2, to=1-3]
        \arrow["#3", from=1-3, to=1-4]
        \arrow[from=1-4, to=1-5]
\end{tikzcd}
}{ 
 \begin{tikzcd}[ampersand replacement=\&]
        0 \& {#4}
        \&  {#5}
        \& {#6}
        \& {0 #7}
        \arrow[from=1-1, to=1-2]
        \arrow["#2", from=1-2, to=1-3]
        \arrow["#3", from=1-3, to=1-4]
        \arrow[from=1-4, to=1-5]
\end{tikzcd}
}}

\numberwithin{equation}{section}
\AddToHook{env/subequations/after}{%
}
\newtheorem{thm}[equation]{Theorem}
\newtheorem{lem}[equation]{Lemma}
\newtheorem{cor}[equation]{Corollary}
\newtheorem{prop}[equation]{Proposition}
\newtheorem{defnthm}[equation]{Definition/Theorem}
\theoremstyle{definition}
\newtheorem{exm}[equation]{Example}
\newtheorem{defn}[equation]{Definition}
\newtheorem{conj}[equation]{Conjecture}

\theoremstyle{remark}
\newtheorem{rem}[equation]{Remark}
\AddToHook{env/thm/begin}{\crefalias{equation}{thm}}
\AddToHook{env/lem/begin}{\crefalias{equation}{lem}}
\AddToHook{env/prop/begin}{\crefalias{equation}{prop}}
\AddToHook{env/cor/begin}{\crefalias{equation}{cor}}
\AddToHook{env/exm/begin}{\crefalias{equation}{exm}}
\AddToHook{env/rem/begin}{\crefalias{equation}{rem}}
\AddToHook{env/defn/begin}{\crefalias{equation}{defn}}
\AddToHook{env/defnthm/begin}{\crefalias{equation}{defnthm}}
\AddToHook{env/conj/begin}{\crefalias{equation}{conj}}
\AddToHook{env/ques/begin}{\crefalias{equation}{ques}}

\crefname{thm}{Theorem}{Theorems}
\crefname{lem}{Lemma}{Lemmas}
\crefname{cor}{Corollary}{Corollaries}
\crefname{prop}{Proposition}{Propositions}
\crefname{exm}{Example}{Examples}
\crefname{defn}{Definition}{Definitions}
\crefname{claim}{Claim}{Claims}
\crefname{rem}{Remark}{Remarks}
\crefname{defnthm}{Definition/Theorem}{Definitions/Theorems}

\usepackage[mathscr]{eucal}
\makeatletter
\def\instring#1#2{TT\fi\begingroup
  \edef\x{\endgroup\noexpand\in@{#1}{#2}}\x\ifin@}
\def\isuppercase#1{%
  \instring{#1}{ABCDEFGHIJKLMNOPQRSTUVWXYZ}%
}%
\newcommand{\C@lIfUpper}[1]{
 \if\isuppercase{#1}\mathscr{#1}%
 \else #1%
 \fi
}
\newcommand{\cat}[1]{\mathit{\@tfor\next:=#1\do{\C@lIfUpper{\next}}}}
\makeatother
\newcommand{\newcat}[2]{\newcommand{#1}{\cat{#2}}}
 \newcat{\fC}{C}
 \newcat{\fD}{D}
 \newcat{\Grp}{Grp} 
 \newcat{\Mod}{Mod} 
 \newcat{\Alg}{Alg} 
 \newcat{\Vect}{Vect} 
 \newcat{\sVect}{sVect} 
 \def\Top{\cat{Top}} 
 \newcat{\Fun}{Fun}
\newcommand{\op}{^{\cat{op}}}

\newcommand{\ko}{\mathit{ko}}
\newcommand{\ku}{\mathit{ku}}
\newcommand{\KO}{\mathit{KO}}
\newcommand{\G}{\mathbb G}
\newcommand{\Sq}{\mathrm{Sq}}
\newcommand{\pt}{\mathrm{pt}}
\newcommand{\E}{\mathbb E}

\newcommand{\SH}{\mathit{SH}}
\newcommand{\Map}{\mathrm{Map}}
\newcommand{\cE}{\mathcal E}
\newcommand{\cD}{\cat D}
\newcommand{\cQ}{\mathcal Q}
\newcommand{\sE}{\mathit{s\mathbb{E}}}

\newcommand{\ev}{\mathit{ev}}
\newcommand{\Ad}{\mathrm{Ad}}
\newcommand{\PO}{\mathrm{PO}}
\AtBeginDocument{
}

\newcommand{\term}{\emph}

\DeclareRobustCommand*{\RaiseBoxByDepth}{%
    \raisebox{\depth}%
}
\newcommand{\uQ}{\RaiseBoxByDepth{\protect\rotatebox{180}{$Q$}}}
\newcommand{\sAlg}{\cat{sAlg}}
\newcommand{\Pic}{\mathrm{Pic}}
\newcommand{\tmf}{\mathit{tmf}}
\newcommand{\PL}{\mathrm{PL}}
\newcommand{\Fil}{\cat{Fi\ell}}

\newcommand{\Post}{\mathrm{AH}}
\newcommand{\AdFilt}{\mathrm{Ad}}

\newcommand{\ksc}{\mathit{ksc}}

\usepackage[dvipsnames]{xcolor}

\definecolor{darkblue}{rgb}{0.,0.,0.4}
\definecolor{darkred}{rgb}{0.5,0.,0.}
\definecolor{BlueViolet}{RGB}{138,43,226}
\definecolor{SkyBlue}{RGB}{30,144,255}
\definecolor{DarkGreen}{RGB}{0,100,0}

\usepackage{tikz}
\usepackage{booktabs}

\theoremstyle{plain}
\newtheorem*{maintheorem}{\cref{map_of_filtered_spectra}}
\newtheorem*{maindefinition}{\cref{def_of_g}}
\newtheorem*{sseqmap}{\cref{sseq_map_cor}}
\newtheorem*{hass_construction}{\cref{cQthm}}
\crefname{conj}{Conjecture}{Conjectures}
\Crefname{conj}{Conjecture}{Conjectures}

\title[An AHSS to Adams map and the HASS for supercohomology]{An Atiyah--Hirzebruch to Adams map and the hastened Adams spectral sequence for supercohomology}

\author{Arun Debray}
\address{Department of Mathematics, University of Kentucky,
719 Patterson Office Tower, Lexington, KY 40506-0027}
\email{\href{mailto:a.debray@uky.edu}{a.debray@uky.edu}}

\author{Weicheng Ye}
\address{Department of Physics and Astronomy, and Stewart Blusson Quantum Matter Institute,
University of British Columbia, Vancouver, BC, Canada V6T 1Z1}
\address{Department of Physics, The Chinese University of Hong Kong, Shatin, New Territories, Hong Kong}
\email{\href{mailto:victoryeofphysics@gmail.com}{victoryeofphysics@gmail.com}}

\author{Matthew Yu}
\address{Mathematical Institute, Andrew Wiles Building, Woodstock Road, Oxford, OX2 6GG}
\email{\href{mailto:yumatthew70@gmail.com}{yumatthew70@gmail.com}}

\date{\today}
\subjclass[2020]{\href{https://mathscinet.ams.org/msc/msc2020.html?t=55T15&btn=Current}{55T15}, \href{https://mathscinet.ams.org/msc/msc2020.html?t=55N20&btn=Current}{55N20}}

\begin{document}

\begin{abstract}
We construct a map of filtered spectra inducing a map from the Atiyah--Hirzebruch spectral sequence (AHSS) to the Baker--Lazarev Adams spectral sequence of a $\ko$-module $M$. Using this and a theorem of Bruner--Rognes, we construct a hastened Adams spectral sequence (HASS) for the Postnikov truncations of $\ko$ and the connective Adams summand.
%
Then, motivated by questions in quantum field theory, topological phases of matter, and tensor category theory, we apply the HASS to \term{supercohomology}, the cohomology theory $I_{\C^\times}(\tau_{\le 2}\ko)$, and observe that in examples of interest it is easier to calculate with than either the AHSS or the usual Adams spectral sequence. We include a User's Guide to the HASS and several worked examples, as well as hastened Adams charts for many commonly encountered $\cA(1)$-modules.
\end{abstract}

\maketitle
\tableofcontents

\setcounter{section}{-1}
\section{Introduction}
	\label{s:intro}

As our modern understanding of quantum phases and quantum anomalies has deepened, physicists and mathematicians have been led to explore stable homotopy calculations involving more exotic generalized cohomology theories \cite{2025PhRvX..15b1062M,GJF:2017zba,Dec24,JF25,Debray:2025kfg,DHJNPPRY24}. \term{Supercohomology}\footnote{Unfortunately, ``supercohomology'' is commonly used to refer to two closely related, but distinct, things: \term{restricted supercohomology}, as constructed by Freed~\cite{Fre08} and Gu--Wen~\cite{GW14}, and \term{extended supercohomology}, as constructed by Wang--Gu~\cite{WG18, WG20} and Kapustin--Thorngren~\cite{KT17}. In this paper, ``supercohomology'' and $\SH$ always refer to extended supercohomology; when we briefly discuss restricted supercohomology in \S\ref{ss:restricted}, we will use the notation $\mathit{rSH}$.} \cite{WG18, KT17, WG20} is a prominent example: it has been applied to understanding modular extensions \cite{JF25} and quantum anomalies \cite{Debray:2025kfg,DY25,Wan:2025ymd,ZJSNew} in fermionic topological orders, mathematically modeled using tensor category theory.


For generalized cohomology theories (GCTs) like supercohomology, 
there are two standard tools available for making computations:
\begin{itemize}
    \item the \term{Atiyah--Hirzebruch spectral sequence} (AHSS), which approximates $\SH^*(X)$ in terms of $H^*(X;\SH^*(\mathrm{pt}))$.
    \item The \term{Adams spectral sequence}, which uses some extra structure on $H^*(X;\Z/2)$ to compute $\SH^*(X)$ (up to $2$-completion).
\end{itemize}
In the mathematical physics literature, these tools are both routinely applied to a closely related GCT, denoted $\mho_\Spin^*$ and called the \term{Anderson dual of spin bordism}; both tools are relatively easy to use and quite powerful. See, for example, \cite{Cam17, BC18}. The Adams computations make use of a change-of-rings theorem, substantially simplifying the homological algebra calculations.
For supercohomology, however, \emph{neither} of these tools is easy to compute with. Specifically, the $\SH^*$ Adams spectral sequence generically has more differentials that mix with the ``native'' $\ko$-Adams differentials in subtle ways, even after applying the change-of-rings simplification to the Adams spectral sequence.\footnote{We will give an example of this behavior in \cref{non_modified_Adams}.}\textsuperscript{,}\footnote{As we discuss further in \S\ref{s:guide}, it is more convenient to run the Adams spectral sequence for the Pontrjagin dual to supercohomology, then dualize afterwards. This does not change our qualitative conclusions about this spectral sequence in the introduction. We will also implicitly dualize in the HASS that we discuss below.} 

In this paper, we introduce another technique, called the \term{hastened Adams spectral sequence} (HASS), for supercohomology and some related theories. We show in examples that the HASS is tractable in examples of interest in physics and tensor category theory, and generally easier to use than the AHSS or the Adams spectral sequence.

Heuristically, the HASS is constructed by quotienting the $\ko$-homology Adams spectral sequence by the corresponding spectral sequence for a Postnikov cover of $\ko$, in such a way that taking the quotient eliminates these new, difficult differentials. It is not automatic that such a quotient exists or is well-behaved, but Bruner--Rognes \cite{BR21} give a sufficient criterion. We prove their criterion in our setting by constructing a map from the Atiyah--Hirzebruch spectral sequence to the Adams spectral sequence for a $\ko$-module. On the AHSS side, the analogue of Bruner--Rognes' criterion is vacuously true, and using our map we carry that over to the Adams side. Though our main goal in this project was to construct and use the HASS, we believe this AHSS-to-Adams map could be of independent interest.

Constructing a map between spectral sequences is a lot like constructing a natural transformation between GCTs: one must show compatibility with a large amount of algebraic data. For this reason, when constructing maps between GCTs, researchers often instead construct maps between their representing spectra: the category of spectra provides a homotopy theory for GCTs, allowing tools first developed in the homotopy theory of topological spaces to be applied to GCTs. This simplifies proofs and constructions in most cases. By analogy, we would like to construct our AHSS-to-Adams map in a well-behaved category describing some sort of homotopy theory of spectral sequences.

The category $\Fil$ of \term{filtered spectra} fills this need.\footnote{We will define filtered spectra in \cref{fildefn}. The reader does not need to know the definition in order to read the introduction.}
The perspective to use $\Fil$ to prove theorems about spectral sequences has been around in the homotopy theory, if somewhat implicitly, for a while, but recently has become quite prominent: see for example~\cite{BMS19, Hed20, Ari21, GIKR22, HRW22, Pst23a, Pst23, Ait24, Ant24, AR24, BJM24, CD24, CDN24, Ant25a, Ant25, AP25, AP25a, Car25, CDvN25, Cre25, HM25, QZ25, vN25, BHS26, Emm26, Pie26, PP26}.


We construct the HASS in two steps. First, we construct a map $\AtoA^\ko$ from the filtered spectrum representing the AHSS of a $\ko$-module $M$ to the filtered spectrum representing Baker--Lazarev's $\ko$-module Adams spectral sequence for $M$~\cite{BL01}. This map requires reindexing the filtration, which means that the induced map of spectral sequences is not degree-preserving. We give analogous results for the connective Adams summand $\ell$, which specializes to (the $2$-localization of) $\ku$ at the prime $2$.


The second step is to use $\Phi^\ko$ to construct the HASS. Hastened Adams spectral sequences were introduced by Behrens--Hill--Hopkins--Mahowald~\cite{BHHM08}, then studied in great generality by Bruner--Rognes~\cite[\S 12.6]{BR21}. A HASS is associated to a cofiber sequence of spectra
\begin{subequations}
\begin{equation}
    X \xrightarrow p Y \xrightarrow q Z,
\end{equation}
where $p$ has positive Adams filtration. The spectral sequence converges to the homotopy of $Z$. We apply the HASS to the Postnikov fiber sequences
\begin{equation}
    \tau_{\ge(n+1)}\ko \longrightarrow \ko \longrightarrow \tau_{\le n}\ko,
\end{equation}
\end{subequations}
smashed over $\ko$ with a $\ko$-module $M$. We use our Atiyah--Hirzebruch-to-Adams map $\AtoA^\ko$ to prove a key technical criterion which, by~\cite[Proposition 12.33]{BR21}, gives a usable description of the $E_2$-page. In examples, the HASS is about as difficult as the ordinary Adams spectral sequence for $\ko$-homology. We also give analogous results with $\ell$ in place of $\ko$. Supercohomology specifically is represented by the Pontrjagin dual of $\tau_{\le 2}\ko$ (see \cref{SHdefn}); thus the $\tau_{\le 2}\ko$-homology and the supercohomology of a space or spectrum straightforwardly determine each other (\cref{SH_as_hom}), so in a slight abuse of notation we refer to the $\tau_{\le 2}\ko$ HASS as the supercohomology HASS.


Finally, we run the HASS to compute (twisted) supercohomology in several physically motivated examples. This includes the nontrivial twist of supercohomology over $B^2\Z/2$, which by work of Johnson-Freyd~\cite{johnson2024minimal} computes the sixteen-fold way of modular extensions of $\cat{sVect}$~\cite{Kit09,BGHNPRW17, LKW17, johnson2024minimal}. We included these examples to demonstrate both how to use the HASS and that it is useful and tractable in applications.


\subsection*{Main results}

One of our main results is the construction of the HASS for a class of theories including twisted $\tau_{\le 2}\ko$-homology, stated precisely in \cref{cQthm} below. Before reading these results, though, the reader may want to look at our ``user's guide'' to the HASS in \S\ref{s:guide}. This guide includes the following.
\begin{itemize}
    \item First, in \S\ref{subsection:procedure}, we summarize the steps of the procedure to use the HASS to compute $\SH^*(X; a, b)$. The steps are similar to the computation of the corresponding twisted $\ko$-homology groups of $X$, as described in \cite{BC18} and~\cite[\S 3.1]{DY23}, but instead of using Ext over $\cA(1)$, one uses the functor $\cQ$ defined in \cref{cQthm}.
    \item In \S\ref{subsection:ExtCharts}, we give $\cQ(M)$ for thirteen frequently encountered $\cA(1)$-modules $M$.
    \item In \S\ref{ss:users_examples}, we work out three examples of the HASS. In \S\ref{sss:cyclic}, we compute $\SH^*(B\Z/2n)$ and $\SH^*(B\Z/2n, 0, y)$, where $y\in H^2(B\Z/2n;\Z/2)$ is the unique nonzero class. In \S\ref{sss:16}, we compute the twisted supercohomology of $K(\Z/2, 2)$, and relate it to Kitaev's sixteen-fold way~\cite[\S 9]{Kit06}. In \S\ref{subsection:Q8}, we compute $\SH^*(BQ_8)$.
\end{itemize}

Now let us discuss the specifics of our main theorems.
We have two main results: the construction of the map from the Atiyah--Hirzebruch to the Adams spectral sequence for $\ko$- and $\ell$-modules (where $\ell$ is the connective Adams summand, see \cref{lem:Adams_summand}), and the existence and formula for the $E_2$-page of the HASS for the Postnikov truncation $\tau_{\le n}\ko$- and $\tau_{\le n}\ell$-homology (see \cref{defn:tau}). Dualizing the result for $\tau_{\le 2}\ko$, the corresponding HASS calculates supercohomology.

Before stating the main theorem, we need the following function to shift the degrees of the Atiyah--Hirzebruch filtration and the Adams filtration.

\begin{maindefinition}
For $x\in\Z$, write $x = 8n+k$ for $0\le k<8$. Then define
\begin{subequations}
\begin{equation}
	g(x) \coloneqq \begin{cases}
		4n+k, &0\le k\le 3\\
		4n+3, &k = 4\\
		4n+4, &5\le k\le 7.
	\end{cases}
\end{equation}
Likewise, for each prime number $p$, define $g^c(p,\bl)\colon\Z\to\Z$ by
\begin{equation}
	g^c(p, x) \coloneqq \left\lceil \frac{x}{2p-2} \right\rceil.
\end{equation}
\end{subequations}
\end{maindefinition}
(In the main text, we give a different, but equivalent, definition of $g$ and $g^c$.)

Then, we have the following promised ``Atiyah--Hirzebruch to Adams'' map for the two corresponding filtrations, defined in \cref{AH_filt_defn,filterR}, respectively. 

For $R$ an $E_\infty$-ring spectrum, $\Fil_R$ denotes the $\infty$-category of filtered $R$-module spectra. For example, the $R$-module Atiyah--Hirzebruch, resp.\ Adams filtrations (\cref{AH_filt_defn}, resp.\ \ref{filterR}) define functors $\Post^R_\bullet, \AdFilt^R_\bullet\colon \Mod_R\to\Fil_R$. The colimit of a filtered spectrum $X$ is denoted $X_{-\infty}$.
\begin{maintheorem}
\hfill
\begin{subequations}
\begin{enumerate}
	\item There is a natural transformation of functors $\Mod_\ko\to \Fil_\ko$
	\begin{equation}
		\AtoA^\ko\colon \Post^\ko_\bullet \Longrightarrow \AdFilt_{g(\bullet)}^{\ko}
	\end{equation}
	which, when composed with $(\bl)_{-\infty}$, is homotopy equivalent to $2$-completion.
	\item Let $p$ be a prime. There is a natural transformation of functors $\Mod_\ell\to\Fil_\ell$
	\begin{equation}
		\AtoA^\ell\colon \Post^\ell_\bullet \Longrightarrow \AdFilt_{g^c(p,\bullet)}^{\ell}
	\end{equation}
	which, when composed with $(\bl)_{-\infty}$, is homotopy equivalent to $p$-completion.
\end{enumerate}
\end{subequations}
\end{maintheorem}

From this, we have a map of spectral sequences, albeit one that shifts the grading.

Our proof of \cref{map_of_filtered_spectra} is fairly specific to $\ell$ and $\ko$. In \S\ref{ss:comparison}, we indicate how it generalizes to any $E_\infty$-ring spectrum which is \term{rising} (\cref{rising_defn}), but this property does not seem to be especially common. (See, however, \cref{ksc_conj}.)

\begin{sseqmap}
Let $X$ be a space or spectrum. After regrading the Adams spectral sequence, there is a map of spectral sequences
from the Atiyah--Hirzebruch spectral sequence for $\ko_*(X)$ to the Adams spectral sequence for $\ko_*(X)_2^\wedge$
(and likewise for $\ell$ at all primes); this map converges to the $2$-completion map $\ko_*(X)_2^\wedge$.
\end{sseqmap}
Computations in examples suggest that this map
explains a curious phenomenon that difficult
differentials in the 2-primary Adams spectral sequence for $\ko_*(X)$ tend to correspond to difficult
differentials in the AHSS for $\ko_*(X)$. We plan to study this map in detail in future work.

Afterwards, we build on the construction of Behrens--Hill--Hopkins--Mahowald~\cite{BHHM08} and Bruner--Rognes \cite{BR21} (see \cref{BHHM} and \cref{coh_MASS}), to construct the hastened Adams resolution and the associated HASS corresponding to the truncation of $\ko$. The following theorem gives an explicit construction of the HASS for supercohomology.

Let $C_k\coloneqq \Sigma^{-k}H_\ko^*(\tau_{\ge k}\ko)$, which we describe explicitly in \cref{A1_module_of_connective_covers}.
\begin{hass_construction}\hfill
\begin{enumerate}
    \item
    There is a map of $\Z^2$-graded $\Ext(\Z/2)$-modules
    \begin{equation}
        f_{n-1}\colon\Ext_{\cA(1)}^{s,t}(C_{n}) \longrightarrow
            \Ext_{\cA(1)}^{s+g(n), t+g(n)+n}(\Z/2),
    \end{equation}
    which is induced from the Postnikov cover map $\tau_{\ge n}\colon \tau_{\ge n}\ko\to\ko$.
    \item 
    There is a functor $\mathcal Q_n^{*,*}$ from $\cA(1)$-modules to $\Z^2$-graded $\Ext(\Z/2)$-modules which commutes with direct sums and such that for all $\cA(1)$-modules $M$, there is a long exact sequence
    \begin{equation}
    \dotsb\to \Ext_{\cA(1)}^{s,t}(C_{n}\otimes M) \xrightarrow{f_{n-1}}
            \Ext_{\cA(1)}^{s+g(n), t+g(n)+n}(M)
    \xrightarrow{\partial_{n-1}}
    \mathcal Q_n^{s,t}(M) \to
    \Ext_{\cA(1)}^{s+1,t}(C_{n}\otimes M)\xrightarrow{f_{n-1}}\dotsb
    \end{equation}
    \item
    Let $M$ be a finite type $\ko$-module.\footnote{The finite-type assumption appears only for technical reasons, and is true in all examples we consider in this paper.}
    Then the HASS for $\tau_{\le n}\ko\wedge_\ko M$ converges strongly and has signature
    \begin{equation}
        E_2^{s,t} = \mathcal Q_n^{s,t}(H^*_\ko(M)) \Longrightarrow \pi_{t-s}(\tau_{\le n}\ko\wedge_\ko M)_2^\wedge.
    \end{equation}
    The map $\pi_*(M)\to \pi_*(\tau_{\le n}\ko\wedge_\ko M)$ lifts to a map from the ordinary Adams spectral sequence to the HASS which on the $E_2$-page is the map $\partial_n$ from~\eqref{cQLES}.
    \end{enumerate}
\end{hass_construction}
An analogous result is true for truncations of $\ell$ at an arbitrary prime, and the proof is identical.

Though we prove the general result above, in applications we focus on the case $n = 2$, as $\tau_{\le 2}\ko$-homology is dual to supercohomology. As such, we will use $\cQ$ to denote $\cQ_2$. In~\cite[Theorem E.6]{Debray:2025kfg}, we stated this special case of \cref{cQthm} without proof\footnote{In the arXiv v1 and v2 of that paper, there is a typo in the degree of the map $f_2$, which is called $g_4$ there; instead of shifting $(s,t)\mapsto(s+3,t+2)$, it should shift $(s,t)\mapsto (s+3,t+7)$. This typo did not propagate to affect any of our computations.} and used it to make supercohomology calculations in (\textit{ibid.}, Propositions E.20, E.23, and F.3). Thus, our proof of \cref{cQthm} in this paper supports those computations.

In practice, the HASS works quite similarly to the usual Adams spectral sequence, except that it replaces the Ext module in the Adams spectral sequence with the module given by the functor $\mathcal Q$ of the theorem. In the users' guide in \cref{s:guide}, we compute $\mathcal Q$ for many common $\cA(1)$-modules, then use them to give several example HASS calculations, motivated by problems arising in physics. In particular, we rederive the 16-fold way~\cite{Kit09,BGHNPRW17, LKW17, johnson2024minimal} from the calculation of (twisted) supercohomology, after Johnson-Freyd's original proof~\cite{johnson2024minimal}. We learned from these examples that running the HASS to compute $\tau_{\le 2}\ko$-homology is generally no more difficult than running the corresponding $\ko$-homology Adams spectral sequence. 

Along the way to our main results and examples, we proved the following results, which may be of independent interest.
\begin{itemize}
    \item In Appendix~\ref{s:supercoh}, we show that the definition of (twisted) supercohomology that arises from the study of fermionic fusion 2-categories agrees with a homotopical definition of supercohomology given as a twisted generalized cohomology theory à la the appoach of Ando--Blumberg--Gepner--Hopkins--Rezk~\cite{ABGHR14a, ABGHR14b}. 
    \item Sanford~\cite[Theorem 5.10]{San26} proves that the kernel of the canonical map from the Witt group of $\cat{Vect}_\R$ to the Witt group of $\cat{sVect}_\R$ is isomorphic to $\Z/2\times\Z/2$, which is an $\R$-analogue of the sixteenfold way over $\C$. Following~\cite[\S 4.4]{johnson2024minimal}, we reinterpret this in \cref{sanford_cor} as saying a certain twisted $\abs{\Pic(\sAlg_\R)}$-cohomology group of $B^2\Z/2$ is isomorphic to $\Z/2\times\Z/2$. We then prove this directly, giving another proof of Sanford's result. The bulk of this proof is identifying the homotopy type of $\abs{\Pic(\sAlg_\R)}$ (\cref{real_pic}), which was a folklore result. Our proof of \cref{sanford_cor} does not directly use the HASS, but does use our computations from \S\ref{sss:16}.
    \item In \cref{Q_ext_calc} we compute $\cQ(M)$ for several commonly occurring $\cA(1)$-modules $M$. We used these calculations to draw the $\cQ$ charts in \cref{Q_ext_charts_1}. Along the way to proving \cref{Q_ext_calc}, we make some $\Ext_{\cA(1)}$ computations that may be of independent interest. For example, in \cref{E5_ext} we determine the algebra structure on $\Ext_{\cA(1)}(H^*((B\Z/2)^{2\sigma-2};\Z/2), \Z/2)$, which could be useful for studying the ring structure on $\Spin\times_{\set{\Z/2}}\Z/4$ bordism (see~\cite[\S 7.8]{Cam17}). 
\end{itemize}

In \cref{s:further} we conclude with some conjectures and further directions regarding the HASSes in this paper.
\begin{itemize}
    \item In \S\ref{ss:restricted} we give some potential applications of the $\tau_{\le 1}\ko$-homology HASS (i.e.\ for restricted supercohomology).
    \item In \S\ref{ss:conjectures} we conjecture that $\AtoA^\ko$ and $\AtoA^\ell$ are lax monoidal and that the HASS is multiplicative.
    \item In \cref{ksc_conj,equivariant_AHSS_to_Adams,,motivic_AHSS_to_Adams}, we speculate on generalizations to connective self-conjugate $K$-theory and equivariant and motivic versions of $\ko$ and $\ku$.
\end{itemize}

\subsection*{Outline}
\S \ref{s:filts} is devoted to proving \cref{map_of_filtered_spectra}, which constructs the Atiyah--Hirzebruch-to-Adams maps $\AtoA^\ko$ and $\AtoA^\ell$. First, in \S\ref{ss:filtered_background}, we introduce our point of view spectral sequences from the point of view of filtered spectra, and describe the Atiyah--Hirzebruch and Adams filtrations and the induced spectral sequences. We finish this subsection by stating \cref{map_of_filtered_spectra} and discussing what we learn, at the qualitative level, from the existence of the AHSS-to-Adams map. In \S\ref{ss:complex}, we prove \cref{map_of_filtered_spectra} for $\ell$ by comparing the $H\Z_{(p)}$- and $H\Z/p$-based $\ell$-module Adams spectral sequences. In \S\ref{ss:comparison}, we prove the $\ko$ part of \cref{map_of_filtered_spectra}.

In \S\ref{s:MASS}, we prove our other main theorem,  \cref{cQthm}, constructing the HASS for $\tau_{\le n}\ko$- and $\tau_{\le n}\ell$-modules. First, in \S \ref{s:bakerlazarev}, we work out the Baker--Lazarev $\ko$-module Adams spectral sequence for $\tau_{\le 2}\ko$. In \S \ref{subsection:HASS}, we review work of Behrens--Hill--Hopkins--Mahowald~\cite{BHHM08} and Bruner--Rognes~\cite[\S 12]{BR21} on hastened Adams spectral sequences in general. In \S\ref{ss:apply_HASS_to_ko}, we state and prove \cref{map_of_filtered_spectra}, applying the general setup of \S\ref{subsection:HASS} to truncations of $\ko$ and $\ell$. We finish this section by working out the HASS for $\tau_{\le 2}\ko$ as a module over itself and observing that this spectral sequence is easier than the ordinary Adams spectral sequence from \S\ref{s:bakerlazarev}.

\S\ref{s:guide} is our user's guide. In \S\ref{subsection:procedure} we walk through the steps of the method of computing twisted supercohomology with the HASS. Then, in \S\ref{subsection:ExtCharts}, we give the charts for $\cQ(\bl)$, the input to the HASS, for many commonly encountered $\cA(1)$-modules. In  \S\ref{ss:users_examples}, we give several examples of the HASS: $B\Z/2n$ in \S\ref{sss:cyclic}, $B^2\Z/2$ and the sixteen-fold way in \S\ref{sss:16}, and $BQ_8$ in \S\ref{subsection:Q8}. In addition, we use our computations to study Sanford's real fourfold way~\cite[Theorem 5.10]{San26} in \S\ref{ss:real_fourfold}.

In \S\ref{s:further}, we speculate on further directions and conjectures. In \S \ref{ss:restricted}, we list three potential applications of the HASS for $\tau_{\le 1}\ko$-homology. In \S\ref{ss:conjectures}, we discuss generalizations of our work in a few different directions.

We conclude with two appendices. In Appendix~\ref{s:supercoh}, we prove \cref{SH_is_SH}, that two different definitions of twisted supercohomology are naturally isomorphic. The first definition, expressed using category theory and more closely tied to applications, appears in \S\ref{ss:cattwist}, and the second, defined with homotopy theory and better for computations, is given in \S\ref{ss:homtwist}. Then, in \S\ref{ss:sametwist}, we prove \cref{SH_is_SH}. Appendix~\ref{ext_Q_computations} contains all of the $\Ext_{\cA(1)}$ and $\cQ$ computations needed to support the $\cQ$ charts in \S\ref{subsection:ExtCharts}.


%

\subsection*{Acknowledgements}
It is a pleasure to thank
Prerna Dhankhar,
Mike Hopkins,
and Morgan Opie
for helpful conversations. MY is supported by the EPSRC Open Fellowship EP/X01276X/1 and the Fund for Natural Science and the Adler Family Fund. 
No authors have competing
interests to declare that are relevant to the content of this article. 

\subsection*{Statement on AI use} 
We did not use AI tools to develop any ideas, proofs, or examples in this paper. We completed the main results and examples in 2025, before the most advanced LLM tools became available.


\section{Compatibility of the Atiyah--Hirzebruch and Adams filtrations}
	\label{s:filts}

  \subsection{Dramatis personæ}
\label{ss:filtered_background}
Here we review the language of filtered objects in an $\infty$-category and the examples that enter into the
construction of the Atiyah--Hirzebruch and Adams spectral sequences. In particular, we prove that there is a natural map from the Atiyah-Hirzebruch spectral sequence to the Adams spectral sequence for $\ko$ or the connective Adams summand $\ell$, in the sense of \cref{sseq_map_cor}.

First, we introduce the Atiyah--Hirzebruch and Adams filtrations, which induce the corresponding spectral sequences in a sense we discuss below. 

\begin{defn}[{\cite[Definition 1.2.2.9]{HA}}]
\label{fildefn}
Let $\fC$ be an $\infty$-category.
\begin{enumerate}
    \item A \term{filtered object} $C_\bullet$ of $\fC$ is a contravariant functor from
the poset category $(\Z, \le)$ to $\fC$. Equivalently, it is data of an object $C_s\in\fC$ for each $s\in\Z$ and
$\fC$-morphisms $f_s\colon C_s\to C_{s-1}$.
    \item Assume $\fC$ is stable. The \term{associated graded} of a filtered object $C_\bullet$ is the filtered spectrum
    $\mathrm{gr}(C)_\bullet$ defined by
    \begin{equation}
        \mathrm{gr}(C)_s\coloneqq
        \mathrm{cofib}(\mathrm f_{s+1}\colon C_{s+1}\to C_s).
    \end{equation}
    Its structure maps $\mathrm{gr}(C)_{s+1}\to\mathrm{gr}(C)_s$ are zero.
    \item The \term{abutment} of $C$ is $C_{-\infty}\coloneqq \operatorname{colim}_s C_s$ (which may or may not exist).
\end{enumerate}
\end{defn}
We will let $\Fil(\fC)\coloneqq \Fun(\Z\op, \fC)$ denote the $\infty$-category of filtered objects in $\fC$. For the most part, our filtered spectra
will have $C_s = 0$ for $s<0$, and we will let $\fC$ be the category $\Mod_R$ of modules over an $E_\infty$-ring spectrum $R$. In this case we abbreviate $\Fil(\Mod_R)$ as $\Fil_R$.

For $X_\bullet\colon\Z\op\to\fC$ a filtered spectrum, we will let $X_s\coloneqq X_\bullet(s)$. If $t\le s$ we will
use lowercase $\mathrm f_{s,t}\colon X_s\to X_{t}$ to denote the image under $X_\bullet$ of the unique morphism
$t\to s$ in the category $\Z$ and $\mathrm f_s\coloneqq \mathrm f_{s,s-1}$.

For $s\in\Z$, $\ev_s\colon\Fil(\fC)\to\fC$ denotes the ``evaluation at $s$'' functor $X_\bullet\mapsto
X_s$.
\begin{rem}[Reindexing]
\label{reindexing}
Let $C_\bullet$ be a filtered object of an $\infty$-category $\cat C$. We will sometimes need to reindex the filtration on $C_\bullet$, e.g.\ changing filtration $s$ to filtration $s+k$ for some $k$. Unwinding the definitions, reindexing amounts to precomposing the functor $C\colon\Z\op\to\cat C$ with a functor $f\colon \Z\op\to\Z\op$. Since $\Hom_{\Z\op}(m, n)$ is empty if $m<n$ and a singleton if $m\ge n$, the data of such a functor $f$ is precisely a nondecreasing function $f\colon\Z\to\Z$.

Suppose our nondecreasing function $f$ is unbounded above. Then we may define $f^\dashv,f^\vdash\colon \Z\to\Z$ by
\begin{equation}
\label{poset_adjoint}
    f^\dashv(m) \coloneqq \max\{n\mid m\ge f(n)\},\qquad\qquad
    f^\vdash(m)\coloneqq \min\{n\mid f(n)\ge m\}.
\end{equation}
One may directly check that $f^\dashv$ and $f^\vdash$ are nondecreasing, so that they define functors $\Z\op\to\Z\op$, and that $f^\dashv$, resp.\ $f^\vdash$, is left, resp.\ right adjoint to $f$. 
\end{rem}
\begin{defn}\label{defn:tau}
Let $\tau_{\le n}\cat{Sp}$ be the full sub-$\infty$-category of spectra whose homotopy groups vanish in degrees $n+1$ and above. The inclusion $\tau_{\le n}\cat{Sp}\inj \cat{Sp}$ has a left adjoint $\tau_{\le n}\colon\cat{Sp}\to\tau_{\le n}\cat{Sp}$, which we call \term{Postnikov $n$-truncation}.

Dually, let $\tau_{\ge n}\cat{Sp}$ be the full sub-$\infty$-category of spectra whose homotopy groups vanish in degrees $n-1$ and below. The inclusion $\tau_{\ge n}\cat{Sp}\inj\cat{Sp}$ admits a right adjoint $\tau_{\ge n}\colon\cat{Sp}\to\tau_{\ge n}\cat{Sp}$, which we call the \term{Postnikov $(n-1)$-connected cover}.
\end{defn}
See~\cite[\S 1.2.1]{HA} for a discussion in the generality of stable $\infty$-categories, and (\textit{ibid.}, Propositions 1.4.3.4 and 1.4.3.6) for why this applies to $\cat{Sp}$ in particular.

For any spectrum $E$, the unit and counit of these adjunctions, resp., give rise to maps $\tau_{\le n}\colon E\to\tau_{\le n}E$ and $\tau_{\ge n}\colon \tau_{\ge n}E\to E$. By~\cite[Remark 1.2.1.8]{HA} these maps participate in a fiber sequence $\tau_{\ge (n+1)}E\to E\to \tau_{\le n}E$. Furthermore, by \cite[Proposition 7.1.1.13]{HA}, $\tau_{\le n}$ and $\tau_{\ge n}$ send $R$-modules to $R$-modules.
\begin{defn}\label{AH_filt_defn}
Let $R$ be an $E_\infty$-ring spectrum and $M$ be an $R$-module. The \term{($R$-module) Atiyah--Hirzebruch filtration} of $M$, denoted $\Post_\bullet^R(M)$, is the functor sending $n\mapsto\tau_{\ge n}(R)\wedge_R M$ and the unique map $m\to n$, $m\le n$, to $\tau_{\ge n}(R)\wedge_R \id\colon\tau_{\ge n}(R)\wedge_R M\to\tau_{\ge m}(R)\wedge_R M$.
\end{defn}
\begin{rem}

We call this the ``Atiyah--Hirzebruch filtration'' because it gives rise to a generalization of the Atiyah--Hirzebruch spectral sequence in the category of $R$-module spectra; see~\cite[Theorem IV.3.7]{EKMM97}. If $M = R\wedge X$ for a spectrum $X$, the spectral sequence obtained from the filtered spectrum $\Post_\bullet^R(M)$ is naturally isomorphic to the usual Atiyah--Hirzebruch spectral sequence for the $R$-homology of $X$.\footnote{Atiyah--Hirzebruch~\cite[\S 2]{AH61} defined their spectral sequence using the skeletal filtration on $X$. Maunder~\cite[Theorem 3.3]{Mau63} showed that using the Postnikov filtration on $R$, as we do in \cref{AH_filt_defn}, leads to an isomorphic spectral sequence from the $E^2$-page onward. Antieau~\cite[Proposition 9.2]{Ant24} lifts Maunder's theorem to a natural equivalence of filtered spectra.}
\end{rem}

\begin{defn}\label{filterR}
Let $\theta\colon R\to S$ be a map of $E_\infty$-ring spectra, $\overline S$ be the cofiber of $\theta$, and $\alpha\colon \Sigma^{-1} \overline S\to R$ be the fiber of $\theta$. The \term{$S$-based $R$-module Adams filtration} of an $R$-module $M$, denoted $\Ad^{R,S}_\bullet(M)$, is the functor sending $n\mapsto (\Sigma^{-1}\overline{S})^{\wedge_R n}\wedge_R M$, and the map $n\rightarrow n+1$ to
\begin{equation}
    \alpha\wedge\id_{\Ad^{R,S}_{n}(M)}\colon \Sigma^{-1}\overline S\wedge_R (\Sigma^{-1}\overline{S})^{\wedge_R n}\wedge_R M \to
    R\wedge_R (\Sigma^{-1}\overline{S})^{\wedge_R n}\wedge_R M \simeq (\Sigma^{-1}\overline{S})^{\wedge_R n}\wedge_R M.
\end{equation}
\end{defn}
Often, $S = H\Z/p$ for a prime $p$ and $\theta$ is the Postnikov $0$-truncation of $R$ followed by modulo $p$ reduction. In this setting, we will write $\Ad_\bullet^R(M)$ instead of $\Ad_\bullet^{R,H\Z/p}(M)$. However, we will need the more general version of \cref{filterR} in \S\ref{ss:complex} and \S\ref{ss:comparison}.

\begin{rem}
Though \cref{filterR} is essentially a restatement of Baker--Lazarev's construction of their version of the Adams spectral sequence in~\cite[\S 2]{BL01}, often one sees the Adams filtration defined differently, as the Bousfield--Kan filtration~\cite{BK72, BK73b, BK73a, Bou03} associated to the cosimplicial $R$-module $S^{\wedge_R \bullet}\wedge_R M$, e.g.\ in van Nigtevecht~\cite[\S 2.5.2]{vN25}. These two approaches coincide, and by Mathew--Naumann--Noel~\cite[Proposition 2.14]{MNN17}, this fact can be lifted to a natural equivalence of filtered $R$-modules between the two definitions of the Adams filtration. See van Nigtevecht~\cite[\S 2.5]{vN25} for further discussion.
\end{rem}

\begin{defnthm}[{\cite[Construction 2.31]{vN25}}]
\label{con:exact-couple}
Let $X$ be a filtered spectrum.
\begin{subequations}
\begin{enumerate}
    \item The \emph{associated exact couple} of $X$ is the data of the bigraded abelian groups
    \begin{equation}
        A_{n,s}(X) := \pi_n X_s
        \qquad \text{and} \qquad
        E_{n,s}(X) := \pi_n \mathrm{gr}(X)_s,
    \end{equation}
    together with the usual maps in the long exact sequence relating them:
    \begin{equation}
        \begin{tikzcd}[column sep=large]
            \pi_n X_s
                \arrow[rr, "{(0,\,-1)}"]
            & &
            \pi_n X_s
                \arrow[dl, "{(0,\,0)}"]
            \\
            &
            \pi_n \mathrm{gr}(X)_s
                \arrow[ul, "{(-1,\,1)}"]
        \end{tikzcd}
    \end{equation}
    (Each map is labeled with its $(n,s)$-bidegree.)

    \item The \emph{underlying spectral sequence} of $X$, written $\{{E}^r_{n,s}(X),\, d^r\}_{r\ge 1}$, is the spectral sequence arising from the exact couple associated with $X$. By definition, this spectral sequence is of the form
    \begin{equation}
        {E}^1_{n,s}(X) := \pi_n \mathrm{gr}(X)_s \implies \pi_n X_{-\infty}
    \end{equation}
    and the differential $d^r$ has $(n,s)$-bidegree $(-1, r)$ for $r \ge 1$.
   
\end{enumerate}
\end{subequations}
\end{defnthm}
The ``theorem'' part of \cref{con:exact-couple} is the claim that $(A_{n,s}(X), E_{n,s}(X))$ is an exact couple, and that $(E_{n,s}^r, d^r)$ is a spectral sequence. 

In many but not all cases, this spectral sequence converges to $X_{-\infty}$.
\begin{rem}
Though we will use the $\infty$-category of filtered spectra as our setting to study the homotopy theory of spectral sequences, it is worth mentioning another perspective. Work of
Brotherston,
Circi,
Egas Santander,
Fu,
Guan,
Livernet,
Muro,
Rotzheim,
and
Whitehouse~\cite{MR19, CESLW20, FGLW22, Bro24, LW24, LW24a, Bro25, CLW25} constructs and studies for each $r$ model structures on filtered (multi)complexes in which weak equivalences are isomorphisms on the $E_r$-pages of the induced spectral sequences. See Aitken~\cite[\S 3.2]{Ait24} for a comparison to the $\infty$-categorical point of view on filtered spectra.
\end{rem}
\begin{exm}[{Maunder~\cite{Mau63}}]
Choose a space $X$. Let's see what \cref{con:exact-couple} looks like for the Atiyah--Hirzebruch filtration for an $E_\infty$-ring spectrum $R$ and the $R$-module $R\wedge X_+$. The cofiber of $\tau_{\ge(s+1)}R\to\tau_{\ge s}R$ is $\tau_{\le s}\tau_{\ge s}R\simeq \Sigma^n H\pi_s(R)$. Smash with $X$ to conclude that the associated graded of $\mathrm{AH}_\bullet^R(R\wedge X)$ has $s^{\mathrm{th}}$ object
\begin{equation}
    \mathrm{gr}(\mathrm{AH}^R(R\wedge X))_s \simeq \Sigma^s H\pi_s(R)\wedge X,
\end{equation}
whose homotopy groups are $H_*(X; \pi_s(R))$ shifted in degree by $s$. That is, in the language of \cref{con:exact-couple}, the associated spectral sequence is
\begin{equation}
\label{surprise_E1}
    E^1_{n,s} \cong \widetilde H_n(X; \pi_s(R)).
\end{equation}
Now this looks like the $E^2$-page, not the $E^1$-page, of the Atiyah--Hirzebruch spectral sequence. What happened? The answer is that we may reindex the filtration, or equivalently reindex the spectral sequence, such that~\eqref{surprise_E1} happens on $E^2$, rather than $E^1$. See Antieau~\cite[\S 7, \S 9]{Ant24} for further information.

In general, when $R$ and $M$ are bounded below, the spectral sequence associated to $\Post_\bullet^R(M)$ converges to its abutment, which is $M$: see~\cite[\S 6]{Ant24}.
\end{exm}



\begin{defn}[{\cite[Definition IV.1.7]{EKMM97}}]
Let $p$ be a prime number and $R$ be an $E_\infty$-ring spectrum with an $E_\infty$-ring map $\theta\colon R\to H\Z/p$. The \term{$R$-module homology and cohomology} of an $R$-module $M$ are
\begin{equation}
    H^R_*(M) \coloneqq \pi_*(M\wedge_R H\Z/p),\qquad\qquad
    H_R^*(M) \coloneqq \pi_{-*}(\mathrm{Map}_R(M, H\Z/p)).
\end{equation}
We write $H^R_*H$ for $H^R_*(H\Z/p)$, and analogously for $H_R^*H$.
\end{defn}
Though $p$ and $\theta$ are implicit in this notation, they will be clear from context wherever we use $H_R^*$ and $H^R_*$ in this paper.
\begin{defn}
\label{adams_filtration}
Let $f\colon M\to N$ be an $R$-module map.
\begin{enumerate}
	\item We say $f$ has \term{($R$-module) Adams filtration $\ge s$} if there is a sequence of $R$-module maps
	\begin{equation}
		M = X_0\overset{t_1}{\longrightarrow} X_1\overset{t_2}{\longrightarrow}\dotsb\overset{t_s}{\longrightarrow}
		X_s = N
	\end{equation}
	such that $(t_i)_*\colon H^R_*(X_{i-1})\to H^R_*(X_i)$ is the zero map.
	\item If $f$ has filtration $\ge s$ but not $\ge(s+1)$, we say $f$ has \term{($R$-module) Adams filtration $s$}.
\end{enumerate}
\end{defn}


We have now defined two Adams filtrations: on an $R$-module (\cref{filterR}) and on the abelian group of homotopy classes of maps between two $R$-modules (\cref{adams_filtration}). When the domain is $R$, this filters $\pi_*(N)$. The content of the next lemma is that $\pi_*$ of the filtration on $N$ is the filtration on $\pi_*(N)$.

\begin{lem}[{Adams~\cite[Theorem III.15.1]{Ada74}, Baker--Lazarev~\cite[\S 2]{BL01}}]
Let $p$ be a prime.
The filtration induced on $\pi_*(M)\cong\pi_*(\Map_R(R, M))$ by applying $\pi_*$ to the filtered spectrum $\Ad_n^{R, H\Z/p}(M)$ defined in \cref{filterR} equals the Adams filtration on $\pi_*(M)$ defined in \cref{adams_filtration}.
\end{lem}
Suppose $\fC$ is now a stable $\infty$-category with a $t$-structure and containing sequential colimits, such as
$\Mod_R$. If $X_\bullet$ is a filtered object of $\fC$, then it gives rise to a spectral sequence in a manner
generalizing the spectral sequence associated to a filtered chain complex; see~\cite[\S 1.2.2]{HA}.
\begin{defn}
The \term{Baker--Lazarev Adams spectral sequence} for an $R$-module $M$ is the spectral sequence associated to the
canonical $2$-primary Adams filtration of $M$.
\end{defn}
\begin{thm}[{Baker--Lazarev~\cite[Theorems 2.1 and 2.3]{BL01}}]\label{thm:basic_BL_ss}
Let $R$ be as above and $M$ be an $R$-module.
\begin{enumerate}
	\item\label{BLA_conv} Suppose $R$ and $M$ are bounded-below, cellular spectra of finite type.
	Then the Baker--Lazarev Adams spectral sequence converges to
	$\pi_*(M)_p^\wedge$ with the $R$-module Adams filtration.
	\item\label{BLA_E2} The $E_2$-page of the Baker--Lazarev Adams spectral sequence
	is\footnote{Here, ``$\Ext$'' means the derived functor of the functor of \emph{co}module homomorphisms. Some
	references will therefore write $\mathrm{coExt}$ or $\mathrm{Coext}$ to emphasize this.}
	\begin{equation}
		E_2^{s,t} = \Ext_{H^R_*H}^{s,t}(\Z/p, H^R_*(M)).
	\end{equation}
\end{enumerate}
\end{thm}
Baker--Lazarev (\textit{ibid.}, \S 2) also discuss convergence of their spectral sequence. Using the criteria they provide, one can show that, for $R = \ko$ or $\ell$ and $M$ bounded below, their spectral sequence converges strongly to the abutment of $\Ad^R_\bullet(M)$.
\begin{rem}
Baker--Lazarev work in greater generality than we do, and therefore their version of~\eqref{BLA_E2} requires a
flatness assumption. This assumption is always true in the setting we work in. In this generality, their spectral
sequences converges to the homotopy groups of the ``$R$-module $H\Z/p$-nilpotent completion'' of $M$, which can be
identified with the $p$-completion (see~\cite[\S 2.2]{DY23}, which follows~\cite[Example 1.16]{Rav84}).
\end{rem}

We would like to specialize to real and complex connective $K$-theory. At odd primes, $\ku$ is decomposable, so we replace it with one of its summands.
\begin{lem}[{Adams~\cite[Lecture 4, Corollary 8]{Ada69}, Johnson--Wilson~\cite[Proposition 2.7]{JW73}}]\label{lem:Adams_summand}
For any prime $p$, there is a $p$-local spectrum $\ell$ and a $p$-local homotopy equivalence
\begin{equation}
    \ku\simeq \bigvee_{i=1}^{p-1} \Sigma^{2i-2}\ell.
\end{equation}
\end{lem}
The spectrum $\ell$ admits an essentially unique $E_\infty$-ring structure~\cite[Corollary 1.4]{BR08} and is called the \term{connective Adams summand}. It is sometimes denoted $\mathit{BP}\ang 1$. When we refer to $\ell$, we have implicitly chosen a prime $p$ to localize at.

At $p = 2$, $\ku\simeq\ell$. When $p$ is odd, $\ko$ is also a sum of shifts of $\ell$.
\begin{rem}
Adams~\cite[Lecture 4, Corollary 8]{Ada69} proves a version of \cref{lem:Adams_summand} for periodic $K$-theory; Johnson--Wilson~\cite[Proposition 2.7]{JW73} explicitly prove the connective case, which they say was a folk theorem.
\end{rem}
In this paper, the  base ring spectrum $R$ will either be $\ko$ or $\ell$. For $R = \ko$, we will restrict to $p = 2$.
We first give some basic properties of these spectra that we will use many times.
\begin{prop}
\label{homotopy_rings}
There are isomorphisms of graded rings
\begin{subequations}
  \begin{align}
    \label{ko-ring}
    \pi_*(\ko) &\cong \Z[\eta, \alpha, \beta]/(2\eta, \eta\alpha, \alpha^2 - 4\beta)\\
    \label{ellring}
    \pi_*(\ell) &\cong \Z_{(p)}[v_1],
  \end{align}
\end{subequations}
where $\abs\eta = 1$, $\abs\alpha = 4$, $\abs\beta = 8$, and $\abs{v_1} = 2p-2$.
\end{prop}
The isomorphisms in \cref{homotopy_rings} are not unique, but one can construct a preferred choice of isomorphism using the map from spin bordism (for $\ko$), resp.\ (the connective cover of) Adams' choice of splitting in \cref{lem:Adams_summand} (for $\ell$).
\begin{lem}[{Baker~\cite[Lemma 5.1]{Bak20}}]
\label{HH_Baker}
\hfill
\begin{enumerate}
    \item\label{A1_ko_H} Let $\cA(1)$ denote the subalgebra of the mod $2$ Steenrod algebra $\cA$ generated by $\Sq^1$ and $\Sq^2$. Then there is a unique isomorphism $H_\ko^*H\cong\cA(1)$ identifying the forgetful map $H_\ko^*\to H_{\Sph}^*H\cong\cA$ with the inclusion $\cA(1)\inj\cA$.
    \item For any prime $p$, let $\cE(1)$ denote the subalgebra of the mod $p$ Steenrod algebra $\cA_p$ generated by the Bockstein $\beta$ and the operation $Q_1\coloneqq \beta\mathcal P^1 - \mathcal P^1\beta$. Then there is a unique isomorphism $H_\ell^*H\cong\cE(1)$ identifying the forgetful map $H_\ell^*H\to H_{\Sph}^*H\cong\cA_p$ with the inclusion $\cE(1)\inj\cA_p$.
    \item If $\cA(1)_*$, resp.\ $\cE(1)_*$ denote the dual Hopf algebras to $\cA(1)$ and $\cE(1)$, there are unique isomorphisms $H_*^\ko H\cong\cA(1)_*$ and $H_*^\ell H\cong\cE(1)_*$ such that the base change maps $H^\Sph_*H\to H^\ko_*H$, resp.\ $H^\Sph_*H\to H^\ell_*H$ are identified with the usual quotient maps $\cA_*\to\cA(1)_*$, resp.\ $(\cA_p)_*\to\cE(1)_*$.
\end{enumerate}
\end{lem}
At $p = 2$, $Q_1 = \Sq^1\Sq^2 + \Sq^2\Sq^1$. 

The following lemma is standard, following for example from the analogous facts for the spaces $B\mathrm O$ and $B\mathrm U$.
\begin{lem}
$\ko$ and $\ell$ satisfy the hypotheses on $R$ in item~\eqref{BLA_conv} of \cref{thm:basic_BL_ss}, i.e., they are bounded-below, cellular spectra of finite type.
\end{lem}

\begin{defn}\label{def_of_g}
Let $R$ be a connective spectrum. Define the function $g^R\colon \Z\to\Z$ as follows: $g^R(n) = m$ means that, counting upwards from zero, $\pi_m(R)$ is the $n^{\mathrm{th}}$ nonzero homotopy group of $R$.

We also let $g\coloneqq g^\ko$ and $g^c(p,\bl)\coloneqq g^\ell$.
\end{defn}
Explicitly,
for $x\in\Z$, write $x = 8n+k$ for $0\le k<8$. Then
\begin{subequations}
\begin{equation}
	g(x) \coloneqq \begin{cases}
		4n+k, &0\le k\le 3\\
		4n+3, &k = 4\\
		4n+4, &5\le k\le 7.
	\end{cases}
\end{equation}
Likewise, for each prime number $p$,
\begin{equation}
	g^c(p, x) \coloneqq \left\lceil \frac{x}{2p-2} \right\rceil.
\end{equation}
\end{subequations}
If $x\ge 0$, then in the $E_2 = E_\infty$-page of the $2$-primary Adams spectral sequence for $\ko$, the element of
minimal Adams filtration in topological degrees $\ge t-s$ has filtration $g(t-s)$. The analogous fact is true for
the $p$-primary Adams spectral sequence for $\ell$ and $g^c(p,\bl)$.

As mentioned in the introduction, we aim to produce a map from the Atiyah--Hirzebruch to the Adams filtrations. 
\begin{thm}\hfill
\label{map_of_filtered_spectra}
\begin{subequations}
\begin{enumerate}
	\item There is a natural transformation of functors $\Mod_\ko\to \Fil_\ko$
	\begin{equation}
		\AtoA^\ko\colon \Post^\ko_\bullet \Longrightarrow \AdFilt_{g(\bullet)}^{\ko}
	\end{equation}
	which, when composed with $(\bl)_{-\infty}$, is homotopy equivalent to $2$-completion.
	\item\label{cpx_case} Let $p$ be a prime. There is a natural transformation of functors $\Mod_\ell\to\Fil_\ell$
	\begin{equation}
		\AtoA^\ell\colon \Post^\ell_\bullet \Longrightarrow \AdFilt_{g^c(p,\bullet)}^{\ell}
	\end{equation}
	which, when composed with $(\bl)_{-\infty}$, is homotopy equivalent to $p$-completion.
\end{enumerate}
\end{subequations}
\end{thm}
This yields a somewhat curious corollary.
\begin{cor}
\label{sseq_map_cor}
Let $X$ be a space or spectrum. After regrading the Adams spectral sequence, there is a map of spectral sequences
from the Atiyah--Hirzebruch spectral sequence for $\ko_*(X)$ to the Adams spectral sequence for $\ko_*(X)_2^\wedge$
(and likewise for $\ell$ at all primes); this map converges to the $2$-completion map $\ko_*(X)_2^\wedge$ (resp.\ to $p$-completion).
\end{cor}
We plan to investigate this map in future work. Specifically, we believe this explains a phenomenon we noticed
studying Adams spectral sequences for computing twisted $\ko$-homology: ``difficult differentials'' in the $2$-primary
Adams spectral sequence for $\ko_*(X)$ tend to correspond to difficult differentials, sometimes of different
lengths, in the AHSS for $\ko_*(X)$, in that one vanishes if and only if the other does. Two examples of 
this phenomenon occur in~\cite[\S 9]{Ped17} and~\cite[\S 5.3]{Deb21}.

\Cref{sseq_map_cor} and the aforementioned examples offer a concrete answer to a frequently asked question in this
area: \emph{why, and to what extent, is the Adams spectral sequence so much more effective than the AHSS for
$\ko$-homology?} We have learned three things towards the answer.
\begin{enumerate}
	\item The Adams filtration makes multiplication by $2$ visible in most cases thanks to the algebraic structure
	on the $E_\infty$-page, but the AHSS doesn't. Thus, the Adams spectral sequence is much better at addressing
	extensions by $2$.
	\item The comparison map from \cref{sseq_map_cor} suggests that Adams and Atiyah--Hirzebruch differentials could obey similar formulas. This is also something we see in examples: for example, higher differentials in both spectral sequences can be described using higher cohomology operations. See for example~\cite{BJ04, AS06}.
    \item However, we do not expect $\AtoA^\ko$ to identify all differentials between these two spectral sequences. For example, the May--Milgram theorem~\cite{MM81} calculates differentials between $h_0$-towers in Adams spectral sequences. The corresponding elements in the AHSS collapse into a single bidegree, so there is not in general an analogue of the May--Milgram theorem available in the AHSS. In this case, we think these differentials will not be in the image of $\AtoA^\ko$.
\end{enumerate}

Moreover, some difficult differentials of the AHSS for $\ko_*(X)$ in low degrees have been identified in \cite{WG18,2025arXiv251225069N}, etc. It is interesting to see how these results about the AHSS can help with the differentials of the Adams spectral sequence.

\begin{rem}
\label{reindexing_2}
\Cref{reindexing} implies that we can rephrase the natural transformation $\Phi^\ko$ from \cref{map_of_filtered_spectra} using the adjunction $g\dashv g^\vdash$ to obtain a natural transformation
\begin{equation}
\label{other_regrade}
    (\AtoA^\ko)^\vdash\colon \mathrm{AH}_{g^\vdash(\bullet)}^\ko \Longrightarrow \Ad_{\bullet}^\ko,
\end{equation}
and analogously for $\ell$ using $(g^c)^\vdash$. That is, one can regrade the Atiyah--Hirzebruch filtration instead if one wants, though in practice we found the result to be cleaner when we lengthen the Adams filtration rather than compress the Atiyah--Hirzebruch filtration.

To obtain~\eqref{other_regrade}, we used a standard if slightly unintuitive fact in category theory, that given an adjunction $f\colon \cat C_1\rightleftarrows \cat C_2\colon g$ and another category $\cat D$, the precomposition functor
\begin{equation}
    \bl\circ g\colon\cat{Fun}(\cat C_2, \cat D)\longrightarrow\cat{Fun}(\cat C_1, \cat D)
\end{equation}
is \emph{left} adjoint to $\bl\circ f$, as $\cat{Fun}(\bl, \cat D)$ is contravariant.
\end{rem}

\subsection{Proof in the complex case}
\label{ss:complex}
In this subsection, we first prove part~\eqref{cpx_case} of \cref{map_of_filtered_spectra}: the compatibility between the Atiyah--Hirzebruch and $p$-primary Adams filtrations for the Adams summand $\ell$. When $p = 2$, this specializes to $\ku$, which is why we call this the ``complex case''. 

Recall that the Postnikov $0$-truncation map $\tau_{\le 0}\colon\ell\to H\Z_{(p)}$ realizes $H\Z_{(p)}$ as an $E_\infty$-$\ell$-algebra. Therefore, following Baker--Lazarev~\cite{BL01}, we can consider the $H\Z_{(p)}$-based, $\ell$-module Adams spectral sequence. Unlike the usual $H\Z_{(p)}$-based, $\Sph$-module Adams spectral sequence, this will be easy to compute.
\begin{thm}
\label{post_ku_HZ}
Let $\Ad^{\ell,H\Z_{(p)}}_\bullet$ denote the $H\Z_{(p)}$-based, $\ell$-module Adams filtration. Then there is a natural isomorphism $\widetilde\AtoA{}^\ell\colon \Post_\bullet^\ell\simeq \Ad_{g^c(p,\bullet)}^{\ell,H\Z_{(p)}}$
of functors $\cat{Mod}_\ell\to\Fil_\ell$.
\end{thm}
We believe this is not a new theorem, but we were unable to find a statement or proof in the literature.
\begin{proof}
Let $M$ be an $\ell$-module. The functors $\Post_\bullet^\ell$ and $\Ad_{g^c(p,\bullet)}^{\ell, H\Z_{(p)}}$ are both defined by constructing the respective filtrations on $\ell$, then levelwise smashing with $M$; therefore it suffices to construct an equivalence between the filtered $\ell$-modules $\Post_\bullet^\ell(\ell)$ and $\Ad_{g^c(p,\bullet)}^{\ell, H\Z_{(p)}}(\ell)$. We will describe these filtrations and see that they are equivalent.

Recall from \cref{homotopy_rings} that $\ell_*\cong\Z_{(p)}[v_1]$ with $\abs{v_1} = 2p-2$. There is an $\ell$-module equivalence $\Sigma^{2p-2}\ell\simeq \tau_{\ge(2p-2)}\ell$ (this is Bott periodicity if $p = 2$), and if $n\not\equiv 0\bmod 2p-2$, $\pi_n(\ell) = 0$. Therefore the Atiyah--Hirzebruch filtration is easy to describe.
\begin{itemize}
    \item If $n\not\equiv 0\bmod 2p-2$, then $\tau_{\ge n}\colon\tau_{\ge n}\ell\to\tau_{\ge(n-1)}\ell$ is an $\ell$-module equivalence.
    \item If $n\equiv 0\bmod 2p-2$, then $\tau_{\ge n}\colon\tau_{\ge n}\ell\to\tau_{\ge(n-1)}\ell$ can be identified via the above equivalences with $\tau_{\ge n}\colon \tau_{\ge n}\ell\to\tau_{\ge(n-2p-2)}\ell$, which by $v_1$-periodicity is $v_1\colon \Sigma^n\ell\to\Sigma^{n-2p+2}\ell$.
\end{itemize}
The $H\Z_{(p)}$-based Adams filtration is not much harder. The zeroth piece is $\ell$, by definition, and the first is the fiber of $\tau_{\le 0}\colon \ell\to H\Z_{(p)}$. By inspecting homotopy groups, the fiber is the Postnikov $(2p-3)$-connected cover, which, as mentioned above, is the map $v_1\colon\Sigma^{2p-2}\ell\to\ell$. Thus, the first piece of the filtration is $\Sigma^{2p-2}\ell$ and the whole story repeats to imply the $n^{\mathrm{th}}$ step of the filtration is $v_1\colon \Sigma^{n(2p-2)}\ell\to\Sigma^{(n-1)(2p-2)}\ell$. Thus, after reindexing the filtration with $g^c(p,\bullet)$, we obtain the Atiyah--Hirzebruch filtration.
\end{proof}
\begin{proof}[Proof of \cref{map_of_filtered_spectra}, part~\eqref{cpx_case}]
For an $E_\infty$-ring spectrum $R$ and a map of $E_\infty$-$R$-algebra spectra $\psi\colon S_1\to S_2$, Baker--Lazarev's construction of the $S_i$-based, $R$-module Adams filtration is manifestly functorial in the bases $S_i$. Apply this with $R = \ell$, $S_1 = H\Z_{(p)}$, $S_2 = H\Z/p$, and $\psi = (\mathrm{mod}\, p)$ to obtain a natural transformation
\begin{equation}
    \AtoA^\ell\colon \Post_\bullet^\ell \xrightarrow[\eqref{post_ku_HZ}]{\widetilde\AtoA{}^\ell} \Ad_{g^c(p,\bullet)}^{\ell, H\Z_{(p)}} \xrightarrow{(\mathrm{mod}\, p)_*} \Ad_{g^c(p,\bullet)}^{\ell}.
    \qedhere
\end{equation}
\end{proof}
\begin{rem}
\Cref{post_ku_HZ} also implies that the $H\Z_{(p)}$-based $\ell$-module Adams spectral sequence for $\ell\wedge X$
has a tractable $E_2$-page. This may be surprising, as the analogous fact over the sphere spectrum is well-known to be false (see for example~\cite[\S 1.1]{BP25}), and when $p = 2$ it is also false over $\ko$ (see below).

One can verify that the $E_2$-page of the $H\Z_{(p)}$-based $\ell$-module Adams spectral sequence is tractable
directly in the standard way, by showing that the algebra $\pi_{-*}\mathrm{Map}_\ell(H\Z_{(p)}, H\Z_{(p)})$ of
$\ell$-module $H\Z_{(p)}$-cohomology operations is a flat $\Z_{(p)}$-module (see Baker--Lazarev~\cite[Theorem 2.1]{BL01}, who remark that their proof follows that of Adams~\cite[Theorem III.15.1]{Ada74}). To check this, use the universal coefficient spectral sequence~\cite[Theorem IV.4.1]{EKMM97}
\begin{equation}
\label{univcoeffs}
    E_2^{*,*} = \Ext_{\Z_{(p)}[v_1]}^{*,*}(\Z_{(p)}, \Z_{(p)}) \Longrightarrow \pi_{-*}\mathrm{Map}_\ell(H\Z_{(p)}, H\Z_{(p)}).
\end{equation}
Using Koszul duality,
\begin{equation}
    \Ext_{\Z_{(p)}[v_1]}^{*,*}(\Z_{(p)}, \Z_{(p)})\cong \Z_{(p)}[\varpi]/(\varpi^2),\qquad\varpi\in\Ext^{1, 2p^n-2},
\end{equation}
so the spectral sequence~\eqref{univcoeffs} collapses, implying the algebra of $\ell$-module $H\Z_{(p)}$ cohomology operations is exterior on a single generator, hence a free $\Z_{(p)}$-module as promised.

The story does not generalize to $2$-local $\ko$: Baker~\cite[Example 5.5]{Bak20} exhibits torsion in the algebra
$\pi_{-*}\Map_{\ko_{(2)}}(H\Z_{(2)}, H\Z_{(2)})$.
\end{rem}
\begin{rem}
Let $k(n)$ denote the connective Morava $K$-theory spectrum at the implicit prime $p$; then $\pi_*(k(n))\cong\Z/p[v_n]$, with $\abs{v_n} = 2p^n-2$, and the fiber of the Postnikov $0$-truncation map $\tau_{\le 0}\colon k(n)\to H\Z/p$ can be identified with $v_n\colon \Sigma^{2p^n-2}k(n)\to k(n)$, in a higher analogue of Bott periodicity. Thus the proof of \cref{post_ku_HZ} generalizes straightforwardly to prove that the Atiyah--Hirzebruch and $H\Z/p$-Adams filtrations of a $k(n)$-module are naturally equivalent, hence that the induced Atiyah--Hirzebruch and $H\Z/p$-Adams spectral sequences are isomorphic after a regrading. This is a known theorem: see the discussion in~\cite{MO_Kn}, which inspired our proof of \cref{post_ku_HZ}.
\end{rem}

\subsection{Proof in the real case}
\label{ss:comparison}
In this subsection, we finish the proof of \cref{map_of_filtered_spectra} by proving the case of $\ko$ at the prime $2$. We will prove this as a corollary of a more general result.
\begin{defn}
\label{rising_defn}
Let $R$ be a connective $E_\infty$-ring spectrum and $S$ be an $E_\infty$-$R$-algebra spectrum. We call $R$ \term{rising with respect to $S$} if for all $n>0$, either $\pi_n(R) = 0$ or, after applying $\bl\wedge_R S$, the Postnikov cover map $\tau_{\ge n}\colon\tau_{\ge n}R\to\tau_{\ge(n-1)}R$ is null-homotopic as an $S$-module map.
\end{defn}

\begin{prop}
\label{it_works_for_ko}
$\ko$ is rising with respect to $H\Z/2$.
\end{prop}

The \textit{rising} property is essential in our proof in the following sense.

\begin{thm}
\label{0_homology_filt}
Let $R$ be a connective $E_\infty$-ring spectrum and $S$ be an $E_\infty$-$R$-algebra spectrum. If $R$ is rising with respect to $S$, there is a natural transformation of functors $\Mod_R\to\Fil_R$
\begin{equation}
    \AtoA^R\colon \mathrm{AH}^R_\bullet \Longrightarrow \Ad_{g^R(\bullet)}^{R, S}
\end{equation}
which is the $S$-nilpotent completion~\cite[\S 5]{Bou79} on abutments, where $g^R\colon\Z\to\Z$ is the degree of the $n^{\mathrm{th}}$ nonzero homotopy group of $R$, given in \cref{def_of_g}.
\end{thm}
When $S = H\Z/p$, the $H\Z/p$-nilpotent completion is the usual $p$-completion~\cite[Example 1.16]{Rav84}.

\begin{lem}[{Baker~\cite[\S 5]{Bak20}}]
\label{A1_module_of_connective_covers}
\hfill
\begin{enumerate}
    \item[(0)] $H_\ko^*(\tau_{\ge 0}\ko) = H_\ko^*(\ko)\cong\Z/2$. Thus by Bott periodicity, $H_\ko^*(\tau_{\ge 8}\ko)\cong\Sigma^8 \Z/2$.
    \item[(1)] $H_\ko^*(\tau_{\ge 1}\ko)\cong\Sigma Q$, where $Q\coloneqq \cA(1)/(\Sq^2)$. 
    \item[(2)] $H_\ko^*(\tau_{\ge 2}\ko)\cong \Sigma^2 J$, where $J\coloneqq\cA(1)/(\Sq^3)$. Since $\pi_3(\ko) = 0$, this is also true for $\tau_{\ge 3}\ko$.
    \item[(4)] $H_\ko^*(\tau_{\ge 4}\ko)\cong \Sigma^4 \uQ$, where $\uQ\coloneqq\cA(1)/(\Sq^1, \Sq^2\Sq^3)$. Since $\pi_k(\ko)$ vanishes for $k = 5,6,7$, this is also true for $\tau_{\ge k}\ko$ for those $k$.
\end{enumerate}
Thus, in particular, for all $n\ge 0$, $H_\ko^*(\tau_{\ge n}\ko)$ is a cyclic $\cA(1)$-module.
\end{lem}
See \cref{Q_ext_F2,Q_ext_uQ,,Q_ext_J,,Q_ext_Q} for pictures of these $\cA(1)$-modules.

We will need the following analogue of the universal coefficient theorem.
\begin{lem}[Relative universal coefficient theorem]
\label{UCT_rel_thm}
Let $R$ be an $E_\infty$-ring spectrum with an $E_\infty$-ring map $\theta\colon R\to H\Z/p$. Then for an $R$-module $M$, there is a natural isomorphism
\begin{equation}
\label{UCT_rel}
    H_R^*(M) \overset\cong\longrightarrow \Hom_{\Z/p}(H^R_*(M), \Z/p).
\end{equation}
\end{lem}
Like for the ordinary universal coefficient theorem, naturality implies~\eqref{UCT_rel} is an isomorphism of $H_R^*H$-modules, where the $H_R^*H$-module structure on $\Hom_{\Z/p}(H^R_*(M), \Z/p)$ is the one defined to be dual to the $H^R_*H$-comodule structure on $H^R_*(M)$.
\begin{proof}
The universal coefficient spectral sequence~\cite[Theorem IV.4.1]{EKMM97}
\begin{equation}
    E_2^{*,*} = \Ext_{\pi_*(H\Z/p)}^{*,*}(H^R_*(M), \pi_*(H\Z/p)) \Longrightarrow H_R^*(M)
\end{equation}
(see also~\cite[(3.1)]{BL01}) is concentrated in degree $0$ because $\pi_*(H\Z/p)$ is a field. Thus it collapses to imply the isomorphism in the lemma statement.
\end{proof}
The proof of the following lemma is likewise straightforward from the relative Künneth spectral sequence~\cite[Theorem 4.1]{EKMM97}.
\begin{lem}[Relative Künneth formula]
\label{rel_kunneth}
With $p$, $R$, and $M$ as in \cref{UCT_rel_thm}, and $N$ another $R$-module, there is a natural isomorphism
\begin{equation}
\label{kunneth_rel_isom}
    H_R^*(M\wedge_R N) \overset\cong\to H_R^*(M)\otimes_{\Z/p} H_R^*(N).
\end{equation}
\end{lem}
Using the Hopf algebra structure on $H_R^*H$, \eqref{kunneth_rel_isom} upgrades to an $H_R^*H$-module isomorphism.
\begin{proof}[Proof of \cref{it_works_for_ko}]
We will show that $\tau_{\ge n}$ is the zero map on $H_*^\ko$, which implies that $\id_{H\Z/2}\wedge_\ko\tau_{\ge n}\simeq 0$, because a map of $H\Z/2$-modules is determined up to homotopy by its behavior on homotopy groups. 
By \cref{UCT_rel_thm} we may equivalently use $H_\ko^*$ instead of $H^\ko_*$. Real Bott periodicity produces a $\ko$-module equivalence $\Sigma^8\ko\simeq\tau_{\ge 8}\ko$ which is commutative with respect to the multiplication by the Bott class; thus for all $k\ge 0$ and $n>0$, $\tau_{\ge(8k+n)}\colon\tau_{\ge(8k+n)}\ko\to\tau_{\ge(8k+n-1)}\ko$ is $\Sigma^{8k}$ of $\tau_{\ge n}\colon\tau_{\ge n}\ko\to\tau_{\ge(n-1)}\ko$. Therefore, it suffices to prove the theorem for $n\le 8$. And since $\pi_n\ko$ vanishes for $n = 3,5,6,7$, we need only consider the four cases $n = 1,2,4,8$.

Recall from \cref{A1_module_of_connective_covers} the $\cA(1)$-modules $H_\ko^*(\tau_{\ge n}\ko)$ for these values of $n$.
All four of these $\cA(1)$-modules are cyclic, so the $\cA(1)$-module map $(\tau_{\ge n})^*\colon H_\ko^*(\tau_{\ge(n-1)}\ko)\to H_\ko^*(\tau_{\ge n}\ko)$ is determined by its value on the generator of the domain. But in all four cases we need, namely $n = 1,2,4,8$, $H_\ko^{n-1}(\tau_{\ge n}\ko) = 0$, so the map must vanish.
\end{proof}

\begin{proof}[Proof of \cref{0_homology_filt}]
As in the proof of \cref{post_ku_HZ}, it suffices to construct a map between the respective filtrations of $R$ itself, as the filtrations on a more general $R$-module $M$ are defined by smashing the filtrations on $R$ levelwise with $M$.

We induct on the filtration degree $n$; the base case is the observation that in degree $0$, both filtrations are just $R$. Thus, assume that we have constructed a map $\phi_{n-1}\colon \tau_{\ge (n-1)}R\to \Ad_{g^R(n-1)}^{R,S}(R)$. We would like to construct a map $\phi_{n}\colon \tau_{\ge n}R\to \Ad_{g^R(n)}^{R,S}(R)$ such that the following diagram commutes.
\begin{equation}
\label{try_to_lift}
\begin{tikzcd}
	{\tau_{\ge n}R} & {\Ad_{g^R(n)}^{R,S}(R)} \\
	{\tau_{\ge(n-1)}R} & {\Ad_{g^R(n-1)}^{R,S}(R).}
	\arrow["{\phi_n}", dashed, from=1-1, to=1-2]
	\arrow["{\tau_{\ge n}}"', from=1-1, to=2-1]
	\arrow["\alpha", from=1-2, to=2-2]
	\arrow["{\phi_{n-1}}", from=2-1, to=2-2]
\end{tikzcd}\end{equation}
If $\pi_n(R)\cong 0$, this is vacuous: $g^R(n) = g^R(n-1)$, so the vertical maps in~\eqref{try_to_lift} are both the identity, up to the natural isomorphism in the Postnikov cover. Thus, assume $\pi_n(R)\not\cong 0$. In this case, to show that a map $\phi_n$ as desired exists, we could equivalently show that $\phi_{n-1}\circ\tau_n$ vanishes after mapping to the cofiber of $\alpha$, which is by definition $S\wedge_R \Ad_{g^R(n-1)}^{R,S}R$. This is base change along $R\to S$, hence also factors through the base changes of $\tau_{\ge n}$ and $\phi_{n-1}$ along $R\to S$:
\begin{equation}\begin{tikzcd}
	& {\tau_{\ge n}R} & {\Ad_{g^R(n)}^{R,S}(R)} \\
	{\textcolor{Fuchsia}{S\wedge_R \tau_{\ge n}R}} & {\tau_{\ge(n-1)}R} & {\Ad_{g^R(n-1)}^{R,S}(R)} \\
	& {\textcolor{Fuchsia}{S\wedge_R \tau_{\ge(n-1)}R}} & {S\wedge_R\Ad_{g^R(n-1)}^{R,S}(R)}
	\arrow["{\phi_n}", dashed, from=1-2, to=1-3]
	\arrow["{\id\wedge_R S}"', from=1-2, to=2-1, Fuchsia]
	\arrow["{\tau_{\ge n}}", from=1-2, to=2-2]
	\arrow["\alpha", from=1-3, to=2-3]
	\arrow["0"', from=2-1, to=3-2, Fuchsia]
	\arrow["{\phi_{n-1}}", from=2-2, to=2-3]
	\arrow["{\id\wedge_R S}", from=2-2, to=3-2]
	\arrow["{\mathrm{cofib}(\alpha)}", from=2-3, to=3-3]
	\arrow["{\id_S\wedge \phi_{n-1}}"', from=3-2, to=3-3, Fuchsia]
\end{tikzcd}\end{equation}
Since $R$ is rising with respect to $S$, $\id_S\wedge_R \tau_{\ge n}\simeq 0$, so the map to the cofiber vanishes, allowing us to choose a lift $\phi_n$ as desired.
\end{proof}

One can prove that $\ell$ is rising with respect to $H\Z/p$ in essentially the same way as for $\ko$, giving another proof of the complex case of \cref{map_of_filtered_spectra}.


\section{The hastened Adams spectral sequence}
	\label{s:MASS}
The main result of this section is \cref{cQthm}, which calculates the $E_2$-page of the HASS for all Postnikov truncations $\tau_{\le n}\ko$ of $\ko$. As we discuss in \cref{SHdefn}, $\tau_{\le 2}\ko$ is Pontrjagin dual to supercohomology, so it is straightforward to infer twisted supercohomology groups from the $\tau_{\le 2}\ko$-homology HASS (\cref{SH_as_hom}). For this reason, whenever we refer to the ``HASS for supercohomology,'' we mean the $\tau_{\le 2}\ko$-homology HASS, leaving the duality implicit.

As an appetizer for why we use the HASS, in \S\ref{s:bakerlazarev} we explain in an example the complications due to differentials that arise when looking at truncations of $\ko$. In \S\ref{subsection:HASS} we introduce the concept of a hastened Adams resolution. From this resolution we obtain a corresponding HASS, and we then present the form of its $E_2$-page.

\subsection{The Baker--Lazarev Adams spectral sequence for supercohomology}
\label{s:bakerlazarev}
In \S\ref{ss:filtered_background}, we introduced the Baker--Lazarev Adams filtration (\cref{filterR}) for $\ko$-modules and briefly discussed the induced spectral sequence. In this subsection, we work through this spectral sequence for the $\tau_{\le 2}\ko$-homology of a point. As mentioned above, and as we discuss in more detail in \cref{SH_as_hom}, this is Pontrjagin dual to $\SH^*(\pt)$.

We discussed the Baker--Lazarev Adams spectral sequence in the abstract in \S\ref{ss:filtered_background}, so let us start by making it more concrete. Combining \cref{thm:basic_BL_ss,HH_Baker}, we learn that the $E_2$-page of the Baker--Lazarev Adams spectral sequence for a $\ko$-module $M$ is
\begin{subequations}
\begin{equation}
    E_2^{s,t} \cong \Ext_{\cA(1)_*}^{s,t}(\Z/2, H^\ko_*(M)),
\end{equation}
since $\cA(1)_* \cong H^\ko_*H$. We find it more intuitive to dualize and work with $H_\ko^*H\cong\cA(1)$ and the $\cA(1)$-module structure on $H_\ko^*(M)$ rather than the $\cA(1)_*$-comodule structure on $H^\ko_*(M)$. If $M$ is finite type, as is the case for all examples in this paper and in~\cite{Debray:2025kfg, DYY3}, we can invoke \cref{UCT_rel_thm} to do this: see~\cite{Ada69, Sal26}. Thus in this case the $E_2$-page is
\begin{equation}
\label{koBL}
    E_2^{s,t} \cong \Ext_{\cA(1)}^{s,t}(H_\ko^*(M), \Z/2).
\end{equation}
\end{subequations}
In this form, the $\ko$-module Baker--Lazarev Adams spectral sequence may be familiar from prior work employing it, such as in~\cite{DY23, Kur25, CDDH26, BottSpiral, DYY25}; more generally, as discussed in~\cite[\S 2.1]{DY23}, if $M = \ko\wedge X$, \eqref{koBL} coincides with the $E_2$-page of the usual Adams spectral sequence for $M$ via the change-of-rings theorem, and there are many examples of such computations in the literature, including but not limited to~\cite{ABP69, Gia73b, Gia73a, Wil73, Gia76, Per78, Bay94, Sto94, Yu95, BB96, Fra11, Pea14, Cam17, Ped17, BC18, Guo18, WW19, DL20, GOPWW20, KPMT20, BR21, Deb21, Davighi:2020uab, DL20c, FH21, LOT21, Lee:2020ewl, LY22, Davighi:2023luh, BGR24, DDHM24, BDDHM25, DW25, Kne25, She25, Dav26, GLWZ26, KKSW26}. We point the reader to Beaudry--Campbell's guide~\cite{BC18} for an excellent introduction to this spectral sequence; in what follows we will assume familiarity with the material in that paper.

For $M = \ko$, $H_\ko^*(\ko)\cong\Z/2$ (\cref{A1_module_of_connective_covers}), and we get the following $E_2$-page.
%
\begin{thm}[{Liulevicius~\cite[Theorem 3]{Liu62}}]
\label{liuthm}
There is an isomorphism of bigraded $\Z/2$-algebras
\begin{equation}
\label{extZ2}
	\E\coloneqq \Ext_{\cA(1)}(\Z/2, \Z/2)\cong\Z/2[h_0, h_1, v, w]/(h_0h_1, h_1^3, vh_1, h_0^2w + v^2),
\end{equation}
where $h_0\in\Ext^{1,1}$, $h_1\in\Ext^{1,2}$, $v\in\Ext^{3,7}$, and $w\in\Ext^{4,12}$.
\end{thm}
The algebra $\E$ acts naturally on all Baker--Lazarev Adams spectral sequences over the ring $\ko$, as the filtering of the $\ko$-action on any $\ko$-module. The class $h_0\in\E$ detects products with $2\in\ko_0$, $h_1$ detects $\eta\in\ko_1$, $v$ detects $\alpha\in\ko_4$, and
$w$ detects the Bott class $\beta\in\ko_8$ (see \cref{homotopy_rings}).


We will also repeatedly use the following computation of $\Ext_{\cA(1)}(\uQ)$, where $\uQ$ is the $\cA(1)$-module we defined in \cref{A1_module_of_connective_covers}, ultimately because $\tau_{\le 2}\ko$ is the cofiber of the Postnikov cover map $\tau_{\ge 4}\ko\to\ko$ and $H_\ko^*(\tau_{\ge 4}\ko)\cong\Sigma^4\uQ$ (\cref{A1_module_of_connective_covers}).

\begin{lem}[{\cite[Figure 2.5]{BR21}}]
\label{ext_ups_Q}
There is an $\E$-module isomorphism\footnote{Adams--Priddy~\cite[\S 3]{AP76} compute $\Ext(\textcolor{MidnightBlue}{\uQ})$ as a bigraded vector space, and the $\E$-module structure is implicit in their method of proof.}
\begin{equation}
    \Ext_{\cA(1)}(\uQ)\cong \E\{\delta, \epsilon\}/(h_1\delta, v\delta + h_0^2\epsilon, w\delta + v\epsilon),
\end{equation}
with $\delta\in\Ext^{0,0}$ and $\epsilon\in\Ext^{1,5}$.
\end{lem}

We will use this to show:

\begin{prop}\hfill
\label{non_modified_Adams}
\begin{enumerate}
    \item\label{E2calc} The $E_2$-page of the Baker--Lazarev $\ko$-module Adams spectral sequence for $\tau_{\le 2}\ko$ is isomorphic as $\E$-modules to
    \begin{equation}
        \textcolor{BrickRed}{\E\cdot 1}\oplus \textcolor{MidnightBlue}{\E\cdot\set{\delta, \epsilon}/(h_1\delta, v\delta + h_0^2\epsilon, w\delta + v\epsilon)},
    \end{equation}
    where $1\in\Ext^{0,0}$, $\delta\in\Ext^{0,5}$, and $\epsilon\in\Ext^{1,10}$.
    \item\label{runsseq} If $r\ge 2$ and $r\ne 3$, $d_r = 0$. On the $E_3$-page, $d_3(1) = 0$, $d_3(\delta) = v\cdot 1$, and $d_3(\epsilon) = w\cdot 1$.
\end{enumerate}
\end{prop}

The classes $\delta$ and $\epsilon$ are responsible for the new $d_3$-differentials that are absent in the ordinary $\ko$-Adams spectral sequence. This is a serious problem for our goal of computing $(\tau_{\le 2}\ko)_*(BG, a, b)$ for various groups $G$ and twisting data $(a, b)$: if there are differentials in the Baker--Lazarev Adams spectral sequence for $\ko_*(BG, a, b)$, they interact with the $d_3$s we found in \cref{non_modified_Adams} in ways that significantly complicate the computation. There are plenty of small groups $G$ and twisting data $(a, b)$ such that the Adams spectral sequence for $\ko_*(BG, a, b)$ has differentials in low degrees, including several examples in~\cite{Bay94, BB96, Fra11, Cam17, Ped17, GOPWW20, Deb21, Davighi:2020uab, Davighi:2022icj, LY22, BGR24, DDHM24, BDDHM25, Debray:2025kfg, Kne25, Dha26}.

We turn to the HASS in this section to solve this problem. In particular, we will modify this family of Adams spectral sequences in a way that cleanly separates the
differentials coming from $\ko_*(BG, a, b)$ from those differentials coming from $\tau_{\le 2}\ko$. The resulting
spectral sequences will be simpler in general. 

\begin{figure}[htbp]
\begin{subfigure}[c]{0.45\textwidth}
    \includegraphics[width=\linewidth]{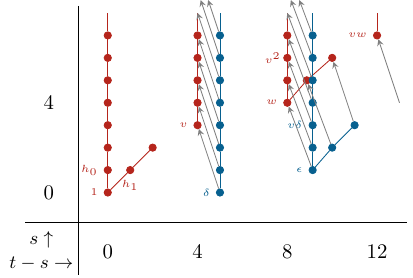}
\end{subfigure}%
\qquad\qquad
\begin{subfigure}[c]{0.45\textwidth}
       \includegraphics[width=\linewidth]{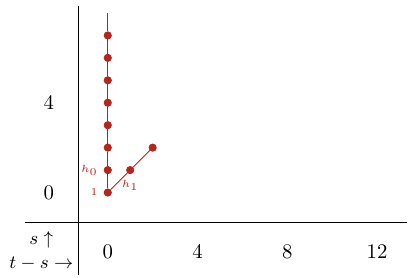}
\end{subfigure}
\caption{Left: the $E_2 = E_3$-page of the Baker--Lazarev Adams spectral sequence for $\tau_{\le 2}\ko$, from
\cref{non_modified_Adams}. Right: the $E_4= E_\infty$-page.}
\label{non_modified_fig}
\end{figure}

The broad outline of \cref{non_modified_Adams}, and its proof, applies in the same way to $\tau_{\le n}\ko$ for $n>1$, as well as (slightly more broadly) to $\tau_{\le n}\ell$.
As a first step, we must calculate $H_\ko^*(\tau_{\le 2}\ko)$.
\begin{lem}
\label{coh_of_superhom}
There is an $\cA(1)$-module isomorphism $H_\ko^*(\tau_{\le 2}\ko)\cong \textcolor{BrickRed}{\Z/2} \oplus \textcolor{MidnightBlue}{\Sigma^5\uQ}$.
\end{lem}
\begin{proof}
Essentially by
definition, $H_\ko^*(\ko)\cong\Z/2$, which necessarily has trivial $\cA(1)$-action. By \cref{A1_module_of_connective_covers}, there
is an isomorphism of $\cA(1)$-modules $H_\ko^*(\tau_{\ge 3}\ko)\cong\Sigma^4\uQ$. The cofiber sequence
\begin{equation}
	\tau_{\ge 3}\ko\longrightarrow \ko \overset{\tau_{\le 2}}{\longrightarrow} \tau_{\le 2}\ko
\end{equation}
induces a long exact sequence on $H_\ko^*$; exactness uniquely characterizes the $\cA(1)$-module structure on $H_\ko^*(\tau_{\le 2}\ko)$.
\end{proof}
\begin{proof}[Proof of \cref{non_modified_Adams}]
\Cref{coh_of_superhom} implies that to prove part~\eqref{E2calc} of the proposition statement, it suffices to compute $\Ext_{\cA(1)}(\textcolor{BrickRed}{\Z/2} \oplus \textcolor{MidnightBlue}{\Sigma^5\uQ})$ as an $\E$-module. $\Ext$ commutes with direct sums, and $\Ext_{\cA(1)}^{s,t}(\Sigma^n M)\cong\Ext_{\cA(1)}^{s,t+n}(M)$ for an $\cA(1)$-module
$M$, so to compute the $E_2$-page of the spectral sequence, it suffices to know
$\Ext_{\cA(1)}(\textcolor{BrickRed}{\Z/2})$ and $\Ext_{\cA(1)}(\textcolor{MidnightBlue}{\uQ})$, which we gave in \cref{liuthm,ext_ups_Q} respectively.

On to part~\eqref{runsseq}.
For all $r\ge 2$, $d_r(1) = 0$ for degree
reasons; since differentials commute with the $\E$-action on the $E_r$-page, to figure out the rest of the spectral
sequence all we have to do is compute $d_r(\delta)$ and $d_r(\epsilon)$.

This spectral sequence converges to
$\pi_*(\tau_{\le 2}\ko)_2^\wedge$, so $\Z_2$ in degree $0$, $\Z/2$ in degrees $1$ and $2$, and $0$ in all other
degrees. In particular, all classes with $t-s\ge 2$ must vanish before the $E_\infty$-page. This uniquely
determines all differentials: the class $v\in E_2^{3,7}$ must be in the image of a differential, but the only way
for this to happen is for $d_3(\delta) = v$. $\E$-linearity of differentials then implies $d_3(v\delta) = v^2$;
since $v\delta = h_0^2\epsilon$ and $v^2 = h_0^2 w$, if $d_3(\epsilon) = x$ then $h_0^2x = h_0^2w$. The map
$h_0^2\colon E_2^{4,12}\to E_2^{6,14}$ is injective, so $x = w$. To summarize, the spectral sequence is completely
determined by the $\E$-action and the differentials $d_3(\delta) = v$, $d_3(\epsilon) = w$, and $d_r(1)= 0$ for all
$r$. We draw the $E_3$-page with these differentials in \cref{non_modified_fig}, left, and the $E_4 =
E_\infty$-page in \cref{non_modified_fig}, right. 
\end{proof}

\subsection{Hastened Adams spectral sequences}\label{subsection:HASS}
%

The computation in \cref{non_modified_Adams} suggests that the ordinary Baker--Lazarev Adams filtration is not well adapted to truncated theories. Given an $R$-module map $f\colon M\to N$, let $Cf$ denote the cofiber of $f$. 
 We will use the idea of Behrens–Hill–Hopkins–Mahowald \cite{BHHM08} to modify the Adams resolution on $Cf$ so that the information coming from the truncation is directly input into the resolution. This will be done in \cref{def:hassadamsresolution}. The resulting spectral sequence converges to $\pi_*(Cf)_2^\wedge$, and we calculate the $E_2$-page in \cref{BHHM}.

Let $R$ be an $E_\infty$-ring spectrum with an $E_\infty$-ring map $R\to H\Z/p$ for some prime $p$. Recall from \cref{filterR} the notation $\Ad_\bullet^R(M)$ for the $H\Z/p$-based, $R$-module Adams filtration on an $R$-module $M$. For $s\ge 0$, we will let $\alpha_M\colon \Ad_{s+1}^R(M)\to\Ad_s^R(M)$ denote the map between adjacent entries in this filtration. The values of $s$ and $p$ will be clear from context when needed, and if $M$ is clear from context, we will abbreviate $\alpha_M$ by $\alpha$.
\begin{defn}[{Behrens--Hill--Hopkins--Mahowald~\cite[\S 3]{BHHM08}, Bruner--Rognes~\cite[Definition 12.31]{BR21}}]
\label{def:hassadamsresolution}
Given an $R$-module map $f\colon M\to N$, let $Cf$ denote the cofiber of $f$. Suppose that there exists a lift $g\colon M\to \Ad^R_m(N)$ of $f$, i.e.\ there is a homotopy of $R$-module maps $\alpha_N^m\circ g\simeq f$.
The \term{hastened Adams resolution} of $Cf$ is the filtered $R$-module $(W_\bullet,
\phi_\bullet)$ defined by
\begin{equation}
	W_s \coloneqq \begin{cases}
		\mathrm{cofib}(\alpha_N^{m-s}\circ g), &s\le m\\
		\Ad_s^R(Cf), &s\ge m
	\end{cases}.
\end{equation}
For $s = m$ these cases agree.

The map $\phi_s$ equals $\alpha_{Cf}$ for $s\ge m$; for $s\le m$, define $\phi_s$ to be the map
$C(\alpha^{m-s}_N\circ g)\to C(\alpha^{m-s-1}_N\circ g)$ induced by $\alpha_N$. (See~\cite[(12.6)]{BR21}.)
\end{defn}
A priori, the hastened Adams resolution depends on both $f$ and $g$; in fact, the
dependence on $g$ is very weak, but this construction very much depends on $f$.
\begin{defn}[{Hastened Adams spectral sequence~\cite[Proposition 3.1]{BHHM08}, \cite[Definition 12.31]{BR21}}]
\label{hastened_defn}
The \term{hastened Adams spectral sequence} for $Cf$ is the spectral sequence associated to the hastened Adams
resolution.
\end{defn}
Let $\cD(H^R_*H)$ denote the derived $\infty$-category of graded $H^R_*H$-comodules; the objects of this $\infty$-category
are chain complexes of graded $H^R_*H$-comodules, and the homotopy equivalences are the grading-preserving
quasi-isomorphisms. For an object $L_\bullet$ of $\cD(H^R_*H)$, we let $\Sigma^t L_\bullet$ denote the chain
complex with each entry $(\Sigma^t L)_i$ equal to the graded $H^R_*H$-module $\Sigma^t (L_i)$, and let
$L_\bullet[s]$ denote the shift with respect to the chain complex grading.
\begin{lem}\label{lem:underg}
The map $\alpha^m\colon \textcolor{MidnightBlue}{\Ad_m^R(N)}\to \textcolor{BrickRed}{N}$ induces a natural homotopy class of morphisms in $\cD(H_*^RH)$ of the form
\begin{equation}
    \underline\alpha^m\colon \textcolor{MidnightBlue}{\Sigma^m H^R_*(\Ad_m^R(N))[-m]} \longrightarrow \textcolor{BrickRed}{H^R_*(N)}.
\end{equation}
\end{lem}
\begin{proof}
We prove the lemma for the case $m = 1$ and then iterate to get the general case. Consider the long exact sequence on $H^R_*$ induced by the fiber sequence $\textcolor{MidnightBlue}{\Ad_1^R(N)}\to \textcolor{BrickRed}{N}\to \textcolor{Green}{N\wedge_R H}$:
\begin{subequations}
\begin{equation}
\label{LEStoSES}
    \dotsb\to H^R_{k+1}(\Ad_1^R(N)) \xrightarrow 0 \textcolor{BrickRed}{H^R_{k+1}(N)} \to \textcolor{Green}{H^R_{k+1}(N\wedge_R H)} \xrightarrow\partial \textcolor{MidnightBlue}{H^R_k(\Ad_1^R(N))} \xrightarrow 0 H^R_k(N) \to\dotsb
\end{equation}
Since $\Ad_1^R(N)\to N$ is the zero map on $H^R_*$, \eqref{LEStoSES} degenerates into a short exact sequence of the form
\begin{equation}
\label{SES_loop}
    \shortexact{\textcolor{BrickRed}{H_{*+1}^R(N)}}{\textcolor{Green}{H_{*+1}^R(N\wedge_R H)}}{\textcolor{MidnightBlue}{H_*^R(\Ad_1^R(N))}}.
\end{equation}
\end{subequations}
That is, we have produced an $H^R_*H$-comodule extension of $\textcolor{MidnightBlue}{H_*^R(\Ad_1^R(N))}$ by $\textcolor{BrickRed}{H_{*+1}^R(N)}\cong \textcolor{BrickRed}{\Sigma^{-1}H_*^R(N)}$, or equivalently of $\textcolor{MidnightBlue}{\Sigma H_*^R(\Ad_1^R(N))}$ by $\textcolor{BrickRed}{H_*^R(N)}$, which is classified up to homotopy by its class in
\begin{equation}
    \Ext_{H^R_*H}^1(\textcolor{MidnightBlue}{\Sigma H^R_*(\Ad_1^R(N))}, \textcolor{BrickRed}{H^R_*(N)})\cong 
    \Hom_{\cD(H_*^RH)}(\textcolor{MidnightBlue}{\Sigma H^R_*(\Ad_1^R(N))[-1]}, \textcolor{BrickRed}{H^R_*(N)}).
    \qedhere
\end{equation}
\end{proof}
At this point, assume $g_*\colon H^R_*(M)\to H^R_*(\Ad_m^R(N))$ is injective. Consider the map $\underline g\coloneqq \underline\alpha^m\circ g$, which is a morphism in $\cD(H^R_*H)$ of the form
\begin{equation}\label{overline_g}
	\underline g\coloneqq \underline\alpha^m\circ g\colon \Sigma^m H^R_*(M)[-m] \longrightarrow H^R_*(N).
\end{equation}
\begin{thm}[{Behrens--Hill--Hopkins--Mahowald~\cite[\S 3]{BHHM08}, Bruner--Rognes~\cite[Proposition 12.33]{BR21}}]
\label{BHHM}
Assume that $g_*\colon H^R_*(M)\to H^R_*(\Ad_m^R(N))$ is injective and let $\mathcal Q_\bullet$ denote the cofiber in $\cD(H^R_*H)$
of the map $\underline g$ from~\eqref{overline_g}. Then, the $E_2$-page of the hastened Adams spectral sequence from
\cref{hastened_defn} is
\begin{equation}
E_2^{s,t} = \Ext_{\cD(H^R_*H)}^{s,t}(\Z/p, \mathcal Q_\bullet),
\end{equation}
and this spectral sequence converges strongly to $\pi_*(Cf)_p^\wedge$ with the $R$-module Adams filtration.
\end{thm}
As always, with finite type hypotheses, we may pass to cohomology. The derived $\infty$-category of
$H_R^*H$-modules, $\cD(H_R^*H)$, is defined analogously to $\cD(H^R_*H)$. 
Analogously to the map in \eqref{overline_g}, we get a morphism in $\cD(H_R^*H)$
\begin{equation}\label{coh_over_g}
	\overline g\colon H_R^*(N)\longrightarrow \Sigma^{-m} H_R^*(M)[m].
\end{equation}
\begin{cor}\label{coh_MASS}
Assume $M$ and $N$ are finite type and that $g^*\colon H_R^*(N)\to H_R^*(M)$ is surjective. Let $\mathcal Q^\bullet$ denote
the cofiber in $\cD(H_R^*H)$ of the map $\overline g$ from~\eqref{coh_over_g}; then the $E_2$-page
of the hastened Adams spectral sequence from \cref{hastened_defn} is
\begin{equation}
	E_2^{s,t} = \Ext_{\cD(H_R^*H)}^{s,t}(\mathcal Q^\bullet, \Z/p).
\end{equation}
\end{cor}


\begin{rem}[The synthetic point of view]
\label{synthetic}
A great deal of recent work on the Adams spectral sequence uses what have come to be called ``synthetic techniques;'' in this remark, we briefly explain the connection to the HASS.

For $E$ an \term{Adams-type ring spectrum} (see~\cite[\S 1.1]{Pst23} and~\cite[Condition III.13.3]{Ada74}; this includes $H\Z/p$ for all primes $p$), Pstrągowski~\cite{Pst23} constructs a symmetric monoidal $\infty$-category $\cat{Syn}_E$ and lax monoidal functors
\begin{equation}
    \cat{Sp}\xrightarrow\nu \cat{Syn}_E \xrightarrow\sigma \Fil(\cat{Sp})
\end{equation}
whose composite is the $E$-based, $\Sph$-module Adams filtration (\cref{filterR}). Just as filtered spectra refine spectral sequences, objects of $\cat{Syn}_E$ refine and generalize $E$-based Adams spectral sequences, and $\cat{Syn}_E$ provides a useful setting for proving things about these Adams spectral sequences.

The functor $\nu$ does not preserve cofiber sequences (see~\cite[Proposition 4.9(2)]{vN25}). Therefore, if $f\colon X\to Y$ is a map of spectra and $C\in\cat{Syn}_E$ is the cofiber of $\nu(f)$, $\sigma(C)$ is a filtered spectrum that, heuristically speaking, ought to correspond to something like an Adams spectral sequence, but in general will not be the Adams spectral sequence for the cofiber $Cf$ of $f$. Instead, as explained in~\cite[Remark 3.21]{LWX24} (see also~\cite[\S 4.3.4]{vN25} and~\cite[\S 2.2]{CD24}), the spectral sequence associated to the filtered spectrum $\sigma(C)$ is precisely the HASS of $Cf$!

To connect this approach to the HASSes appearing in our paper, one would want a generalization of $\cat{Syn}_E$ and its connection to the Adams filtration to a category of synthetic $R$-module spectra affiliated with Baker--Lazarev's $R$-module Adams filtration. To our knowledge, this has not appeared in the literature, though we believe it is likely that the approach for $R = \Sph$ generalizes. Patchkoria--Pstrągowski~\cite{PP26} construct synthetic analogues in a very general setting including module spectra.

Suppose such a generalization existed, in the form of a stable $\infty$-category $\cat{Syn}_{H\Z/2}^\ko$ and lax symmetric monoidal functors $\nu\colon\Mod_\ko\to\cat{Syn}_{H\Z/2}^\ko$ and $\sigma\colon\cat{Syn}_{H\Z/2}^\ko\to\Fil_\ko$ such that $\sigma\circ\nu\simeq\Ad^\ko$ and, like over $\Sph$, $\sigma$ admits a right adjoint $\rho$.\footnote{Since $\nu$ is lax symmetric monoidal, for any $\ko$-module $M$, $\nu(M)$ has a natural $\nu(\ko)$-module structure, suggesting that $\cat{Syn}_{H\Z/2}^\ko$ could be the $\infty$-category of $\nu(\ko)$-modules in $\cat{Syn}_{H\Z/2}$. However, this appears to be the wrong definition: it lifts the action of the Adams spectral sequence for $\ko$ on the ordinary Adams spectral sequence for $M$. In general, the usual and Baker--Lazarev Adams spectral sequences do not coincide, e.g., for $\tmf_{(3)}$ as a module over itself: compare Hill's computation~\cite[\S 2]{Hil07} of the latter with Culver's computation~\cite{Cul21} of the former.}
Then it would be possible to lift \cref{map_of_filtered_spectra} to synthetic spectra as follows. Use \cref{reindexing_2} to rephrase the AHSS-to-Adams map $\AtoA^\ko$ as a natural transformation $(\AtoA^\ko)^\vdash\colon\mathrm{AH}_{g^\vdash(\bullet)}^\ko\Rightarrow\Ad_\bullet^\ko$, i.e.
\begin{subequations}
\begin{equation}
    (\AtoA^\ko)^\vdash\in \Hom_{\cat{Fun}(\Mod_{\ko}, \Fil_\ko)}(\mathrm{AH}^\ko\circ g^\vdash, \Ad^\ko).
\end{equation}
Since $\Ad^\ko\simeq \sigma\circ\nu$, we may use the $\sigma\dashv\rho$ adjunction to port $(\AtoA^\ko)^\vdash$ to a natural transformation
\begin{equation}
    \AtoA^\ko_{\mathrm{Syn}}\in \Hom_{\cat{Fun}(\Mod_{\ko}, \cat{Syn}_{H\Z/2}^\ko)}(\rho(\mathrm{AH}^\ko\circ g^\vdash), \nu),
\end{equation}
\end{subequations}
that is, a natural transformation $\AtoA^\ko_{\mathrm{Syn}}\colon \rho(\mathrm{AH}^\ko_{g^\vdash(\bullet)}) \Rightarrow\nu$ in the setting of synthetic $\ko$-modules. It would be interesting to give a direct construction of $\widetilde\Psi$ and see whether this could provide an another conceptual explanation for the existence of the HASSes we constructed in this paper, or whether such a construction could shed light on \cref{monoidality_conjecture}.
\end{rem}
\subsection{Applying the HASS to truncations of \texorpdfstring{$\ko$}{ko}}
\label{ss:apply_HASS_to_ko}
In this section, we prove the following main theorem about the constructions of the HASS for truncations of $\ko$. Let $C_k\coloneqq \Sigma^{-k}H_\ko^*(\tau_{\ge k}\ko)$. We describe these $\cA(1)$-modules explicitly in \cref{A1_module_of_connective_covers}; in particular, the desuspension by $k$ puts the lowest-degree nonzero element of $C_k$ in degree $0$.
\begin{thm}\hfill
\label{cQthm}
\begin{enumerate}
    \item\label{g4_defn}
    There is a map of $\Z^2$-graded $\Ext(\Z/2)$-modules
    \begin{equation}
        f_{n-1}\colon\Ext_{\cA(1)}^{s,t}(C_{n}) \longrightarrow
            \Ext_{\cA(1)}^{s+g(n), t+g(n)+n}(\Z/2),
    \end{equation}
    which is induced from the Postnikov cover map $\tau_{\ge n}\colon \tau_{\ge n}\ko\to\ko$.
    \item \label{cQ_defn}
    There is a functor $\mathcal Q_n^{*,*}$ from $\cA(1)$-modules to $\Z^2$-graded $\Ext_{\cA(1)}(\Z/2)$-modules which commutes with direct sums and such that for all $\cA(1)$-modules $M$, there is a long exact sequence
    \begin{equation}
    \label{cQLES}
    \dotsb\to \Ext_{\cA(1)}^{s,t}(C_{n}\otimes M) \xrightarrow{f_{n-1}}
            \Ext_{\cA(1)}^{s+g(n), t+g(n)+n}(M)
    \xrightarrow{\partial_{n-1}}
    \mathcal Q_n^{s,t}(M) \to
    \Ext_{\cA(1)}^{s+1,t}(C_{n}\otimes M)\xrightarrow{f_{n-1}}\dotsb
    \end{equation}
    \item\label{HASS_defn}
    Let $M$ be a finite type, bounded-below $\ko$-module.\footnote{The finite-type assumption appears only for technical reasons, and is true in all examples we consider in this paper. Likewise, all examples we consider in this paper are bounded below.}
    Then the HASS for $\tau_{\le n}\ko\wedge_\ko M$ converges strongly and has signature
    \begin{equation}
        E_2^{s,t} = \mathcal Q_n^{s,t}(H^*_\ko(M)) \Longrightarrow \pi_{t-s}(\tau_{\le n}\ko\wedge_\ko M)_2^\wedge.
    \end{equation}
    The map $\pi_*(M)\to \pi_*(\tau_{\le n}\ko\wedge_\ko M)$ lifts to a map from the ordinary Adams spectral sequence to the HASS which on the $E_2$-page is the map $\partial_n$ from~\eqref{cQLES}.
\end{enumerate}
\end{thm}
Notice the shift in $n$ in part~\eqref{HASS_defn} -- to calculate the $E_2$-page of the HASS for $\tau_{\le n}\ko\wedge_\ko M$, we use the map $f_n$ with domain $\Ext(C_{n+1})$!

An analogous result is true for truncations of $\ell$ at an arbitrary prime, and the proof is identical. In both cases, the proof amounts to verifying the hypotheses in \cref{coh_MASS} and then making explicit the details of the calculation of the $E_2$-page of the HASS.
\begin{rem}
The reader may be wondering -- where on Earth did $g(n)$ and $n+g(n)$ come from in \cref{cQthm}? To answer this question, we first clarify what we mean by ``[$f_{n-1}$] is induced from the Postnikov cover map'' in \cref{HASS_defn}. Specifically, $f_{n-1}$ is the composition of the following maps.
\begin{subequations}
\begin{enumerate}
    \item First, the map obtained by applying $\Ext_{\cA(1)}$ to the Atiyah--Hirzebruch-to-Adams map
    \begin{equation}
    \label{Ext_filtmap}
        \tau_{\ge n}\ko \longrightarrow \Ad_{g(n)}^\ko(\ko)
    \end{equation}
    that we built in \cref{map_of_filtered_spectra}. This map preserves $(s,t)$-bidegrees of Ext classes.
    \item The second map uses the fact that the Ext of the $m^{\mathrm{th}}$ Adams cover of $M$ gives rise to the filtration-$m$ elements of the Adams spectral sequence for $M$. Specifically, we get a map
    \begin{equation}
        \label{adams_map_on_ext}
        \Ext_{\cA(1)}^{s,t}(H_\ko^*(\Ad_m^\ko(M))) \longrightarrow
        \Ext_{\cA(1)}^{s+m,t+m}(H_\ko^*(M)).
    \end{equation}
\end{enumerate}
\end{subequations}
Now we put these together. By definition, $H_\ko^*(\tau_{\ge n}\ko)\cong \Sigma^{n}C_{n}$, so we have an isomorphism
\begin{subequations}
\begin{equation}\label{1stcomp}
    \Ext_{\cA(1)}^{s,t}(C_{n+1}) \xrightarrow\cong
    \Ext_{\cA(1)}^{s,t+n}(H_\ko^*(\tau_{\ge n}\ko)).
\end{equation}
Now apply $\Ext_{\cA(1)}^{s,t+n}$ to~\eqref{Ext_filtmap}:
\begin{equation}\label{2ndcomp}
    \Ext_{\cA(1)}^{s,t+n}(H_\ko^*(\tau_{\ge n}\ko)) \to \Ext_{\cA(1)}^{s,t+n}(H_\ko^*(\Ad_{g(n)}^\ko\ko)).
\end{equation}
Now apply~\eqref{adams_map_on_ext} with $m = g(n)$ and again shifted by $n$ in the $t$-direction:
\begin{equation}\label{3rdcomp}
    \Ext_{\cA(1)}^{s,t+n}(H_\ko^*(\Ad_m^\ko(\ko))) \longrightarrow
        \Ext_{\cA(1)}^{s+g(n),t+n+g(n)}(H_\ko^*(\ko)),
\end{equation}
\end{subequations}
and $H_\ko^*(\ko)\cong\Z/2$ (\cref{A1_module_of_connective_covers}). The map $f_{n-1}$ is the composition of~\eqref{1stcomp}, \eqref{2ndcomp}, and~\eqref{3rdcomp}.
\end{rem}

The construction of \cref{HASS_defn} follows the discussion in the previous section closely. Specializing to truncations of $\ko$, this reduces to verifying the following two hypotheses in \cref{BHHM}, as well as their analogues with $\ell$ in place of $\ko$.
\begin{enumerate}
    \item\label{from_filtmap} For all $n\ge 0$, the Postnikov $n$-connected covering map $\tau_{\ge n}\ko\to\ko$ lifts to the $g(n)^{\mathrm{th}}$ Adams cover $\Ad_{g(n)}^\ko(\ko)$ of $\ko$. 
    \item\label{verify_surj} The induced map $H_\ko^*(\Ad_{g(n)}^\ko(\ko))\to H_\ko^*(\tau_{\ge n}\ko)$ is surjective.
\end{enumerate}

In fact, we already took care of item~\eqref{from_filtmap} in \cref{map_of_filtered_spectra}. Thus, we will focus on part~\eqref{verify_surj} in this subsection, which we restate as the following proposition.
\begin{prop}
\label{prop_surj_on_Hko}
Let $\gamma_n\colon \tau_{\ge n}\ko\to\Ad^\ko_{g(n)}(\ko)$ be any $\ko$-module map whose composition with the $g(n)$-fold Adams cover $\alpha_n\colon \Ad_{g(n)}^\ko(\ko)\to\ko$ is the Postnikov $n$-connected covering map $\tau_{\ge n}\colon\tau_{\ge n}\ko\to\ko$. Then $\gamma_n^*\colon H_\ko^*(\Ad^\ko_{g(n)}(\ko))\to H_\ko^*(\tau_{\ge n}\ko)$ is surjective.
\end{prop}
To reiterate, by \cref{BHHM}, this and \cref{map_of_filtered_spectra} prove \cref{cQthm}.

\begin{figure}[htbp]
\includegraphics{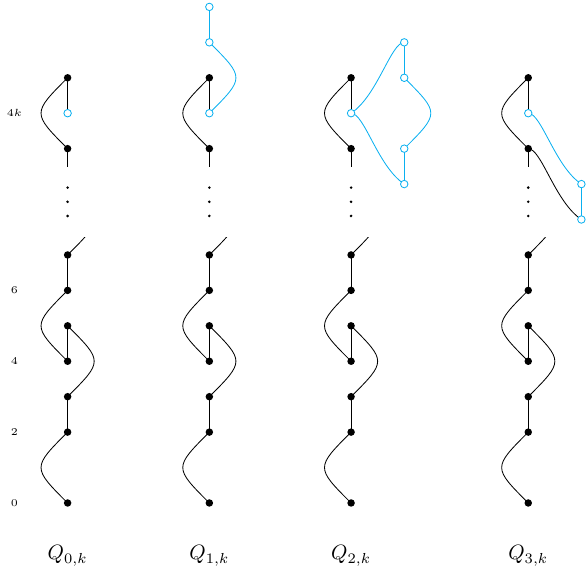}
\caption{The Milgram modules $Q_{j,k}$ (\cref{milgram_defn}). Each module is an extension of the cyclic $\cA(1)$-module $C_{j,k}$ (in blue, empty vertices) by the seagull module $\Upsilon_k$ (in black, filled vertices).}
\label{milgram_modules_picture}
\end{figure}

To get started, we should identify the $\cA(1)$-modules that $\gamma_n^*$ maps between. We have already seen $H_\ko^*(\tau_{\ge n}\ko)$ in \cref{A1_module_of_connective_covers}: if $n = 4k+\ell$, with $0\le\ell<4$, $H_\ko^*(\tau_{\ge n}\ko)$ is isomorphic to $\textcolor{BrickRed}{\Sigma^{8k}\Z/2}$ for $\ell = 0$, $\textcolor{Green}{\Sigma^{8k+1}Q}$ for $\ell = 1$, $\textcolor{MidnightBlue}{\Sigma^{8k+2}J}$ for $\ell = 2$, and $\textcolor{Fuchsia}{\Sigma^{8k+4}\uQ}$ for $\ell = 3$.
\begin{defn}
\label{seagull_defn}
Inductively define the \term{seagull module} $\Upsilon_k$ as follows: let $\Upsilon_1\coloneqq\cA(1)\otimes_{\cA(0)}\Z/2$. Then one can directly verify that, up to isomorphism of extensions, there is a unique nonsplit $\cA(1)$-module extension
\begin{equation}
    \shortexact{\Sigma^{4k-4}\Upsilon_1}{\Upsilon_k}{\Upsilon_{k-1}}.
\end{equation}
\end{defn}
Adamyk~\cite{Ada23} introduced the name ``seagulls'' for these modules, as well as the seagull-like notation $\Upsilon_k$; see (\textit{ibid.}, Figure 1) for pictures of $\Upsilon_k$ for some small $k$.
\begin{defn}[{Milgram~\cite[Definitions 1.6, 1.8, and 1.10]{Mil75}}]
\label{milgram_defn}
Let $k\ge 0$ and $j\in\Z/4$ and define the $\cA(1)$-module
\begin{equation}
    C_{j,k}\coloneqq\begin{cases}
        \textcolor{BrickRed}{\Sigma^{4k}\Z/2}, &j = 0\\
        \textcolor{Green}{\Sigma^{4k}\uQ}, &j = 1\\
        \textcolor{MidnightBlue}{\Sigma^{4k-2}J}, &j = 2\\
        \textcolor{Fuchsia}{\Sigma^{4k-3}Q}, &j = 3.
    \end{cases}
\end{equation}
Then, the reader may verify that, up to isomorphism of extensions, there is a unique nonsplit $\cA(1)$-module extension
\begin{equation}
    \shortexact{\Upsilon_k}{Q_{j,k}}{C_{j,k}}.
\end{equation}
The middle term of this extension is called a \term{Milgram module}. We refer to the submodule $\Upsilon_n$ as the \term{seagull submodule} and the quotient $C_{j,k}$ as the \term{cyclic quotient} of $Q_{j,k}$.
\end{defn}
We draw the Milgram modules, their seagull submodules, and their cyclic quotients in \cref{milgram_modules_picture}.
\begin{prop}[{Milgram~\cite[Corollary 1.15]{Mil75}}]
\label{milgram_tensor}
For all $j_1,k_1,j_2,k_2$, there is a finite-rank free $\cA(1)$-module $F_{j_1,k_1,j_2,k_2}$ and an $\cA(1)$-module isomorphism
\begin{equation}
        Q_{j_1, k_1} \otimes Q_{j_2, k_2} \cong Q_{j_1+j_2, k_1+k_2+\delta} \oplus F_{j_1,k_1,j_2,k_2},
\end{equation}
where $j_1+j_2$ is interpreted modulo $4$ and
\begin{equation}
        \delta = \begin{cases}
                -1 &\text{if $\set 3\subseteq\set{j_1,j_2}\subseteq \set{2,3}$}\\
                \phantom{-}1 &\text{if $\set 1\subseteq\set{j_1,j_2}\subseteq \set{1,2}$}\\
                \phantom{-}0, &\text{else.}
        \end{cases}
\end{equation}
\end{prop}
\begin{cor}\label{which_milgram}
For a positive integer $k$, write $k = 4\ell+m$ for $m\in\set{1,2,3,4}$ and $\ell\in\Z$. Let $\delta_{m, 4} = 1$ if $m =
4$ and $\delta_{m, 4} = 0$ otherwise. Then there is a finite-rank, free $\cA(1)$-module $F_k$ and an $\cA(1)$-module isomorphism
\begin{equation}
\label{its_milgram}
 H_\ko^*(\Ad_k^\ko(\ko)) \cong Q_{4-k, 2\ell + \delta_{k,4}+1} \oplus F_k.
\end{equation}
\end{cor}
See also~\cite[Theorem 2.4]{Mil75} and~\cite[Lemma 3.13]{DGM81}, which are closely related.
\begin{proof}
Since $\Ad_k^\ko(\ko)$ is a smash product over $\ko$ of $k$ copies of $\Sigma^{-1}\overline H$, the $\ko$-module Künneth theorem (\cref{rel_kunneth}) allows us to induct. The base case is the fact that $H_\ko^*(\Sigma^{-1}\overline H)\cong Q_{3,1}$, which can be verified using the long exact sequence in $H_\ko^*$ associated to the fiber sequence $\Sigma^{-1}\overline H\to \ko\to H$. The inductive step is \cref{milgram_tensor}, together with the fact that if $F$ is a free $\cA(1)$-module and $M$ is any $\cA(1)$-module, $F\otimes M$ is free.
\end{proof}
For convenience, we let $\vec v(n)\coloneqq (4- g(n), 2\ell + \delta_{g(n), 4}+1)$, where $g(n) = 4\ell+m$, $\ell\in\Z$, and $1\le m\le 4$. Then $H_\ko^*(\Ad_{g(n)}^\ko(\ko))\cong Q_{\vec v(n)}\oplus F_{g(n)}$.

By directly comparing \cref{A1_module_of_connective_covers,milgram_defn}, we observe:
\begin{lem}
\label{cyclic_is_cover}
The cyclic quotient $C_{\vec v(n)}$ of $Q_{\vec v(n)}$ is \emph{graded} isomorphic to $H_\ko^*(\tau_{\ge n}\ko)$.
\end{lem}
This is not a coincidence.

\begin{lem}\label{zero_or_surj}
For any choice of the isomorphism $H_\ko^*(\Ad_{g(n)}^\ko(\ko))\cong Q_{\vec v(n)}\oplus F_{g(n)}$ from \cref{which_milgram} and any grading-preserving $\cA(1)$-module map $\gamma^*\colon H_\ko^*(\Ad_{g(n)}^\ko(\ko))\to H_\ko^*(\tau_{\ge n}\ko)$, the restriction of $\gamma^*$ to $Q_{\vec v(n)}$ is either zero or surjective.
\end{lem}
\begin{proof}
Let $\psi(n)\coloneqq 2\ell+\delta_{g(n), 4}+1$, where $\ell = \lfloor g(n)/4\rfloor$. As an $\cA(1)$-module, $\Upsilon_{\psi(n)}$ is generated by classes in degrees $4\psi(n)-4$ and below. By \cref{A1_module_of_connective_covers}, $H_\ko^*(\tau_{\ge n}\ko)$ is $(4\psi(n)-4)$-connected, so the restriction of $\gamma$ to the seagull submodule $\Upsilon_{\psi(n)}\subset Q_{\vec v(n)}$  vanishes and $\gamma|_{Q_{\vec v(n)}}$ descends to the cyclic quotient $C_{4-{g(n)}, 2\ell + \delta_{g(n), 4}+1}$. By \cref{cyclic_is_cover}, $C_{\vec v(n)}\cong H_\ko^*(\tau_{\ge n}\ko)$, so by inspection the only two options for a grading-preserving $\cA(1)$-module map out of $C_{\vec v(n)}$ are $0$ or an isomorphism. Pulling back to $Q_{\vec v(n)}$, these options become $0$ or surjectivity. 
\end{proof}
By Margolis' theorem for $\ko$-modules~\cite[Theorem 5.25]{BottSpiral} (a sharpening of the original Margolis theorem~\cite{Mar74}), any $\cA(1)$-module splitting of the form~\eqref{its_milgram} lifts to a $\ko$-module splitting
\begin{equation}
    \Ad_{g(n)}^\ko(\ko)\xrightarrow{p} \overline{\Ad}{}_{g(n)}^\ko(\ko)\vee M_n,
\end{equation}
where $H_\ko^*(\overline{\Ad}{}_{g(n)}^\ko(\ko))\cong Q_{\vec v(n)}$, $M_n$ is $\ko$-module equivalent to a sum of shifts of $H\Z/2$, and $H_\ko^*(M_n)\cong F_{g(n)}$, such that $H_\ko^*(p)$ is the splitting of $H_\ko^*(\Ad_{g(n)}^\ko(\ko))$ chosen. Thus, in what follows, choose such a splitting, and let $p_i$ denote the projection from $\Ad_{g(n)}^\ko(\ko)$ onto the $i^{\mathrm{th}}$ factor of the direct summand.
\begin{lem}
\label{no_Z2_factor}
Let $c\in\pi_n(\ko)$ be nonzero. Then $c$, interpreted as a $\ko$-module map $c\colon\Sigma^{n}\ko\to\ko$, does not factor as a $\ko$-module map through any $\ko$-module $M$ which is $\ko$-module equivalent to a sum of shifts of $H\Z/2$.
\end{lem}
\begin{proof}
Since $c$ is a $\ko$-module map, it commutes with the Bott class $\beta$, so in particular $\beta c\colon\Sigma^{n+8}\ko\to \Sigma^8\ko$ is nonzero on $\pi_*$ (as, tautologically, all nonzero elements of $\pi_n\ko$ induce nonzero maps $\pi_{*+n}\ko\to\pi_*(\ko)$). However, $\beta$ is the zero map on $\pi_*(H\Z/2)$, hence also on $\pi_*(M)$, so if $c$ factors through $M$, $\beta c$ is the zero map on $\pi_*$. And a map $\Sigma^{n}\ko\to\ko$ is nonzero only if it is nonzero on homotopy groups.
\end{proof}
Because $\Sigma^{n}\ko$ is $(n-1)$-connected, $c$ factors as a $\ko$-module map as a composition
\begin{equation}
    \widetilde c\colon \Sigma^{n}\ko \longrightarrow \tau_{\ge n}\ko
\end{equation}
followed by the Postnikov cover $\tau_{\ge n}\colon\tau_{\ge n}\ko\to\ko$. This factorization is unique up to homotopy of $\ko$-module maps. And as we showed in \cref{map_of_filtered_spectra}, $\tau_{\ge n}$ factors through the Adams cover $\alpha_{g(n)}\colon \Ad_{g(n)}^\ko(\ko)\to\ko$.
\begin{cor}
\label{no_Z2_factor_part_2}
The composition
\begin{equation}\begin{tikzcd}
	{\Sigma^{n}\ko} & {\tau_{\ge n}\ko} & {\Ad_{g(n)}^\ko(\ko)} & {\overline{\Ad}{}_{g(n)}^\ko(\ko)\vee M_n} & \ko,
	\arrow["{\widetilde c}", from=1-1, to=1-2]
	\arrow["{\gamma_n}", from=1-2, to=1-3]
	\arrow["p", from=1-3, to=1-4]
	\arrow["{\alpha_{g(n)}}", from=1-4, to=1-5]
\end{tikzcd}\end{equation}
which is $\ko$-module equivalent to the map $c$ as discussed above, factors though the $\overline\Ad{}_{g(n)}^\ko(\ko)$ summand of $\Ad_{g(n)}^\ko(\ko)\simeq \overline{\Ad}{}_{g(n)}^\ko(\ko)\vee M_n$.
\end{cor}
\begin{proof}
For $i=1,2$, let $f_i\coloneqq p_i\circ \gamma_n\circ \widetilde c$, so that $p\circ \gamma_n\circ\widetilde c \simeq f_1\vee f_2$. We wish to prove that $c\simeq \alpha_n|_{\overline{\Ad}{}_{g(n)}^\ko(\ko)}\circ f_1$, so assume this is \emph{not} true. Then $c - \alpha_n|_{\overline{\Ad}{}_{g(n)}^\ko(\ko)}\circ f_1$ is a nontrivial $\ko$-module map $\Sigma^{n}\ko\to\ko$ that factors through a sum of shifts of $H\Z/2$, which by \cref{no_Z2_factor} cannot happen.
\end{proof}

\begin{lem}
\label{ko_classes_filt}
For all $k\ge 0$, the classes $\textcolor{BrickRed}{\beta^k\in\pi_{8k}\ko}$, $\textcolor{Green}{\beta^k\eta\in\pi_{8k+1}\ko}$, $\textcolor{MidnightBlue}{\beta^k\eta^2\in\pi_{8k+2}\ko}$, and $\textcolor{Fuchsia}{\beta^k a\in\pi_{8k+4}\ko}$ have Adams filtrations equal to $\textcolor{BrickRed}{4k}$, $\textcolor{Green}{4k+1}$, $\textcolor{MidnightBlue}{4k+2}$, and $\textcolor{Fuchsia}{4k+3}$ respectively. In particular, of these classes, the class in degree $n$ has Adams filtration $g(n)$.
\end{lem}
\begin{proof}[Proof sketch]
This is standard: it follows from the fact that they are detected by the classes $\textcolor{BrickRed}{w^k}$, $\textcolor{Green}{w^kh_1}$, $\textcolor{MidnightBlue}{w^kh_1^2}$, and $\textcolor{Fuchsia}{w^kv}$ on the Adams $E_\infty$-page, which have the corresponding $s$-degrees.
\end{proof}
This means, for example, it is possible to write $\textcolor{Green}{\beta\eta}$ as a composition of five $\ko$-module maps which are all $0$ on $H_\ko^*$, but it is not possible to write it as a composition of six such $\ko$-module maps.
\begin{prop}
\label{gotta_be_nonzero}
$(p_1\circ \gamma_n)^*\colon H_\ko^*(\overline{\Ad}{}_{g(n)}^\ko(\ko))\to H_\ko^*(\tau_{\ge n}\ko)$ is nonzero.
\end{prop}
\begin{proof}
Suppose instead that $(p_1\circ \gamma_n)^* = 0$. By construction, $\alpha_{g(n)}$ is a composition of $g(n)$ $\ko$-module maps $h_{g(n)}\circ\dotsb\circ h_1$, each of which induces the zero map on $H_\ko^*$; thus the same is true if we restrict $h_1$ to $\overline\Ad{}_{g(n)}^\ko(\ko)$.

Let $c$ be a generator of $\pi_{n}(\ko)$; thus, if $n = \textcolor{BrickRed}{8k}$, $\textcolor{Green}{8k+1}$, $\textcolor{MidnightBlue}{8k+2}$, or $\textcolor{Fuchsia}{8k+4}$, respectively, $c$ is $\textcolor{BrickRed}{\pm\beta^k}$, $\textcolor{Green}{\beta^k\eta}$, $\textcolor{MidnightBlue}{\beta^k\eta^2}$, or $\textcolor{Fuchsia}{\pm\beta^k a}$, and by \cref{ko_classes_filt}, $c$ has Adams filtration $g(n)$. (Since $\pm1$ are indistinguishable in any $\Z/2$-vector space, for any map $f$, $f$ and $-f$ have equal Adams filtrations.) By \cref{no_Z2_factor_part_2}, $c$ factors as the following composition of $\ko$-module maps:
\begin{equation}
    c\simeq h_{g(n)}\circ\dotsb\circ h_1|_{\overline\Ad{}_{g(n)}^\ko(\ko)} \circ (p_1\circ \gamma_n\circ\widetilde c),
\end{equation}
and by hypothesis, each of these $g(n)+1$ maps is the zero map on $H_\ko^*$. This implies that $c$ has Adams filtration at least $g(n)+1$, which contradicts \cref{ko_classes_filt}. 
\end{proof}
\begin{proof}[Proof of \cref{prop_surj_on_Hko}]
By \cref{zero_or_surj,gotta_be_nonzero}, $\gamma_n^*$ is surjective when restricted to a summand of $\Ad_{g(n)}^\ko(\ko)$. Thus, the entire map is also surjective.
\end{proof}
Thus, we have finished the proof of \cref{cQthm}.
\begin{exm}
\label{HASS_exm}
We use the HASS to compute $\tau_{\le 2}\ko_*(\pt)$ and compare with the usual Baker--Lazarev Adams spectral sequence computation, which we did in \cref{non_modified_Adams}. The difference is that the map $f_2$ from \cref{cQthm} replaces the $d_3$ differentials -- in other words, we have isolated the information those differentials contained from the rest of the spectral sequence computation.

We will use \cref{cQthm} with $n = 2$. Throughout this example, and the next section, $\cQ^{s,t}\coloneqq\cQ_2^{s,t}$. We need to calculate $\cQ(H_\ko^*(\ko\wedge_\ko\pt_+))$, i.e.\ $\cQ(\Z/2)$.
\begin{lem}
\label{sE_defn}\hfill
\begin{enumerate}
\item
The $\E$-module map
\begin{equation}
    f_2\colon \Ext_{\cA(1)}^{s,t}(\uQ) \longrightarrow \Ext_{\cA(1)}^{s+3, t+7}(\Z/2).
\end{equation}
sends $\delta\mapsto v$ and $\epsilon\mapsto w$, where $\delta,\epsilon$ are as in \cref{ext_ups_Q}.
\item
Let $s\E\coloneqq\cQ(\Z/2)$.
$s\E\cong\Z/2[h_0,h_1]/(h_0h_1, h_1^3)$, with $h_0\in\cQ^{1,1}$ and $h_1\in\cQ^{1,2}$.
\end{enumerate}
\end{lem}
\begin{proof}
Recall the $\Z^2$-graded $\Z/2$-algebra $\E\coloneqq\Ext_{\cA(1)}^{*,*}(\Z/2)$ from \cref{liuthm}. By part~\eqref{cQ_defn} of \cref{cQthm}, $\cQ(\Z/2)$ fits into a long exact sequence~\eqref{cQLES} with $f_2$.
Here we used that $g(2) = 2$ (\cref{def_of_g}) and $C_3\cong \Sigma^4\uQ$ (\cref{A1_module_of_connective_covers}). 

Recall from \cref{ext_ups_Q} that, as an $\E$-module, $\Ext(\uQ)$ is generated by two classes $\delta\in\Ext^{0,0}$ and $\epsilon\in\Ext^{1,5}$, so $f_2$ is uniquely determined by its values on $\delta$ and $\epsilon$. In particular, $f_2(\delta)\in\Ext^{3,7}(\Z/2)$ and $f_2(\epsilon)\in\Ext^{4,12}(\Z/2)$; looking at the structure of $\Ext(\Z/2)$ given in \cref{liuthm}, we see that $f_2(\delta)$ is either $0$ or $v$, and $f_2(\epsilon)$ is either $0$ or $w$. The relations $v\delta = h_0^2\epsilon\ne 0$ and $v^2 = h_0^2w\ne 0$ in $\Ext(\uQ)$ (\cref{ext_ups_Q}) and $\Ext(\Z/2)$ respectively imply that $f_2(\delta) = v$ iff $f_2(\epsilon) = w$.

It is traditional to draw long exact sequences of bigraded Ext groups as Adams charts (so, graded by $(t-s,s)$ such that the boundary map has bidegree $(-1, 1)$ in the horizontal, resp.\ vertical directions. See~\cite[Figures 24, 26, 33, 36, and 41]{BC18} for some examples. In our case, we would like to think of $f_2$ as the boundary map, so that the ``middle term'' of the long exact sequence~\eqref{cQLES} is $\cQ(\Z/2)$. To do this, we should shift $\Ext(\cQ)$ five units to the right and two upwards. We draw $f_2$ in this way in \cref{modified_fig}, left. 

The reader can check that if $f_2(\delta) = v$, so that $f_2(\epsilon) = w$, then exactness of~\eqref{cQLES} calculates $\cQ(\Z/2)$ in agreement with the lemma statement. Thus we are done if we can derive a contradiction from the assumption $f_2(\delta) = 0$. Assuming $f_2(\delta) = 0$, then $f_2(\epsilon) = 0$ and therefore the entire map $f_2$ vanishes, implying an $\E$-module isomorphism
\begin{equation}
\label{fake_zero}
    \cQ(\Z/2)\cong \E\oplus \Sigma^{2,7}\Ext(\uQ).
\end{equation}
The bigraded suspension uses the bigrading $(s,t)$ (so that, like above, we shifted $\Ext(\uQ)$ five units to the right and two upwards in the $(t-s,s)$ bigrading). \eqref{fake_zero} and our knowledge of $\E$ and $\Ext(\uQ)$ from \cref{liuthm,ext_ups_Q}  respectively would then imply that $\cQ(\Z/2)$ contains $h_0$-towers in topological degrees $4$ and $5$, generated by $v$ and $\delta$ in filtrations $3$ and $2$ respectively, and there are no other nonzero classes in topological degrees $3$, $4$, and $5$.

Now run the HASS. By part~\eqref{HASS_defn} of \cref{cQthm}, the $E_2$-page is $\cQ(\Z/2)$, and this spectral sequence converges to $\pi_*(\tau_{\le 2}\ko)_2^\wedge$. Since $\pi_4(\tau_{\le 2}\ko) \cong 0$, all classes in the $h_0$-tower generated by $v$ must be the sources or targets of differentials. For $r\ge 2$, $E_2^{3-r, 8-r} = \cQ^{3-r,8-r}(\Z/2)\cong 0$, so $d_r\colon E_r^{3-r,8-r}\to E_r^{3,7}$ is the zero map and cannot kill $v$. Likewise, $d_r(v)\in E_r^{3+r, 6+r}$, but $E_2^{3+r,6+r} = \cQ^{3+r,6+r}(\Z/2)\cong 0$, so $v$ is a permanent cycle, giving the desired contradiction.
\end{proof}

\begin{figure}[htbp]
\begin{subfigure}[c]{0.45\textwidth}
\includegraphics[width=\linewidth]{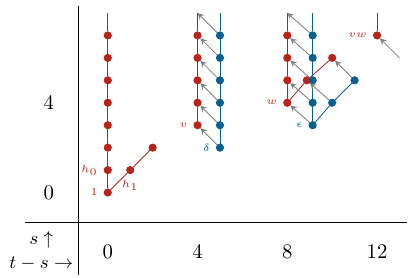}
\end{subfigure}%
\hfill
\begin{subfigure}[c]{0.45\textwidth}
\includegraphics[width=\linewidth]{figures/non_modified_E4.pdf}
\end{subfigure}
\caption{Left: the map $f_2$ from~\eqref{cQthm}, part~\eqref{g4_defn}. We calculate this map, in particular showing $f_2(\delta) = v$ and $f_2(\epsilon) = w$, as part of the proof of \cref{sE_defn}. Right: $\sE\coloneqq\cQ(\Z/2)$, calculated in \cref{sE_defn}. This is the $E_2 = E_\infty$-page of the HASS for $\tau_{\le 2}\ko_*(\pt)$; we run this spectral sequence in \cref{HASS_exm}. Compare these two charts to the ones in \cref{non_modified_fig} to compare the HASS (here) with the ordinary Baker--Lazarev Adams spectral sequence (there).}
\label{modified_fig}
\end{figure}

We draw $\sE$ in \cref{modified_fig}, right. By part~\eqref{HASS_defn} of \cref{cQthm}, this is the $E_2$-page of the HASS for $\tau_{\le 2}\ko_*(\pt)$. This spectral sequence immediately collapses for degree reasons, just as the usual Adams spectral sequence for $\ko_*(\pt)$ does, and there are no hidden extensions.
\end{exm}

\section{User's guide}
	\label{s:guide}

In this section, we outline the general procedure for using the HASS to compute twisted supercohomology groups and give several examples. To start with, we need to define supercohomology.
\begin{defn}[Brown--Comenentz~\cite{BC76}]
\label{BCdefn}
Let $X$ be a spectrum. The \term{Pontrjagin dual of the sphere}, denoted $I_{\C^\times}$, is the spectrum representing the generalized cohomology theory
\begin{equation}
    X\longmapsto \Hom(\pi_{-*}(X), \C^\times).
\end{equation}
The \term{Pontrjagin dual} of $X$ is $I_{\C^\times}X\coloneqq \mathrm{Map}(X, I_{\C^\times})$.
\end{defn}
Brown--Comenetz used $\R/\Z$ instead of $\C^\times$; the resulting theory is mostly the same.
\begin{defn}[Wang--Gu~\cite{WG18, WG20}, Kapistin--Thorngren~\cite{KT17}]
\label{SHdefn}
\term{Supercohomology} is the spectrum $\SH\coloneqq I_{\C^\times}(\tau_{\le 2}\ko)$.
\end{defn}
Thus $\SH$ defines a generalized cohomology theory on spaces $X$ in the usual way. We will also need to consider twisted supercohomology, with the twist given by $a\in H^1(X;\Z/2)$ and $b\in H^2(X;\Z/2)$. We define twisted supercohomology in Appendix~\ref{s:supercoh}. For applying the HASS, what one needs to know is the following lemma.
\begin{lem}
\label{SH_as_hom}
For any spectrum $X$, there is a natural isomorphism
\begin{equation}
    \pi_{-n}(I_{\C^\times}X) \overset\cong\longrightarrow
    \Hom(\pi_n(X), \C^\times).
\end{equation}
Thus, in particular, $\SH^n(X, a, b)$ is the Pontrjagin dual group to $(\tau_{\le 2}\ko)_n(X, a, b)$.
\end{lem}
This is a standard lemma that follows directly from \cref{BCdefn}; see~\cite[\S 1]{BC76}. Because of \cref{SH_as_hom}, we will use the HASS to compute twisted $\tau_{\le 2}\ko$-homology, then trust the reader to apply Pontrjagin duality to obtain the corresponding twisted supercohomology groups.

\subsection{Summary of the procedure}\label{subsection:procedure}

The purpose of this section is to show how to use the HASS that we constructed in the previous section.
The procedure follows the calculation of the Adams spectral sequence over $\ko$ closely.

The computation naturally separates into two pieces: at the prime $2$ and away from $2$.\footnote{See for example \cite[\S 10.2]{DDHM24}.} The latter is easier, so we do it first.
\begin{lem}
\label{away_from_2}
For any space $X$, $a\in H^1(X;\Z/2)$, and $b\in H^2(X;\Z/2)$, the map $H^*(X;\C^\times_a)\to \SH^*(X, a, b)$ induces an isomorphism on odd-primary torsion subgroups.
\end{lem}
Here $\C^\times_a$ is the abelian group $\C^\times$ with $\pi_1(X)$-action given by the cohomology class $a$, interpreted as a homomorphism $a\colon \pi_1(X)\to\{\pm 1\}$.
\begin{proof}
Applying \cref{SH_as_hom}, it suffices to show that $\tau_{\le 0}\colon \tau_{\le 2}\ko\to H\Z$ induces an isomorphism on odd-primary subgroups for $(X, a, b)$-twisted homology. That is, we will show that $\tau_{\le 0}$ is an equivalence away from the prime $2$. This is because its fiber is the spectrum $\tau_{\le 2}\tau_{\ge 1}\ko$, whose homotopy groups are $0$ away from degrees $1$ and $2$, where they are isomorphic to $\Z/2$. Thus after inverting $2$, these homotopy groups are zeroed out, so away from the prime $2$, $\tau_{\le 0}$ has contractible fiber and thus is an equivalence of spectra, finishing the proof. 
\end{proof}

At the prime 2, we will use the HASS to calculate $(\tau_{\leq 2}\ko)_n (X, a, b)_{(2)}$ or, dually, $\SH^n(X, a, b)_{(2)}$ for a space $X$. As with the usual Adams spectral sequence, the starting data is a module over the Steenrod algebra $\cA(1)$.

\begin{defn}\label{defn:twistedcoh}
Let $X$ be a space, $a\in H^1(X;\Z/2)$, and $b\in H^2(X;\Z/2)$. Let $H_{a,b}^*(X;\Z/2)$ be the $\cA(1)$-module which, as a $\Z$-graded $\Z/2$-vector space, is isomorphic to a free $H^*(X;\Z/2)$-module on a single degree-$0$ generator $U$, but with $\cA(1)$-action given by
\begin{equation}
\label{twisted_Steenrod}
\begin{aligned}
    \Sq^1(Ux) &= U(ax + \Sq^1(x))\\
    \Sq^2(Ux) &= U(bx + \Sq^1(a)\Sq^1(x) + \Sq^2(x)).
\end{aligned}
\end{equation}
\end{defn}
It is not a priori obvious that the action we defined in~\eqref{twisted_Steenrod} satisfies the Adem relations and therefore actually defines an $\cA(1)$-action, but this is true and is proven in~\cite[Lemma 2.38(3)]{DY23}. If there is a vector bundle $V\to X$ such that $a = w_1(V)$ and $b = w_2(V)$, this follows from the Wu formula and the Thom isomorphism, but such a $V$ does not in general exist.
\begin{lem}[{\cite[Theorem 2.39(3)]{DY23}}]
\label{fake_Wu_formula}
With $(X, a, b)$ in \cref{defn:twistedcoh}, let $M^\ko h_{a,b}$ be the $\ko$-module Thom spectrum of the map $h_{a,b}\colon X\to B\GL_1(\ko)$ defined in \cref{stop_at_ko}. Then there is a natural $\cA(1)$-module isomorphism
\begin{equation}
    H_\ko^*(M^\ko h_{a,b}) \xrightarrow\cong H_{a,b}^*(X;\Z/2).
\end{equation}
\end{lem}
By \cref{SH_is_SH}, $\SH^*(X, a, b) = \pi_{-*}(I_{\C^\times} (\tau_{\le 2}\ko\wedge_\ko M^\ko h_{a,b}))$. Thus, by \cref{cQthm}, part~\eqref{HASS_defn}, the $E_2$-page of the HASS computing $\SH^*(X, a, b)$ is $\cQ(H_\ko^*(M^\ko h_{a,b}))$, and therefore \cref{fake_Wu_formula} calculates the $\cA(1)$-module we want to apply $\cQ$ to.

The steps of the HASS are as follows.
\begin{enumerate}
    \item Using \cref{fake_Wu_formula}, compute the $\cA(1)$-module structure on $H_\ko^*(M^\ko h_{a,b})$. In practice, this is usually a direct sum of suspensions of modules appearing in the list in \cref{names}.
    \item Draw the $E_2$-page of the HASS: for each $\cA(1)$-module summand $M$, determine $\cQ(M)$. Then direct sum all of these $\cQ(M)$ charts together. 
%
These $\cQ$ charts play the same role for the HASS that Ext charts play for the ordinary Adams spectral sequence. 

\item\label{the_hard_part} Compute the differentials and resolve any hidden extensions.
\item Using \cref{SH_as_hom}, $\SH^k(X, a, b) \cong \Hom(\tau_{\le 2}\ko_k(X, a, b), \C^\times)$.
\end{enumerate}
Part~\eqref{the_hard_part} may be trivial or difficult, depending on the specific problem the HASS has been applied to. In practice, it tends to be about as difficult as the corresponding Adams spectral sequence for twisted $\ko$-homology, which is usually tractable in examples.

%
%


\subsection{Ext charts}\label{subsection:ExtCharts}
In this section, we describe the Ext charts $\cQ(M)$ for frequently encountered $\cA(1)$-modules $M$. As described
in the previous subsection, these charts give the input to the hastened Adams spectral sequence computing
$\tau_{\le 2}\ko$-homology.

This section contains the following information:
\begin{enumerate}
	\item A theorem describing each $\cQ(M)$ as a module over $\sE\coloneqq \cQ(\Z/2)\cong\Z/2[h_0,h_1]/(h_0h_1,
	h_1^3)$ (see \cref{sE_defn}).
	\item A figure displaying $M$ and $\cQ(M)$.
\end{enumerate}
We defer proofs to Appendix~\ref{ext_Q_computations}.

We begin by defining and naming the $\cA(1)$-modules we will consider. Names will mostly follow~\cite{BC18,Deb21}. 
\begin{defn}[Names of common $\cA(1)$-modules]\label{names}
\hfill
\begin{enumerate}
	\item $B\coloneqq \cA(1)\otimes_{\cE(1)}\Z/2$ (called the \term{bow} in~\cite{BGR24}).
	\item $N_1\coloneqq \cA(1)/(\Sq^2, \Sq^2\Sq^1)$.
	\item $\uQ\coloneqq\cA(1)/(\Sq^1, \Sq^2\Sq^3)$ as in \cref{A1_module_of_connective_covers} (often called the \term{upside-down
	question mark}).
	\item $J\coloneqq\cA/(\Sq^3)$ as in \cref{A1_module_of_connective_covers} (called the \term{Joker} by Adams).
	\item $Q\coloneqq\cA/(\Sq^2)$ as in \cref{A1_module_of_connective_covers} (often called the \term{question mark}).
	\item $R_2\coloneqq Q_{3,1}$ (called the \term{elephant} in~\cite{BM23}).
	\item $\Upsilon_1\coloneqq \cA(1)\otimes_{\cA(0)} \Z/2$ as in \cref{seagull_defn} (called the \term{seagull} in~\cite{Ada23}).
	\item $P\coloneqq \Sigma^{-1}\widetilde H^*(B\Z/2;\Z/2)$.
	\item $R_0\coloneqq H^*((B\Z/2)^{1-\sigma};\Z/2)$, where $\sigma\to B\Z/2$ denotes the tautological line
	bundle.
	\item $R_1\coloneqq H^*((B\Z/2)^{2\sigma-2};\Z/2)$.
	\item $R_3$ is the unique non-free $\cA(1)$-module summand in $P\otimes P$.
	\item $R_5$ is the unique non-free $\cA(1)$-module summand in $P\otimes P\otimes P$.
	\item $R_6\coloneqq P\otimes C\eta$.
\end{enumerate}
\end{defn}
\begin{rem}
\label{coalg_remark}
If $R$ is an $E_\infty$-$\ko$-algebra, there is an induced $\Z/2$-coalgebra structure on $H_\ko^*(R)$. This applies to $B\cong H_\ko^*(\ku)$~\cite[(5.12)]{Bak20}, $\Upsilon_1\cong H_\ko^*(H\Z)$~\cite[\S 5]{Bak20}, and $R_1\cong H_\ko^*(\ko\wedge (B\Z/2)^{2\sigma-2})$. Thus in these cases one obtains $\E$-algebra structures on $\Ext_{\cA(1)}$ and $s\E$-algebra structures on $\cQ$ of these $\cA(1)$-modules.
\end{rem}

\begin{thm}
\label{heuristic_Q_calc}
Let $M$ be one of the $\cA(1)$-modules named in \cref{names}. Then $\cQ(M)$ is described in \cref{Q_ext_charts_1,Q_ext_charts2}.
\end{thm}
We give more detail in \cref{Q_ext_calc}, including explicit presentations of $\cQ(M)$ for these $M$ and a description of the image of $\partial_2\colon\Ext_{\cA(1)}(M)\to\cQ(M)$. In particular, \cref{heuristic_Q_calc} follows from the proof of \cref{Q_ext_calc} that we give in Appendix~\ref{ext_Q_computations}.

\begin{figure}[htbp]
\begin{subfigure}[c]{0.45\textwidth}
\begin{subfigure}[c]{0.22\linewidth}
\includegraphics{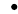}
\end{subfigure}
\begin{subfigure}[c]{0.73\linewidth}
\includegraphics{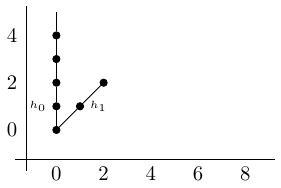}
\end{subfigure}
\caption{$\cQ(\Z/2)$}
\label{Q_ext_F2}
\end{subfigure}
\begin{subfigure}[c]{0.45\textwidth}
\begin{subfigure}[c]{0.22\linewidth}
\includegraphics{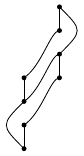}
\end{subfigure}
\begin{subfigure}[c]{0.73\linewidth}
\includegraphics{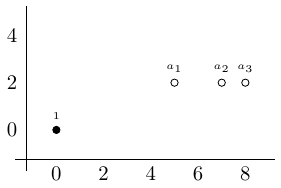}
\end{subfigure}

\caption{$\cQ(\cA(1))$}
\label{Qext_A1}
\end{subfigure}
\begin{subfigure}[c]{0.45\textwidth}
\begin{subfigure}[c]{0.22\linewidth}
\includegraphics{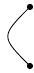}
\end{subfigure}
\begin{subfigure}[c]{0.73\linewidth}
\includegraphics{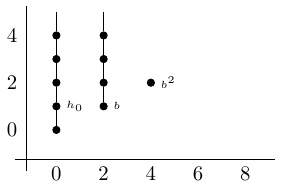}
\end{subfigure}

\caption{$\cQ(B)$}
\label{Q_ext_B}
\end{subfigure}%
\begin{subfigure}[c]{0.45\textwidth}
\begin{subfigure}[c]{0.22\linewidth}
\includegraphics{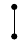}
\end{subfigure}
\begin{subfigure}[c]{0.73\linewidth}
\includegraphics{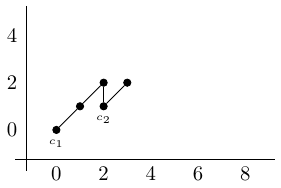}
\end{subfigure}
\caption{$\cQ(N_1)$}
\label{Q_N1_fig}
\end{subfigure}
\begin{subfigure}[c]{0.45\textwidth}
\begin{subfigure}[c]{0.22\linewidth}
\includegraphics{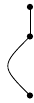}
\end{subfigure}
\begin{subfigure}[c]{0.73\linewidth}
\includegraphics{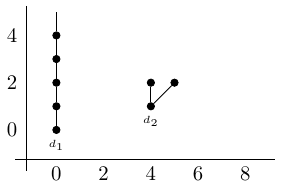}
\end{subfigure}
\caption{$\cQ(\uQ)$}
\label{Q_ext_uQ}
\end{subfigure}%
\begin{subfigure}[c]{0.45\textwidth}
\begin{subfigure}[c]{0.22\linewidth}
\includegraphics{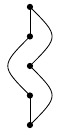}
\end{subfigure}
\begin{subfigure}[c]{0.73\linewidth}
\includegraphics{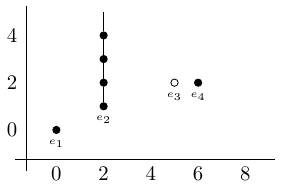}
\end{subfigure}

\caption{$\cQ(J)$}
\label{Q_ext_J}
\end{subfigure}
\begin{subfigure}[c]{0.45\textwidth}
\begin{subfigure}[c]{0.22\linewidth}
\includegraphics{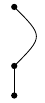}
\end{subfigure}
\begin{subfigure}[c]{0.73\linewidth}
\includegraphics{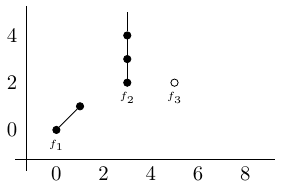}
\end{subfigure}
\caption{$\cQ(Q)$}
\label{Q_ext_Q}
\end{subfigure}
\begin{subfigure}[c]{0.45\textwidth}
\begin{subfigure}[c]{0.22\linewidth}
\includegraphics{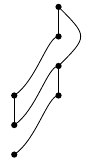}
\end{subfigure}
\begin{subfigure}[c]{0.73\linewidth}
\includegraphics{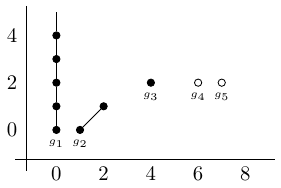}
\end{subfigure}

\caption{$\cQ(R_2)$}
\end{subfigure}
\caption{Ext charts $\cQ(M)$ for various $\cA(1)$-modules $M$, with $M$ displayed on the left for reference. The filled-in classes are the image of the map $\partial_2\colon \Ext(M)\to\cQ(M)$. These
charts are calculated in \cref{Q_ext_calc}.}
\label{Q_ext_charts_1}
\end{figure}
\begin{figure}[htbp]
\ContinuedFloat
\captionsetup[subfigure]{skip=2pt}

\begin{subfigure}[c]{0.45\textwidth}
\begin{subfigure}[c]{0.22\linewidth}
\includegraphics{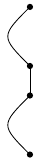}
\end{subfigure}
\begin{subfigure}[c]{0.73\linewidth}
\includegraphics{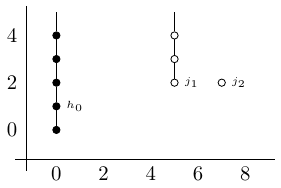}
\end{subfigure}

\caption{$\cQ(\Upsilon_1)$}
\label{Qext_M0}
\end{subfigure}
\begin{subfigure}[c]{0.45\textwidth}
\begin{subfigure}[c]{0.22\linewidth}
\includegraphics{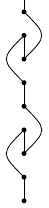}
\end{subfigure}
\begin{subfigure}[c]{0.73\linewidth}
%
\includegraphics{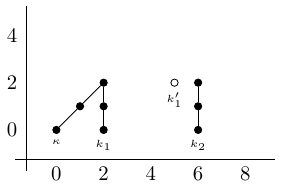}
\end{subfigure}

\caption{$\cQ(P)$}
\label{Q_ext_P}
\end{subfigure}
\begin{subfigure}[c]{0.45\textwidth}
\begin{subfigure}[c]{0.22\linewidth}
\includegraphics{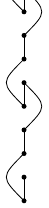}
\end{subfigure}
\begin{subfigure}[c]{0.73\linewidth}
%
%
\includegraphics{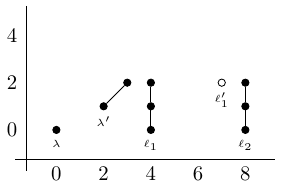}
\end{subfigure}

\caption{$\cQ(R_0)$}
\end{subfigure}
\begin{subfigure}[c]{0.45\textwidth}
\begin{subfigure}[c]{0.22\linewidth}
\includegraphics{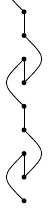}
\end{subfigure}
\begin{subfigure}[c]{0.73\linewidth}
%
\includegraphics{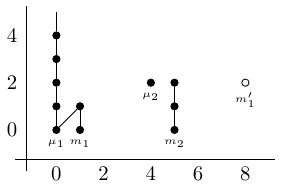}
\end{subfigure}

\caption{$\cQ(R_1)$}
\label{Q_ext_R1}
\end{subfigure}
\begin{subfigure}[c]{0.45\textwidth}
\begin{subfigure}[c]{0.22\linewidth}
\includegraphics{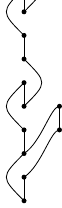}
\end{subfigure}
\begin{subfigure}[c]{0.73\linewidth}
\includegraphics{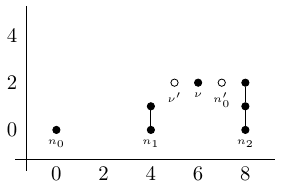}
\end{subfigure}

\caption{$\cQ(R_3)$}
\end{subfigure}
\begin{subfigure}[c]{0.45\textwidth}
\begin{subfigure}[c]{0.22\linewidth}
\includegraphics{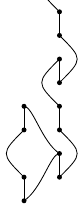}
\end{subfigure}
\begin{subfigure}[c]{0.73\linewidth}
\includegraphics{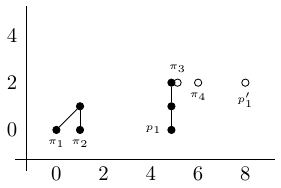}
\end{subfigure}
\caption{$\cQ(R_5)$}
\end{subfigure}

\begin{subfigure}[c]{0.45\textwidth}
\begin{subfigure}[c]{0.22\linewidth}
\includegraphics{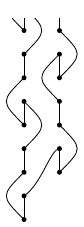}
\end{subfigure}
\begin{subfigure}[c]{0.73\linewidth}
\includegraphics{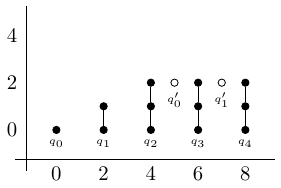}
\end{subfigure}
\caption{$\cQ(R_6)$}
\end{subfigure}

\caption{Ext charts $\cQ(M)$ for various $\cA(1)$-modules $M$, with $M$ displayed on the left for reference. The filled-in classes are the image of the map $\partial_2\colon \Ext(M)\to\cQ(M)$. These
charts are calculated in \cref{Q_ext_calc}.}
\label{Q_ext_charts2}
\end{figure}

\subsection{Examples}
\label{ss:users_examples}
We now use the $\cQ$ charts worked out above to give some example computations of twisted supercohomology.
\subsubsection{Cyclic groups}
\label{sss:cyclic}
Recall that if $n$ is even, $H^*(B\Z/n;\Z/2)\cong\Z/2[x, y]/(x^2 + ky)$, where $\abs x = 1$, $\abs y = 2$, and $k = 1$ if $n = 2$ and $k = 0$ otherwise.

\begin{figure}[htbp]
\begin{subfigure}[c]{0.16\textwidth}
\includegraphics{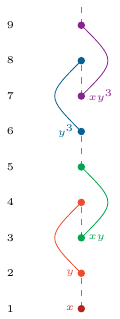}
\end{subfigure}%
\hfill
\begin{subfigure}[c]{0.41\textwidth}
\includegraphics[width=\linewidth]{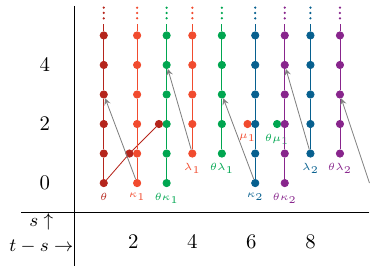}
\end{subfigure}%
\hfill
\begin{subfigure}[c]{0.41\textwidth}
\includegraphics[width=\linewidth]{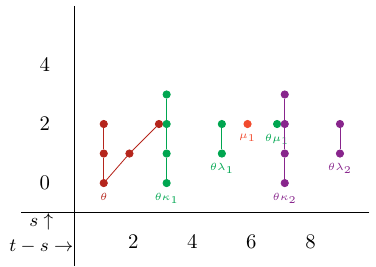}
\end{subfigure}%
\\[1em]
\begin{subfigure}[c]{0.16\textwidth}
\includegraphics{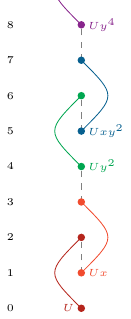}
\end{subfigure}%
\hfill
\begin{subfigure}[c]{0.41\textwidth}
\includegraphics[width=\linewidth]{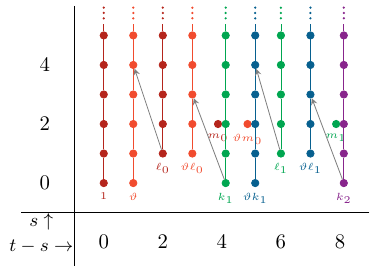}
\end{subfigure}%
\hfill
\begin{subfigure}[c]{0.41\textwidth}
\includegraphics[width=\linewidth]{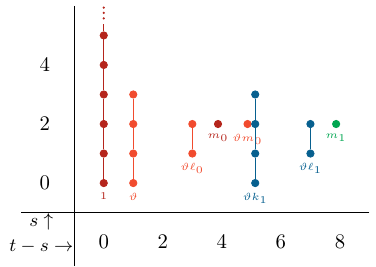}
\end{subfigure}%

\caption{Upper left: the $\cA(1)$-module structure on $\widetilde H^*(B\Z/n;\Z/2)$, when $n\equiv 0\bmod 4$. Let $2^k$ be the largest power of $2$ dividing $n$. The dashed lines represent nonzero $2^k$-Bockstein actions. Figure adapted from~\cite[Figure 7.5]{Cam17}. Upper center: the $E_2 = E_k$-page of the HASS computing $\tau_{\le 2}\ko_*(B\Z/n)$ (here $k = 3$). Notation follows~\cite[Theorem 13.36]{BDDHM25}. We discuss the pictured differentials in \cref{cycdiffs}. Upper right: the $E_{k+1} = E_\infty$-page.
Lower left: the $\cA(1)$-module structure on $H^*((B\Z/n)^{L-2};\Z/2)$, where $L$ is the complex line bundle induced from the rotation representation of $\Z/n$. The dashed lines again represent nonzero $2^k$-Bockstein operations. Lower center: the $E_2 = E_k$-page of the HASS computing $\tau_{\le 2}\ko_*(B\Z/n, 0, y)$, including the differentials we compute in the proof of \cref{cycdiffs} (here $k = 3$). Lower right: the $E_{k+1} = E_\infty$-page.}
\label{cyclic_adams}
\end{figure}

We compute $\SH^*(B\Z/n)$ and $\SH^*(B\Z/n, 0, y)$ for all even $n$.
\begin{prop}
\label{prop_cyclic}
Let $n$ be even. Then the untwisted and $(0,y)$-twisted supercohomology groups of
$B\Z/n$ are given in \cref{cyclic_supercoh}.
\end{prop}
In low degrees, these computations have been done before by the following authors, usually using the Atiyah--Hirzebruch spectral sequence.
\begin{itemize}
    \item $\SH^k(B\Z/2^n)$: Wang--Gu~\cite[Table II]{WG18} for $k\le 3$, Décoppet~\cite[Example 4.13]{Dec24} for $k = 4$, and~\cite[Propositions E.11 and E.16]{Debray:2025kfg} for $k = 5$.
    \item $\SH^k(B\Z/2^n, 0, y)$: Wang--Gu~\cite[Table III]{WG20} and Zhang--Wang--Yang--Qi--Gu~\cite{ZWYQG20} for $k\le 3$ and \cite[Propositions E.20 and E.23]{Debray:2025kfg} for $k = 4,5$.
\end{itemize}
\begin{proof}[Proof of \cref{prop_cyclic}]
As described in \S\ref{subsection:procedure}, we will compute the $\tau_{\le 2}\ko$-homology of $B\Z/n$, as well as of the $(0,y)$ twist thereof. The (twisted) supercohomology groups are then $\SH^k(X, a, b) \cong \Hom(\tau_{\le 2}\ko_k(X, a, b), \C^\times)$. Thus in the rest of this proof we focus only on (twisted) $\tau_{\le 2}\ko$-homology. \Cref{away_from_2} takes care of the calculation away from the prime $2$, so in the rest of the proof, we work at $2$, where we use the HASS.

In the untwisted case, split $\tau_{\le 2}\ko_*(X)$ into the known summand $\tau_{\le 2}\ko_*(\pt)$ and the reduced $\tau_{\le 2}\ko$-homology of $X$. We compute the latter using the reduced HASS, which differs only in taking the $\cA(1)$-module structure on reduced cohomology as input.

Let $k$ be the largest number such that $2^k$ divides $n$. The inclusion $\Z/2^k\to\Z/n$ induces a $2$-local stable equivalence $B\Z/2^k\to B\Z/n$, so we may without loss of generality assume $n = 2^k$.

As described in the user's guide, our first step in these computations is to determine the $\cA(1)$-modules $\widetilde H^*(B\Z/2^k;\Z/2)$ and $H_{0,y}^*(B\Z/2^k; \Z/2)$ for all $k\ge 1$. These can be looked up: there are $\cA(1)$-module isomorphisms
\begin{subequations}
\begin{align}
    \label{cyccoh1}
    \widetilde H^*(B\Z/2; \Z/2) &\cong \Sigma P,\\
    \label{cyccoh2}
    H_{0,y}^*(B\Z/2;\Z/2) &\cong R_1,
\intertext{and for all $k\ge 1$,}
\label{cyccoh3}
    \widetilde H^*(B\Z/2^k;\Z/2) &\cong \Sigma \Z/2\oplus\bigoplus_{i\ge 1}^\infty \Sigma^{4i+2} B \oplus \Sigma^{4i+3}B\\
    \label{cyccoh4}
    H_{0,y}^*(B\Z/2^k;\Z/2) &\cong \bigoplus_{i\ge 0}^\infty \Sigma^{4i} B \oplus \Sigma^{4i+1}B.
\end{align}
\end{subequations}
For~\eqref{cyccoh1} and~\eqref{cyccoh2}, this is essentially \cref{names}; see~\cite[\S 7.8]{Cam17} and~\cite[Figure 4]{BC18}. For~\eqref{cyccoh3} and~\eqref{cyccoh4}, see~\cite[Figure 7.5]{Cam17}, \cite[Figure 17]{DDHM24}, \cite[Lemma 3.8]{BGR24}, \cite[Proposition 13.20]{BDDHM25}, or~\cite[Corollary 2.10]{Dha26}. We draw $P$ and $R_1$ in \cref{Q_ext_P}, left, resp.\ \cref{Q_ext_R1}, left, and we draw the decompositions~\eqref{cyccoh3} and~\eqref{cyccoh4} in~\cref{cyclic_adams}, top left resp.\ bottom left.

Next, using the $\cQ^{s,t}(\bl)$ charts for $\Z/2$, $B$, $P$, and $R_1$ we computed in~\cref{Q_ext_F2,Q_ext_B,Q_ext_P,Q_ext_R1}, respectively, we obtain the $E_2$-page of the HASS for the $\tau_{\le 2}\ko$-homology, twisted by $(0,0)$ and $(0,y)$, of $B\Z/2^k$. For $B\Z/2$, untwisted, this is \cref{Q_ext_P} shifted right by $1$; for $B\Z/2$, twisted, the $E_2$-page is exactly \cref{Q_ext_R1}. 

For $k>1$ we give the $E_2$-pages for the untwisted and twisted computations in \cref{cyclic_adams}, top center resp.\ bottom center, as well as a brief generators-and-relations summary here.
\begin{itemize}
    \item First $(B\Z/2^k, 0, 0)$. By \cref{Q_ext_calc}, part~\eqref{is_surj}, the map $\partial_2\colon\Ext(\widetilde H^*(B\Z/2^k;\Z/2), \Z/2)\to\cQ(\widetilde H^*(B\Z/2^k;\Z/2))$ is surjective, so we can reuse notation from the presentation of $\Ext(\widetilde H^*(B\Z/2^k;\Z/2), \Z/2)$ given in~\cite[Theorem 13.36]{BDDHM25}.
    \begin{itemize}
        \item $\theta\in\cQ^{0,1}$, no relations.
        \item For each $i\ge 1$, classes $\kappa_i\in \cQ^{0,4i-2}$ and $\theta\kappa_i\in\cQ^{0,4i-1}$, with the relations $h_1\kappa_i = 0$ and $h_1\theta\kappa_i = 0$.\footnote{\label{noprod}There is an algebra structure on $\Ext_{\cA(1)}(H^*(B\Z/2^k;\Z/2), \Z/2)$ coming from the group multiplication $\Z/2^k\times\Z/2^k\to\Z/2^k$. With respect to this algebra structure, $\theta\kappa_i$ is indeed the product of $\theta$ and $\kappa_i$. Since we have not proven that the HASS for $\tau_{\le 2}\ko$-homology is multiplicative, we cannot say that $\theta\kappa_i\in\cQ^{0,4i-1}$ is necessarily a product in the same way, and the notation is optimistic. Nevertheless, we expect that the HASS \emph{is} multiplicative (see \cref{monoidality_conjecture}), which would make $\theta\kappa_i = \theta\cdot\kappa_i$ in $\cQ(H^*(B\Z/2^k;\Z/2))$ as well. Similar considerations apply to $\lambda_i$ and $\theta\lambda_i$, and to $\mu_i$ and $\theta\mu_i$.}
        \item For each $i\ge 1$, classes $\lambda_i\in\cQ^{1,4i+1}$ and $\theta\lambda_i\in\cQ^{1,4i+2}$, with $h_1\lambda_i = 0$ and $h_1\theta\lambda_i= 0$.
        \item For each $i\ge 1$, classes $\mu_i\in\cQ^{2,4i+4}$ and $\theta\mu_i\in\cQ^{2,4i+5}$, each annihilated by both $h_0$ and $h_1$.
    \end{itemize}
    \item For $(B\Z/2^k, 0, y)$:\footnote{An analogue of Footnote~\ref{noprod} applies here, except that to our knowledge the algebra structure on $\Ext(H_{0,y}^*(B\Z/2^k;\allowbreak \Z/2))$ has not been worked out. Our notation is based on partial computations of this algebra structure, and we do not claim it implies products in $\cQ$.}
    \begin{itemize}
        \item A class $1\in\cQ^{0,0}$, with no relations.
        \item A class $\vartheta\in\cQ^{0,1}$, with $h_1\vartheta = 0$.
        \item Classes $k_i\in\cQ^{0,4i}$ and $\vartheta k_i\in\cQ^{0,4i+1}$, $i\ge 1$, killed by $h_1$.
        \item Classes $\ell_i\in\cQ^{1,4i-1}$ and $\vartheta\ell_i\in\cQ^{1,4i}$ for all $i\ge 0$, killed by $h_1$.
        \item Classes $m_i\in\cQ^{2,4i+2}$ and $\vartheta m_i\in\cQ^{2,4i+3}$ for all $i\ge 0$, killed by both $h_0$ and $h_1$.
    \end{itemize}
\end{itemize}

When $k = 1$, in both the twisted and untwisted cases, almost all differentials vanish for degree reasons or because this is forced by commutativity with the $\sE$-action. The exceptions are:
\begin{itemize}
    \item For the untwisted computation, for all $i\ge 2$, $d_2\colon E_2^{0,4i-2}\to E_2^{2,4i-1}$ could a priori send $k_i\mapsto k_{i-1}'$. For all $i\ge 1$, $k_i$ is in the image of the map from $\partial_2\colon \Ext(\Sigma P)\to\cQ(\Sigma P)$, as we saw in \cref{cokernel_info}. By \cref{cQthm}, part~\eqref{HASS_defn}, $\partial_2$ is the map on $E_2$-pages from the ordinary Adams spectral sequence for $\widetilde{\ko}_*(B\Z/2)$ to the HASS for $\tau_{\le 2}\widetilde{\ko}_*(B\Z/2)$. Therefore we may compute this differential in the ordinary Adams spectral sequence and push it forward using $f_2$. But in the ordinary Adams spectral sequence, this differential vanishes for degree reasons. Thus this $d_2$ vanishes.
    \item In the twisted computation, for $i\ge 1$, $d_2\colon E_2^{0,4i}\to E_2^{2,4i-1}$ may a priori send $\ell_i$ to $h_1\lambda'$ ($i = 1$) or $\ell_{i-1}'$ ($i>1$). Again, \cref{cokernel_info} tells us that $\ell_i$ is in the image of $\partial_2$, so to show the differential vanishes in the HASS it suffices to show it vanishes in the usual $\ko$-Adams spectral sequence, which was shown in~\cite[\S 7.8]{Cam17}.
\end{itemize}
Thus in both cases the Adams spectral sequence collapses. On the $E_\infty$-page, all extension questions are detected by the action of $h_0$, finishing the computation and giving the groups in \cref{cyclic_supercoh}. Thus in the rest of the proof we restrict to $k>1$.

Inspecting \cref{cyclic_adams}, top and bottom center, there are many differentials which could be nonzero. As in the $k = 1$ case, \cref{Q_ext_calc} shows their sources and targets are in $\mathrm{Im}(\partial_2)$, so we may compute these differentials in the corresponding $\ko$-homology Adams spectral sequences, then push them forward. Though these $\ko$-homology differentials have been studied in several different references, including~\cite{Sie13, Cam17, Davighi:2020uab, BGR24, DDHM24, BDDHM25, Dha26}, a combination of complicated arguments and some errors in the literature (see \cref{err}) make it hard to apply their arguments in this proof. As such, we give a direct account.

As with all of the cited calculations, we use a technique inspired by the May--Milgram theorem~\cite{MM81} to compute all of the HASS differentials. Specifically, consider the Postnikov truncation map $\tau_{\le 0}\colon\ko\to H\Z$, a map of $E_\infty$-ring spectra. For any $\ko$-module $M$, we get an induced map from the $\ko$-module Baker--Lazarev Adams spectral sequence of $M$ to the corresponding $H\Z$-module spectral sequence for $H\Z\wedge_\ko M$. On $E_2$-pages, this map is precisely the map $(\tau_{\le 0})_*\colon \Ext_{\cA(1)}\to \Ext_{\cA(0)}$ induced by the inclusion $\cA(0) \coloneqq \ang{\Sq^1}\inj\cA(1)$~\cite[Proposition 13.1]{BDDHM25}.\footnote{In \textit{loc.\ cit.} $M$ is assumed to be of the form $\ko\wedge X$ for a spectrum $X$, but the more general version we gave here is true, and follows from the naturality of the Baker--Lazarev Adams spectral sequence in the ring spectrum used as a base.}
Thus we have a diagram
\begin{equation}
\label{ext_span}
\begin{tikzcd}[column sep=1em]
	& {\Ext_{\cA(1)}(M)} & \\
	{\Ext_{\cA(0)}(M)} && {\cQ(M).}
	\arrow["{(\tau_{\le 0})_*}"', from=1-2, to=2-1]
	\arrow["{\partial_2}", from=1-2, to=2-3]
\end{tikzcd}
\end{equation}
We will pull differentials back along $(\tau_{\le 0})_*$, then push them forward along $\partial_2$.
\begin{rem}
\label{err}
Let $M = \ko\wedge (B\Z/2^k)$. In this case, \cite[Corollary 13.30]{BDDHM25} calculates the map $(\tau_{\le 0})_*$, and shows it is far from injective. This means that the image of $d_r(x)$ in the $\ko$-module spectral sequence under $(\tau_{\le 0})_*$ does not uniquely determine $d_r(x)$. As such, the statement and proof of~\cite[Proposition 13.38]{BDDHM25} are not correct: see~\cite[Remark 4.7]{Dha26}. For example, when $k = 3$, this method of proof cannot tell whether $d_3(\kappa_2) = h_0^2\theta\lambda_1$ or $d_3(\kappa_2) = h_0^2\theta\lambda_1 + v\theta$.
This error does not affect any results in that paper other than Proposition 13.38. We thank Prerna Dhankhar for pointing this error out to us.

Dhankhar~\cite[Proposition 4.6]{Dha26} determines the correct values of these differentials in degrees $14$ and below, as well as in all degrees if $k = 2$. 
\end{rem}
Fortunately, this ambiguity in the values of the differentials disappears in the HASS.
\begin{defn}
Let $\mathbb M$ be an $\E$- or $\sE$-module.
A class $x\in\mathbb M$ is part of an \term{$h_0$-tower} if for all $j\ge 0$, $h_0^kx\ne 0$.
\end{defn}
\begin{lem}
\label{how_port_diffs}
Let $M\in\{\ko\wedge (B\Z/2^k), M^\ko h_{0,y}\}$ and let $x\in\Ext_{\cA(1)}(H_\ko^*(M))$ be part of an $h_0$-tower. If $x\in\ker((\tau_{\le 0})_*)$, then for $j\gg 0$, $h_0^jx\in\ker(\partial_2)$.
\end{lem}
\begin{proof}
For $M = \ko\wedge (B\Z/2^k)$, $\ker((\tau_{\le 0})_*)$ is calculated in~\cite[Corollary 13.30]{BDDHM25}, and $\ker(\partial_2)$ can be calculated using the $\cA(1)$-module decomposition~\eqref{cyccoh3} and the computation of $\partial_2$ on each summand given in \cref{Q_ext_calc}. Comparing the two finishes the proof for this choice of $M$. The proof for the other choice of $M$ is analogous.
\end{proof}
\begin{cor}
\label{cycdiffs}
The differentials in the HASSes computing $\SH^*(B\Z/2^k)$ and $\SH^*(B\Z/2^k, 0, y)$ vanish except on the $E_k$-page, where they are uniquely determined by $h_0$- and $h_1$-linearity and the following values on generators.
\begin{equation}
\label{HASS_cyc_diffs}
    \begin{alignedat}{2}
        d_k(\kappa_1) &= h_0^k \theta & d_k(\ell_0) &= h_0^{k+1}\vartheta\\
        d_k(\kappa_i) &= h_0^{k-1}\theta\lambda_{i-1},\; i>1 \qquad\qquad& d_k(k_i) &= h_0^{k-1}\vartheta\ell_{i-1}\\
        d_k(\lambda_i) &= h_0^{k+1}\theta\kappa_i & d_k(\ell_i) &= h_0^{k+1}\vartheta k_i\; i\ge 1,\\
    \end{alignedat}
\end{equation}
and $d_k$ vanishes on the classes $1$, $\theta$, $\vartheta$, $\mu_i$, $m_i$, $\theta\kappa_i$, $\vartheta k_i$, $\theta\lambda_i$, $\vartheta\ell_i$, $\theta\mu_i$, and $\vartheta m_i$.
\end{cor}
\begin{proof}
Differentials vanish on $\mu_i$, $m_i$, $\theta\mu_i$, and $\vartheta m_i$ because that is the only possibility compatible with $h_0$-linearity of differentials. For the remaining claimed differentials, use May--Milgram's identification of Adams differentials over $\cA(0)$ with differentials in the $\Z/2$-Bockstein spectral sequence for homology with $\Z_{(2)}$ coefficients~\cite{MM81}. The latter differentials are uniquely determined by the values of $H_*(B\Z/2^k;\Z_{(2)})$, and therefore so are the differentials in the $H\Z$-module Adams spectral sequence for $H\Z\wedge_\ko M$. \Cref{how_port_diffs} ports these differentials over to the HASS, determining all differentials on classes $h_0^jx$ for $j\gg 0$. This determines all differentials modulo the kernel of $h_0$. For degree reasons, the only differentials that could be in the kernel of $h_0$ are $d_2(\theta\kappa_i)\in\{0,\mu_{i-1}\}$ for $i\ge 2$ and $d_2(\vartheta k_i)\in\{0, m_0\}$ for $i\ge 1$. These differentials can be shown to vanish because the corresponding differentials in the $\ko$-module Adams spectral sequence for $M$ vanish, which follows from the Leibniz rule.
\end{proof}


Thus we can draw the $E_\infty$-page in \cref{cyclic_adams}, top right (untwisted case), resp.\ bottom right (twisted case). Despite the presence of tricky extension questions in the Adams spectral sequence for $\widetilde{\ko}_*(B\Z/2^k)$ and $\ko_*(B\Z/2^k, 0, y)$ (see, for example, \cite[\S 13.2]{DDHM24} and~\cite[\S 13.2]{BDDHM25}), in the HASSes for $\tau_{\le 2}\ko$-homology, the corresponding extension questions for multiplication by $2$ are completely solved by the $h_0$-action on the $E_\infty$-page.
Therefore we obtain the supercohomology groups claimed in \cref{cyclic_supercoh}.
\end{proof}
\begin{rem}[Eventual periodicity]
Looking at \cref{cyclic_supercoh}, the twisted supercohomology groups of $B\Z/n$ are periodic in degrees $3$ and above. This is an instance of a quite general fact: if $E$ is a ring spectrum with finitely many nonzero homotopy groups and $X$ is a space with periodic cohomology, one can use the Atiyah--Hirzebruch spectral sequence to show that $E^*(X)$ is periodic in degrees $N$ and above for some $N\gg 0$. Thus this applies to $\tau_{\le n}\ko_*(B\Z/n)$, $\tau_{\le n}\ku_*(B\Z/n)$, and twists thereof, even though the (ordinary and twisted) $\ko$- and $\ku$-(co)homology groups of $B\Z/n$ are not periodic. See for example~\cite[\S 2.2 and \S 3.4]{BG03} and~\cite[\S 6.3 and \S 7.3]{BG10}.
\end{rem}
\begin{rem}[Gysin sequences]
\Cref{cyclic_supercoh} has another pattern: for $k$ big enough, $\SH^k(B\Z/n)\cong\SH^{k+2}(B\Z/n, 0, y)$ and $\SH^k(B\Z/n, 0, y)\cong\SH^{k+2}(B\Z/n)$. This is also not a coincidence. Following the discussion in~\cite[Example 7.35]{MathSmith}, but with $\SH$ in place of $\MTSpin$, we have Gysin sequences for the complex line bundle $L\to B\Z/n$ associated to the representation as the $n^{\mathrm{th}}$ roots of unity:
\begin{subequations}
\label{SH_gysin}
\begin{gather}
    \dotsb\to\SH^k(B\Z/n) \xrightarrow{e^\SH(L)} \SH^{k+2}(B\Z/n, 0, y) \longrightarrow \SH^{k+2}(S^1) \longrightarrow \SH^{k+1}(B\Z/n) \to\dotsb\\
    \dotsb\to\SH^k(B\Z/n, 0, y) \xrightarrow{e^\SH(L)} \SH^{k+2}(B\Z/n) \longrightarrow \SH^{k+2}(S^1) \longrightarrow \SH^{k+1}(B\Z/n) \to\dotsb
\end{gather}
\end{subequations}
Here $e^\SH(V)$ denotes the product with the Euler class in (twisted) $\tau_{\le 2}\ko$-cohomology. See~\cite[\S 4.1]{MathSmith} for more information.

A quick Atiyah--Hirzebruch spectral sequence argument shows $\SH^k(S^1) = 0$ for $k>3$, implying that the Euler class maps in~\eqref{SH_gysin} are isomorphisms for those $k$.
\end{rem}


\begin{table}[htbp]
\begin{tabular}{c c c c c}
\toprule
& \multicolumn{2}{c}{$n\equiv 2\bmod 4$} &\multicolumn{2}{c}{$n\equiv 0\bmod 4$}\\
\cmidrule(r){2-3}
\cmidrule(r){4-5}
$k$ & $\SH^k(B\Z/n)$ & $\SH^k(B\Z/n, 0, y)$ & $\SH^k(B\Z/n)$ & $\SH^k(B\Z/n, 0, y)$ \\
\midrule
$0$ & $\C^\times$ & $\C^\times$ & $\C^\times$ & $\C^\times$\\
$1$ & $\Z/n\oplus\Z/2$ & $\Z/2n$ & $\Z/2n$ & $\Z/2n$\\
$2$ & $\Z/2\oplus\Z/2$ & $0$ & $\Z/2$ & $0$\\
$4j$, $j>0$ & $0$ & $\Z/2$ & $0$ & $\Z/2$\\
$4j+1$, $j>0$ & $\Z/(n/2)$ & $\Z/4n$ & $\Z/(n/2)$ & $\Z/2n\oplus\Z/2$\\
$4j+2$, $j>0$ &$0$ & $0$ & $\Z/2$ & $0$\\
$4j+3$, $j>0$ &$\Z/4n$ & $\Z/(n/2)$ & $\Z/2n\oplus\Z/2$ & $\Z/(n/2)$\\
\bottomrule
\end{tabular}
\caption{Twisted supercohomology groups of cyclic groups. We prove these isomorphisms in \cref{prop_cyclic}.} 
\label{cyclic_supercoh}
\end{table}

\subsubsection{The $16$-fold way and ``the $2$-fold way''}
\label{sss:16}
Let $B\in H^2(B^2\Z/2;\Z/2)$ denote the tautological class. In this subsubsection, we will prove $\SH^4(B^2\Z/2, 0, B)\cong\Z/16$, then identify this with the ``$16$-fold way'' classification of \cite{Kit06,johnson2024minimal}. Similarly, we will prove $\SH^5(B^2\Z/2, 0, B)\cong\Z/2$, then identify this with the two different ``fermionic $(3+1)$D topological orders'' $S$ theory and $T$ theory in \cite{JF25}.

\begin{prop}[{Johnson-Freyd~\cite[Remark 2.2.5]{JF25}}]
\label{16fold_computation}
There are isomorphisms
\begin{equation}
    \SH^k(B^2\Z/2, 0, B) = \begin{cases}
        \Z/2, &k = 5\\
        \Z/16, &k = 4\\
        0, &k = 1,2,3\\
        \C^\times, &k=0.
    \end{cases}
\end{equation}
\end{prop}
Johnson-Freyd uses the Atiyah--Hirzebruch spectral sequence to make this computation; he only discusses $\SH^4$ and $\SH^5$ explicitly, but the lower-degree groups follow clearly from his proof. This proof requires addressing a tricky extension question in $\SH^4$. We will see that the corresponding extension question in the HASS is much easier to solve.
\begin{proof}
Since $B^2\Z/2$ is trivial at odd primes, we focus on $p = 2$, where we use the HASS. Recall Serre's computation of $H^*(B^k\Z/2;\Z/2)$~\cite[Théorème 2]{Ser53}; specialized to $k = 2$, we have
\begin{equation}
    H^*(B^2\Z/2;\Z/2)\cong\Z/2[B, C, E, \dotsb],
\end{equation}
where $\abs B = 2$, $\abs C = 3$, $\abs E =5$, $\Sq(B) = B + C + B^2$, $\Sq(C) = C + E + C^2$, and $\Sq(E) = E + C^2 + \Sq^4(E) + E^2$.

Using \cref{defn:twistedcoh} we have an $\cA(1)$-module isomorphism
\begin{equation}
\label{16_fold_A1}
    H_{0,B}^*(B^2\Z/2; \Z/2)\cong
        \textcolor{Green}{\Upsilon_1} \oplus
        \textcolor{MidnightBlue}{\Sigma^4 \Upsilon_1} \oplus
        \textcolor{Fuchsia}{\Sigma^5\cA(1)} \oplus L,
\end{equation}
where $L$ is $5$-connected (and therefore cannot affect the computation of $\SH^4$). We draw the isomorphism~\eqref{16_fold_A1} in \cref{16_figure}, left. Then, using the $\cQ(\bl)$ charts we computed in \cref{Qext_A1,Qext_M0}, we draw the $E_2$-page of the HASS for $\tau_{\le 2}\ko_*(B^2\Z/2, 0, B)$ in \cref{16_figure}, right. Already we see that $\SH^0\cong\C^\times$ and that $\SH^1$, $\SH^2$ and $\SH^3$ vanish, so for the rest of the proof we focus on $\SH^4$ and $\SH^5$.

\begin{figure}[htbp]
\begin{subfigure}[c]{0.35\textwidth}
    \includegraphics{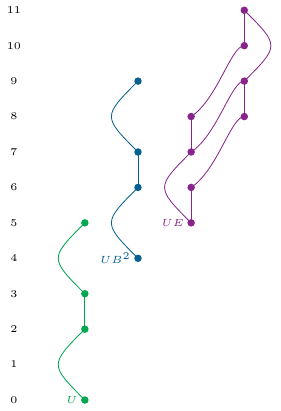}
\end{subfigure}
\begin{subfigure}[c]{0.38\textwidth}
\includegraphics{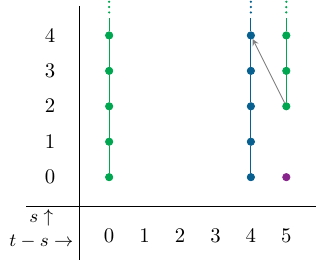}
\end{subfigure}
\caption{Left: the $\cA(1)$-module structure on $H_{0,B}^*(B^2\Z/2; \Z/2)$ in low degrees, as in~\eqref{16_fold_A1}. Right: the $E_2$-page of the HASS computing $(\tau_{\le 2}\ko)_*(B^2\Z/2, 0, B)$. We use this in the proof of \cref{16fold_computation}, where we show that the $d_2$ depicted is nonzero. 
}
\label{16_figure}
\end{figure}

In topological degrees $5$ and below, all differentials vanish for degree reasons, except for the differentials $d_r$ from the $h_0$-tower in topological degree $5$ to the $h_0$-tower in topological degree $4$. Since the $h_0$-action on $E_\infty$ lifts to multiplication by $2$, by inspection of \cref{16_figure}, right, one sees that the proposition statement is equivalent to $d_r\ne 0$ for $r = 2$, and that is what we will prove.

First, it is not too hard to show that for \emph{some} $r\ge 2$, $d_r\ne 0$: if not, $(\tau_{\le 2}\ko)_4(B^2\Z/2, 0, B)\otimes\Q$ would be nonzero, but $\tau_{\le 2}\ko$ is rationally equivalent to $H\Q$ (since it has $\pi_0\cong\Z$ and all other homotopy groups torsion), and $H_4(B^2\Z/2; \Q)\cong 0$. If $d_r\ne 0$, then by inspection $(\tau_{\le 2}\ko)_4(B^2\Z/2, 0, B)\cong\Z/2^{r+2}$ (since $h_0$ lifts to multiplication by $2$), with a minimum order of $16$ if $r = 2$. Then, from the $(r+1)$-page, there is only a single $\Z/2$ in topological degree $5$, implying $\SH^5\cong\Z/2$. Thus it remains only to compute $\SH^4$.

The $E^2$-page of the Atiyah--Hirzebruch spectral sequence for $(\tau_{\le 2}\ko)_*(B^2\Z/2, 0, B)$ provides an upper bound for the orders of these twisted $\tau_{\le 2}\ko$-homology groups, as differentials lower the order. In total degree $4$, we have a $\Z/4$ and two copies of $\Z/2$, giving an upper bound of order $16$. Putting these two together, we conclude that $(\tau_{\le 2}\ko)_4(B^2\Z/2, 0, B)\cong\Z/16$.
\end{proof}
The $\Z/16$ in $\SH^4(B^2\Z/2, 0, B)$ is the same $16$ in Kitaev's $16$-fold way~\cite[\S 9]{Kit06} (see also~\cite{BGHNPRW17, LKW17, johnson2024minimal}), as well as in the fermionic $16$-fold way of~\cite{KPY26}. In the rest of this subsubsection, we describe this connection, as well as the meaning of $\SH^k(B^2\Z/2, 0, B)$ for $0\le k\le 5$.
\begin{defn}\hfill
\label{classsp}
\begin{enumerate}
\item
A \term{Picard $n$-groupoid} is a symmetric monoidal $n$-category which is an $n$-groupoid (all $k$-morphisms are invertible, for $1\le k\le n$), such that all objects are invertible with respect to the tensor product.
\item If $\cat G$ is a Picard $n$-groupoid, the geometric realization of the nerve of $\cat G$ is a grouplike $E_\infty$-space (grouplike because all objects of $\cat G$ have inverses, $E_\infty$ via the monoidal product), and therefore $\Omega^\infty$ of a spectrum $\abs{\cat G}$, which we call the \term{classifying spectrum} of $\cat G$.
\end{enumerate}
\end{defn}
 By construction, $\pi_0(\abs{\cat G})$ is isomorphic to the abelian group of isomorphism classes of objects of $\cat G$, $\pi_1(\abs{\cat G})\cong\Aut_{\cat G}(\boldsymbol 1)$, $\pi_2(\abs{\cat G})\cong \Aut_{\Aut_{\cat G}(\boldsymbol 1)}(\id_{\boldsymbol 1})$, and so on; in particular, $\pi_k(\abs{\cat G}) = 0$ if $k\not\in [0,n]$.
\begin{thm}[{Stable homotopy hypothesis~\cite{MOPSV22}}]
\label{stab_htpy_hyp}
The assignment from a Picard $n$-groupoid to its classifying spectrum defines an equivalence of $\infty$-categories from Picard $n$-groupoids to spectra whose homotopy groups are concentrated in degrees $[0,n]$.
\end{thm}
\begin{rem}
For $n\le 2$, \cref{stab_htpy_hyp} was proven earlier. The case $n = 1$ is a folk theorem, with proofs or sketches given in~\cite{BCC93, HS05, Dri06, Pat12, JO12, GK14, Sha23}. The case $n = 2$ is a theorem of Gurski--Johnson--Osorno~\cite{GJO19}.
\end{rem}
A priori, the $16$-fold way is about tensor (1-)categories. In fact, as explained in \cite{JF25,johnson2024minimal}, it lives in the world of fusion $2$-categories, and the explicit way it is relevant is discussed below. Since braided fusion $2$-categories are in particular $E_2$-monoidal $2$-categories, they belong to a Morita $5$-category by work of Johnson-Freyd--Scheimbauer~\cite[Theorem/Definition 1.12]{JFS17}.

The connection to supercohomology passes through the following fact.
\begin{prop}[{Johnson-Freyd--Reutter~\cite{johnson2024minimal}, Johnson-Freyd~\cite{JF25}}]
\label{morita_16fold}
Let $\cat F$ be the Picard $5$-groupoid of nondegenerate, braided fermionic fusion $2$-categories with $\cat{sVect}$ as 1-morphisms. Then for $0\le k\le 5$ there are injective homomorphisms
\begin{equation}
\label{SH_to_F}
    \SH^k(B^2\Z/2, 0, B)\hookrightarrow \pi_{5-k}(\abs{\cat F}).
\end{equation}
For $k = 0,5$ this map is an isomorphism.
\end{prop}
For context, the Picard spectrum $\cat F$ is constructed out of the core of the (unpointed) Morita 5-category of nondegenerate, braided fermionic fusion 2-categories with $\sVect$ as 1-morphisms. The monoidal structure on objects is the Deligne tensor product of fusion 2-categories~\cite{Dec24a}, and the monoidal structure on higher morphisms is given by composition of bimodules.

Combining \cref{16fold_computation,morita_16fold}, we recover some known facts in and around the $16$-fold way.
\begin{itemize}
    \item On $\pi_0$, we recover the discussion in~\cite[Remark 3.16]{DY25}. $\pi_0(\abs{\cat F})\cong \SH^5(B^2\Z/2, 0, B)\cong \Z/2$, meaning that the abelian group of Morita classes of nondegenerate braided fermionic strongly fusion $2$-categories is $\Z/2$. The two Morita classes $0, 1\in\Z/2$ correspond to the theories $\mathcal S\coloneqq Z(\cat{2sVect})$, resp.\ $\mathcal T$, defined by Johnson-Freyd~\cite[\S 2.2]{JF25}; see also Johnson-Freyd--Reutter~\cite[\S 3.1]{johnson2024minimal}.
    \item On $\pi_1$, we recover Johnson-Freyd's interpretation of the $16$-fold
    way~\cite[\S 2]{JF25}:
    \[
        \SH^4(B^2\Z/2, 0, B)\cong\Aut(\mathcal S)\cong\Z/16
    \]
    (see \textit{ibid.}, Remarks 2.2.5 and 2.3.9). In more detail, the 16-fold way is the group parameterizing the  minimal nondegenerate extension of $\cat{sVect}$. In \cite{johnson2024minimal}, Johnson-Freyd--Reutter proved the existence of minimal nondegenerate extensions of $\cat{sVect}$ using the theory $\mathcal{S}$, and found that choices of automorphisms for  $\mathcal{S}$ corresponded to choices of minimal nondegenerate extension for $\cat{sVect}$. Since  nondegenerate braided fermionic strongly fusion 2-categories are classified by $\SH^5(B^2\Z/2, 0, B)$, their automorphisms are classified by $\SH^4(B^2\Z/2, 0, B)$ which is the group we compute.
    \item $\pi_5(\abs{\cat F})\cong\SH^0(B^2\Z/2, 0, B)\cong\C^\times$ is a standard fact: $\Omega^4\cat F\simeq \cat{sLine}_\C$, the Picard groupoid of invertible super vector spaces, so this $\C^\times$ is the nonzero scalars.
\end{itemize}
\begin{rem}
Barkeshli--Hsin--Kobayashi~\cite{Barkeshli:2023bta} also studied the automorphisms of $\mathcal{S}$, in the context of quantum codes on the boundary. For their applications, they use the subgroup $\Z/8$ which exchanges the two ``magnetic'' objects in $\mathcal{S}$, and this group corresponds to the 8 Ising type braided fusion categories among the 16 minimal nondegenerate extensions of $\sVect$.
\end{rem}
It is natural to guess that~\eqref{SH_to_F} lifts to a map of spectra.
\begin{conj}
\label{morita_conj}
There is a map of spectra $\tau_{\ge 0}(\Sigma^5 I_{\C^\times} M^{\tau_{\le 2}\ko}(f_B))\to\abs{\cat F}$ inducing the maps~\eqref{SH_to_F} on homotopy groups.
\end{conj}
Here $f_B\colon B^2\Z/2\to B\GL_1(\tau_{\le 2}\ko)$ is the map $B^2\Z/2\to B\rO/B\Spin$ from \cref{AB_defn}, part~\eqref{B_defn} followed by the map $(\tau_{\le 2})_*\circ T\colon B\rO/B\Spin\to B\GL_1(\tau_{\le 2}\ko)$ from \cref{htp_SH}.

One can show that the Postnikov invariants of $\tau_{\ge 0}(\Sigma^5 I_{\C^\times} M^{\tau_{\le 2}\ko}(f_B))$ vanish, so that it is a sum of shifts of Eilenberg--Mac Lane spectra. Thus an affirmative answer to \cref{morita_conj} would be progress towards computing the $k$-invariants of $\abs{\cat F}$, which have category-theoretic meaning (see~\cite{JO12, GJOS17} for lower-categorical analogues).
\subsubsection{A fourfold way over $\R$}
\label{ss:real_fourfold}
Recently, Sanford~\cite{San26} proved the following analogue of the sixteen-fold way over $\R$.
\begin{thm}[{Sanford~\cite[Theorem 5.10]{San26}}]
\label{sanford_thm}
The kernel of the map from the Witt group of $\cat{Vect}_\R$ to the Witt group of $\cat{sVect}_\R$ is isomorphic to $\Z/2\times\Z/2$.
\end{thm}
Analogously to the reformulation in~\cite[\S 4.4]{johnson2024minimal}, \cref{sanford_thm} can be rephrased as follows.
\begin{cor}
\label{sanford_cor}
Let $\cat{sAlg}_\R$ denote the Morita category of real superalgebras and $\abs{\Pic(\cat{sAlg}_\R)}$ denote the classifying spectrum of its Picard $2$-groupoid of units (see \cref{classsp}). Then there is an isomorphism
\begin{equation}
    (\Sigma^{-2}\abs{\Pic(\cat{sAlg}_\R)})^5(B^2\Z/2, 0, B) \cong\Z/2\times\Z/2.
\end{equation}
\end{cor}
We will prove this directly, using our calculations from the proof of \cref{16fold_computation}. Just like in that proof, the Atiyah--Hirzebruch spectral sequence approach to \cref{sanford_cor} is conceptually simpler but involves some more difficult computations.
\begin{rem}
To make sense of $(0, B)$-twisted supercohomology, we require $\abs{\Pic(\cat{sAlg}_\R)}$ to be a $\ko$-module. We define its $\ko$-module structure as follows: since Postnikov truncation is lax symmetric monoidal (see~\cite[Example/Proposition 4.15]{HNP25}) and the homotopy groups of $\abs{\Pic(\cat{sAlg}_\R)}$ are concentrated between degrees $0$ and $2$ by construction, $\abs{\Pic(\cat{sAlg}_\R)}$ acquires a canonical $\tau_{\le 2}\Sph$-module structure. Since the unit map $\Sph\to\ko$ is $2$-connected, it induces an equivalence $\tau_{\le 2}\Sph\xrightarrow\simeq\tau_{\le 2}\ko$, giving $\abs{\Pic(\cat{sAlg}_\R)}$ the structure of a $\tau_{\le 2}\ko$-module. Finally, precomposing with the algebra map $\ko\to\tau_{\le 2}\ko$, we get a $\ko$-module structure on $\abs{\Pic(\cat{sAlg}_\R)}$.
\end{rem}
In order to prove \cref{sanford_cor}, we need a simpler description of $\abs{\Pic(\cat{sAlg}_\R)}$. We will give it in \cref{real_pic} after a few lemmas.

The following lemma is well-known.
\begin{lem}
\label{eta_and_htpy1}
Up to $\ko$-module equivalence, there is a unique $\ko$-module $M$ such that
\begin{enumerate}
    \item $\pi_0(M)\cong \pi_1(M)\cong \Z/2$ ,
    \item $\pi_n(M) = 0$ for $n\ne 0,1$, and
    \item $\eta\colon \pi_0(M)\to\pi_1(M)$ is nonzero.
\end{enumerate}
\end{lem}
\begin{proof}
Existence is straightforward: $\Sigma^{-1}\tau_{1:2}\ko$ meets the conditions in the lemma statement. On to uniqueness.

Consider the Postnikov fiber sequence
\begin{subequations}
\begin{equation}
\label{firstpost}
    \tau_{\ge 1}M \longrightarrow M \overset{\tau_{\le 0}}\longrightarrow \tau_{\le 0}M,
\end{equation}
which is a fiber sequence of $\ko$-modules.
$\tau_{\ge 1}M$ has a single nonzero homotopy group, which is $\pi_1\cong \Z/2$, so $\tau_{\ge 1}M\simeq H\Z/2$. By a similar argument, $\tau_{\le 0}M\simeq H\Z/2$. If one takes the cofiber of $\tau_{\le 0}$ in~\eqref{firstpost}, the result is a cofiber sequence
\begin{equation}
\label{secondpost}
 M \overset{\tau_{\le 0}}{\longrightarrow} H\Z/2 \overset k\longrightarrow \Sigma^2 H\Z/2,
\end{equation}
\end{subequations}
for some $\ko$-module map $k$. In particular, the homotopy class of $k$ in $\pi_0\Map_\ko(H\Z/2, \Sigma^2 H\Z/2)\cong H_\ko^2(H\Z/2)$ determines $M$ up to $\ko$-module equivalence, as the fiber of $k$. By \cref{HH_Baker}, part~\eqref{A1_ko_H}, $H_\ko^*(H\Z/2)\cong\cA(1)$, whose subgroup of homogeneous, degree-$2$ elements is isomorphic to $\Z/2$, so there are exactly two homotopy classes of $\ko$-modules belonging to a fiber sequence of the form~\eqref{secondpost}.

If $k \simeq 0$, \eqref{secondpost} splits and $M\simeq H\Z/2\vee\Sigma H\Z/2$. Since $\eta$ acts trivially on $\pi_*(H\Z/2)$, $\eta$ acts trivially on this module; therefore there is at most $\ko$-module up to homotopy equivalence meeting the three conditions in the lemma statement.
\end{proof}
\begin{lem}
\label{EM_Hko}
There are $\cA(1)$-module isomorphisms
\begin{subequations}
\begin{align}
    H_\ko^*(\tau_{1:2}\ko) &\cong \Sigma Q\oplus \Sigma^5\uQ\\
    H_\ko^*(H\Z/8) &\cong \Upsilon_1
    \oplus \Sigma\Upsilon_1.
\end{align}
\end{subequations}
\end{lem}
The $\ko$-module structure on $H\Z/8$ is induced from the $H\Z$-module structure on $H\Z/8$ by the Postnikov quotient $\tau_{\le 0}\colon \ko\to H\Z$.
Names of $\cA(1)$-modules are as in \cref{names}.
\begin{proof}[Proof sketch]
Apply $H_\ko^*$ to the cofiber sequences of $\ko$-modules
\begin{subequations}
  \begin{gather}
        \tau_{\ge 4}\ko \overset{\tau_{\ge 3}}{\longrightarrow} \tau_{\ge 1}\ko \overset{\tau_{\le 2}}{\longrightarrow} \tau_{1:2}\ko\\
        H\Z\overset 8\longrightarrow H\Z\longrightarrow H\Z/8
  \end{gather}
\end{subequations}
to obtain long exact sequences of $\cA(1)$-modules. By \cref{A1_module_of_connective_covers}, $H_\ko^*(\tau_{\ge 4}\ko)\cong \Sigma^4 \uQ$ and $H_\ko^*(\tau_{\ge 1}\ko)\cong \Sigma Q$, and by~\cite[\S 5]{Bak20}, $H_\ko^*(H\Z)\cong \Upsilon_1$. The lemma statement then follows directly from exactness of these two long exact sequences.
\end{proof}
\begin{cor}
\label{maps_cyclic}
$\pi_0\Map_\ko(\Sigma^{-1}\tau_{1:2}\ko, \Sigma^3 H\Z/8)$ is isomorphic to either $\Z/4$ or $\Z/8$.
\end{cor}
\begin{proof}
Since $H\Z/8$ is $2$-local, we apply the Baker--Lazarev Adams spectral sequence. The $E_2$-page is
\begin{subequations}
\begin{equation}
    \begin{aligned}
        E_2^{s,t} &\cong \Ext_{\cA(1)}^{s,t}(H_\ko^*(\Sigma^3 H\Z/8), H_\ko^*(\Sigma^{-1}\tau_{1:2}\ko))\\
        &\cong \Ext^{s,t}_{\cA(1)}(\Sigma^3\Upsilon_1 \oplus \Sigma^4 \Upsilon_1, Q\oplus \Sigma^4\uQ)
    \end{aligned}
\end{equation}
by \cref{EM_Hko}. Since $\Upsilon_1\cong\cA(1)\otimes_{\cA(0)}\Z/2$ (\cref{seagull_defn}), the change-of-rings theorem (see for example \cite[\S 4.5]{BC18}) simplifies the $E_2$-page to
\begin{equation}
\label{A0_ext}
    \Ext_{\cA(0)}^{s,t}(\Sigma^3 \Z/2 \oplus \Sigma^4 \Z/2, Q\oplus \Sigma^4 \uQ).
\end{equation}
\end{subequations}
One can directly verify that there are $\cA(0)$-module isomorphisms $Q\cong \cA(0)\oplus \Sigma^3 \Z/2$ and $\uQ\cong \Z/2\oplus \Sigma^2 \cA(0)$ (see, e.g., \cite[Figure 31]{DDHM24}). Ext commutes with direct sums in both arguments, reducing the computation of~\eqref{A0_ext} to that of $\Ext_{\cA(0)}(\Z/2, \cA(0))$ and $\Ext_{\cA(0)}(\Z/2, \Z/2)$. The former consists of a single $\Z/2$ in bidegree $(s,t) = (0, -1)$, which follows from graded self-injectivity of $\cA(0)$ (see~\cite[Proposition 12.8, Theorem 12.9]{Mar83}) or can be worked out directly, and the latter is $\Z/2[h_0]$ with $h_0$ in bidegree $(1, 1)$~\cite[Remark 4.5.4]{BC18}.

Putting it all together, the $E_2$-page contains four $\Z/2[h_0]$ summands (``$h_0$-towers'') in topological degrees $-1$, $0$, $0$, and $1$, together with $\Z/2$ summands in filtration $0$ and topological degrees $-4$, $-3$, $2$, and $3$. Since $H\Z/8$ vanishes rationally, $\pi_*\Map_\ko(\Sigma^{-1}\tau_{1:2}\ko, \Sigma^3 H\Z/8)\otimes \Q\cong 0$. Therefore differentials must kill all but finitely many of the elements of the four $h_0$-towers. This means that one of the $h_0$-towers in degree $0$ supports a $d_{r_1}$ differential to the $h_0$-tower in degree $-1$, so this entire tower vanishes after the $E_{r_1}$-page, and the other $h_0$-tower in degree $0$ must receive a $d_{r_2}$ differential from the $h_0$-tower in degree $1$, for some $r_2\ge 2$. Thus, on the $E_\infty$-page, in topological degree $0$ we have $\Z/2[h_0]/(h_0^{r_2+1})$; since $h_0$ lifts to multiplication by $2$, this implies $\pi_0\Map_\ko(\Sigma^{-1}\tau_{1:2}\ko, \Sigma^3 H\Z/8)\cong\Z/2^{r_2}$.

Finally, we need to show $r_2\le 3$. If it were larger, a group of homotopy classes of maps into $H\Z/8$ would contain elements on which $8\ne 0$, but this cannot occur.
\end{proof}

In what follows, we will make use of the \term{Anderson dual of the sphere}, $I_\Z$~\cite{And69}; see~\cite[Appendix B]{FMS07} or~\cite[\S 5.3]{FH21} for a review.

\begin{lem}
\label{eta_and_htpy2}
Up to $\ko$-module equivalence, there is a unique $\ko$-module $N$ such that
\begin{enumerate}
    \item\label{gps__} $\pi_0(N)\cong\Z/8$, $\pi_1(N)\cong\Z/2$, and  $\pi_2(N)\cong\Z/2$,
    \item\label{2type} $\pi_n(N) = 0$ for $n\not\in\{0,1,2\}$, and
    \item\label{eta2} $\eta^2\colon \pi_0(N)\to\pi_2(N)$ is nonzero.
\end{enumerate}
\end{lem}
\begin{proof}
The universal property of Anderson duality implies that it is equivalent to show that there is a unique $\ko$-module $N'$ up to $\ko$-module equivalence satisfying conditions~\eqref{2type} and~\eqref{eta2}, but with~\eqref{gps__} replaced with $\pi_0(N')\cong\Z/2$, $\pi_1(N')\cong\Z/2$, and $\pi_2(N')\cong\Z/8$. This is because, if $N'$ satisfies these desiderata, then $N\coloneqq \Sigma^3 I_\Z N'$ satisfies the conditions in the lemma statement, and vice versa, so uniqueness of $N'$ implies uniqueness of $N$. Thus in the rest of the proof we will focus on $N'$.

For existence, take $N'\coloneqq \tau_{\le 2}(\ko\wedge (B\Z/2)^{\sigma-1})$. The facts we need about its homotopy groups and $\eta^2$-action are standard, and follow from Mahowald's~\cite[Lemma 7.3]{Mah82}, resp.\ Bruner--Mira--Stanley--Snaith's ~\cite[\S 3.1]{BMSS15} computation of $\ko_*(B\Z/2)$ as a graded abelian group, resp.\ $\ko_*$-module, after using the Smith isomorphism $\Sigma^{-1}\Sigma^\infty B\Z/2\simeq (B\Z/2)^{\sigma-1}$ to pass from $\Sigma^{-1} \ko\wedge B\Z/2$ to $\ko\wedge (B\Z/2)^{\sigma-1}$: see, e.g., \cite[Lemma 2.6.5]{Koc96}.

Thus we turn to uniqueness. 
By a similar line of reasoning as we used to build~\eqref{secondpost} during the proof of \cref{eta_and_htpy1}, there is a fiber sequence of $\ko$-modules
\begin{equation}
\label{thirdpost}
    N' \overset{\tau_{\le 1}}{\longrightarrow} M \overset k\longrightarrow \Sigma^3 H\Z/8,
\end{equation}
where $M$ is a $\ko$-module satisfying the conditions of \cref{eta_and_htpy1}. Thus, as we proved in \cref{eta_and_htpy1}, $M\simeq \Sigma^{-1}\tau_{1:2}\ko$ as $\ko$-modules, so by \cref{maps_cyclic}, the group of homotopy classes of possible maps $k$ in~\eqref{thirdpost} is isomorphic to $\Z/4$ or $\Z/8$.
%

Precomposing with $\tau_{\ge 1}\colon \Sigma H\Z/2\simeq\tau_{\ge 1}M\to M$ defines a homomorphism
\begin{equation}
    \tau_{\ge 1}^*\colon \pi_0\Map_\ko(M, \Sigma^3 H\Z/8) \longrightarrow \pi_0\Map_\ko(\Sigma H\Z/2, \Sigma^3 H\Z/8)
\end{equation}
sending $k$ to the Postnikov invariant of $\tau_{\ge 1}N'$. Since $2 = 0$ in $\pi_0\Map_\ko(\Sigma H\Z/2, \Sigma^3 H\Z/8)$, then if $k$ is even, then $\tau_{\ge 1}^*(k) = 0$, implying $\tau_{\ge 1}N'\simeq \Sigma H\Z/2\vee \Sigma^2 H\Z/8$. This is incompatible with our requirement that $\eta^2\colon\pi_0(N')\to\pi_2(N')$ is nonzero.

Thus $k$ is an odd element of $\Z/4$ or $\Z/8$.
%
In this case, there is an automorphism of $\Sigma^3 H\Z/8$ sending $k\mapsto 1$. Thus there is at most one homotopy class of $\ko$-modules meeting the conditions in the lemma statement.
%
\end{proof}

\begin{prop}
\label{real_pic}
There is a $\ko$-module equivalence \begin{equation}
\abs{\Pic(\cat{sAlg}_\R)}\simeq \Sigma^2 I_\Z(\tau_{\le 2}(\ko\wedge (B\Z/2)^{\sigma-1})).\end{equation}
\end{prop}
We believe \cref{real_pic} is known to some experts in the field, but are unaware of an explicit reference in the literature. See~\cite[\S 6.4]{GJF:2017zba} for a closely related result.
\begin{proof}[Proof sketch]
We will show both of these $\ko$-modules meet the hypotheses of \cref{eta_and_htpy2}, so that the conclusion of that lemma implies they are $\ko$-module equivalent. In fact, for $\Sigma^2 I_\Z\tau_{\ge 0}(\ko\wedge (B\Z/2)^{\sigma-1})$, this follows from our analysis of $\tau_{\le 2}(\ko\wedge (B\Z/2)^{\sigma-1})$ during the proof of \cref{eta_and_htpy2} together with the universal property of Anderson duality.

Thus we focus on $\abs{\Pic(\cat{sAlg}_\R)}$. Freed~\cite[(1.40)]{Fre12} showed that $\pi_*(\abs{\Pic(\cat{sAlg}_\R)})$ satisfies conditions~\eqref{gps__} and~\eqref{2type} of \cref{eta_and_htpy2}. That $\eta^2\ne 0$ as maps $\pi_0(\abs{\Pic(\cat{sAlg}_\R)})\to\pi_2(\abs{\Pic(\cat{sAlg}_\R)})$ can be shown in the same way as is shown for $\cat{sAlg}_\C$ in~\cite[\S 4.3.5]{DG18}.
\end{proof}
\begin{cor}
\label{reduce_to_ko_BZ2}
Let $h_{0,B}\colon B^2\Z/2\to B\GL_1(\ko)$ be the twist induced by $(0, B)$ as defined in \cref{stop_at_ko}. There is an isomorphism
\begin{equation}
\label{dual_not_can}
(\Sigma^{-2}\abs{\Pic(\cat{sAlg}_\R)})^5(B^2\Z/2, 0, B) \cong \pi_4(\ko\wedge (B\Z/2)^{\sigma-1}\wedge_\ko M^\ko h_{0,B}).
\end{equation}
\end{cor}
\begin{proof}[Proof sketch]
Apply \cref{real_pic} to replace $\Sigma^{-2}\abs{\Pic(\cat{sAlg}_\R)}$ with $\Sigma^4 I_\Z(\tau_{\le 2}(\ko\wedge (B\Z/2)^{\sigma-1}))$. But since $\ko_n((B\Z/2)^{\sigma-1}) = 0$ for $n = 3,4,5$ (see for example \cite[Theorem 6.4]{Cam17} or~\cite[Example 5.31]{BottSpiral}), truncation defines a $\ko$-module equivalence
\begin{equation}
\tau_{\le 2}^*\colon
\Sigma^4 I_\Z(\tau_{\le 2}(\ko\wedge (B\Z/2)^{\sigma-1}))\overset\simeq\longrightarrow
\Sigma^4 I_\Z(\tau_{\le 5}(\ko\wedge (B\Z/2)^{\sigma-1})).
\end{equation}
Next apply the universal property of Anderson duality to identify the degree-$5$ twisted $\Sigma^4 I_\Z(\tau_{\le 5}(\ko\wedge (B\Z/2)^{\sigma-1}))$-cohomology of $B^2\Z/2$ with the degree-$4$ twisted $\tau_{\le 5}(\ko\wedge (B\Z/2)^{\sigma-1})$-homology of $B^2\Z/2$. Finally, since $\tau_{\le 5}$ is $5$-connected, it induces an isomorphism on (twisted) homology in degrees $5$ and below, so we may remove $\tau_{\le 4}$, giving the right-hand side of~\eqref{dual_not_can}.
\end{proof}
Finally we can return to our reformulation of Sanford's result~\cite[Theorem 5.10]{San26}.
\begin{proof}[Proof of \cref{sanford_cor}]
By \cref{reduce_to_ko_BZ2}, it suffices to show
\begin{equation}
\label{only_left}
\pi_4(\ko\wedge (B\Z/2)^{\sigma-1}\wedge_\ko M^\ko h_{0,B})\cong\Z/2\oplus\Z/2.
\end{equation}
We will do this with the Baker--Lazarev Adams spectral sequence. The $E_2$-page is the Ext of the $\cA(1)$-module
\begin{subequations}
\begin{equation}
    H_\ko^*(\ko\wedge (B\Z/2)^{\sigma-1}\wedge_\ko M^\ko h_{0,B}) \cong H^*((B\Z/2)^{\sigma-1}; \Z/2) \otimes H_\ko^*(M^\ko h_{0,B})
\end{equation}
by \cref{rel_kunneth} (and the identification $H_\ko^*(\ko\wedge X)\cong H^*(X;\Z/2)$~\cite[Lemma 5.16]{BottSpiral}). By~\eqref{16_fold_A1} and the identification $B\Z/2\simeq \Sigma (B\Z/2)^{\sigma-1}$ mentioned above,
\begin{equation}
H_\ko^*(M^\ko h_{0,B})\otimes
    H^*((B\Z/2)^{\sigma-1}; \Z/2) \cong (\textcolor{Green}{\Upsilon_1} \oplus \textcolor{MidnightBlue}{\Sigma^4\Upsilon_1} \oplus \textcolor{Fuchsia}{\Sigma^5\cA(1)} \oplus L)\otimes P,
\end{equation}
\end{subequations}
where $L$ is $5$-connected, and therefore can be ignored for our $\pi_4$ computation.

Since $\Upsilon_1\cong \cA(1)\otimes_{\cA(0)}\Z/2$ (\cref{seagull_defn}), the change-of-rings theorem~\cite[\S 4.5]{BC18} implies
\begin{equation}
    \Ext_{\cA(1)}(\Upsilon_1\otimes P, \Z/2)\cong\Ext_{\cA(0)}(P, \Z/2).
\end{equation}
One can check directly that, as an $\cA(0)$-module, $P$ is free on a homogeneous basis with one element in each nonnegative even degree. Since $\Ext_{\cA(0)}(\cA(0), \Z/2)$ consists of a $\Z/2$ in bidegree $(0, 0)$, we conclude that $\Ext_{\cA(0)}(P, \Z/2)$ consists of a $\Z/2$ in each bidegree $(s, t) = (0, 2k)$ for each $k\ge 0$.

By a similar argument, $\Ext_{\cA(1)}(\cA(1)\otimes P; \Z/2)\cong\Ext_{\Z/2}(P, \Z/2)$. Since $\Z/2$ is a field and $P$ is connective, these Ext groups vanish except when $s = 0$ and $t\ge 0$.

Combining all of these facts, we see that in topological degree $4$, the $E_2$-page of the Baker--Lazarev Adams spectral sequence has $E_2^{0,4}\cong \Z/2\oplus\Z/2$, and for all $t-s\le 5$, the $E_2$-page vanishes for $s>0$. Thus the spectral sequence collapses in this range to imply~\eqref{only_left}.
\end{proof}

\subsubsection{The quaternion group of order 8}\label{subsection:Q8}
When performing the calculations of either $\ko$-homology or supercohomology, it is sometimes helpful to be able to compute the Atiyah--Hirzebruch filtration of a given class $x$. In the mathematical physics literature, this is often called the \term{layer} of $x$~\cite{WG20}. This controls where the image of $x$ on the $E^\infty$-page of the AHSS is.

Specializing to supercohomology, for any space $X$ and twisting data $(a,b)$, there is a map $(\bl)_{\mathrm{Maj}}\colon \SH^k(X, a, b)\to H^{k-2}(X;\Z/2)$, called the \term{Majorana layer} \cite{WG20,Debray:2025kfg}, dual to the $1$-connected cover
\begin{equation}
    \tau_{\ge 2}\colon \Sigma^2 H\Z/2\simeq \tau_{\ge 2}\tau_{\le 2}\ko\longrightarrow\tau_{\le 2}\ko.
\end{equation}
Johnson-Freyd--Yu~\cite[\S 1]{JFY21} asked whether there is a finite group $G$ such that the Majorana layer homomorphism $\SH^4(BG, 0, 0)\to H^2(BG;\Z/2)$ is nonzero, and Décoppet--Johnson-Freyd~\cite[Appendix A]{Dec24} showed that this happens for $G = \Z/4\times\Z/4$. This information helps determine whether the corresponding topological orders are bosonic gauge theories in the usual sense. In~\cite[Appendix C]{Debray:2025kfg} (see also \textit{ibid.}, ``Definition'' 3.15), we defined the \term{$p+ip$ layer} associated to a $\ko$-homology or spin bordism class, which is the part of the Atiyah--Hirzebruch filtration corresponding to $\pi_4\ko\cong\Z$. The $p+ip$ layer acts as an obstruction for an invertible field theory to be the anomaly of a unitary TFT~\cite{Debray:2025kfg,DYY3}.


Using the map $\AtoA^\ko$ from \cref{map_of_filtered_spectra}, which maps from the Atiyah--Hirzebruch to the Adams filtration, we can use the Adams filtration of a class to constrain which layer it lies in.
We will apply this for $G = Q_8$ and show that the Majorana layer homomorphism $\SH^4(BG, 0, 0)\to H^2(BG;\Z/2)$ is in fact zero, using the HASS to calculate $\SH^*(BQ_8)$ and $\AtoA^\ko$ to compute the Majorana layer homomorphism.

First, the supercohomology calculation.
\begin{prop}
\label{SH_Q8}
Let $j\ge 1$. Then there are isomorphisms
\begin{equation}
    \begin{alignedat}{2}
        \SH^0(BQ_8) &\cong \C^\times\qquad\qquad\quad
            & \SH^{4j}(BQ_8) &\cong \Z/2\\
        \SH^1(BQ_8) &\cong (\Z/2)^{\oplus 3}
            & \SH^{4j+1}(BQ_8) &\cong (\Z/2)^{\oplus 2}\\
        \SH^2(BQ_8) &\cong (\Z/2)^{\oplus 3}
            & \SH^{4j+2}(BQ_8) &\cong (\Z/2)^{\oplus 3}\\
        && \SH^{4j-1}(BQ_8) &\cong \Z/8\oplus (\Z/4)^{\oplus 2}.
    \end{alignedat}
\end{equation}
The unique nonzero element of $(\tau_{\le 2}\ko)_{4j}(BQ_8)\cong\Z/2$ has Adams filtration $1$.
\end{prop}
\begin{proof}
We will compute reduced $\tau_{\le 2}\ko$-homology and trust the reader to dualize and add back in $\SH^*(\pt)$ at the end. Mitchell--Priddy~\cite[Theorem D]{MP84} produce a $2$-local stable splitting
\begin{equation}
    \Sigma^\infty BQ_8 \simeq_{(2)} \Sigma^\infty B\SL_2(\F_3) \vee C\vee C,
\end{equation}
where $C$ is the fiber of the map $\Sigma^\infty B\Pin_2^-\to \Sigma^\infty B\Spin_3$ induced by the usual inclusion $\Pin_2^-\inj\Spin_3$. Thus we may study $\widetilde\SH{}^*(B\SL_2(\F_3)))$ and $\SH^*(C)$ separately, then $2$-localize and direct-sum them together.

Bayen--Bruner~\cite[\S 3]{BB96} produce $\cA(1)$-module isomorphisms
\begin{subequations}
\begin{align}
    \widetilde H^*(B\SL_2(\F_3);\Z/2) &\cong \bigoplus_{j\ge 1}^\infty \Sigma^{4j-1}\Z/2 \oplus \Sigma^{4j}\Z/2\\
    H^*(C;\Z/2) &\cong \bigoplus_{j\ge 0}^\infty \Sigma^{4j+1} N_1,
\end{align}
\end{subequations}
using Quillen's calculation~\cite[\S 8]{Qui72} of $H^*(B\SL_2(\F_3);\Z/2)$ (see~\cite[Theorem 2.4]{MP84} for the ring structure), resp.\ Mitchell--Priddy's calculation~\cite[\S 6]{MP84} of $H^*(C;\Z/2)$.

Thus, using the charts for $\cQ(\Z/2)$ and $\cQ(N_1)$ that we computed in \cref{Q_ext_F2,Q_N1_fig}, we draw the $E_2$-pages of the HASSes for $\tau_{\le 2}\ko$ of $B\SL_2(\F_3)$ and $C$ in \cref{two_Q8_HASS}, upper left and lower left, respectively. \Cref{Q_ext_calc} shows that $\partial_2\colon \Ext(M)\to\cQ(M)$ is surjective for $M \in\{\Z/2, N_1\}$, so for both $B\SL_2(\F_3)$ and $C$, the map on $E_2$-pages from the usual Adams spectral sequences for $\ko$-homology to the HASSes for $\tau_{\le 2}\ko$-homology are surjective. Thus the same is true on $E_3$ if there are no differentials in the $\ko$-homology Adams spectral sequence on $E_2$.

Bayen--Bruner~\cite[\S 3]{BB96} compute all differentials in the Adams spectral sequences for $\widetilde\ko_*(B\SL_2(\F_3))$ and $\ko_*(C)$ (see \textit{ibid.}, Figures 1, 2, 3, and 4). As described in the previous paragraph, these imply the following differentials in the HASSes for $B\SL_2(\F_3)$ and $C$.
\begin{itemize}
    \item For $B\SL_2(\F_3)$, Bayen--Bruner's computations together with our calculation of $\partial_2$ in \cref{Q_ext_calc} imply that in the HASS, $d_2 = 0$ and $d_3$ is uniquely determined by the fact that it vanishes except in topological degrees $t-s = 4j$, $j\ge 1$. On all later pages, differentials vanish for degree reasons, so $E_4 = E_\infty$.
    \item For $C$, Bayen--Bruner's computations together with \cref{Q_ext_calc} imply that in the HASS, $d_2$ is nonzero precisely in topological degrees $t-s = 4j+1$, $j\ge 1$, which determines it uniquely. From $E_3$ onward, the HASS collapses for degree reasons.
\end{itemize}
We draw the $E_\infty$-pages of these spectral sequences in \cref{two_Q8_HASS}, top right and bottom right, respectively. For $C$, degree reasons imply there are no extension problems and we are done. For $B\SL_2(\F_3)$, there a priori could be a hidden extension in degrees $4j+1$, $j\ge 1$. These classes are in the image of the map from $\ko$, and Bayen--Bruner~\cite[Theorem 5]{BB96} show that the corresponding extension in $\ko$-homology splits, implying the extension in $\tau_{\le 2}\ko_{4j+1}(BQ_8)$ also splits.

\begin{figure}[htbp]
\begin{subfigure}[c]{0.48\textwidth}
    \includegraphics{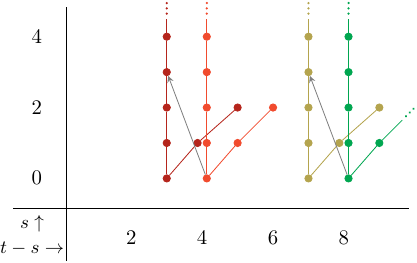}
\end{subfigure}
\begin{subfigure}[c]{0.48\textwidth}
    \includegraphics{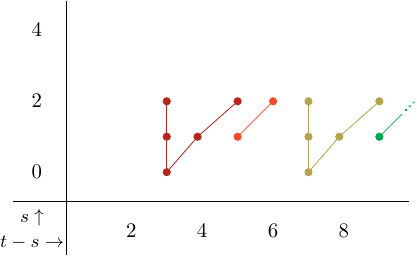}
\end{subfigure}
\\[0.5em]
\begin{subfigure}[c]{0.48\textwidth}
    \includegraphics{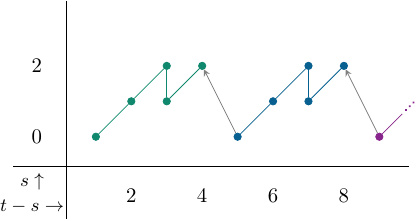}
\end{subfigure}
\begin{subfigure}[c]{0.48\textwidth}
    \includegraphics{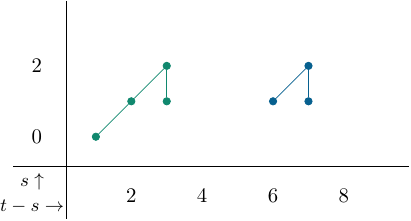}
\end{subfigure}
\caption{The $\tau_{\le 2}\ko$ HASSes for the direct summands of $BQ_8$. Upper left: the $E_2 = E_3$-page of the HASS computing $(\widetilde{\tau_{\le 2}\ko})_*(B\SL_2(\F_3))_2^\wedge$. In the proof of \cref{SH_Q8}, we deduce the pictured differentials from the corresponding differentials in the ordinary Adams spectral sequence for $\widetilde\ko_*(B\SL_2(\F_3))$, which Bayen--Bruner~\cite[\S 3]{BB96} calculated. Upper right: the $E_4 = E_\infty$-page. Lower left: the $E_2$-page of the corresponding HASS for the spectrum $C\coloneqq\mathrm{fib}(\Sigma^\infty B\Pin_2^-\to \Sigma^\infty B\Spin)$, whose differentials follow in the same way. Lower right: the $E_3 = E_\infty$-page of this HASS. We use these computations in \cref{SH_Q8} as key steps in the computation of $\SH^*(BQ_8)$.}
\label{two_Q8_HASS}
\end{figure}

The last claim to verify is the Adams filtration of the unique nonzero class $x_{4j}$ in degrees $4j$, $j\ge 1$. This class is in the $B\SL_2(\F_3)$ direct summand, since we have just shown $\tau_{\le 2}\ko_{4j}(C) = 0$. Inspecting \cref{two_Q8_HASS}, upper right, we see that the image of $x_{4j}$ in the $E_\infty$-page has filtration $1$, so we are done.
\end{proof}
\begin{cor}
For all $j$, the Majorana layer homomorphism
\[
    (\bl)_{\mathrm{Maj}}\colon\SH^{4j}(BQ_8)\to H^{4j-2}(BQ_8;\Z/2)
\]
vanishes.
\end{cor}
\begin{proof}
By the universal property of Anderson duality, it suffices to show that the Anderson dual map $\tau_{\ge 2}\colon H_{4j-2}(BQ_8;\Z/2)\to (\tau_{\le 2}\ko)_{4j}(BQ_8)$ is zero. Equivalently, we want to show that $\tau_{\le 2}\ko_{4j}(BQ_8)$ has no elements $x$ of Atiyah--Hirzebruch filtration $2$ or higher. The map $\Phi^\ko$ between the AHSS and Adams filtrations implies that
any such element $x$ would then have Adams filtration at least $2$, which we proved does not happen in \cref{SH_Q8}.
\end{proof}

\section{Further directions}
\label{s:further}
In this section, we list some computations and conjectures regarding the HASSes we studied in this paper and potential generalizations.

\subsection{The HASS for restricted supercohomology}
\label{ss:restricted}
Though we focused on $\tau_{\le 2}\ko$ and applications to supercohomology in this paper, our main theorems hold for all Postnikov truncations of $\ko$ and $\ku$, and the calculation of the corresponding HASSes proceeds similarly. It would be interesting to study the HASS for other truncations of $\ko$, or of any truncation of $\ku$.

In the specific example of $\tau_{\le 1}\ko$, Pontrjagin dual to \term{restricted supercohomology} $\mathit{rSH}$~\cite{Fre08, GW14}, one could apply HASS computations towards several problems. 
\begin{exm}
Restricted supercohomology classifies levels of spin Chern--Simons theories~\cite{Jen05, FN23} as well as a class of twists of string bordism and topological modular forms: see~\cite{FHT10, JF20, DY23, TY23a, TY23b, TY25b, BLM23}.
\end{exm}
\begin{exm}
Let $J_{\mathrm{PL}}\colon \mathrm{PL}\to \GL_1(\Sph)$ be the
    piecewise linear $J$-homomorphism. By work of Sullivan~\cite{Sul66}, maps from a space $M$ to the cofiber $\GL_1(\Sph)/\mathrm{PL}$ of $J_{\mathrm{PL}}$ are closely related to the classification of piecewise linear manifold structures on $M$. Sullivan~\cite[Theorem 4(i)]{Sul96} (see also Madsen--Milgram~\cite[Theorem 4.34(a)]{MM79}) showed that, $2$-locally,\footnote{The homotopy type of $\GL_1(\Sph)/\PL$ away from $2$ is closely related to that of $\KO$ ~\cite[Theorem 4(ii)]{Sul96}.} $\GL_1(\Sph)/\mathrm{PL}$ splits as a space into a product of Eilenberg--Mac Lane spaces and one factor of $\Omega^{\infty-4}\mathit{rSH}$. Thus to classify the PL structures on $M$, one needs to compute $\mathit{rSH}^4(M)$, which by Anderson duality is determined by $\tau_{\le 1}\ko_i(M)$ for $i = 3,4$.

    Prior work computing $[\bl, \GL_1(\Sph)/\mathrm{PL}]$ for applications to classifying PL structures includes \cite{Wal68, LdM71, Rud74, Rub79, Muk03, Muk04, Muk06, Wan22, Bak26}. These works used a variety of ad hoc methods to resolve extension questions in $\mathit{rSH}$-cohomology of various manifolds of interest, such as projective spaces and Dold manifolds. The $\tau_{\le 1}\ko$-HASS may enable a systematic approach or an extension to more classes of manifolds.
\end{exm}
\begin{exm}
Johnson-Freyd--Treumann~\cite[\S 6]{JFT20} describe how understanding $\mathit{rSH}^*(BQ_{16})$ can be used to calculate fractional Pontrjagin classes of representations of the Conway group $\mathrm{Co}_0$, with applications to conformal field theory. See also the discussion in~\cite[\S 2.3]{JF20}.
\end{exm}
We predict that, like for (extended) supercohomology, the HASS will solve extension problems better than the AHSS (which in this case simplifies to the long exact sequence of~\cite[(1.3)]{Fre08}).
\subsection{Generalizations}
\label{ss:conjectures}
In this section we gather some conjectures on the HASSes we constructed and some possible generalizations.

First, both the $R$-module Adams and Atiyah--Hirzebruch spectral sequences are \emph{multiplicative}: applied to an $R$-algebra spectrum, one obtains a natural product on each page of these spectral sequences, and differentials satisfy a Leibniz rule with respect to this product. More generally, there is an analogue of the cross product relating the spectral sequences of two $R$-modules $M$ and $N$ with that of $M\wedge_R N$. The analogous statement at the level of filtered spectra is the claim that the Atiyah--Hirzebruch and Adams filtrations, as functors from $R$-modules to filtered $R$-modules have lax symmetric monoidal structures.\footnote{See Hedenlund~\cite[Theorem II.3.5]{Hed20}, Antieau~\cite[Theorem 8.9]{Ant24}, and Wu~\cite[Theorem E]{Wu26} for the details of the connection between multiplicative spectral sequences as defined above and lax symmetric monoidal structures on filtrations.} Surprisingly, this does not appear to be in the literature in general.
\begin{itemize}
    \item The symmetric monoidal structure on $\Post^\Sph_\bullet$ is discussed in~\cite[Remark 2.25]{vN25}. The same proof also applies to $\Post_\bullet^R$, because the result in \textit{loc.\ cit.} follows directly from the results in~\cite[\S II.1.4]{Hed20}, which are proven in the generality of a stable $\infty$-category with $t$-structure, such as $\cat{Mod}_R$.
    \item For $\Ad_\bullet^\Sph$, the lax symmetric monoidal structure is due to Pstrągowski~\cite{Pst23}: see Patchkoria--Pstrągowski~\cite[Theorem 5.64 and Remark 5.65]{PP26}, who also removed an Adams type assumption. For more general $E_\infty$-ring spectra $R$, this does not seem to be written down, though it would presumably follow from a good theory of synthetic $R$-module spectra (see \cref{synthetic}).
\end{itemize}
\begin{conj}
\label{monoidality_conjecture}
There is a lax symmetric monoidal structure on $\Ad_\bullet^R$ refining the multiplicative structure on the Baker--Lazarev Adams spectral sequence, such that the natural transformations $\AtoA^\ko\colon \Post_\bullet^\ko\Rightarrow\Ad_{g(\bullet)}^\ko$ and $\AtoA^\ell\colon  \Post_\bullet^\ell\Rightarrow\Ad_{g^c(\bullet)}^\ell$ we constructed in \cref{map_of_filtered_spectra} are lax monoidal.
\end{conj}
Carrick~\cite{Car25} studies a map between the Adams--Novikov and slice filtrations in genuine equivariant homotopy theory, and proves that this map sends $E_\infty$-algebras to $E_\infty$-algebras, so \cref{monoidality_conjecture} suggests a parallel to his result.

\Cref{monoidality_conjecture} would imply that $\AtoA^\ko$ and $\AtoA^\ell$ commute with products, and should also lead to a proof that the HASSes we constructed for $\tau_{\le n}\ko$ and $\tau_{\le n}\ell$ are also multiplicative. This would be helpful for computations.
\begin{rem}
We gave two quite different proofs of \cref{map_of_filtered_spectra} in two cases of $\ko$ and $\ell$. The proof for $\ell$ began with an identification of the $\ell$-module Atiyah--Hirzebruch filtration with the $H\Z_{(p)}$-based, $\ell$-module Adams filtration (\cref{post_ku_HZ}). Showing this is a monoidal natural isomorphism would immediately imply \cref{monoidality_conjecture} in the case of $\ell$.

Unfortunately, this approach does not seem to work for $\ko$. The approach we found in this case seems much harder to make monoidal.
\end{rem}
Our other dream is that there are surely more examples of this Atiyah--Hirzebruch-to-Adams phenomenon. We list a few that we are especially optimistic about.
\begin{conj}
\label{ksc_conj}
Recall (e.g.\ from~\cite[\S 6]{Hys23}) that the \term{connective, self-conjugate $K$-theory spectrum}, $\ksc\coloneqq\ko/\eta^2$, is an $E_\infty$-$\ko$-algebra.
\begin{enumerate}
    \item \label{ksc_1} Let $\mathrm{HAd}_\bullet^{\ksc}$ denote the ($2$-primary) hastened Adams filtration associated to the fiber sequence $\Sigma^2\ko\overset{\eta^2}{\to}\ko\to\ksc$. Then, there is a natural transformation of functors $\Mod_{\ksc}\to\Fil_{\ksc})$
    \begin{equation}
        \AtoA^\ksc\colon \Post_\bullet^{\ksc} \Longrightarrow \mathrm{HAd}_{g^{\ksc}(\bullet)}^\ksc
    \end{equation}
    which, when composed with $(\bl)_{-\infty}$, is $2$-completion. (The function $g^{\ksc}$ is defined in \cref{0_homology_filt}.)
    \item There is a hastened Adams spectral sequence for $\tau_{\le k}\ksc$-homology analogous to the ones we constructed for $\ko$ and $\ell$.
\end{enumerate}
\end{conj}
The hastened Adams filtration mentioned in part~\eqref{ksc_1} of \cref{ksc_conj}, as well as the induced spectral sequence, do not appear in the literature to our knowledge, but are known to some experts in the field. They can be constructed by adapting the argument we gave for $\tau_{\le 2}\ko$, using the fact that $\eta^2\colon\Sigma^2\ko\to\ko$ factors through $\tau_{\le 2}\ko\to\ko$.

Another direction we speculate in is motivic and equivariant analogues of our main results. Motivic and $G$-equivariant stable homotopy theory are generalizations of stable homotopy theory in which, approximately speaking, topological spaces are replaced with smooth projective varieties over a base scheme, resp.\ topological spaces with a $G$-action, where $G$ is typically a compact Lie group. In these settings, there are analogues of $\ko$ and $\ku$ whose Adams and Atiyah--Hirzebruch spectral sequences behave somewhat similarly to the classical case (at least for $G = \Z/2$ in equivariant stable homotopy theory).
\begin{itemize}
    \item In motivic stable homotopy theory, the analogue of $\ku$ is $\mathit{kg\ell}$, the effective cover of $K$-theory, constructed by Voevodsky (see~\cite[\S 6.2]{Voe98} and~\cite[\S 2]{Voe02}). The analogue of $\ko$, denoted $\mathit{kq}$, is the very effective cover of Hermitian $K$-theory.\footnote{Spitzweck--{\O}stvær~\cite{SO12} constructed the very effective cover functor and Bachmann~\cite{Bac17} first applied it to $\mathit{KQ}$ to define $\mathit{kq}$. Ananyevskiy--R\"ondigs--{\O}stvær~\cite{ARO20} argued that $\mathit{kq}$ is a good motivic analogue of $\ko$, answering questions of Hill~\cite[\S 5.2]{Hil11} and Isaksen--Shkembi~\cite{IS11}.}
    \item In $\Z/2$-equivariant stable homotopy theory, the analogues of $\ku$, resp.\ $\ko$ are $\mathit{kr}$, resp.\ $\ko_{\Z/2}$, the equivariant connective covers (see~\cite[Proposition 3.3]{Lew95}) of Atiyah's Real $K$-theory~\cite{Ati66} resp.\ $\Z/2$-equivariant $\KO$-theory.\footnote{We also recommend Bruner--Greenlees~\cite{BG10} and Greenlees~\cite{Gre18} for more on $\mathit{kr}$ and its relatives.}
\end{itemize}
In all of these cases, there is a version of the Adams spectral sequence whose $E_2$-page is isomorphic to $\Ext$ over an algebra analogous to $\cE(1)$ or $\cA(1)$. Examples of such motivic and $\Z/2$-equivariant Adams spectral sequences appear in~\cite{Hil11, IS11, Orm11, OO13, Ghe18, Qui19, GIHR20, CKQ21, CQ21, DFNQY25, LPT25, Mor25, MPT25, Dha26, Mor25a}.

%
%
%
%
\begin{conj}
\label{equivariant_AHSS_to_Adams}
For $R$ an $E_\infty$-$\Z/2$-ring spectrum, let $\cat{Mod}_R$ denote the $\infty$-category of $R$-modules in genuine $\Z/2$-spectra. Then:
\begin{itemize}
    \item Let $\Post_\bullet^R\colon\cat{Mod}_R\to\Fil_R$ be the $\Z/2$-equivariant Atiyah--Hirzebruch filtration, defined analogously to \cref{defn:tau}.
    \item For $R = \mathit{kr}$, $\mathit{ko}_{\Z/2}$, let $\AdFilt^R_\bullet\colon\cat{Mod}_R\to\Fil_R$ be the $\Z/2$-equivariant, $H\underline{\Z/2}$-based, $R$-module Adams filtration, where $\underline{\Z/2}$ denotes the constant Mackey functor on $\Z/2$.\footnote{To construct this we also need maps $\ko_{\Z/2}\to H\underline{\Z/2}$ and $\mathit{kr}\to H\underline{\Z/2}$. For $\mathit{kr}$, Dugger~\cite[\S 6]{Dug05} constructs a map $\mathit{kr}\to H\underline\Z$, and the map we want is that map followed by mod $2$ reduction. For $\ko_{\Z/2}$, precompose with the map $c\colon\ko_{\Z/2}\to\mathit{kr}$ defined in~\cite[Definition 10.11]{GIHR20}.}
\end{itemize}
We conjecture that, for these choices of $R$, there are natural transformations of functors $\AtoA^R_{\mathrm{mot}}\colon \Post_\bullet^R\Rightarrow \Ad_{\widetilde g{}^R(\bullet)}^R$, where $\widetilde g{}^R$ is $g$ applied to the underlying spectrum of $R$, and such that, upon taking underlying spectra, one recovers the natural transformations of \cref{map_of_filtered_spectra}.
\end{conj}
For groups $G$ larger than $\Z/2$, very little is known about $G$-equivariant Adams spectral sequences for $\ko_G$ or $G$-equivariant analogues of $\mathit{kr}$. For this reason, we restrict to $G = \Z/2$ in \cref{equivariant_AHSS_to_Adams}.
\begin{conj}
\label{motivic_AHSS_to_Adams}
Let $F$ be either $\C$ or $\R$, and let $\mathit{kg\ell}$ and $\mathit{kq}$ denote the effective cover of $K$-theory, resp.\ the very effective cover of Hermitian $K$-theory, over the field $F$. When $R$ is one of these $E_\infty$-rings, let $\cat{Mod}_R$ be the $\infty$-category of $R$-modules in $F$-motivic spectra. Then:
\begin{itemize}
    \item Let $\Post_\bullet^R\colon \cat{Mod}_R\to \Fil_R$ be the motivic Atiyah--Hirzebruch filtration, defined similarly to \cref{defn:tau} (see also~\cite{DFNQY25} for some examples).
    \item Let $\AdFilt_\bullet^R\colon\cat{Mod}_R\to\Fil_R$ be the $H\Z/2$-based, $R$-module Adams filtration.
\end{itemize}
We conjecture that, for these choices of $F$ and $R$, there are natural transformations of functors $\AtoA^R_{\mathrm{eq}}\colon \Post_\bullet^R\to\AdFilt_{\widehat g(\bullet)}^R$, where $\widehat g$ is $g$ after Betti realization. We conjecture that after Betti realization $r_F$ to spectra ($F = \C$) or $\Z/2$-spectra ($F = \R$), these natural transformations recover $\AtoA^{r_\C(R)}$ (\cref{map_of_filtered_spectra}), resp.\ $\AtoA^{r_\R(R)}_{\mathrm{mot}}$ (\cref{equivariant_AHSS_to_Adams}).
\end{conj}
Balderrama--Culver--Quigley's motivic Lefschetz principle~\cite[Proposition 5.2.1]{BCQ25} and Dhankhar--Field--Nigam--Quigley--Yang's real motivic Lefschetz principle~\cite[Proposition 5.1]{DFNQY25} strongly suggest that \cref{motivic_AHSS_to_Adams} generalizes at least to other characteristic $0$ algebraically closed and real closed fields. It may generalize further.


\appendix

\section{Two definitions of twisted supercohomology are equivalent}
	\label{s:supercoh}
%

In this appendix, we will characterize twisted supercohomology in two different ways.
\begin{enumerate}
	\item The first way, which we give in \cref{cat_supercoh_defn}, uses the language of $2$-categories. From this
	point of view, supercohomology has been used to study fusion $2$-categories in the work of Johnson-Freyd--Yu~\cite{JFY21, JFY22},
	Décoppet--Yu~\cite{DY23a, DY25}, Décoppet~\cite[\S 4, Appendix]{Dec24},
	Décoppet--\allowbreak Huston--\allowbreak Johnson-Freyd--\allowbreak Nikshych--\allowbreak Penneys--\allowbreak Plavnik--\allowbreak Reutter--\allowbreak Yu~\cite{DHJNPPRY24}, Johnson-Freyd~\cite{JF25},
	Teixeira--Yu~\cite[\S 2]{TY25},
    Gu--Luo--Wang--Zhang~\cite{GLWZ26}, and Kobayashi--Prem--Yu~\cite[\S A.4]{KPY26a}.
	\item The second way we characterize twisted supercohomology, in \cref{htp_SH}, is homotopical, using the
	Ando--Blumberg--Gepner--Hopkins--Rezk~\cite{ABGHR14a, ABGHR14b} approach to twisted generalized cohomology.
	From this point of view, it is easier to compute twisted supercohomology groups.
\end{enumerate}
We will prove these two approaches are equivalent in \cref{SH_is_SH}. In particular, these two definitions give rise to two different twists labeled by $a\in H^1(X;\Z/2)$ and $b\in H^2(X;\Z/2)$ for a space $X$, and we prove that these two different formalisms for twisted supercohomology are the same.
\subsection{The categorical twist}
\label{ss:cattwist}
Throughout this section, when we write ``$2$-category'' we mean a bicategory, i.e.\ we do not require composition to be strictly associative, only up to coherent isomorphism.
%
\begin{defn}
Let $\sAlg$ denote the symmetric monoidal $2$-category whose objects are complex superalgebras, $1$-morphisms are
super bimodules, and $2$-morphisms are grading-preserving bimodule homomorphisms. The symmetric monoidal structure
comes from the tensor product with the Koszul sign rule.
\end{defn}
Given a symmetric monoidal $2$-category $\fC$, let $\Pic(\fC)$ denote the sub-$2$-category of $\fC$ of
$\otimes$-invertible objects and composition-invertible $1$- and $2$-morphisms. This is a \term{Picard
$2$-groupoid}, meaning a symmetric monoidal $2$-groupoid whose objects are all invertible with respect to the
monoidal product.

The following theorem is a consequence of Baez--Dolan's cobordism hypothesis~\cite{BD95}, proven by Lurie~\cite{Lur09}.
\begin{thm}
\label{O2_action}
There is a canonical $\rO_2$-action on $\Pic(\sAlg)$.
\end{thm}
Let $\mathbb G\coloneqq\Aut(\Pic(\sAlg))$; since $\Pic(\sAlg)$ is a Picard $2$-groupoid, $\mathbb G$ is a symmetric $3$-group, i.e.\ a Picard $3$-groupoid with a single isomorphism class of objects. Thus the classifying space $B\mathbb G$ is an $E_\infty$-algebra whose homotopy groups vanish outside of degrees $[1,3]$ (see \cref{stab_htpy_hyp}).

\Cref{O2_action} in particular means that there is a functor $\rho$ from the fundamental $\infty$-groupoid of $B\rO_2$ to $B\mathbb G$, which one can think of as a ``principal $\mathbb G$-bundle,'' and thus obtain an associated bundle $\mathcal X\to B\rO_2$ with ``fiber'' $\Pic(\sAlg)$. This has a few more concrete consequences.
\begin{cor}\hfill
\label{unspool_cobhyp}
\begin{enumerate}
    \item The isomorphism $\pi_1(B\rO_2)\cong\Z/2$ means that $\rho$ defines an involution on $\Pic(\sAlg)$. By~\cite[Example 2.4.12]{Lur09}, this is the duality involution. 

    From the $\Pic(\sAlg)$-bundle point of view, this says that the monodromy of the bundle $\mathcal X\to B\rO_2$ around the nontrivial loop is the duality involution.
    \item $\rho$ sends $\pi_2(B\rO_2)\cong\Z$ to the \term{Serre automorphism} $S$ of $\Pic(\sAlg)$~\cite[Example 2.4.14]{Lur09}.
    \item Since $\Pic(\sAlg)$ is a $2$-category, the action $\rho$ is trivial on homotopy groups in degrees $4$ and above. Since $\pi_3(B\rO_2)\cong 0$, the action is also trivial in degree $3$.
\end{enumerate}
\end{cor}
\begin{prop}[{Gunningham~\cite[\S 2.5, \S 2.6]{Gun16}}]
\label{serre_squared}
The restriction of $\rho$ across the map $r\colon \Spin_2\to\rO_2$, defined to be the double cover of the identity component, is canonically trivializable.
\end{prop}
Heuristically, the pullback of $\mathcal X$ to $B\Spin_2$ is canonically equivalent to the trivial bundle.
\begin{proof}
Since $\Pic(\sAlg)$ is a $2$-category, this question factors through the $3$-skeleton of $B\Spin_2$, which is $S^2$. Thus it suffices to find a trivialization for the generator of $\pi_2(B\Spin_2)\cong\Z$. The map $r\colon B\Spin_2\to B\SO_2$ induces multiplication by $2$ on $\pi_2$, since $\Spin_2\to\SO_2$ is a double cover; the map $B\SO_2\to B\rO_2$ is the universal cover, hence is an isomorphism on $\pi_2$. Thus the generator of $\pi_2(B\Spin_2)$ acts via $\rho$ as the \emph{square} of the Serre automorphism, and it suffices to show this acts trivially.

Let $A$ be a superalgebra. Gunningham~\cite[\S 2.6]{Gun16} shows that trivializing the square of the Serre automorphism acting on $A$ is naturally equivalent to specifying an isomorphism of bimodules from $A^!\coloneqq\Hom_{A^e}(A, A^e)$ to $A^*\coloneqq\Hom_{\C}(A, \C)$, where $A^e\coloneqq A\otimes A^{\mathrm{op}}$ (here all tensor products and Homs are in the super sense: see (\textit{ibid.}), Example 2.5). If $A$ is invertible, its inverse is its dual up to a contractible space of choices, and the dual in the Morita $2$-category of superalgebras is the opposite superalgebra (\textit{ibid.}). Thus for $A$ an invertible superalgebra, $A\otimes A^{\mathrm{op}}$ is canonically equivalent to $\C$ as superalgebras, and therefore $A^!$ is canonically equivalent to $A^*$.
\end{proof}
\begin{prop}[{Lurie~\cite[Example 2.4.15]{Lur09}}]
\label{ijk_defn}
The map $\rho$ factors as a composition
\begin{equation}
    \rho\colon B\rO_2 \xrightarrow i B\rO\xrightarrow J B\GL_1(\Sph) \xrightarrow k B\mathbb G,
\end{equation}
where $i$ is induced from the usual inclusion $\rO_2\inj\rO$, $J$ is the $J$-homomorphism, and both $J$ and $k$ are maps of $E_\infty$-spaces.
\end{prop}
Lurie does not discuss $B\mathbb G$ or $\rho$ specifically; the point is that, for any Picard $n$-groupoid $\cat G$, the $\rO_n$-action on $\cat G$ guaranteed by the cobordism hypothesis factors in an analogous way through the map $J$. See~\cite{DS26} for a detailed discussion of this point,
\begin{lem}
\label{3trunc}
Let $u_*\colon B\GL_1(\Sph)\to B\GL_1(\tau_{\le 2}\ko)$ be the map induced by the unit map $\Sph\to\tau_{\le 2}\ko$. Then the map $k$ from \cref{ijk_defn} factors through $u_*$.
\end{lem}
\begin{proof}
Because $\Pic(\sAlg)$ is a $2$-category, the $\infty$-group $\G$ is in fact a $3$-group: it has no non-identity $n$-morphisms for $n>3$. Thus $\pi_n(B\G) = 0$ for $n>3$, so the map $k\colon B\GL_1(\Sph)\to B\G$ factors through $\tau_{\le 3} B\GL_1(\Sph)$.

There is a natural identification of $\pi_n(B\GL_1(R))$ with $1$ (if $n = 0$), $\pi_0(R)^\times$ (if $n = 1$), or $\pi_{n-1}(R)$ (if $n>1$); applying this, the map $u_*$ is an isomorphism on homotopy groups in degrees $3$ and below, and kills all higher-degree homotopy groups of $B\GL_1(\Sph)$; thus $u_*$ is $3$-truncation, finishing the proof.
\end{proof}
\begin{defn}
\label{iota_q_defn}
Throughout this appendix, let $q\colon B\rO\to B\rO/B\Spin$ denote the cofiber, in the $\infty$-category of $E_\infty$-spaces, of the map $\iota\colon B\Spin\to B\rO$ induced from the map $\Spin_n\to\rO_n$, $n\ge 3$, which is the universal cover of the connected component of the identity.
\end{defn}
\begin{lem}
\label{trivspin}
The space of null-homotopies of the composition
\begin{equation}
    B\Spin \xrightarrow\iota B\rO \xrightarrow{k\circ J} B\mathbb G,
\end{equation}
as a map of $E_\infty$-spaces, is contractible (hence in particular nonempty).
\end{lem}
\begin{cor}
\label{rho_tilde_defn}
Descending $k\circ J$ across $q$, we obtain a map of $E_\infty$-spaces
\begin{equation}
    \widetilde\rho\colon B\rO/B\Spin\longrightarrow B\mathbb G.
\end{equation}
The unique nonzero element of $\pi_1(B\rO/B\Spin)\cong\Z/2$ acts on $\Pic(\sAlg)$ by duality, and the unique nonzero element of $\pi_2(B\rO/B\Spin)\cong\Z/2$ acts by the Serre automorphism, which squares to $0$ by \cref{serre_squared}.
\end{cor}
The actions of the generators follow from the description of the actions of $\pi_1$ and $\pi_2$ of $B\rO_2$ in \cref{unspool_cobhyp}, together with the fact that $\pi_i(B\rO_2)\to \pi_i(B\rO)\to \pi_i(B\rO/B\Spin)$ is surjective for $i = 1,2$.
\begin{proof}[Proof of \cref{trivspin}]
Since $B\Spin$, $B\rO$, and $B\mathbb G$ are grouplike $E_\infty$-spaces, it is equivalent to show that the space of null-homotopies of the  corresponding map of spectra
\begin{equation}
\label{bspin_to_bg}
    \mathit{bspin} \longrightarrow \mathit{bo} \longrightarrow b\mathbb G
\end{equation}
is contractible. Here, under the equivalence between connective spectra and grouplike $E_\infty$-spaces, $\mathit{bspin}$, $\mathit{bo}$, and $b\mathbb G$ are the connective spectra corresponding to $B\Spin$, $B\rO$, and $B\mathbb G$ respectively.\footnote{We warn the reader that some authors use $\mathit{bo}$ to denote the spectrum $\ko$. In our notation, $\mathit{bo}$ is the $0$-connected cover of $\ko$. We apologize for this notational inconsistency and hope that, since we do not use the notation $\mathit{bo}$ outside of the proof of \cref{trivspin}, it does not cause any confusion.}

Since $B\Spin$ is $3$-connected, $\mathit{bspin}$ is also $3$-connected. Since $B\mathbb G$ is $3$-truncated, so is $b\mathbb G$. Thus $[\mathit{bspin}, b\mathbb G] = 0$: all maps of the form~\eqref{bspin_to_bg} are null-homotopic. The space of null-homotopies of any such map is homotopy equivalent to the connective cover of $\mathrm{Map}_{\cat{Sp}}(\mathit{bspin}, \Sigma^{-1}b\mathbb G)$, which is contractible by an analogous connectivity argument.
\end{proof}
\begin{defn}\hfill
\label{AB_defn}
\begin{enumerate}
    \item Let $A\colon B\Z/2\to B\rO/B\Spin$ be the composition of the map $B\Z/2\to B\rO$ given by the tautological line bundle followed by $q\colon B\rO\to B\rO/B\Spin$.
    \item\label{B_defn} Let $B\colon B^2\Z/2\simeq B\SO/B\Spin\to B\rO/B\Spin$ be the map induced from the usual map $B\SO\to B\rO$.
\end{enumerate}
\end{defn}
\begin{prop}[{\cite[Proposition 1.37]{DY23}}]
\label{T_prop}
The map $S\coloneqq (A, B)\colon B\Z/2\times B^2\Z/2\to B\rO/B\Spin$ is a homotopy equivalence of spaces.
\end{prop}
Both $B\rO/B\Spin$ and $B\Z/2\times B^2\Z/2$ are $E_\infty$-spaces; however, they are not equivalent even as $E_1$-spaces (\textit{ibid.}, Lemma 1.41).

We are now ready for our first definition of twisted supercohomology.
\begin{defn}
\label{cat_supercoh_defn}
Let $X$ be a space, $a\in H^1(X;\Z/2)$, and $b\in H^2(X;\Z/2)$. The pair $(a,b)$ defines a map
\begin{subequations}
\begin{equation}
\label{before_assoc}
	f_{a,b}\colon X\xrightarrow{(a,b)} K(\Z/2, 1)\times K(\Z/2, 2) \xrightarrow[\eqref{T_prop}]{S} B\rO/B\Spin \xrightarrow[\eqref{rho_tilde_defn}]{\widetilde\rho} B\mathbb G.
\end{equation}
Let $P_k\to B\mathbb G$ be the associated bundle
\begin{equation}\label{assocbun}
	E\mathbb G \times_{\mathbb G} B^{k-2}\mathrm \lvert \mathrm N(\Pic(\sAlg))\rvert \longrightarrow B\mathbb G,
\end{equation}
\end{subequations}
where $\abs{\mathrm N(\bl)}$ denotes the geometric realization of the nerve. The \term{$(a,b)$-twisted categorical supercohomology of $X$} in
degree $k$, denoted $\SH_c^{k}(X; a, b)$, is the abelian group of homotopy classes of sections of $f_{a,b}^*P_k\to
X$.
\end{defn}
Since $\Pic(\sAlg)$ is symmetric monoidal, its nerve is an $E_\infty$-space, hence may be delooped arbitrarily many times as we do in~\eqref{assocbun}. In particular, for $0\le k<2$ we interpret $B^{k-2}$ as $\Omega^{2-k}$. The mixing construction over the higher group $\mathbb G$ refers to the (geometric realization of the nerve of the) action (higher) groupoid associated to the $\mathbb G$-action.
\begin{lem}
Let $X$ be a space and $b\in H^2(X;\Z/2)$. There is a natural isomorphism from $\SH_c^*(X; 0, b)$ to the $b$-twisted supercohomology of $X$ as defined by Décoppet~\cite[Definition 4.3]{Dec24}.
\end{lem}
This is essentially true by definition: indeed, we built \cref{cat_supercoh_defn} to generalize Décoppet's definition to the case $a\ne 0$.

Let $x\in H^1(B\Z/2;\Z/2)\cong\Z/2$ be the unique nonzero element.
\begin{exm}
\label{pinp_cat}
We will show $\SH_c^1(B\Z/2, x, 0)\cong \Z/2$. This is by definition the space of sections of the pullback $f_{x,0}^* P_1\to B\Z/2$, whose fiber is $\Omega\Pic(\sAlg)\simeq\Pic(\sVect)$. We will use the standard CW structure on $B\Z/2$ to make computations; the results are independent of the choice of CW structure.
The map $(x, 0)\colon B\Z/2\to K(\Z/2, 1)\times K(\Z/2, 2)\cong B\rO/B\Spin$ is an equivalence on the $1$-skeleton, and is trivial on the $2$-cell of $B\Z/2$.
The monodromy around the cells of $B\Z/2$ can be described as follows.
\begin{itemize}
    \item Around the $1$-cell $e^1$ of $B\Z/2$, the monodromy of $f_{x,0}^*P_1$ is the same as the monodromy of $P_1$ around the $1$-cell of $B\rO/B\Spin$, which by \cref{unspool_cobhyp} sends a morphism $\phi$ in $\Pic(\sVect)$ to $\phi^\vee$.
    \item Around the $2$-cell $e^2$ of $B\Z/2$, the monodromy of $f_{x,0}^*P_1$ is the identity.
    \item Since $f_{x,0}^*P_1$ is a bundle of Picard groupoids, which have no higher morphisms, all higher-dimension cells of $B\Z/2$ necessarily carry trivial monodromy.
\end{itemize}
We will build a section $s$ of $f_{x,0}^*P_1\to B\Z/2$ one cell at a time. Since this is a bundle whose fibers are Picard groupoids, our section is data over the $0$- and $1$-cells and a condition on the $2$-cell, then automatically extends to higher degrees.
\begin{enumerate}
    \setcounter{enumi}{-1}
    \item $B\Z/2$ has a single $0$-cell $e^0$, so over this cell, our section chooses a single object $s_0$ of $\Pic(\sVect)$, isomorphic to the even line $\C$ or the odd line $\Pi\C$. We may without loss of generality assume that $s_0$ \emph{is} $\C$ or $\Pi\C$: if not, we can choose an isomorphism $\alpha\colon s_0\overset\cong\to\C$ (or $\Pi\C$) and replace our section $s$ with $\alpha^{-1}\circ s\circ\alpha$, which is in the same homotopy class of sections. Thus we assume $s_0$ is $\C$ or $\Pi\C$ in the rest of this argument.
    \item To extend $s$ over the $1$-skeleton, we must choose an isomorphism $s_1\colon s_0\overset\cong\to s_0^\vee$ (since the monodromy around $e^1$ is the duality involution on $\Pic(\sVect)$). As of yet we can choose $s_1$ freely, giving a $\C^\times$ of choices for either $s_0\cong\C$ or $s_0\cong\Pi\C$ a priori. However, like for $s_0$, different choices can lead to homotopy-equivalent sections. Specifically, if one acts on $s_0$ by $\lambda\in\C^\times$, one necessarily obtains a homotopy-equivalent section, yet $s_1$ is multiplied by $\lambda^2$ (on the domain $s_0$ we act by $\lambda$, and on the codomain $s_0^\vee$ we act by $\lambda^\vee = \lambda^{-1}$; moving the latter action to the domain we get $\lambda^2$). Thus, if $s_1 = \lambda^2 \widetilde s_1$ for another choice $\widetilde s_1$, the data of $s_1$ and $\widetilde s_1$ define homotopy-equivalent sections over the $1$-skeleton. So for both $\C$ and $\Pi\C$, we have \emph{one} choice of homotopy class of section.
    \item The $2$-cell $e^2$ is attached by $\partial e^2 = 2e_1$, which gives the condition that ``squaring'' $s_1$ must equal the monodromy around the $2$-cell. By ``squaring'' we mean that the composition
    \begin{equation}
    \label{so_called_square}
        s_0\xrightarrow{s_1} s_0^\vee \xrightarrow{(s_1^\vee)^{-1}} s_0^{\vee\vee} \xrightarrow{\psi^{-1}} s_0
    \end{equation}
So the composition is the number $s_1^{-1}s_1 = 1$. We saw that the monodromy of $f_{x,0}^*P_1$ is also equal to $1$, so the condition on the $2$-cell is met by both $s_0 = \C$ and $s_0 = \Pi\C$.
\end{enumerate}
To summarize, there are two choices of section up to homotopy, so $\SH_c^1(B\Z/2, x, 0)$ is an abelian group of order $2$, hence is isomorphic to $\Z/2$.
\end{exm}
\begin{exm}
\label{pinm_cat}
We will show $\SH_c^1(B\Z/2, x, x^2)\cong 0$. Most of the story is the same, except that since the degree-$2$ part of the twist is nonzero, around the unique $2$-cell $e^2$ of $B\Z/2$, the monodromy of $f_{x,x^2}^*P_1$ equals the monodromy of the $2$-cell in $K(\Z/2, 2)$, which by \cref{unspool_cobhyp} is the grading operator $(-1)^F$.

Thus when we attempt to build a section $s$, the argument is exactly the same for $s_0$ and $s_1$: $s_0$ is isomorphic to $\C$ or $\Pi\C$, and given $s_0$ there is a $\C^\times$ of choices for $s_1$, all of which are homotopy equivalent. Now, though, the condition on the $2$-cell is that $s_1^{-1}s_1$ equals the monodromy around the $2$-cell. For $s_0$ an even line this is fine, as $(-1)^F$ acts by $1$ on even lines, so we have a unique homotopy class of section with $s_0\cong\C$. But if $s_0$ is an odd line, $(-1)^F$ acts as $-1$, and there is no choice of $s_1$ such that $s_1^{-1}s_1 = -1$. Thus there is a unique homotopy class of section, so $\SH_c^1(B\Z/2, x, x^2)\cong 0$.
\end{exm}
\begin{rem}
\label{other_two_cat_twists}
The reader may enjoy calculating $\SH_c^1(B\Z/2, 0, x^2)\cong\Z/4$ and $\SH_c^1(B\Z/2, 0, 0)\cong\Z/2\times\Z/2$ using arguments similar to those in \cref{pinp_cat,pinm_cat}.
\end{rem}

\subsection{The homotopical twist}
\label{ss:homtwist}
Now we introduce the other twist, using Ando--Blumberg--Gepner--Hopkins--Rezk's definition~\cite{ABGHR14a, ABGHR14b} of twisted generalized cohomology. Given an $E_\infty$-ring spectrum $R$, May~\cite[\S III.2]{MQRT77} defines a grouplike $E_\infty$-space $\GL_1(R)$, called the \term{$\infty$-group of units} of $R$, whose definition is manifestly natural in $R$; see~\cite[\S 1.2]{ABGHR14b} for a review. Ando--Blumberg--Gepner--Hopkins--Rezk~\cite[Definition 2.10, Corollary 2.14]{ABGHR14a} construct a natural transformation $M$ from the fundamental $\infty$-groupoid of $B\GL_1(R)$ to the category of $R$-module spectra.
\begin{defn}[{\cite[Definition 2.20]{ABGHR14a}}]
\label{thomdefn}
For $R$ an $E_\infty$-ring spectrum, the \term{$R$-module Thom spectrum} of a map $f\colon X\to B\GL_1(R)$, denoted $M^R f$, is the colimit of the $X$-shaped diagram
\begin{equation}
    X\xrightarrow f B\GL_1(R) \xrightarrow M \cat{Mod}_R.
\end{equation}
\end{defn}
Here we have identified $X$ and $B\GL_1(R)$ with their respective fundamental $\infty$-groupoids.

Using \cref{thomdefn}, Ando--Blumberg--Gepner~\cite[Definition 1.2]{ABG11} define $f$-twisted $R$-homology to be $\pi_*(M^Rf)$ and $f$-twisted $R$-cohomology to be $\pi_{-*}\mathrm{Map}_R(M^Rf, R)$. This is close to how we define twisted supercohomology, but because supercohomology is not a ring spectrum we modify their approach slightly.

We will work with the Postnikov $2$-truncation $\tau_{\le 2}\ko$ of the connective real $K$-theory spectrum $\ko$. $\tau_{\le 2}\ko$ inherits an $E_\infty$-ring structure from the standard $E_\infty$-ring structure on $\ko$.

\begin{defn}
\label{htp_SH}
Let $T\colon B\rO/B\Spin\to B\GL_1(\ko)$ be the map defined in~\cite[(1.36)]{DY23}. Let $X$ be a space, $a\in H^1(X;\Z/2)$, and $b\in H^2(X;\Z/2)$. Let $g_{a,b}$ be the composition
\begin{subequations}
\begin{equation}
    g_{a,b}\colon X\xrightarrow{(a,b)} K(\Z/2, 1)\times K(\Z/2, 2) \xrightarrow[\eqref{T_prop}]{S} B\rO/B\Spin \xrightarrow[\text{\cite{DY23}}]{T} B\GL_1(\ko) \xrightarrow{(\tau_{\le 2})_*} B\GL_1(\tau_{\le 2}\ko).
\end{equation}
Let $M^{\tau_{\le 2}\ko}g_{a,b}$ denote the $\tau_{\le 2}\ko$-module Thom spectrum of $g_{a,b}$ in the sense of \cref{thomdefn}. The \term{$(a,b)$-twisted homotopical supercohomology of $X$} in degree $k$ is
\begin{equation}
    \SH_h^k(X; a, b) \coloneqq\pi_{-k}(I_{\C^\times}M^{\tau_{\le 2}\ko}g_{a,b}).
\end{equation}
\end{subequations}
\end{defn}
\begin{rem}
\label{stop_at_ko}
With notation as in \cref{htp_SH}, let
\begin{equation}
    h_{a,b}\coloneqq T\circ S\circ(a,b)\colon X \longrightarrow B\GL_1(\ko).
\end{equation}
In other words, we use the twisting data $(a,b)$ to make a $\ko$-module Thom spectrum $M^\ko h_{a,b}$ analogously to how we made a $\tau_{\le 2}\ko$-module Thom spectrum in \cref{htp_SH}. By~\cite[Lemma 1.8]{DY23}, there is a $\tau_{\le 2}\ko$-module equivalence $M^\ko h_{a,b}\wedge_\ko \tau_{\le 2}\ko\simeq M^{\tau_{\le 2}\ko}g_{a,b}$.
\end{rem}
\begin{exm}
\label{pinp_h}
We will show $\SH_h^1(B\Z/2, x, 0)\cong \Z/2$. By definition and the universal property of $I_{\C^\times}$,
\begin{equation}
    \SH_h^1(B\Z/2, x, 0)\cong\Hom(\pi_1(M^{\tau_{\le 2}\ko}g_{x,0}), \C^\times).
\end{equation}
Thus it suffices to show that $\pi_1(M^{\tau_{\le 2}\ko}g_{x,0}) \cong\Z/2$. Because the map $T$ from \cref{htp_SH} factors through the map $B\GL_1(\MTSpin)\to B\GL_1(\ko)$ (see~\cite[(1.36)]{DY23}), $g_{a,b}$ factors through the map on $B\GL_1$ induced by the composition of maps of $E_\infty$-ring spectra $\MTSpin\to \ko\to\tau_{\le 2}\ko$. Thus, by~\cite[Lemma 1.8]{DY23},
\begin{equation}
    M^{\tau_{\le 2}\ko}g_{x,0} \simeq M^{\MTSpin}g_{x,0}\wedge_\MTSpin \tau_{\le 2}\ko,
\end{equation}
where $M^{\MTSpin}g_{x,0}$ is the $\MTSpin$-module Thom spectrum associated to the lift of the map $g_{a,b}$ to $B\GL_1(\MTSpin)$. The map $\MTSpin\to\tau_{\le 2}\ko$ is $2$-connected, because it is the composition of the Atiyah--Bott--Shapiro map $\MTSpin\to\ko$, which Anderson--Brown--Peterson~\cite{ABP67} showed is $7$-connected, with $\tau_{\le 2}$, which is $2$-connected by definition.
Thus the map
\begin{equation}\label{thom_change_base}
    M^{\MTSpin}g_{x,0}\to M^{\MTSpin}g_{x,0}\wedge_\MTSpin \tau_{\le 2}\ko\simeq M^{\tau_{\le 2}\ko}g_{a,b}
\end{equation}
is also $2$-connected (this can be proven, e.g., with the Künneth spectral sequence). Thus~\eqref{thom_change_base} is an isomorphism on $\pi_1$. There is an $\MTSpin$-module equivalence $M^{\MTSpin}g_{x,0}\simeq\MTPin^-$ (in this language, this is~\cite[Example 6.23]{MathSmith}, but the essence of the proof dates to~\cite[\S 7]{Pet68}), so
\begin{equation}
    \pi_1(M^{\tau_{\le 2}\ko}g_{x,0})\cong\pi_1(\MTPin^-)\cong \Omega_1^{\Pin^-},
\end{equation}
which is isomorphic to $\Z/2$~\cite[Theorem 5.1]{ABP69}.
\end{exm}
\begin{exm}
\label{pinm_h}
We will show $\SH_h^1(B\Z/2, x, x^2)\cong 0$. The proof is completely analogous to that of \cref{pinp_h}, but with $M^{\MTSpin}g_{x,x^2}$ in place of $M^{\MTSpin}g_{x,0}$. There is an $\MTSpin$-module equivalence $M^{\MTSpin}g_{x,x^2}\simeq\MTPin^+$ (in this language, this is~\cite[Example 6.23]{MathSmith}, but the real content of the proof is in~\cite[\S 8]{Sto88}), so $\pi_1(M^{\tau_{\le 2}\ko} g_{x,x^2})\cong\Omega_1^{\Pin^+}$, which vanishes~\cite[\S 2]{Gia73a}.
\end{exm}
Compare \cref{pinp_cat,pinm_cat}.
\subsection{The two twists are the same}
\label{ss:sametwist}
Our main result in this appendix is:
\begin{thm}
\label{SH_is_SH}
There is a natural isomorphism $\kappa\colon\SH_c^*(X, a, b)\overset\cong\to \SH_h^*(X, a, b)$.
\end{thm}
We will thus just say ``supercohomology'' and ``$\SH^*(X, a, b)$'' to refer to this object.

Before we prove \cref{SH_is_SH}, we will need a few lemmas.
\begin{defn}
Let $\mathit{pic}(\sAlg)$ denote the classifying spectrum of $\sAlg$ (see \cref{classsp}), i.e.\ a connective spectrum whose infinite loop space is equivalent as $E_\infty$-spaces to $\abs{\mathrm N(\Pic(\sAlg))}$.
\end{defn}
By \cref{stab_htpy_hyp}, this definition characterizes $\mathit{pic}(\sAlg)$ up to a contractible space of choices.
\begin{lem}
\label{what_is_picsalg}
There is a canonical homotopy equivalence of spectra $\mathit{pic}(\sAlg)\overset\simeq\to \Sigma^2 I_{\C^\times}(\tau_{\le 2}\ko)$.
\end{lem}
\begin{proof}
Debray--Gunningham~\cite[\S 4.3.5]{DG18} constructed a map $t\colon \mathit{pic}(\sAlg)\to \Sigma^2 I_{\C^\times}$ and showed that it is the connective cover, giving a canonical equivalence $\mathit{pic}(\sAlg)\to\tau_{\ge 0}\Sigma^2 I_{\C^\times}$. The Hurewicz map $u\colon \Sph\to\tau_{\le 2}\ko$ is an isomorphism on homotopy groups in degrees $2$ and below (this is standard for $\Sph\to\ko$, and to get $\tau_{\le 2}\ko$, use the universal property of Postnikov truncation). Thus by the universal property of $I_{\C^\times}$, the map $\Sigma^2 I_{\C^\times}(u)\colon \Sigma^2 I_{\C^\times}(\tau_{\le 2}\ko)\to \Sigma^2 I_{\C^\times}(\Sph) = \Sigma^2 I_{\C^\times}$ is an isomorphism on nonnegative-degree homotopy groups, hence a homotopy equivalence on connective covers. Let $\widetilde u$ denote this equivalence on connective covers. Finally, since $\pi_i(\tau_{\le 2}\ko)$ vanishes for $i\ge 3$, $\pi_j(\Sigma^2 I_{\C^\times}(\tau_{\le 2}\ko))$ vanishes for $j<0$, so $\Sigma^2 I_{\C^\times}(\tau_{\le 2}\ko)$ is canonically equivalent to its connective cover. Combining these facts, we have that $\widetilde u{}^{-1} \circ t$ defines a homotopy equivalence $\mathit{pic}(\sAlg)\overset\simeq\to \Sigma^2 I_{\C^\times}(\tau_{\le 2}\ko)$.
\end{proof}
\begin{rem}
\label{heres_xi}
Passing the $\tau_{\le 2}\ko$-module structure on $\Sigma^2 I_{\C^\times}(\tau_{\le 2}\ko)$ through
\cref{what_is_picsalg}, we get a canonical $\tau_{\le 2}\ko$-module structure on $\mathit{pic}(\sAlg)$. This gives rise to an $\infty$-group morphism $\xi\colon \GL_1(\tau_{\le 2}\ko)\to\G$ by restricting the action map $\tau_{\le 2}\ko\to\End(\mathrm{pic}(\sAlg))$ to the $\infty$-group of units and observing that the units act by automorphisms, rather than just endomorphisms.
\end{rem}
\begin{cor}
\label{trunc_str_grp}
The following diagram commutes up to homotopy.
\begin{equation}\begin{tikzcd}
	{B\rO/B\Spin} && {B\G.} \\
	& {B\GL_1(\tau_{\le 2}\ko)}
	\arrow["{\widetilde\rho}", "\eqref{rho_tilde_defn}"', from=1-1, to=1-3]
	\arrow["{S\circ (a,b)}"', "\eqref{htp_SH}", from=1-1, to=2-2]
	\arrow["\xi"', "\eqref{heres_xi}", from=2-2, to=1-3]
\end{tikzcd}\end{equation}
\end{cor}
Heuristically, \cref{trunc_str_grp} tells us that if $P\to X$ is a principal $\G$-bundle whose classifying map factors through $B\rO/B\Spin$, $P$ has a canonical reduction of structure group to $\GL_1(\tau_{\le 2}\ko)$. We have not provided definitions to make this heuristic parse, but as an idea, it helped us plot the proof of \cref{SH_is_SH}.
\begin{proof}
Since we constructed $\widetilde\rho$ in \cref{trivspin,rho_tilde_defn} by canonically descending the map $k\circ J\colon B\rO\to B\G$ across $B\rO/B\Spin$, it suffices to prove the corollary with $B\rO$ in place of $B\Spin$, $k\circ J$ in place of $\widetilde\rho$, and $S\circ (a,b)$ precomposed with the quotient $q\colon B\rO\to B\rO/B\Spin$. By~\cite[Theorem 1.13]{DY23} (see also (\textit{ibid.}, Remark 1.14) and ~\cite[Lemma IV.2.6]{MQRT77}) we may replace $S\circ (a,b)\circ q$ with
\begin{equation}
    B\rO\xrightarrow J B\GL_1(\Sph) \xrightarrow{u_*} B\GL_1(\tau_{\le 2}\ko),
\end{equation}
where $u_*$ is the map induced by the unit map $u\colon\Sph\to \tau_{\le 2}\ko$ of $E_\infty$-ring spectra. That is, we want to show that the diagram
\begin{equation}\begin{tikzcd}[row sep=0.5em]
	&& {B\GL_1(\tau_{\le 2}\ko)} \\
	{B\rO} & {B\GL_1(\Sph)} \\
	&& {B\G}
	\arrow["\xi", from=1-3, to=3-3]
	\arrow["J", from=2-1, to=2-2]
	\arrow["{u_*}", from=2-2, to=1-3]
	\arrow["k"', from=2-2, to=3-3]
\end{tikzcd}\end{equation}
commutes, meaning $k$ factors through $B\GL_1(\tau_{\le 2}\ko)$, which we proved in \cref{3trunc}: one can trace through the definition of $\xi$ to see that $k\simeq \xi\circ u_*$.
\end{proof}
To prove \cref{SH_is_SH}, we use May--Sigurdsson's~\cite{MS06} and Ando--Blumberg--Gepner's~\cite{ABG11} formalisms for parametrized stable homotopy theory. In these settings, the analogues of unbased, resp.\ based spaces are spaces over $B$ (i.e.\ spaces $X$ with maps $X\to B$), resp.\ \term{ex-spaces}, diagrams $B\to X\to B$ whose composition is the identity (see~\cite[\S 1.3]{MS06} for conventions). As with ordinary spectra, which are built from based spaces, May--Sigurdsson's parametrized spectra are built from ex-spaces.

For a space $X$, let $\Top_{/X}$ denote the category of CGWH spaces over $X$, $\mathscr K_X$ denote the category of ex-spaces over $X$, and $\mathscr S_X$ denote May--Sigurdsson's category of parametrized orthogonal spectra over $X$, defined in~\cite[\S 11.2]{MS06}. We will let $\cat{Sp}_{/X}$ denote Ando--Blumberg--Gepner's $\infty$-category of parametrized spectra over $X$~\cite[Definition 3.3]{ABG11}, which is the $\infty$-category of functors $X^{\mathrm{op}}\to\cat{Sp}$.\footnote{\label{noop}Here we implicitly identify $X$ with its fundamental $\infty$-groupoid. In particular, $X^{\mathrm{op}}\simeq X$ (see~\cite[Remark 3.2]{ABG11}); we use $X^{\mathrm{op}}$ here following Ando--Blumberg--Gepner~\cite{ABG11}, who parametrize over more general things than spaces.}
\begin{defn}
\label{append_basepoint}
Given a space over $B$ $f\colon X\to B$, let $X_+$ denote the ex-space $X\amalg B$, where the map $B\to X\amalg B$ is the inclusion of the second factor and the map $X\amalg B\to B$ is $f\amalg\id$.
\end{defn}
\begin{rem}
\label{un_or_reduced}
In particular, if $E\in\mathscr S_B$ and $X$ is a space over $B$, $E^*(X)$ denotes the \emph{unreduced} $E$-cohomology of $X$, which is naturally isomorphic to the \emph{reduced} $E$-cohomology of $X_+$ (see \cite[\S 20.1]{MS06}).\footnote{The notation $X_+$ in (\textit{ibid.}, Proof of Theorem 20.1.5) should agree with \cref{append_basepoint}. This is in conflict with (\textit{ibid.}, Notation 5.1.3), but this appears to be a typo: using \cref{append_basepoint} is necessary for the statements to type-check.}
\end{rem}
Let $\Pi$ be a locally compact, well-based topological group, $P\to X$ be a principal $\Pi$-bundle, and $E$ be a spectrum with $\Pi$-action. In (\textit{ibid.}, \S 15.4), May--Sigurdsson construct an associated parametrized spectrum $P\times_{\Pi} E$ to a principal $\Pi$-bundle $P\to X$ and a spectrum $E$ with a $\Pi$-action (i.e.\ a $\Pi$-spectrum indexed on a trivial $\Pi$-universe in the sense of Lewis--May~\cite[Definition I.2.1]{LMSM86}\footnote{Lewis--May assume that the group $G$ that acts is a compact Lie group; as May--Sigurdsson~\cite[\S 15.4]{MS06} note, since we are using the trivial $G$-universe, this assumption can be removed.}). Then, in (\textit{ibid.}, \S 22.1), May--Sigurdsson remark that it follows from the definitions that if $E_P$ denotes the parametrized cohomology theory over $X$ associated to $P\times_\Pi E$ (\cite[Definition 20.2.4]{MS06}), then the (unreduced) cohomology $E_P^0(X)$ is naturally isomorphic to the abelian group of homotopy classes of sections of $P\times_\Pi \Omega^\infty E$. 

Our $\G$ plays the role of $\Pi$, but we do not know whether it is a topological group, as it is merely a grouplike $E_1$-space. We fix this as follows.\footnote{\label{same_BOBSpin}An alternate approach would be to work with the grouplike $E_1$-space $\Omega(B\rO/B\Spin)$, since the map~\eqref{before_assoc} factors through $B\rO/B\Spin$. Following Maycock Parker~\cite{May88}, $\Omega(B\rO/B\Spin)$ is modeled by the projective orthogonal group of the super real Hilbert space $\ell^2\oplus\Pi\ell^2$, where orthogonal transformations send homogeneous elements to homogeneous elements, yet can change their gradings. This group is also not locally compact.}
Let $\mathbbold\Omega\colon\cat{sSet}\to\cat{sGrp}$ denote \term{Kan's loop group functor}~\cite{Kan58a, Kan58b}, where $\cat{sSet}$ and $\cat{sGrp}$ are the categories of simplicial sets and simplicial groups, respectively. If $\abs{\bl}$ denotes geometric realization, then for any simplicial set $X$, there is a natural homotopy equivalence of $E_1$-spaces $\abs{\mathbbold\Omega X}\overset\simeq\to \Omega\abs X$. Moreover, $\abs{\mathbbold\Omega X}$, being the geometric realization of a simplicial group, is a compactly generated, Hausdorff, paracompact topological group. Thus we may replace $\G$ with $\abs{\mathbbold\Omega(\mathrm{Sing}(B\mathbb G))}$, where $\mathrm{Sing}\colon\cat{Top}\to\cat{sSet}$ denotes the singular complex functor, to obtain a topological group equivalent to $\G$ as $E_1$-spaces. In the rest of this section, we will implicitly use this construction to make $\G$ into a topological group, which one can check is well-based.

However, the topological group $\G$ we constructed in this way is not locally compact. To surmount this, we check that for the applications we have in mind, May--Sigurdsson's local compactness hypothesis is unnecessary.
\begin{lem}
\label{CGWH}
Let $\Pi$ be a compactly generated, Hausdorff, paracompact topological group (e.g.\ the geometric realization of a simplicial group). Let $P\to X$ be a principal $\Pi$-bundle of CGWH spaces and $E$ be an orthogonal $\Omega$-spectrum with $\Pi$-action. Then:
\begin{enumerate}
    \item applying the associated bundle construction $P\times_\Pi\bl$ levelwise to $E$ produces a spectrum over $X$ in the sense of~\cite[\S 11.2]{MS06}, generalizing the construction in (\textit{ibid.}, \S 15.4), and
    \item if $E_P^*$ denotes the parametrized cohomology theory over $X$ associated to $P\times_\Pi E$ in the sense of \cite[Definition 20.2.4]{MS06}, then the unreduced cohomology $E_P^0(X)$ is naturally isomorphic to the abelian group of homotopy classes of sections of $P\times_\Pi \Omega^\infty E$.
\end{enumerate}
\end{lem}
\begin{rem}
May--Sigurdsson~\cite[\S 3.3, \S 15.4]{MS06} consider a much more general setting of a group $\Gamma$ containing $\Pi$ as a normal subgroup, $G\coloneqq \Gamma/\Pi$, and a fiber bundle of $G$-spaces $F\to E\to B$. For us, $\Pi=\Gamma$, $G = 1$, and $E\to B$ is $\id\colon B\to B$.
\end{rem}
\begin{proof}[Proof of \cref{CGWH}]
For the most part, the proof amounts to tracing through May--Sigurdsson's definitions and theorems and observing that the ones we need do not actually require $\Pi$ to be locally compact. We highlight the key parts. 
\begin{enumerate}
    \item In~\cite[Proposition 1.6.1]{MS06} and the surrounding discussion, May--Sigurdsson explain that they restrict to locally compact $\Pi$ so as to avoid confusion as to differing definitions of a $\Pi$-space. For us, the assumptions of the lemma hand us all of the $\Pi$-spaces we need as CW complexes, so this ambiguity does not apply.
    \item Paralleling~\cite[\S 3.3]{MS06}, we want to show that $P\times_\Pi\bl$ defines a functor from based $\Pi$-spaces to ex-spaces over $X$, both within $\cat{CGWH}$. This follows from standard properties of quotients of CGWH spaces, along with the fact that $\Pi$ acts freely on $P\times F$ to get Hausdorffness of the quotient.
    \item The spectrum-level associated bundle construction~\cite[Notation 15.4.1]{MS06} goes through in exactly the same way.
\end{enumerate}
Once we have our parametrized spectrum $P\times_\Pi E$ over $X$, the rest is not any different from the case May--Sigurdsson considered, but we will be explicit. Let $[\bl,\bl]_{\fC}$ denote the homotopy classes of maps in a model category $\fC$,\footnote{In all cases we use, May--Sigurdsson construct the model structure on $\fC$. Specifically, see \cite[\S 6.2]{MS06} for $\Top_{/X}$ and $\mathscr K_X$, and (\textit{ibid.}, \S 12.3) for $\mathscr S_X$.} and consider $X$ as a space over itself by $\id\colon X\to X$. Then by~\cite[Definition 20.2.4]{MS06} and \cref{un_or_reduced},
\begin{subequations}
\begin{align}
    E_P^0(X) &\cong [\Sigma_X^\infty X_+, P\times_\Pi E]_{\mathscr S_X},
\intertext{where $\Sigma_X^\infty$ and $\Omega_X^\infty$ are defined in (\textit{ibid.}, Definition 11.24) and $(\bl)_+$ is \cref{append_basepoint}. $\Sigma_X^\infty$ and $\Omega_X^\infty$ are shown to partake in a Quillen adjunction in~\cite[Proposition 12.6.2]{MS06}, so we can rewrite $E_P^0(X)$ as}
    &\cong [X_+, \Omega_X^\infty(P\times_\Pi E)]_{\mathscr K_X}\\
    &\cong [X, P\times_\Pi \Omega^\infty E]_{\Top_{/X}},
    \label{finally_section}
\end{align}
\end{subequations}
by the definition of $\Omega_B^\infty$ as the $0^{\mathrm{th}}$ space of a parametrized $\Omega$-spectrum and the fact that May--Sigurdsson~\cite[Notation 15.41]{MS06} define $P\times_{\Pi}\bl$ levelwise. To finish the proof, observe that~\eqref{finally_section} is by definition the abelian group of homotopy classes of sections of the bundle of spaces $P\times_\Pi\Omega^\infty E\to X$.
\end{proof}
\begin{lem}
\label{who_is_unit}
Let $X$ be a CGWH space and $r\colon X\to\pt$ be the crush map.
\begin{enumerate}
    \item\label{MS_unit} The parametrized spectrum $\Sigma_X^\infty X_+$ is the unit with respect to May--Sigurdsson's symmetric monoidal structure~\cite[Theorem 11.2.5]{MS06} on $\mathscr S_X$.
    \item\label{ABG_unit} $r^*\Sph$ is the unit with respect to Ando--Blumberg--Gepner's symmetric monoidal structure~\cite[\S 6.1]{ABG11} on $\cat{Sp}_{/X}$.
\end{enumerate}
\end{lem}
Given a map $f\colon Y\to X$ and $E\in\cat{Sp}_{/X}$, i.e.\ a functor $X\overset E\to \cat{Sp}$, the pullback $f^*E$ is by definition the composition $Y\overset f\to X\overset E\to\cat{Sp}$. Thus in \cref{who_is_unit}, $r^*\Sph$ is the functor $X\to\cat{Sp}$ sending $x\mapsto \Sigma_+^\infty x$.
\begin{proof} 
Part~\eqref{ABG_unit} follows directly from the fact that $r^*$ is symmetric monoidal, which is proven in \cite[Proposition 6.8]{ABG11}. Thus in the rest of the proof we focus on part~\eqref{MS_unit}.

May--Sigurdsson describe the monoidal unit $S_X$ of $\mathscr S_X$ in~\cite[Definition 11.2.2]{MS06} (see also (\textit{ibid.}, Theorem 11.2.5)).\footnote{In \textit{loc.\ cit.}, May--Sigurdsson actually describe $S_X$ as an \term{$\mathscr I$-space}, i.e.\ a functor from the category $\mathscr I$ defined in~\cite[Example 4.4]{MMSS:2001} or~\cite[Definition II.2.1]{MM02} (set $G = \set 1$) to $\Top$. This is okay: May--Sigurdsson~\cite[Definition 11.2.3]{MS06} define their parametrized spectra to be $\mathscr I$-spaces together with slightly more data. Here the extra data form an evident homotopy equivalence compatible with our identifications.}
So we need to identify $\Sigma_X^\infty X_+$ and $S_X$. Let $V$ be a finite-dimensional inner product space and $S^V$ be the one-point compactification of $V$. Unwinding~\cite[(11.2.1), Definition 11.2.4]{MS06}, it suffices to identify the ex-spaces $S_X^V\wedge_X X_+$ and $S_X^V$, where $S_X^V$ is the ex-space $S^V\times X$ with the map $X\to S^V\times X$ given by $x\mapsto (x, \infty)$ and $S^V\times X\to X$ is projection onto the second factor. At this point, one can directly check that there is a natural homeomorphism of ex-spaces $S_X^V\wedge_X X_+\simeq S_X^V$ using the definition of $\wedge_X$~\cite[Definition 1.3.4]{MS06}.
\end{proof}

\begin{proof}[Proof of \cref{SH_is_SH}]
%
Let $E_k\coloneqq \Sigma^{k-2} I_{\C^\times}(\tau_{\le 2}\ko)$. By definition, $\G$ acts on the Picard $2$-groupoid $\Pic(\sAlg)$, hence also on $\mathit{pic}(\sAlg)$, and we showed $\mathit{pic}(\sAlg)\simeq \Sigma^{2-k}E_k$ in \cref{what_is_picsalg}. Thus we obtain the structure of a spectrum with $\G$-action on $E_k$. Let $\pi\colon f_{a,b}^*E\G\to X$ be pullback of the universal principal $\G$-bundle to $X$ by the map $f_{a,b}$ from~\eqref{before_assoc} and $E_{k,\pi}\coloneqq (f_{a,b}^*E\G)\times_\G E_k$ be the associated bundle we constructed in \cref{CGWH} following May--Sigurdsson~\cite[\S 11.2]{MS06}.


Let $E_{k,\pi}^*$ be the parametrized cohomology theory for spaces over $X$ associated to $E_{k,\pi}$ as in~\cite[\S 22.1]{MS06}. Then, by \cref{CGWH}, $E_{k,\pi}^0(X)$ is naturally isomorphic to the abelian group of homotopy classes of sections of $f_{a,b}^*E\G\times_\G \Omega^\infty E_k$, which by \cref{cat_supercoh_defn} is $\SH_c^k(X, a, b)$.

By~\cite[Definition 20.2.4]{MS06}, $E_{k,\pi}^0(X)\cong \pi_0 r_* F_X(\Sigma_X^\infty X_+, P\times_{\G} E_k)$ where $r\colon X\to\pt$ is the crush map, and $F_X$ and $\Sigma_X^\infty$ are as defined in~\cite[Proof of Theorem 11.2.5, resp.\ Definition 11.2.4]{MS06} (set $G = \set 1$). Ando--Blumberg--Gepner~\cite[Theorem B.2]{ABG11} prove an equivalence of \term{Wirthmüller contexts} (\cite[Definition 2.12(iii)]{FHM03}) between $\cat{Sp}_{/X}$ and $\mathscr S_X$ -- specifically, from Ando--Blumberg--Gepner's Wirthmuller context on $\cat{Sp}_{/X}$ in~\cite[Corollary 6.9]{ABG11} to May--Sigurdsson's on $\mathscr S_X$ (see~\cite[\S 2.2]{MS06}). This equivalence amounts to the data of a natural equivalence of symmetric monoidal $\infty$-categories $W_X\colon\mathscr S_X \overset\simeq\to\cat{Sp}_{/X}$ commuting with natural isomorphisms between the respective base-change functors $f^*$, $f_*$, and $f_!$ for a map of spaces $f\colon X\to Y$, as well as the parametrized mapping spectrum $F_X$. Invoking this equivalence, we now have a natural isomorphism
\begin{equation}
    \SH_c^k(X, a, b) \xrightarrow\cong \pi_0 r_* F_X( W_X(\Sigma_X^\infty X_+), P\times_\G E_k),
\end{equation}
where on the right-hand side, $r_*$, $F_X$, and the associated bundle construction are in the Ando--Blumberg--Gepner sense. We have not discussed the associated bundle in this sense, and we also need to compute $W_X(\Sigma_X^\infty X_+)$; let us do both now.
Recall that $\cat{Sp}_{/X}$ is the $\infty$-category of functors $X\to\cat{Sp}$ (see Footnote~\ref{noop}).
\begin{itemize}
    \item By part~\eqref{MS_unit} of \cref{who_is_unit}, $\Sigma_X^\infty X_+$ is the monoidal unit of $\mathscr S_X$. Since $W_X$ is symmetric monoidal, $W_X(\Sigma_X^\infty X_+)$ is the monoidal unit of $\cat{Sp}_{/X}$, which is $r^*\Sph$ (\cref{who_is_unit}, part~\eqref{ABG_unit}).
    \item The associated bundle $P\times_\G E_k$ is defined as follows: the $\G$-action on $E_k$ defines a functor $B\G\to\cat{Sp}$ whose value at the basepoint of $B\G$ is $E_k$. The associated bundle of spectra to $P$ is the composition $X\to B\G\to\cat{Sp}$, where the first map is the classifying map for $P$. To compare with May--Sigurdsson's associated bundle, which was defined levelwise, consider the monodromy of $W_{B\G}(E\G\times_{\G} E_k)$. This is by definition a map $B\G\to\cat{Sp}$ which agrees with the above description of the associated bundle.
\end{itemize}
By \cref{trunc_str_grp}, we can ``change structure group'' to $\GL_1(\tau_{\le 2}\ko)$. More precisely, we can replace the associated bundle $P\times_\G E_k$ with the functor
\begin{equation}
\label{paramphi}
    \Phi\colon X \xrightarrow{S\circ(a,b)} B\GL_1(\tau_{\le 2}\ko) \to \Mod_{\tau_{\le 2}\ko} \xrightarrow{(\bl)^\vee} \Mod_{\tau_{\le 2}\ko},
\end{equation}
where $(\bl)^\vee$ is $\Map_{\Sph}(\bl, \Sigma^k I_{\C^\times})$. (That is, the composition of $\Phi$ with the forgetful functor from $\tau_{\le 2}\ko$-modules to spectra is naturally equivalent to $P\times_\G E_k$.)

Let $\widetilde\Phi$ be the composition of only the first two arrows in~\eqref{paramphi}, so $\Phi = F_X(\widetilde\Phi, r^*\Sigma^k I_{\C^\times})$, where $r\colon X\to\pt$ is the crush map. Therefore we can invoke the tensor-hom adjunction:
\begin{equation}
\begin{aligned}
    \pi_0 r_* F_X(r^*\Sph, P\times_\G E_k) &\cong \pi_0 r_* F_X(r^*\Sph, F_X(\widetilde\Phi, r^* \Sigma^k I_{\C^\times}))\\
    &\cong \pi_0 r_* F_X(r^*\Sph\wedge_X \widetilde\Phi, r^* \Sigma^k I_{\C^\times})\\
    &\cong \pi_0 r_* F_X(\widetilde\Phi, r^*\Sigma^k I_{\C^\times}),
\end{aligned}
\end{equation}
since $r^*\Sph$ is the monoidal unit (\cref{who_is_unit}, part~\eqref{ABG_unit}). Now we can invoke the $(r_!, r^*)$ adjunction~\cite[\S 5.4]{ABG11}:
\begin{equation}
\label{finally_hSH}
    \pi_0 r_* F_X(\widetilde\Phi, r^*\Sigma^k I_{\C^\times})\cong \pi_0\Map(r_!\widetilde\Phi, \Sigma^k I_{\C^\times}).
\end{equation}
By~\cite[Proposition 4.23]{ABG11} (where $f$ is the trivial family with fiber $\pt$), $r_!\widetilde\Phi$ is the $\tau_{\le 2}\ko$-module Thom spectrum of the map $S\circ(a,b)\colon X\to B\GL_1(\tau_{\le 2}\ko)$, so~\eqref{finally_hSH} is naturally isomorphic to $\SH_h^k(X, a, b)$ by \cref{htp_SH}.
\end{proof}
\begin{rem}
Many applications of twisted supercohomology only need twists that factor through $K(\Z/2, 2)\simeq B\SO/B\Spin\to B\rO/B\Spin$, i.e.\ for which the degree-$1$ part of the twist is zero. This includes the applications to fusion $2$-categories discussed in~\cite{JFY21, JFY22, DY23a, Dec24, DHJNPPRY24, JF25, TY25, Debray:2025kfg,DY25, DYY3, GLWZ26, KPY26a}, though see~\cite[\S 5]{Debray:2025kfg}.

The theorems in this appendix are quite a bit easier restricted to this class of twists because one can take advantage of uniqueness theorems similar to Antieau--Gepner--Gómez's~\cite[Theorem 1.1]{AGG14}, showing there is a unique nontrivial map $K(\Z/2, 2)\to B\GL_1(\KO)$ up to homotopy. By contrast, there is a $\Z/2$-torsor of ways to extend to twists over $B\rO/B\Spin$. This torsor manifests in several different ways, including but not limited to:
\begin{enumerate}
    \item The choice of whether to put a $(BG_b, a, b)$-twisted spin structure on the tangent bundle or stable normal bundle of a manifold, hence whether to consider tangential or normal Thom spectra $\mathit{MTH}$, resp.\ $\mathit{MH}$.
    \item The choice of whether a $(BG_b, a, b)$-twisted spin structure on a vector bundle $V$ includes a trivialization of $w_2(V) + b$ or of $w_2(V) + a^2 + b$.
    \item The choice of $\Pin^+$ or $\Pin^-$, and relatedly one's sign conventions for Clifford algebras. 
    \item The choice of homotopy equivalence $B\rO/B\Spin\simeq B\Z/2\times B^2\Z/2$ (see~\cite[Remark 1.40]{DY23}).
    \item The set of homotopy classes of $E_1$-equivalences $\PO(\cH)\simeq \Omega(B\rO/B\Spin)$ (see Footnote~\ref{same_BOBSpin}).
\end{enumerate}
See also~\cite[\S 1]{DS26} for some related conventions.

So why did we work at this level of generality in this paper? Though we did not need degree-one twists for the main results of~\cite{Debray:2025kfg}, we conjectured there that our results should extend to this more general class of twists, corresponding to not-necessarily-unitary symmetries of fermionic QFTs. The tools we built in this paper will hopefully be useful in studying that conjecture.
\end{rem}

\section{Ext and \texorpdfstring{$\cQ$}{Q} computations}
\label{ext_Q_computations}

In this appendix, we compute $\cQ(M)$ for each $\cA(1)$-module $M$ in \cref{names}. Our main theorem in this section is \cref{Q_ext_calc}, which we stated heuristically in the main text as \cref{heuristic_Q_calc}.

Before we prove that theorem, we prove some lemmas computing Ext groups of some $\cA(1)$-modules that we will need. Neither the statements nor the proofs of these lemmas are necessary to understand the statements of the Ext charts, just to give a complete account of the calculations.

The following lemma is well-known, and in any case can be checked directly.
\begin{lem}
\label{bow_is_tens}
There is a unique isomorphism $f\colon B\to\cA(1)\otimes_{\cE(1)} \Z/2$ of graded vector spaces. It is also an $\cA(1)$-module isomorphism and a $\Z/2$-coalgebra isomorphism.
\end{lem}
The following lemma is also standard; for example, it follows by combining the $\E$-module and $\Z/2$-algebra structures on $\Ext_{\cA(1)}(B)$, described in~\cite[\S 4.5]{Bay94} and~\cite[Example 4.5.6]{BC18}, respectively.
\begin{lem}
\label{bow_alg}
There is an $\E$-algebra isomorphism $\Ext_{\cA(1)}(B)\cong\E[b]/(h_1, v+h_0b^2, w+b^4)$, where $b\in\Ext^{1,3}$.
\end{lem}
By \cref{bow_is_tens}, the algebra structure can be thought of as arising either from the coalgebra structure on $B$ or from $\Ext_{\cE(1)}(\Z/2)$ by the change-of-rings theorem (see, e.g., \cite[\S 4.5]{BC18}).
\begin{rem}
\label{B_tensor_stuff}
We can use \cref{bow_is_tens,bow_alg} to compute the $\E$-module structure on $\Ext_{\cA(1)}(B\otimes M)$ for an $\cA(1)$-module $M$ in terms of the $\Ext_{\cE(1)}(\Z/2)\cong\Z/2[h_0, b]$-module structure on $\Ext_{\cE(1)}(M)$. Specifically, $\Ext_{\cA(1)}(B\otimes M)$ is naturally isomorphic to the $\E$-module with underlying bigraded vector space $\Ext_{\cE(1)}(M)$ and $\E$-action induced from the $\E$-algebra structure on $\Ext_{\cE(1)}(\Z/2)$ from \cref{bow_alg}.

We apply this to the following examples.
\begin{itemize}
    \item $R_6\coloneqq P\otimes B$, so the $\E$-module structure on $\Ext_{\cA(1)}(R_6)$ follows from the $\Ext_{\cE(1)}(\Z/2)$-module structure on $\Ext_{\cE(1)}(P)$ calculated in~\cite[Proposition 4.47]{Deb21}.
    \item $B\otimes\uQ$, using the $\Ext_{\cE(1)}(\Z/2)$-module structure on $\Ext_{\cE(1)}(\uQ)$ calculated by Adams--Priddy~\cite[\S 3]{AP76}.
    \qedhere
\end{itemize}
\end{rem}

Next we calculate $\Ext(R_1)$. Since $R_1$ is a coalgebra (\cref{coalg_remark}), $\Ext(R_1)$ is an $\E$-algebra. The bigraded vector space structure, together with the actions of $h_0$ and $h_1$, appears in Gitler--Mahowald--Milgram~\cite[\S 2]{GMM68}; the rest of the module structure is implicit in their computation.
\begin{lem}
\label{E5_ext}
There is an $\E$-algebra isomorphism
\begin{equation}
    \Ext_{\cA(1)}(R_1) \overset\cong\longrightarrow
    \E[q, p_i: i\ge 0]/(
        h_0p_0 + h_1, h_0q + v, q^2+w, h_1p_i, qp_i + h_0^2 p_{i+1}, p_ip_j\text{ for all }i,j\ge 0
    ),
\end{equation}
where $q\in\Ext^{2,6}$ and $p_i\in\Ext^{0,4i+1}$.
\end{lem}
\begin{proof}
The key observation is that the claimed algebra structure is the unique $\Z^2$-graded $\E$-algebra structure
compatible with the underlying $\Z^2$-graded $\E$-module structure, so it suffices to compute the $\E$-module
structure. Gitler--Mahowald--Milgram~\cite[\S 2]{GMM68} compute $\Ext_{\cA(1)}^{*,*}(R_1)$ as an
$\Z/2[h_0,h_1]/(h_0h_1, h_1^3)$-module; to determine the actions of $v$ and $w$, one may follow the arguments of
Campbell~\cite[\S 7.8]{Cam17} and Beaudry--Campbell~\cite[Figures 25 and 26]{BC18}, who study $\Ext_{\cA(1)}(R_1)$ using the long exact sequence associated to a short exact sequence of $\cA(1)$-modules. This long exact sequence commutes with the action of elements of $\E$, and this suffices to prove the claimed $v$- and
$w$-actions.
\end{proof}

\begin{figure}[htbp]
\includegraphics{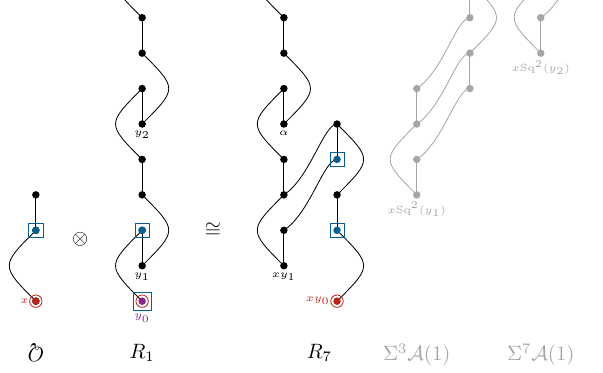}
\caption{A picture proof of \cref{R6}, part~\eqref{R6_tenspart}: that $\uQ\otimes R_1\simeq R_7$. Here
$\alpha\coloneqq \Sq^2(x)\Sq^2(y_1) + xy_2$. The circled classes generate the $Q_0$-Margolis homology of each
$\cA(1)$-module, and the squared classes generate the $Q_1$-Margolis homology.}
\label{R6_fig}
\end{figure}

\begin{lem}\label{R6}\hfill
\begin{enumerate}
	\item\label{R6_unique} Up to isomorphism, there is a unique nonsplit $\cA(1)$-module extension of $\uQ$ by
	$\Sigma R_3$. Call this $\cA(1)$-module $R_7$.
	\item\label{R6_tenspart} There is an $\cA(1)$-module isomorphism
	\begin{equation}
	\label{R6_as_tensor}
		\uQ\otimes R_1 \cong R_7 \oplus \bigoplus_{i\ge 0} \Sigma^{4i+3}\cA(1).
	\end{equation}
	\end{enumerate}
\end{lem}
See \cref{R6_fig} for a picture of $R_7$ and of the isomorphism~\eqref{R6_as_tensor}.
\begin{proof}[Proof sketch]
Part~\eqref{R6_unique} can be checked directly using the Adem relations. For part~\eqref{R6_tenspart}, we need the technology of \term{Margolis homology}~\cite{AM71, Mar83}, which we briefly review. Let $Q_0\coloneqq\Sq^1$ and $Q_1\coloneqq\Sq^1\Sq^2 + \Sq^2\Sq^1$; then for $i = 0,1$, $Q_i$ squares to $0$ and acts naturally on any $\cA(1)$-module, so we may regard an $\cA(1)$-module $M$ as a chain complex with differential $Q_i$. The resulting homology groups are called the \term{Margolis homology} of $M$.

Back to the show. To prove~\eqref{R6_tenspart}, one can
embed $R_7$ as a submodule of $\uQ\otimes R_1$, then verify that this embedding is an isomorphism on $Q_0$- and
$Q_1$-Margolis homology. We draw this embedding and Margolis homology computation in \cref{R6_fig}. By a theorem of Wall~\cite[Theorem 2]{Wal62} (see also~\cite[Theorem 4.1]{AM71}),
this implies the embedding $R_7\inj \uQ\otimes R_1$ is an $\cA(1)$-module isomorphism after quotienting out all free summands. To finish the proof of~\eqref{R6_as_tensor}, we need to find the locations of these free summands. This is a standard Poincaré polynomial argument similar to those in Freed--Hopkins~\cite[\S\S D.4--D.6]{FH21}.
\end{proof}

\begin{figure}[htbp]
\includegraphics{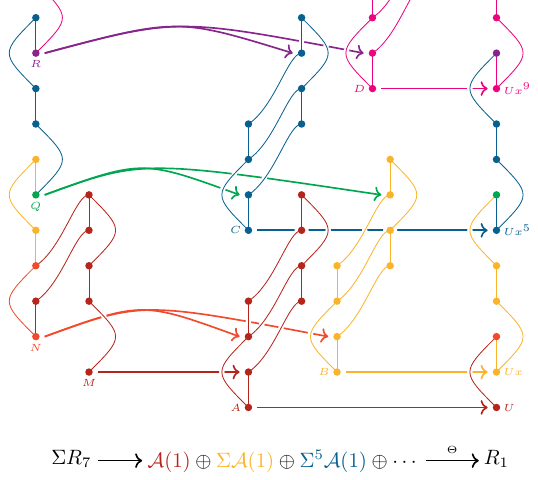}
\caption{A picture proof of \cref{R6ker}: that the kernel of the map $\Theta\colon \cA(1)\oplus \Sigma\cA(1) \oplus
\Sigma^5\cA(1)\oplus\dotsb\to R_1$ is isomorphic to $R_7$. The inclusion of the kernel sends
$\textcolor{BrickRed}{M}\mapsto \textcolor{BrickRed}{\Sq^1(A)}$, $\textcolor{RedOrange}{N}\mapsto
\textcolor{BrickRed}{\Sq^2(A)} + \textcolor{Dandelion}{\Sq^1(B)}$, $\textcolor{Green}{Q}\mapsto
\textcolor{Dandelion}{\Sq^2\Sq^3(B)} + \textcolor{MidnightBlue}{\Sq^1(C)}$, $\textcolor{Fuchsia}{R}\mapsto
\textcolor{MidnightBlue}{\Sq^2\Sq^3(C)} + \textcolor{RubineRed}{\Sq^1(D)}$, etc.}
\label{R6_as_kernel}
\end{figure}

Recall from \cref{names} that $R_1$ was defined to be the cohomology of the Thom spectrum $(B\Z/2)^{2\sigma-2}$.
Thus, by the Thom isomorphism theorem, $R_1$ is one-dimensional over $\Z/2$ in each nonnegative degree $n$,
generated by the class $Ux^n$, where $U\in H^0$ is the Thom class and $x\in H^1(B\Z/2;\Z/2)$ is the unique nonzero
class.
\begin{lem}\label{R6ker}
Let $\Theta\colon\cA(1)\oplus\bigoplus_{i\ge 0}\Sigma^{4i+1}\cA(1)\to R_1$ be the $\cA(1)$-module map sending the
generator of the $\Sigma^k\cA(1)$ summand to $Ux^k$. Then $\Theta$ is surjective with kernel isomorphic to $\Sigma
R_7$.
\end{lem}
This is proven by direct calculation. See \cref{R6_as_kernel} for a picture illustrating the proof.

\begin{cor}
\label{ExtR6}
For $s,t\ge 0$, there is an isomorphism $\Ext_{\cA(1)}^{s,t}(R_7)\cong \Ext_{\cA(1)}^{s+1,t+1}(R_1)$ compatible
with the $\E$-actions on both sides.
\end{cor}
\begin{proof}
The map $\Theta$ from \cref{R6ker} is the first step in a minimal
resolution $P_\bullet$ of $R_1$, so $P_{\bullet-1}$ (for $\bullet\ge 1$) forms a minimal resolution for
$\ker(\Theta)\cong\Sigma R_7$, which implies the corollary statement on Ext groups.
\end{proof}
We have just one more auxiliary lemma before we can give the Adams charts: we have to calculate the $\cA(1)$-modules $M\otimes\uQ$ for the $\cA(1)$-modules $M$ we are interested in.
\begin{lem}
\label{tensor_Q}
There are $\cA(1)$-module isomorphisms
\begin{subequations}
\begin{alignat}{2}
    \Z/2\otimes\uQ &\cong \uQ \qquad\qquad\qquad\qquad & \cA(1)\otimes\uQ &\cong \cA(1) \oplus \Sigma^2\cA(1) \oplus \Sigma^3\cA(1)\\
    \uQ\otimes \uQ &\cong Q_{2,1} & J\otimes\uQ &\cong \cA(1)\oplus \Sigma^2 R_2\\
    Q\otimes\uQ &\cong \Sigma^3\Z/2\oplus \cA(1) & R_2\otimes\uQ &\cong Q_{0,1} \oplus \Sigma\cA(1) \oplus\Sigma^2\cA(1)\\
    \Upsilon_1 \otimes \uQ &\cong \Upsilon_1 \oplus \Sigma^2 \cA(1) & P\otimes\uQ &\cong \Sigma^2 R_0\oplus\bigoplus_{i\ge 0}\Sigma^{4i}\cA(1)\\
    R_0\otimes \uQ &\cong R_3\oplus \bigoplus_{i\ge 0}\Sigma^{4i+2}\cA(1) &
    R_3\otimes\uQ &\cong \Sigma^3 R_5\oplus \cA(1) \oplus \bigoplus_{i\ge 0} \Sigma^{4i+2}\cA(1)
\end{alignat}
\begin{equation}
    R_5\otimes\uQ \cong \Sigma^3 P \oplus \cA(1) \oplus \Sigma\cA(1) \oplus \bigoplus_{i\ge 0} \Sigma^{4i+3}\cA(1).
\end{equation}
\end{subequations}
\end{lem}
\begin{proof}
For $\Z/2$ this is tautological. For $\cA(1)$, $\cA(1)\otimes M$ is a free $\cA(1)$-module on the graded vector space $M$. For $\uQ$, $Q$, and $R_2$, use \cref{milgram_tensor} with $(j_2,n_2) = (1, 0)$ and $(j_1,n_1)$ equal to $(1, 0)$, $(3, 0)$, and $(3, 1)$ respectively; the isomorphism type of the free summand of the tensor product then follows from a dimension-counting argument. For $J$, this is~\cite[Lemma 3.4(2)]{Mil24}. For $\Upsilon_1$, use $\Upsilon_1\cong\cA(1)\otimes_{\cA(0)}\Z/2$ and the fact that, as an $\cA(0)$-module, $\uQ\cong \Z/2\oplus \Sigma^2\cA(0)$.

For $M$ equal to $P$, $R_0$, $R_3$, and $R_5$, use the fact that $\uQ$ is \term{stably invertible}~\cite[\S 3]{AP76} and $M$ is \term{$Q_1$-locally stably invertible}~\cite{Yu95} to deduce that $M\otimes\uQ$ is also $Q_1$-locally stably invertible and is isomorphic to $M\otimes \Sigma^{-1}R_0 \otimes \uQ$ modulo free summands~\cite[Theorem 5.2]{Bru14}. $R_0\otimes\uQ\cong R_3$ modulo free summands~\cite[(6.45)]{BottSpiral}, and Yu~\cite[Lemma 2.5]{Yu95} calculates the tensor product of any two of $P$, $R_0$, $R_3$, and $R_5$ modulo free summands. Finally, put the free summands back with a Poincaré polynomial argument like the one in the proof of \cref{R6}.
\end{proof}

Below, bidegrees are of the form $(t-s,s)$.
\begin{thm}
\label{Q_ext_calc}
Recall from \cref{sE_defn} that $\sE\cong \Z/2[h_0,h_1]/(h_0h_1, h_1^3)$, with $\abs{h_0} = (0, 1)$ and $\abs{h_1} = (1,
1)$.
\begin{enumerate}
\item
We have the following $\sE$-module isomorphisms.
\begingroup
\postdisplaypenalty=10000
\predisplaypenalty=0
\begin{subequations}
\label{lists_of_gens}
\begin{align}
	\cQ(\Z/2) &\cong \sE\\
	\cQ(\cA(1)) &\cong \sE[a_1, a_2, a_3]\big/ \notag\\
	&\quad (h_0, h_1, a_i^2, a_ia_j\text{ for all } i,j),\\
\shortintertext{where $\abs{a_1} = (5, 2)$, $\abs{a_2} = (7, 2)$, and $\abs{a_3} = (8, 2)$.}
	\cQ(B) &\cong \sE[b]/(h_1, b^3, h_0b^2),\qquad\abs{b} = (2, 1).\\
	\cQ(N_1) &\cong \sE\set{c_1,c_2}/(h_0c_1, h_1^2c_1 + h_0c_2, h_1^2c_2),\notag\\
	&\quad \abs{c_1} = (0, 0),\quad \abs{c_2} = (2, 1).\label{Q_N1}\\
	\cQ(\uQ) &\cong \sE\set{d_1, d_2}/(h_1d_1, h_0^2d_2, h_1^2d_2),\notag\\
	&\quad \abs{d_1} = (0, 0),\quad \abs{d_2} = (4, 1).\\
	\cQ(J) &\cong \sE\set{e_1, e_2, e_3, e_4}\big/ \notag\\
	&\quad (h_0e_1, h_1e_1, h_1e_2, h_0 e_3, h_1 e_3, h_0e_4, h_1 e_4),\\
\shortintertext{where $\abs{e_1} = (0, 0)$, $\abs{e_2} = (2, 1)$, $\abs{e_3} = (5, 2)$, and $\abs{e_4} = (6, 2)$.}
	\cQ(Q) &\cong \sE\set{f_1, f_2, f_3}\big/ \notag\\
	&\quad (h_0f_1, h_1^2 f_1, h_1f_2, h_0f_3, h_1f_3),\\
\shortintertext{where $\abs{f_1} = (0, 0)$, $\abs{f_2} = (3, 2)$, and $\abs{f_3} = (5, 2)$.}
	\cQ(R_2) &\cong \sE\set{g_1, g_2, g_3, g_4, g_5}\big/ \notag\\
	&\quad (h_1 g_1, h_0g_2, h_1^2 g_2, h_0g_3, h_1g_3,\notag\\
	&\qquad h_0 g_4, h_1 g_4, h_0g_5, h_1g_5),\\
\shortintertext{where $\abs{g_1} = (0, 0)$, $\abs{g_2} = (1, 0)$, $\abs{g_3} = (4, 2)$, $\abs{g_4} = (6, 2)$, and $\abs{g_5} = (7, 2)$.}
	\cQ(\Upsilon_1) &\cong \sE[j_1, j_2]/(h_1, j_1^2, j_2^2, j_1j_2, h_0j_2),\notag\\
	&\quad \abs{j_1} = (5, 2),\quad \abs{j_2} = (7, 2).\\
	\cQ(P) &\cong \sE\set{\kappa, k_i, k_i': i\ge 1}\big/ \notag\\
	&\quad (h_0\kappa, h_0^3 k_i, h_1k_i, h_0k_i', h_1k_i', h_1^2\kappa + h_0^2k_1),\\
\shortintertext{where $\abs{\kappa} = (0, 0)$, $\abs{k_i} = (4i-2, 0)$, and $\abs{k_i'} = (4i+1, 2)$.}
	\cQ(R_0) &\cong \sE\set{\lambda, \lambda', \ell_i, \ell_i': i\ge 1}\big/ \notag\\
	&\quad (h_0\lambda, h_1\lambda, h_0\lambda',h_1^2\lambda',\notag\\
	&\qquad h_0^3 \ell_i, h_1\ell_i, h_0\ell_i', h_1\ell_i'),\\
\shortintertext{where $\abs{\lambda} = (0, 0)$, $\abs{\lambda'} = (2, 1)$, $\abs{\ell_i} = (4i, 0)$, and $\abs{\ell_i'} = (4i+3, 2)$.}
	\cQ(R_1) &\cong \sE\set{\mu_1, \mu_2, m_i, m_i': i\ge 1}\big/ \notag\\
	&\quad (h_1\mu_1 + h_0m_1, h_0\mu_2, h_1\mu_2,\notag\\
	&\qquad h_0^3 m_i, h_0m_i', h_1m_i, h_1m_i'),\\
\shortintertext{where $\abs{\mu_1} = (0, 0)$, $\abs{\mu_2} = (4, 2)$, $\abs{m_i} = (4i+1, 0)$, and $\abs{m_i'} = (4i+8, 2)$.}
	\cQ(R_3) &\cong (\sE/h_1)\set{\nu, \nu', n_i, n_i': i\ge 0}\big/ \notag\\
	&\quad (h_0\nu, h_0\nu', h_0 n_0, h_0 n_i', h_0^2 n_1,\notag\\
	&\qquad h_0^3 n_i\text{ for } i\ge 2),\\
\shortintertext{where $\abs{n_i} = (4i, 0)$, $\abs{n_i'} = (4i+7, 2)$, $\abs{\nu'} = (5, 2)$, and $\abs{\nu} = (6, 2)$.}
	\cQ(R_5) &\cong \sE\set{\pi_1,\pi_2,\pi_3,\pi_4, p_i, p_i': i\ge 1}\big/ \notag\\
	&\quad (h_0\pi_1, h_1\pi_1 + h_0\pi_2, h_1\pi_2,\notag\\
	&\qquad h_0\pi_3, h_1\pi_3, h_0\pi_4, h_1\pi_4,\notag\\
	&\qquad h_0^3 p_i, h_1 p_i, h_0p_i', h_1p_i'),\\
\shortintertext{where $\abs{\pi_1} = (0, 0)$, $\abs{\pi_2} = (1, 0)$, $\abs{\pi_3} = (5, 2)$, $\abs{\pi_4} = (6, 2)$, $\abs{p_i} = (4i+1, 0)$, and $\abs{p_i'} = (4i+4, 2)$.}
	\cQ(R_6) &\cong (\sE/h_1)\set{q_i, q_i': i\ge 0}\big/ \notag\\
	&\quad (h_0^3 q_i, h_0q_i'\text{ for } i\ge 0,\notag\\
	&\qquad h_0q_0, h_0^2 q_1)
\end{align}
\end{subequations}
where $\abs{q_i} = (2i, 0)$ and $\abs{q_i'} = (2i+5, 2)$.
\endgroup

In the cases of $\Z/2$, $B$, and $\Upsilon_1$, these are also equivalences of $\sE$-algebras.
\item\label{is_surj} For $M\in\set{\Z/2, B, N_1, \uQ}$, the map $\partial_2\colon \Ext_{\cA(1)}(M)\to\cQ(M)$ is surjective. For the remaining $\cA(1)$-modules we consider, $\coker(\partial_2)$ is spanned by the classes listed in \cref{cokernel_info}.
\end{enumerate}
\end{thm}
For all $\cA(1)$-modules $M$ appearing in this theorem, we draw $M$ and $\cQ(M)$ in \cref{Q_ext_charts_1}.

\begin{table}[htbp]
\begin{tabular}{c c}
\toprule
$M$ & generators of $\coker(\partial_2)$\\
\midrule
    $\cA(1)$ & $a_1$, $a_2$, $a_3$\\
    $J$ & $e_3$\\
    $Q$ & $f_3$\\
    $R_2$ & $g_4$, $g_5$\\
    $\Upsilon_1$ & $j_1$, $j_2$\\
    $P$ & $k_i'$, $i\ge 1$\\
    $R_0$ & $\ell_i'$, $i\ge 1$\\
    $R_1$ & $m_i'$, $i\ge 1$\\
    $R_3$ & $\nu'$, $n_i'$, $i\ge 0$\\
    $R_5$ & $\pi_3$, $\pi_4$, $p_i'$, $i\ge 1$\\
    $R_6$ & $q_i'$, $i\ge 0$\\
\bottomrule
\end{tabular}
\caption{For each $\cA(1)$-module $M$ in the table, we give an $\sE$-module generating set for the cokernel of the map $\partial_2\colon\Ext(M)\to\cQ(M)$ from~\eqref{cQLES}. For $M = \Z/2$, $B$, $N_1$, and $\uQ$, $\coker(\partial_2) = \set 0$. In all cases, the image of $\partial_2$ is spanned as an $\sE$-module by the generators we gave in~\eqref{lists_of_gens} that are not present in this table. This table is part of \cref{Q_ext_calc}.}
\label{cokernel_info}
\end{table}

In the proof of \cref{Q_ext_calc}, we will need to know $\Ext_{\cA(1)}(M)$ as an $\E$-module for lots of $\cA(1)$-modules $M$. For the most part, the Ext groups we need have been worked out in the literature, together with the actions of $h_0$ and $h_1$, giving the module structure over $\Z/2[h_0, h_1]/(h_0h_1, h_1^3)$. Most of these references leave the rest of the $\E$-module structure implicit. Thus in \cref{Ext_refs} we gather what is to our knowledge the original references for these computations.

\begin{table}[htbp]
\begin{tabular}{c c c}
\toprule
$M$ & $\Ext(M)$ (partial $\E$-mod.\ str.) & $\Ext(M)$ (full $\E$-mod.\ str.)\\
\midrule
$\Z/2$ & --- & \cite[Theorem 3]{Liu62}\\
$B$ & --- & \cite[\S 4.5]{Bay94} \\
$N_1$ & --- & \cite[\S VII.2]{Per78}\\
$\uQ$ & \cite[\S 3]{AP76} & \cite[\S VII.1]{Per78}\\
$J$ & \cite[\S 3]{AP76} & \cite[\S VII.1]{Per78}\\
$Q$ & \cite[\S 3]{AP76} & \cite[\S VII.1]{Per78}\\
$R_2$ & \cite[\S 3]{AP76} & \cite[\S VII.1]{Per78}\\
$\Upsilon_1$ & --- & \cite[Lemma 5.8]{Mil75}\\
$P$ & \cite[\S 2]{GMM68} & \cite[Figure A.9.17]{BG10}\\
$R_0$ & \cite[\S 2]{GMM68} & \cite[Example 5.40]{BottSpiral}\\
$R_1$ & \cite[\S 2]{GMM68} & \cref{E5_ext}\\
$R_3$ & \cite[Theorem 3.1]{Yu95} & \cite[Example 5.37]{BottSpiral}\\
$R_5$ & \cite[Theorem 3.1]{Yu95} & \cite[Example 5.34]{BottSpiral}\\
$R_6$ & \cite[Theorem A.0.19]{Sie13} & \cref{B_tensor_stuff}\\
$B\otimes\uQ$ & \cite[Appendix A]{DL20} & \cref{B_tensor_stuff}\\
$N_1\otimes\uQ$ & \cite[Figure 6]{Pea14} & \cite[Proposition 5.99]{BottSpiral}\\
$Q_{2,1}$ & \cite[\S 3]{AP76} & \cite[\S VII.1]{Per78}\\
$Q_{0,1}$ & \cite[\S 3]{AP76} & \cite[\S VII.1]{Per78}\\
$R_7$ & --- & \cref{ExtR6}\\
\bottomrule
\end{tabular}
\caption{In the proof of \cref{Q_ext_calc}, we will need to know the $\E$-module structures on $\Ext(M)$ for various $\cA(1)$-modules $M$. In this table, we list the first reference we know of explicitly giving each $\E$-module structure. In many cases, an earlier reference described $\Ext(M)$ either as merely a bigraded vector space or as a module over $\Z/2[h_0, h_1]/(h_0h_1, h_1^3)$, leaving the rest of the $\E$-module structure implicit. In those cases, we give both references for completeness. The reader may also find the Ext calculations in~\cite[Appendix A]{BG10}, \cite[\S\S 4.4--4.6, 5.2, 5.3, 5.5]{BC18}, \cite[\S 6]{BR21}, and~\cite[Proposition 14.11]{BDDHM25} helpful.}
\label{Ext_refs}
\end{table}

\begin{proof}[Proof sketch]
The strategy is exactly the same in all sixteen cases; therefore we walk through the case $M = N_1$ in detail, then discuss the few differences that occur in the remaining fifteen cases.

We will compute the long exact sequence from \cref{cQthm} in order to find $\cQ(N_1)$. We draw this computation in \cref{N1_Q_calc}. Consulting \cref{Ext_refs}, we have $\E$-module isomorphisms
\begin{subequations}
\begin{align}
    \Ext(N_1) &\cong \E\{c_1, c_2\}/(h_0c_1, h_1^2c_1 + h_0c_2, vc_1 + h_1^2c_2, vc_2)\\
    \Ext(\uQ\otimes N_1) &\cong\E\{p,q,r\}/
    (h_0p, h_1p, vp, wp = h_1^2r, h_0q, h_1^2q + h_0r, vq + h_1^2r, vr),
\end{align}
\end{subequations}
where the $(t-s,s)$ bidegrees of the generators are $\abs{c_1} = (0, 0)$, $\abs{c_2} = (2, 1)$, $\abs{p} = (0, 0)$, $\abs{q} = (4, 1)$, and $\abs r = (6, 2)$. The map $f_2\colon \Ext^{s,t}(\uQ\otimes N_1)\to\Ext^{s+3,t+7}(N_1)$ in the long exact sequence~\eqref{cQLES} is $\E$-linear, hence is determined by its values on $p$, $q$, and $r$. For degree reasons, $f_2(p)$ is either $0$ or $vc_1$, $f_2(q)$ is $0$ or $wc_1$, and $f_2(r)$ is $0$ or $wc_2$. Moreover, $\E$-linearity
implies all three of these boundary maps are zero or all three are nonzero: $w$, $h_1^2$, and $v$ respectively act injectively on the domain and codomain of $\partial$ for $p$, $q$, and $r$ respectively, so all three are detected by whether $f_2(wp)$ is $0$ or $vwc_1$. 

\begin{figure}[htbp]
\begin{subfigure}[c]{0.45\textwidth}
\includegraphics[width=\linewidth]{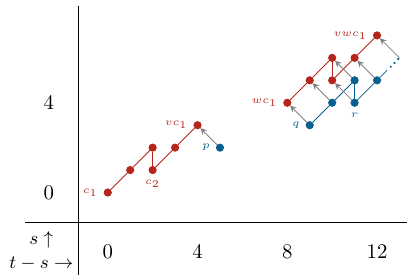}
%
\end{subfigure}%
\hfill
\begin{subfigure}[c]{0.45\textwidth}
    \includegraphics[width=\linewidth]{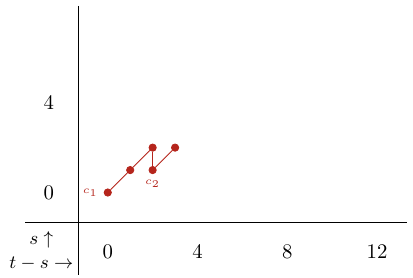}
\end{subfigure}
\caption{The long exact sequence from \cref{cQthm}, used to compute $\cQ(N_1)$ in the proof of \cref{Q_ext_calc}. Left: the map $f_2\colon\Ext^{s, t}(N_1\otimes\uQ)\to\Ext^{s+3,t+7}(N_1)$. This map is injective in all degrees, which determines the entire long exact sequence; in particular, $\cQ(N_1)\cong\coker(f_2)$ and $\partial_2\colon\Ext(N_1)\to\cQ(N_1)$ is surjective. Right: $\cQ(N_1)$.}
\label{N1_Q_calc}
\end{figure}

Thus we need to find $f_2(p)$.  For degree reasons, the Künneth map $\kappa\colon \Ext(\uQ)\otimes\Ext(N_1)\to\Ext(\uQ\otimes N_1)$ must send $\delta\otimes c_1\mapsto p$. Recall from \cref{sE_defn} that, if $N_1$ is replaced with $\Z/2$, $f_2(\delta) = v$. Thus
\begin{equation}
    f_2(p) = f_2(\kappa(\delta \otimes c_1)) = f_2(\delta) f_2(c_1) = vc_1,
\end{equation}
using that $\E$ acts naturally on the entire long exact sequence~\eqref{cQLES}.
As stated above, this forces $\partial(q) = wc_1$ and $\partial(r) = wc_2$.

After computing the boundary map, we see that all classes in topological degrees $4$ and above vanish, leaving only the vector space spanned by $c_1$, $h_1c_1$, $h_1^2c_1 = h_0c_2$, and $c_2$, which is the $\sE$-module claimed in the theorem statement.

For the remaining fifteen $\cA(1)$-modules $M$ appearing in the theorem statement, one calculates $\cQ(M)$ in the same way, using the long exact sequence in \cref{cQthm}. As input, we need to compute the $\cA(1)$-module structure on $\uQ\otimes M$ and the $\E$-module structures on $\Ext(M)$ and $\Ext(\uQ\otimes M)$. We calculated $\uQ\otimes M$ in most cases in \cref{tensor_Q}; see \cref{Ext_refs} for all Ext calculations needed in this proof. This leaves open $\uQ\otimes M$ for $M\in\set{B, R_1, R_6}$, which we did not discuss in \cref{tensor_Q}. We computed $\uQ\otimes R_1$ in \cref{R6}, part~\eqref{R6_tenspart}. We calculate $\Ext(\uQ\otimes B)$ in \cref{B_tensor_stuff}. Finally, since $R_6 \coloneqq P\otimes B$, $\uQ\otimes R_6\cong B\otimes(\uQ\otimes P)$. \Cref{B_tensor_stuff} thus reduces the computation of $\Ext_{\cA(1)}(\uQ\otimes R_6)$ to that of $\Ext_{\cE(1)}(\uQ\otimes P)$. \Cref{tensor_Q} reduces to finding $\Ext_{\cE(1)}$ of $R_0$ and $\cA(1)$. There is an $\cE(1)$-module isomorphism $\cA(1)\cong\cE(1)\oplus\Sigma^2\cE(1)$ (\cite[\S 4.7.2]{Bay94}); $\Ext_{\cE(1)}(\cE(1))$ is tautologically $\Z/2$ concentrated in bidegree $(0, 0)$, and $\Ext_{\cE(1)}(R_0)$ is given in~\cite[Proposition 4.47]{Deb21}.
\end{proof}
\begin{rem}
As we discussed in \cref{coalg_remark}, $R_1$ is a coalgebra, so $\cQ(R_1)$ is an $\sE$-algebra. In \cref{Q_ext_calc}, we only computed the $\sE$-module structure on this algebra. This is because, unlike for $\E$ acting on $\Ext(R_1)$ (c.f.\ \cref{E5_ext}), our method of proof does not uniquely determine the $\sE$-algebra structure. 
\end{rem}

\bibliographystyle{alpha}
\bibliography{references}

@incollection{BC18,
    AUTHOR = {Beaudry, Agnès and Campbell, Jonathan A.},
     TITLE = {A guide for computing stable homotopy groups},
 BOOKTITLE = {Topology and quantum theory in interaction},
    SERIES = {Contemp. Math.},
    VOLUME = {718},
     PAGES = {89--136},
 PUBLISHER = {Amer. Math. Soc., [Providence], RI},
      YEAR = {2018},
   MRCLASS = {55Q10 (55R40 55S10 55T15)},
  MRNUMBER = {3869642},
	DOI = {10.1090/conm/718/14476},
	note = {\href{https://arxiv.org/abs/1801.07530}{arXiv:\allowbreak 1801.07530}},
    eprint = {1801.07530},
    archivePrefix = {arXiv},
}

@article{BL01,
   AUTHOR = {Baker, Andrew and Lazarev, Andrey},
     TITLE = {On the {A}dams spectral sequence for {$R$}-modules},
   JOURNAL = {Algebr. Geom. Topol.},
  FJOURNAL = {Algebraic \& Geometric Topology},
    VOLUME = {1},
      YEAR = {2001},
     PAGES = {173--199},
      ISSN = {1472-2747},
   MRCLASS = {55P42 (55N20 55P43 55T15)},
  MRNUMBER = {1823498},
MRREVIEWER = {Mark Hovey},
       DOI = {10.2140/agt.2001.1.173},
       note = {\href{https://arxiv.org/abs/math/0105079}{arXiv:\allowbreak math/0105079}},
    eprint = {math/0105079},
    archivePrefix = {arXiv},
}

@misc{Cam17,
	 author = {Jonathan A. Campbell},
	 title = {Homotopy Theoretic Classification of Symmetry Protected Phases},
	 year = {2017},
	 note = {\href{https://arxiv.org/abs/1708.04264}{arXiv:\allowbreak 1708.04264}},
    eprint = {1708.04264},
    archivePrefix = {arXiv},
}

@incollection{Mil75,
    AUTHOR = {Milgram, R. James},
     TITLE = {The {S}teenrod algebra and its dual for connective
		{$K$}-theory},
 BOOKTITLE = {Conference on homotopy theory ({E}vanston, {I}ll., 1974)},
    SERIES = {Notas Mat. Simpos.},
    VOLUME = {1},
     PAGES = {127--158},
 PUBLISHER = {Soc. Mat. Mexicana, México},
      YEAR = {1975},
   MRCLASS = {55S10 (55N15 55S25)},
  MRNUMBER = {761725},
}

@article{DGM81,
    AUTHOR = {Davis, Donald M. and Gitler, Sam and Mahowald, Mark},
     TITLE = {The stable geometric dimension of vector bundles over real
		projective spaces},
   JOURNAL = {Trans. Amer. Math. Soc.},
  FJOURNAL = {Transactions of the American Mathematical Society},
    VOLUME = {268},
      YEAR = {1981},
    NUMBER = {1},
     PAGES = {39--61},
      ISSN = {0002-9947},
   MRCLASS = {55N15 (55R25)},
  MRNUMBER = {628445},
MRREVIEWER = {Harvey Margolis},
	DOI = {10.2307/1998336},
	URL = {https://doi.org/10.2307/1998336},
}

@misc{Deb21,
	 author = {Arun Debray},
	 title = {Invertible phases for mixed spatial symmetries and the fermionic crystalline equivalence principle},
	 year = {2021},
	 note = {\href{https://arxiv.org/abs/2102.02941}{arXiv:\allowbreak 2102.02941}},
    eprint = {2102.02941},
    archivePrefix = {arXiv},
}

@article{KPMT20,
    AUTHOR = {Kaidi, Justin and Parra-Martinez, Julio and Tachikawa, Yuji},
     TITLE = {Topological superconductors on superstring worldsheets},
      NOTE = {With a mathematical appendix by Arun Debray. \href{https://arxiv.org/abs/1911.11780}{arXiv:\allowbreak 1911.11780}},
   JOURNAL = {SciPost Phys.},
  FJOURNAL = {SciPost Physics},
    VOLUME = {9},
      YEAR = {2020},
    NUMBER = {1},
     PAGES = {010},
   MRCLASS = {81T30},
  MRNUMBER = {4196335},
	DOI = {10.21468/SciPostPhys.9.1.010},
	URL = {https://doi.org/10.21468/SciPostPhys.9.1.010},
    eprint = {1911.11780},
    archivePrefix = {arXiv},
    pagetotal = {70},
}

@mastersthesis{Ped17,
  author	= {Riccardo Pedrotti},
  title = {Stable classification of certain families of four-manifolds},
  school	= {Max Planck Institute for Mathematics},
  year   = {2017}
}

@article{GOPWW20,
    AUTHOR = {Guo, Meng and Ohmori, Kantaro and Putrov, Pavel and Wan,
		Zheyan and Wang, Juven},
     TITLE = {Fermionic finite-group gauge theories and interacting
		symmetric/crystalline orders via cobordisms},
   JOURNAL = {Commun. Math. Phys.},
  FJOURNAL = {Communications in Mathematical Physics},
    VOLUME = {376},
      YEAR = {2020},
    NUMBER = {2},
     PAGES = {1073--1154},
      ISSN = {0010-3616},
   MRCLASS = {81T45 (57R56)},
  MRNUMBER = {4103964},
	DOI = {10.1007/s00220-019-03671-6},
	 note = {\href{https://arxiv.org/abs/1812.11959}{arXiv:\allowbreak 1812.11959}},
    eprint = {1812.11959},
    archivePrefix = {arXiv},
}

@article{LY22,
    AUTHOR = {Lee, Yasunori and Yonekura, Kazuya},
     TITLE = {Global anomalies in 8d supergravity},
   JOURNAL = {J. High Energy Phys.},
  FJOURNAL = {Journal of High Energy Physics},
    VOLUME = {2022},
      YEAR = {2022},
    NUMBER = {7},
     PAGES = {125},
      ISSN = {1126-6708},
   MRCLASS = {83E50},
  MRNUMBER = {4458304},
	DOI = {10.1007/jhep07(2022)125},
	note = {\href{https://arxiv.org/abs/2203.12631}{arXiv:\allowbreak 2203.12631}},
    eprint = {2203.12631},
    archivePrefix = {arXiv},
}

@incollection{FHT10,
    AUTHOR = {Freed, Daniel S. and Hopkins, Michael J. and Teleman,
		Constantin},
     TITLE = {Consistent orientation of moduli spaces},
 BOOKTITLE = {The many facets of geometry},
     PAGES = {395--419},
 PUBLISHER = {Oxford Univ. Press, Oxford},
      YEAR = {2010},
   MRCLASS = {57R56 (19L50 58D29)},
  MRNUMBER = {2681705},
MRREVIEWER = {Andrey Yu. Lazarev},
	DOI = {10.1093/acprof:oso/9780199534920.003.0019},
	URL = {https://doi.org/10.1093/acprof:oso/9780199534920.003.0019},
}

@article{FH21,
    AUTHOR = {Freed, Daniel S. and Hopkins, Michael J.},
     TITLE = {Reflection positivity and invertible topological phases},
   JOURNAL = {Geom. Topol.},
  FJOURNAL = {Geometry \& Topology},
    VOLUME = {25},
      YEAR = {2021},
    NUMBER = {3},
     PAGES = {1165--1330},
      ISSN = {1465-3060},
   MRCLASS = {81T45 (55N22 57R90 81T50 82D03)},
  MRNUMBER = {4268163},
	DOI = {10.2140/gt.2021.25.1165},
	note = {\href{https://arxiv.org/abs/1604.06527}{arXiv:\allowbreak 1604.06527}},
    eprint = {1604.06527},
    archivePrefix = {arXiv},
}

@article{Lee:2020ewl,
    author = "Lee, Yasunori and Tachikawa, Yuji",
    title = "{Some comments on 6D global gauge anomalies}",
    eprint = "2012.11622",
    archivePrefix = "arXiv",
    primaryClass = "hep-th",
    doi = "10.1093/ptep/ptab015",
    journal = {Prog. Theor. Exp. Phys.},
    volume = "2021",
    number = "8",
    pages = "08B103",
    year = "2021",
	 note = {\href{https://arxiv.org/abs/2012.11622}{arXiv:\allowbreak 2012.11622}},
    fjournal = {Progress of Theoretical and Experimental Physics},
}

@article{Davighi:2020uab,
    author = "Davighi, Joe and Lohitsiri, Nakarin",
    title = {The algebra of anomaly interplay},
    eprint = "2011.10102",
    archivePrefix = "arXiv",
    primaryClass = "hep-th",
    doi = "10.21468/SciPostPhys.10.3.074",
    journal = "SciPost Phys.",
    volume = "10",
    number = "3",
    pages = "074",
    year = "2021",
	 note = {\href{https://arxiv.org/abs/2011.10102}{arXiv:\allowbreak 2011.10102}}
}

@article{Davighi:2022icj,
    AUTHOR = {Davighi, Joe and Gripaios, Ben and Lohitsiri, Nakarin},
     TITLE = {Anomalies of non-{A}belian finite groups via cobordism},
   JOURNAL = {J. High Energy Phys.},
  FJOURNAL = {Journal of High Energy Physics},
    VOLUME = {2022},
      YEAR = {2022},
    NUMBER = {9},
     PAGES = {147},
      ISSN = {1126-6708},
   MRCLASS = {81T45 (55 81T30 81T50)},
  MRNUMBER = {4486791},
	DOI = {10.1007/jhep09(2022)147},
    note = {\href{https://arxiv.org/abs/2207.10700}{arXiv:\allowbreak 2207.10700}},
    eprint = {2207.10700},
    archivePrefix = {arXiv},
}

@article{Davighi:2023luh,
    author = "Davighi, Joe and Lohitsiri, Nakarin",
    title = "{Toric 2-group anomalies via cobordism}",
    eprint = "2302.12853",
    archivePrefix = "arXiv",
    primaryClass = "hep-th",
    doi = "10.1007/JHEP07(2023)019",
    journal = {J. High Energy Phys.},
    volume = {2023},
    pages = "019",
    year = "2023",
    note = {With an appendix by Arun Debray. \href{https://arxiv.org/abs/2302.12853}{arXiv:\allowbreak 2302.12853}},
    number = {7},
    fjournal = {Journal of High Energy Physics},
}

@article{WW19,
  title = {Higher anomalies, higher symmetries, and cobordisms {I}: classification of higher-symmetry-protected topological states and their boundary fermionic/bosonic anomalies via a generalized cobordism theory},
  author = {Zheyan Wan and Juven Wang},
  journal = {Ann. Math. Sci. Appl.},
  volume = {4},
  number = {2},
  pages = {107--311},
  year = {2019},
  publisher = {International Press},
	 note = {\href{https://arxiv.org/abs/1812.11967}{arXiv:\allowbreak 1812.11967}},
    eprint = {1812.11967},
    archivePrefix = {arXiv},
}

@article{DL20c,
	 author = {Joe Davighi and Nakarin Lohitsiri},
	 title = {Omega {\textit{vs}.} pi, and 6d anomaly cancellation},
	 note = {\href{https://arxiv.org/abs/2012.11693}{arXiv:\allowbreak 2012.11693}},
	 volume = {2021},
	 number = {5},
	 year = {2021},
	 journal = {J. High Energy Phys.},
    doi = {10.1007/JHEP05(2021)267},
    pages = {267},
    fjournal = {Journal of High Energy Physics},
    eprint = {2012.11693},
    archivePrefix = {arXiv},
}

@article{DL20,
    AUTHOR = {Davighi, Joe and Lohitsiri, Nakarin},
     TITLE = {Anomaly interplay in {$\rm U(2)$} gauge theories},
   JOURNAL = {J. High Energy Phys.},
  FJOURNAL = {Journal of High Energy Physics},
    VOLUME = {2020},
      YEAR = {2020},
    NUMBER = {5},
     PAGES = {098},
      ISSN = {1126-6708},
   MRCLASS = {81T50 (81T13)},
  MRNUMBER = {4112345},
	DOI = {10.1007/jhep05(2020)098},
	 note = {\href{https://arxiv.org/abs/2001.07731}{arXiv:\allowbreak 2001.07731}},
    pagetotal = {20},
    eprint = {2001.07731},
    archivePrefix = {arXiv},
}

@article{LOT21,
	 title={Revisiting {Wess-Zumino-Witten} terms},
	 author={Yasunori Lee and Kantaro Ohmori and Yuji Tachikawa},
	 journal={SciPost Phys.},
	 volume={10},
	 number={3},
	 pages={61},
	 year={2021},
	 publisher={SciPost},
	 doi={10.21468/SciPostPhys.10.3.061},
		  note = {\href{https://arxiv.org/abs/2009.00033}{arXiv:\allowbreak 2009.00033}},
    eprint = {2009.00033},
    archivePrefix = {arXiv},
}

@misc{JF20,
	 author = {Theo Johnson-Freyd},
	 title = {Topological {M}athieu Moonshine},
	 year = {2020},
	 note = {\href{https://arxiv.org/abs/2006.02922}{arXiv:\allowbreak 2006.02922}},
    eprint = {2006.02922},
    archivePrefix = {arXiv},
}

@misc{And69,
	 author = {D.W. Anderson},
	 title = {Universal coefficient theorems for {$K$}-theory},
	 year = {1969},
	 note = {\url{https://faculty.tcu.edu/gfriedman/notes/Anderson-UCT.pdf}}
}

@article{Ser53,
    AUTHOR = {Serre, Jean-Pierre},
     TITLE = {Cohomologie modulo {$2$} des complexes
		d'{E}ilenberg-{M}ac{L}ane},
   JOURNAL = {Comment. Math. Helv.},
  FJOURNAL = {Commentarii Mathematici Helvetici},
    VOLUME = {27},
      YEAR = {1953},
     PAGES = {198--232},
      ISSN = {0010-2571},
   MRCLASS = {56.0X},
  MRNUMBER = {60234},
MRREVIEWER = {W. S. Massey},
	DOI = {10.1007/BF02564562},
	URL = {https://doi.org/10.1007/BF02564562},
}

@article{Barkeshli:2023bta,
    author = "Barkeshli, Maissam and Hsin, Po-Shen and Kobayashi, Ryohei",
    title = {Higher-group symmetry of {(3+1)D} fermionic $\mathbb{Z}_2$ gauge theory: logical {CCZ}, {CS}, and {T} gates from higher symmetry},
    eprint = "2311.05674",
    archivePrefix = "arXiv",
    primaryClass = "cond-mat.str-el",
    doi = "10.21468/SciPostPhys.16.5.122",
    journal = "SciPost Phys.",
    volume = "16",
    number = "5",
    pages = "122",
    year = "2024",
    note = {\href{https://arxiv.org/abs/2311.05674}{arXiv:\allowbreak 2311.05674}}
}

@book{BR21,
    AUTHOR = {Bruner, Robert R. and Rognes, John},
     TITLE = {The {A}dams spectral sequence for topological modular forms},
    SERIES = {Mathematical Surveys and Monographs},
    VOLUME = {253},
 PUBLISHER = {American Mathematical Society, Providence, RI},
      YEAR = {2021},
     PAGES = {xix+690},
      ISBN = {978-1-4704-5674-0},
   MRCLASS = {55N34 (18G40 55P43 55Q45 55T15)},
  MRNUMBER = {4284897},
MRREVIEWER = {Tilman Bauer},
       DOI = {10.1090/surv/253},
       URL = {https://doi.org/10.1090/surv/253},
}

@article{Hil07,
    AUTHOR = {Hill, Michael A.},
     TITLE = {The 3-local {${\rm tmf}$}-homology of {$B\Sigma_3$}},
   JOURNAL = {Proc. Amer. Math. Soc.},
  FJOURNAL = {Proceedings of the American Mathematical Society},
    VOLUME = {135},
      YEAR = {2007},
    NUMBER = {12},
     PAGES = {4075--4086},
      ISSN = {0002-9939},
   MRCLASS = {55N20 (55N34 55T15)},
  MRNUMBER = {2341960},
MRREVIEWER = {Vincent Giambalvo},
       DOI = {10.1090/S0002-9939-07-08937-X},
	note = {\href{https://arxiv.org/abs/math/0511649}{arXiv:\allowbreak math/0511649}},
    eprint = {math/0511649},
    archivePrefix = {arXiv},
}

@article{ABGHR14a,
    AUTHOR = {Ando, Matthew and Blumberg, Andrew J. and Gepner, David and
	      Hopkins, Michael J. and Rezk, Charles},
     TITLE = {An {$\infty$}-categorical approach to {$R$}-line bundles,
	      {$R$}-module {T}hom spectra, and twisted {$R$}-homology},
   JOURNAL = {J. Topol.},
  FJOURNAL = {Journal of Topology},
    VOLUME = {7},
      YEAR = {2014},
    NUMBER = {3},
     PAGES = {869--893},
      ISSN = {1753-8416},
   MRCLASS = {55P43 (55N20 55U40)},
  MRNUMBER = {3252967},
MRREVIEWER = {Tyler D. Lawson},
       DOI = {10.1112/jtopol/jtt035},
	note = {\href{https://arxiv.org/abs/1403.4325}{arXiv:\allowbreak 1403.4325}},
    eprint = {1403.4325},
    archivePrefix = {arXiv},
}

@article{ABGHR14b,
    AUTHOR = {Ando, Matthew and Blumberg, Andrew J. and Gepner, David and
	      Hopkins, Michael J. and Rezk, Charles},
     TITLE = {Units of ring spectra, orientations and {T}hom spectra via
	      rigid infinite loop space theory},
   JOURNAL = {J. Topol.},
  FJOURNAL = {Journal of Topology},
    VOLUME = {7},
      YEAR = {2014},
    NUMBER = {4},
     PAGES = {1077--1117},
      ISSN = {1753-8416},
   MRCLASS = {55P43 (55N34 55R65)},
  MRNUMBER = {3286898},
MRREVIEWER = {Steven R. Costenoble},
       DOI = {10.1112/jtopol/jtu009},
	note = {\href{https://arxiv.org/abs/1403.4320}{arXiv:\allowbreak 1403.4320}},
    eprint = {1403.4320},
    archivePrefix = {arXiv},
}

@book{LMSM86,
    AUTHOR = {Lewis, Jr., L. G. and May, J. P. and Steinberger, M. and
	      McClure, J. E.},
     TITLE = {Equivariant stable homotopy theory},
    SERIES = {Lecture Notes in Mathematics},
    VOLUME = {1213},
      NOTE = {\url{http://www.math.uchicago.edu/~may/BOOKS/equi.pdf}},
 PUBLISHER = {Springer-Verlag, Berlin},
      YEAR = {1986},
     PAGES = {x+538},
      ISBN = {3-540-16820-6},
   MRCLASS = {55-02 (55Nxx 55Pxx 57S99)},
  MRNUMBER = {866482},
MRREVIEWER = {T. tom Dieck},
       DOI = {10.1007/BFb0075778},
       URL = {https://doi.org/10.1007/BFb0075778},
}

@misc{HA,
note = {\url{https://www.math.ias.edu/~lurie/papers/HA.pdf}},
author = {Jacob Lurie},
title = {Higher Algebra},
year = {2017}
}

@article{ABG11,
    AUTHOR = {Ando, Matthew and Blumberg, Andrew J. and Gepner, David},
     TITLE = {Parametrized spectra, multiplicative {T}hom spectra and the
	      twisted {U}mkehr map},
   JOURNAL = {Geom. Topol.},
  FJOURNAL = {Geometry \& Topology},
    VOLUME = {22},
      YEAR = {2018},
    NUMBER = {7},
     PAGES = {3761--3825},
      ISSN = {1465-3060},
   MRCLASS = {55P43 (55R70)},
  MRNUMBER = {3890766},
MRREVIEWER = {Birgit Richter},
       DOI = {10.2140/gt.2018.22.3761},
	note = {\href{https://arxiv.org/abs/1112.2203}{arXiv:\allowbreak 1112.2203}},
    eprint = {1112.2203},
    archivePrefix = {arXiv},
}

@article{MMSS:2001,
    author = {Mandell, M. A. and May, J. P. and Schwede, S. and Shipley, B.},
    title = "Model Categories of Diagram Spectra",
    journal = {Proc. Lond. Math. Soc. (3)},
    volume = {82},
    number = {2},
    pages = {441--512},
    year = {2001},
    issn = {0024-6115},
    doi = {10.1112/S0024611501012692},
	note = "\url{http://www.math.uchicago.edu/~may/PAPERS/mmssLMSDec30.pdf}"
}

@article{GJF:2017zba,
    author = "Gaiotto, Davide and Johnson-Freyd, Theo",
    title = {Symmetry-protected topological phases and generalized cohomology},
    eprint = "1712.07950",
    archivePrefix = "arXiv",
    primaryClass = "hep-th",
    doi = "10.1007/JHEP05(2019)007",
    journal = {J. High Energy Phys.},
    volume = {2019},
    pages = "007",
    year = "2019",
    note = {\href{https://arxiv.org/abs/1712.07950}{arXiv:\allowbreak 1712.07950}},
    number = {5},
    fjournal = {Journal of High Energy Physics},
}

@article{WG20,
  title = {Construction and Classification of Symmetry-Protected Topological Phases in Interacting Fermion Systems},
  author = {Wang, Qing-Rui and Gu, Zheng-Cheng},
  journal = {Phys. Rev. X},
  volume = {10},
  number = {3},
  pages = {031055},
  numpages = {64},
  year = {2020},
  publisher = {American Physical Society},
  doi = {10.1103/PhysRevX.10.031055},
	note = {\href{https://arxiv.org/abs/1811.00536}{arXiv:\allowbreak 1811.00536}},
    eprint = {1811.00536},
    archivePrefix = {arXiv},
}

@article{KT17,
	author="Kapustin, Anton
	and Thorngren, Ryan",
	title="Fermionic {SPT} phases in higher dimensions and bosonization",
	journal={J. High Energy Phys.},
	year="2017",
	volume="2017",
	number="10",
	pages={080},
	issn="1029-8479",
	doi="10.1007/JHEP10(2017)080",
	note = {\href{https://arxiv.org/abs/1701.08264}{arXiv:\allowbreak 1701.08264}},
    fjournal = {Journal of High Energy Physics},
    eprint = {1701.08264},
    archivePrefix = {arXiv},
}

@article{JFT20,
    AUTHOR = {Johnson-Freyd, Theo and Treumann, David},
     TITLE = {{${\rm H}^4({\rm Co}_0; {\bf Z}) = {\bf Z}/24$}},
   JOURNAL = {Int. Math. Res. Not. IMRN},
  FJOURNAL = {International Mathematics Research Notices. IMRN},
    VOLUME = {2020},
      YEAR = {2020},
    NUMBER = {21},
     PAGES = {7873--7907},
      ISSN = {1073-7928},
   MRCLASS = {20D08 (20J06 55R40)},
  MRNUMBER = {4176841},
MRREVIEWER = {David Benson},
       DOI = {10.1093/imrn/rny219},
    note = {\href{https://arxiv.org/abs/1707.07587}{arXiv:\allowbreak 1707.07587}},
    eprint = {1707.07587},
    archivePrefix = {arXiv},
}

@book{Pet68,
  author    = {F. P. Peterson},
  title     = {Lectures on Cobordism Theory},
  publisher = {Kinokuniya Book Store Co., Ltd.},
  year      = {1968},
  series    = {Lectures in Mathematics},
}

@article{ABP69,
    AUTHOR = {Anderson, D. W. and Brown, Jr., E. H. and Peterson, F. P.},
     TITLE = {Pin cobordism and related topics},
   JOURNAL = {Comment. Math. Helv.},
  FJOURNAL = {Commentarii Mathematici Helvetici},
    VOLUME = {44},
      YEAR = {1969},
     PAGES = {462--468},
      ISSN = {0010-2571},
   MRCLASS = {57.10},
  MRNUMBER = {261613},
MRREVIEWER = {P. S. Landweber},
       DOI = {10.1007/BF02564545},
       URL = {https://doi.org/10.1007/BF02564545},
}

@article{Sto88,
    AUTHOR = {Stolz, Stephan},
     TITLE = {Exotic structures on {$4$}-manifolds detected by spectral
              invariants},
   JOURNAL = {Invent. Math.},
  FJOURNAL = {Inventiones Mathematicae},
    VOLUME = {94},
      YEAR = {1988},
    NUMBER = {1},
     PAGES = {147--162},
      ISSN = {0020-9910},
   MRCLASS = {57R55 (57R85 58G25)},
  MRNUMBER = {958593},
MRREVIEWER = {J\"{u}rgen Eichhorn},
       DOI = {10.1007/BF01394348},
       URL = {https://doi.org/10.1007/BF01394348},
}

@article{Gia73a,
    AUTHOR = {Giambalvo, V.},
     TITLE = {Pin and {$\mathrm{Pin}'$} cobordism},
   JOURNAL = {Proc. Amer. Math. Soc.},
  FJOURNAL = {Proceedings of the American Mathematical Society},
    VOLUME = {39},
      YEAR = {1973},
     PAGES = {395--401},
      ISSN = {0002-9939},
   MRCLASS = {57D90},
  MRNUMBER = {321123},
MRREVIEWER = {J. A. Schafer},
       DOI = {10.2307/2039653},
       URL = {https://doi.org/10.2307/2039653},
}

@book{EKMM97,
    AUTHOR = {Elmendorf, A. D. and Kriz, I. and Mandell, M. A. and May, J.
              P.},
     TITLE = {Rings, modules, and algebras in stable homotopy theory},
    SERIES = {Mathematical Surveys and Monographs},
    VOLUME = {47},
      NOTE = {With an appendix by M. Cole},
 PUBLISHER = {American Mathematical Society, Providence, RI},
      YEAR = {1997},
     PAGES = {xii+249},
      ISBN = {0-8218-0638-6},
   MRCLASS = {55N20 (19D10 19D55 55P42 55T25)},
  MRNUMBER = {1417719},
MRREVIEWER = {Donald M. Davis},
       DOI = {10.1090/surv/047},
       URL = {https://doi.org/10.1090/surv/047},
}

@article{MM02,
    AUTHOR = {Mandell, M. A. and May, J. P.},
     TITLE = {Equivariant orthogonal spectra and {$S$}-modules},
   JOURNAL = {Mem. Amer. Math. Soc.},
  FJOURNAL = {Memoirs of the American Mathematical Society},
    VOLUME = {159},
      YEAR = {2002},
    NUMBER = {755},
     PAGES = {x+108},
      ISSN = {0065-9266},
   MRCLASS = {55P91 (18E30 55P42 55P43 55P48)},
  MRNUMBER = {1922205},
MRREVIEWER = {J. P. C. Greenlees},
       DOI = {10.1090/memo/0755},
	note = "\url{http://www.math.uchicago.edu/~may/PAPERS/MMMFinal.pdf}"
}

@misc{Jen05,
	author = {Jerome A. Jenquin},
	title = {Classical {C}hern-{S}imons on manifolds with spin structure},
	year = {2005},
	note = {\href{https://arxiv.org/abs/math/0504524}{arXiv:\allowbreak math/0504524}},
    eprint = {math/0504524},
    archivePrefix = {arXiv},
}

@article{GW14,
  title = {Symmetry-protected topological orders for interacting fermions: Fermionic topological nonlinear $\ensuremath{\sigma}$ models and a special group supercohomology theory},
  author = {Gu, Zheng-Cheng and Wen, Xiao-Gang},
  journal = {Phys. Rev. B},
  volume = {90},
  number = {11},
  pages = {115141},
  numpages = {59},
  year = {2014},
  publisher = {American Physical Society},
  doi = {10.1103/PhysRevB.90.115141},
	note = {\href{https://arxiv.org/abs/1201.2648}{arXiv:\allowbreak 1201.2648}},
    eprint = {1201.2648},
    archivePrefix = {arXiv},
}

@article{Fre08,
    AUTHOR = {Freed, Daniel S.},
     TITLE = {Pions and generalized cohomology},
   JOURNAL = {J. Differential Geom.},
  FJOURNAL = {Journal of Differential Geometry},
    VOLUME = {80},
      YEAR = {2008},
    NUMBER = {1},
     PAGES = {45--77},
      ISSN = {0022-040X},
   MRCLASS = {53C08 (19L50 53C80 55N20 55S45 81T50 81V05)},
  MRNUMBER = {2434259},
MRREVIEWER = {Michael A. Hill},
	note = {\href{https://arxiv.org/abs/hep-th/0607134}{arXiv:\allowbreak hep-th/0607134}},
    eprint = {hep-th/0607134},
    archivePrefix = {arXiv},
}

@article{AGG14,
    AUTHOR = {Antieau, Benjamin and Gepner, David and Gómez, José Manuel},
     TITLE = {Actions of {$K(\pi,n)$} spaces on {$K$}-theory and uniqueness
              of twisted {$K$}-theory},
   JOURNAL = {Trans. Amer. Math. Soc.},
  FJOURNAL = {Transactions of the American Mathematical Society},
    VOLUME = {366},
      YEAR = {2014},
    NUMBER = {7},
     PAGES = {3631--3648},
      ISSN = {0002-9947},
   MRCLASS = {19L50 (55N15)},
  MRNUMBER = {3192610},
MRREVIEWER = {Bin Zhang},
       DOI = {10.1090/S0002-9947-2014-05937-0},
note = {\href{https://arxiv.org/abs/1106.5099}{arXiv:\allowbreak 1106.5099}},
    eprint = {1106.5099},
    archivePrefix = {arXiv},
}

@article{ABP67,
    AUTHOR = {Anderson, D. W. and Brown, Jr., E. H. and Peterson, F. P.},
     TITLE = {The structure of the {S}pin cobordism ring},
   JOURNAL = {Ann. of Math. (2)},
  FJOURNAL = {Annals of Mathematics. Second Series},
    VOLUME = {86},
      YEAR = {1967},
     PAGES = {271--298},
      ISSN = {0003-486X},
   MRCLASS = {57.10},
  MRNUMBER = {219077},
MRREVIEWER = {P. S. Landweber},
       DOI = {10.2307/1970690},
       URL = {https://doi.org/10.2307/1970690},
}

@article{Cul21,
    AUTHOR = {Culver, D.},
     TITLE = {The {A}dams spectral sequence for 3-local {$\mathrm {tmf}$}},
   JOURNAL = {J. Homotopy Relat. Struct.},
  FJOURNAL = {Journal of Homotopy and Related Structures},
    VOLUME = {16},
      YEAR = {2021},
    NUMBER = {1},
     PAGES = {1--40},
      ISSN = {2193-8407},
   MRCLASS = {55T15 (55N34)},
  MRNUMBER = {4225505},
MRREVIEWER = {Tilman Bauer},
       DOI = {10.1007/s40062-020-00271-3},
note = {\href{https://arxiv.org/abs/1902.04230}{arXiv:\allowbreak 1902.04230}},
    eprint = {1902.04230},
    archivePrefix = {arXiv},
}

@book{MQRT77,
    AUTHOR = {May, J. Peter},
     TITLE = {{$E\sb{\infty }$} ring spaces and {$E\sb{\infty }$} ring
              spectra},
    SERIES = {Lecture Notes in Mathematics, Vol. 577},
     NOTE = {With contributions by Frank Quinn, Nigel Ray, and Jørgen
              Tornehave. \url{http://www.math.uchicago.edu/~may/BOOKS/e_infty.pdf}},
 PUBLISHER = {Springer-Verlag, Berlin-New York},
      YEAR = {1977},
     PAGES = {268},
   MRCLASS = {55D35},
  MRNUMBER = {0494077},
MRREVIEWER = {J. Stasheff},
}

@book{MS06,
    AUTHOR = {May, J. P. and Sigurdsson, J.},
     TITLE = {Parametrized homotopy theory},
    SERIES = {Mathematical Surveys and Monographs},
    VOLUME = {132},
 PUBLISHER = {American Mathematical Society, Providence, RI},
      YEAR = {2006},
     PAGES = {x+441},
      ISBN = {978-0-8218-3922-5; 0-8218-3922-5},
   MRCLASS = {55P42 (19L99 55N20 55N22 55P91)},
  MRNUMBER = {2271789},
MRREVIEWER = {A. A. Ranicki},
       DOI = {10.1090/surv/132},
    note = {\url{http://www.math.uchicago.edu/~may/EXTHEORY/MaySig.pdf}}
}

@article{MM81,
    AUTHOR = {May, J. P. and Milgram, R. J.},
     TITLE = {The {B}ockstein and the {A}dams spectral sequences},
   JOURNAL = {Proc. Amer. Math. Soc.},
  FJOURNAL = {Proceedings of the American Mathematical Society},
    VOLUME = {83},
      YEAR = {1981},
    NUMBER = {1},
     PAGES = {128--130},
      ISSN = {0002-9939},
   MRCLASS = {55T15},
  MRNUMBER = {619997},
MRREVIEWER = {J. F. Adams},
       DOI = {10.2307/2043907},
}

@article{Gia76,
    AUTHOR = {Giambalvo, V.},
     TITLE = {Cobordism of spin manifolds with involution},
   JOURNAL = {Quart. J. Math. Oxford Ser. (2)},
  FJOURNAL = {The Quarterly Journal of Mathematics. Oxford. Second Series},
    VOLUME = {27},
      YEAR = {1976},
    NUMBER = {106},
     PAGES = {241--252},
      ISSN = {0033-5606},
   MRCLASS = {57D85},
  MRNUMBER = {400264},
MRREVIEWER = {R. E. Stong},
       DOI = {10.1093/qmath/27.2.241},
       URL = {https://doi.org/10.1093/qmath/27.2.241},
}

@article{Gia73b,
    AUTHOR = {Giambalvo, V.},
     TITLE = {Cobordism of line bundles with a relation},
   JOURNAL = {Illinois J. Math.},
  FJOURNAL = {Illinois Journal of Mathematics},
    VOLUME = {17},
      YEAR = {1973},
     PAGES = {442--449},
      ISSN = {0019-2082},
   MRCLASS = {57D85},
  MRNUMBER = {321121},
MRREVIEWER = {E. C. Turner},
       URL = {http://projecteuclid.org/euclid.ijm/1256051610},
}

@article{Bou79,
    AUTHOR = {Bousfield, A. K.},
     TITLE = {The localization of spectra with respect to homology},
   JOURNAL = {Topology},
  FJOURNAL = {Topology. An International Journal of Mathematics},
    VOLUME = {18},
      YEAR = {1979},
    NUMBER = {4},
     PAGES = {257--281},
      ISSN = {0040-9383},
   MRCLASS = {55N20 (55N15 55P60)},
  MRNUMBER = {551009},
MRREVIEWER = {Willi Meier},
       DOI = {10.1016/0040-9383(79)90018-1},
       URL = {https://doi.org/10.1016/0040-9383(79)90018-1},
}

@article{Rav84,
    AUTHOR = {Ravenel, Douglas C.},
     TITLE = {Localization with respect to certain periodic homology
              theories},
   JOURNAL = {Amer. J. Math.},
  FJOURNAL = {American Journal of Mathematics},
    VOLUME = {106},
      YEAR = {1984},
    NUMBER = {2},
     PAGES = {351--414},
      ISSN = {0002-9327},
   MRCLASS = {55P42},
  MRNUMBER = {737778},
MRREVIEWER = {J. F. Adams},
       DOI = {10.2307/2374308},
       URL = {https://doi.org/10.2307/2374308},
}

@misc{DY23,
        author = {Arun Debray and Matthew Yu},
        title = {Adams spectral sequences for non-vector-bundle {T}hom spectra},
        year = {2023},
        note = {\href{https://arxiv.org/abs/2305.01678}{arXiv:\allowbreak 2305.01678}},
    eprint = {2305.01678},
    archivePrefix = {arXiv},
}

@article{BB96,
    AUTHOR = {Bayen, Dilip and Bruner, Robert R.},
     TITLE = {Real connective {$K$}-theory and the quaternion group},
   JOURNAL = {Trans. Amer. Math. Soc.},
  FJOURNAL = {Transactions of the American Mathematical Society},
    VOLUME = {348},
      YEAR = {1996},
    NUMBER = {6},
     PAGES = {2201--2216},
      ISSN = {0002-9947},
   MRCLASS = {19L41 (19L47 55M05 55N15 55Q91 55R35)},
  MRNUMBER = {1329527},
       DOI = {10.1090/S0002-9947-96-01516-4},
       URL = {https://doi.org/10.1090/S0002-9947-96-01516-4},
}

@phdthesis{Bay94,
    AUTHOR = {Bayen, Dilip Kumar},
     TITLE = {The connective real {K}-homology of finite groups},
      school = {Wayne State University},
      YEAR = {1994},
}

@article{DDHM24,
AUTHOR = {Debray, Arun and Dierigl, Markus and Heckman, Jonathan J. and
              Montero, Miguel},
     TITLE = {The chronicles of {IIB}ordia: dualities, bordisms, and the
              {S}wampland},
   JOURNAL = {Adv. Theor. Math. Phys.},
  FJOURNAL = {Advances in Theoretical and Mathematical Physics},
    VOLUME = {28},
      YEAR = {2024},
    NUMBER = {3},
     PAGES = {805--1025},
      ISSN = {1095-0761},
   MRCLASS = {81 (55)},
  MRNUMBER = {4818749},
        note = {\href{https://arxiv.org/abs/2302.00007}{arXiv:\allowbreak 2302.00007}},
    eprint = {2302.00007},
    archivePrefix = {arXiv},
}

@misc{BGR24,
	author = {Noé Bárcenas and Luis Eduardo García-Hernández and Raphael Reinauer},
	title = {The {Gromov-Lawson-Rosenberg} Conjecture for {$\mathbb{Z}/4\times\mathbb{Z}/4$}},
	year = {2024},
	note = {\href{https://arxiv.org/abs/2408.07895}{arXiv:\allowbreak 2408.07895}},
    eprint = {2408.07895},
    archivePrefix = {arXiv},
}

@article{BHHM08,
    AUTHOR = {Behrens, M. and Hill, M. and Hopkins, M. J. and Mahowald, M.},
     TITLE = {On the existence of a {$v^{32}_2$}-self map on {$M(1,4)$} at
              the prime 2},
   JOURNAL = {Homology Homotopy Appl.},
  FJOURNAL = {Homology, Homotopy and Applications},
    VOLUME = {10},
      YEAR = {2008},
    NUMBER = {3},
     PAGES = {45--84},
      ISSN = {1532-0073},
   MRCLASS = {55Q51 (55Q40 55T15)},
  MRNUMBER = {2475617},
MRREVIEWER = {Donald M. Davis},
       DOI = {10.4310/hha.2008.v10.n3.a4},
	note = {\href{https://arxiv.org/abs/0710.5426}{arXiv:\allowbreak 0710.5426}},
    eprint = {0710.5426},
    archivePrefix = {arXiv},
}

@article{AP76,
    AUTHOR = {Adams, J. F. and Priddy, S. B.},
     TITLE = {Uniqueness of {$B\mathrm{SO}$}},
   JOURNAL = {Math. Proc. Cambridge Philos. Soc.},
  FJOURNAL = {Mathematical Proceedings of the Cambridge Philosophical
              Society},
    VOLUME = {80},
      YEAR = {1976},
    NUMBER = {3},
     PAGES = {475--509},
      ISSN = {0305-0041},
   MRCLASS = {55D35},
  MRNUMBER = {431152},
MRREVIEWER = {V. P. Snaith},
       DOI = {10.1017/S0305004100053111},
       URL = {https://doi-org.ezproxy.lib.utexas.edu/10.1017/S0305004100053111},
}

@article{BM23,
        AUTHOR = {Buchanan, Jonathan and McKean, Stephen},
     TITLE = {K{S}p-characteristic classes determine {S}pin{$^h$} cobordism},
   JOURNAL = {Algebr. Geom. Topol.},
  FJOURNAL = {Algebraic \& Geometric Topology},
    VOLUME = {26},
      YEAR = {2026},
    NUMBER = {2},
     PAGES = {485--551},
      ISSN = {1472-2747,1472-2739},
   MRCLASS = {19L41 (53C27)},
  MRNUMBER = {5034306},
       DOI = {10.2140/agt.2026.26.485},
       URL = {https://doi.org/10.2140/agt.2026.26.485},
        note = {\href{https://arxiv.org/abs/2312.08209}{arXiv:\allowbreak 2312.08209}},
    eprint = {2312.08209},
    archivePrefix = {arXiv},
}

@phdthesis{Yu95,
    AUTHOR = {Yu, Cherng-Yih},
     TITLE = {The connective real {K}-theory of elementary abelian 2-groups},
      YEAR = {1995},
    school = {The University of Notre Dame}
}

@misc{BottSpiral,
        	author = {Arun Debray and Cameron Krulewski and Luuk Stehouwer},
	title = {Unraveling the {B}ott spiral},
	year = {2026},
	note = {\href{https://arxiv.org/abs/2605.00316}{arXiv:\allowbreak 2605.00316}},
    eprint = {2605.00316},
    archivePrefix = {arXiv},
}

@article{GMM68,
    AUTHOR = {Gitler, S. and Mahowald, M. and Milgram, R. James},
     TITLE = {The nonimmersion problem for {$RP^{n}$} and higher-order
              cohomology operations},
   JOURNAL = {Proc. Nat. Acad. Sci. U.S.A.},
  FJOURNAL = {Proceedings of the National Academy of Sciences of the United
              States of America},
    VOLUME = {60},
      YEAR = {1968},
     PAGES = {432--437},
      ISSN = {0027-8424},
   MRCLASS = {57.20},
  MRNUMBER = {227997},
MRREVIEWER = {K. Aoki},
       DOI = {10.1073/pnas.60.2.432},
       URL = {https://doi-org.ezproxy.lib.utexas.edu/10.1073/pnas.60.2.432},
}

@misc{ZJSNew,
        author = {Wan, Zheyan and Wang, Juven and Yau, Shing-Tung},
        title = {Three Families and {P}ontryagin Class Mod 3: Topological Order via Symmetry Extension},
	year = {2026},
	note = {\href{https://arxiv.org/abs/2605.26202}{arXiv:\allowbreak 2605.26202}},
    eprint = {2605.26202},
    archivePrefix = {arXiv},
}

@misc{DYY3,
	author = {Arun Debray and Matthew Yu and Weicheng Ye},
	title = {Symmetric Gapped States and Symmetry-Enforced Gaplessness in 3 Dimensions},
	year = {2026},
	note = {\href{https://arxiv.org/abs/2602.12335}{arXiv:\allowbreak 2602.12335}},
    eprint = {2602.12335},
    archivePrefix = {arXiv},
}

@article{DYY25,
    author = "Debray, Arun and Ye, Weicheng and Yu, Matthew",
    title = "Global Structure in the Presence of a Topological Defect",
    eprint = "2501.18399",
    archivePrefix = "arXiv",
    primaryClass = "math-ph",
    doi = "10.1007/s00220-026-05629-x",
    journal = "Commun. Math. Phys.",
    volume = "407",
    number = "8",
    pages = "159",
    year = "2026",
	note = {\href{https://arxiv.org/abs/2501.18399}{arXiv:\allowbreak 2501.18399}}
}

@book{Mar83,
    AUTHOR = {Margolis, H. R.},
     TITLE = {Spectra and the {S}teenrod algebra},
    SERIES = {North-Holland Mathematical Library},
    VOLUME = {29},
      NOTE = {Modules over the Steenrod algebra and the stable homotopy
              category},
 PUBLISHER = {North-Holland Publishing Co., Amsterdam},
      YEAR = {1983},
     PAGES = {xix+489},
      ISBN = {0-444-86516-0},
   MRCLASS = {55-02 (55Sxx)},
  MRNUMBER = {738973},
MRREVIEWER = {Andrzej Kozlowski},
}

@article{Pea14,
    AUTHOR = {Pearson, Paul Thomas},
     TITLE = {The connective real {$K$}-theory of {B}rown-{G}itler spectra},
   JOURNAL = {Algebr. Geom. Topol.},
  FJOURNAL = {Algebraic \& Geometric Topology},
    VOLUME = {14},
      YEAR = {2014},
    NUMBER = {1},
     PAGES = {597--625},
      ISSN = {1472-2747},
   MRCLASS = {55P42 (19L99 55N15 55T25)},
  MRNUMBER = {3158769},
MRREVIEWER = {Richard John Steiner},
       DOI = {10.2140/agt.2014.14.597},
       URL = {https://doi.org/10.2140/agt.2014.14.597},
}

@article{DY23a,
    AUTHOR = {Décoppet, Thibault D. and Yu, Matthew},
     TITLE = {Gauging noninvertible defects: a 2-categorical perspective},
   JOURNAL = {Lett. Math. Phys.},
  FJOURNAL = {Letters in Mathematical Physics},
    VOLUME = {113},
      YEAR = {2023},
    NUMBER = {2},
     PAGES = {36},
      ISSN = {0377-9017},
   MRCLASS = {18M15 (18M20 18N10)},
  MRNUMBER = {4565887},
MRREVIEWER = {Raffael Stenzel},
       DOI = {10.1007/s11005-023-01655-1},
	note = {\href{https://arxiv.org/abs/2211.08436}{arXiv:\allowbreak 2211.08436}},
    eprint = {2211.08436},
    archivePrefix = {arXiv},
    pagetotal = {42},
}

@article{Dec24,
    AUTHOR = {D{\'{e}}coppet, Thibault Didier},
     TITLE = {Extension theory and fermionic strongly fusion 2-categories},
      NOTE = {With an appendix by Décoppet and Theo Johnson-Freyd. \href{https://arxiv.org/abs/2403.03211}{arXiv:\allowbreak 2403.03211}},
   JOURNAL = {SIGMA Symmetry Integrability Geom. Methods Appl.},
  FJOURNAL = {SIGMA. Symmetry, Integrability and Geometry. Methods and
              Applications},
    VOLUME = {20},
      YEAR = {2024},
     PAGES = {092},
   MRCLASS = {18M20 (18N25)},
  MRNUMBER = {4843381},
       DOI = {10.3842/SIGMA.2024.092},
    eprint = {2403.03211},
    archivePrefix = {arXiv},
    pagetotal = {20},
}

@misc{DHJNPPRY24,
	author = {Thibault D. Décoppet and Peter Huston and Theo Johnson-Freyd and Dmitri Nikshych and David Penneys and Julia Plavnik and David Reutter and Matthew Yu},
	title = {The Classification of Fusion 2-Categories},
	year = {2024},
	note = {\href{https://arxiv.org/abs/2411.05907}{arXiv:\allowbreak 2411.05907}},
    eprint = {2411.05907},
    archivePrefix = {arXiv},
}

@article{MP84,
    AUTHOR = {Mitchell, Stephen A. and Priddy, Stewart B.},
     TITLE = {Symmetric product spectra and splittings of classifying
              spaces},
   JOURNAL = {Amer. J. Math.},
  FJOURNAL = {American Journal of Mathematics},
    VOLUME = {106},
      YEAR = {1984},
    NUMBER = {1},
     PAGES = {219--232},
      ISSN = {0002-9327},
   MRCLASS = {55P42 (20J06)},
  MRNUMBER = {729761},
MRREVIEWER = {Frederick Cohen},
       DOI = {10.2307/2374436},
       URL = {https://doi.org/10.2307/2374436},
}

@article{Mil24,
    AUTHOR = {Mills, Keith},
     TITLE = {The structure of the spin{$^h$} bordism spectrum},
   JOURNAL = {Proc. Amer. Math. Soc.},
  FJOURNAL = {Proceedings of the American Mathematical Society},
    VOLUME = {152},
      YEAR = {2024},
    NUMBER = {8},
     PAGES = {3605--3616},
      ISSN = {0002-9939},
   MRCLASS = {55N22 (55S10 55T15)},
  MRNUMBER = {4767288},
       DOI = {10.1090/proc/16748},
        note = {\href{https://arxiv.org/abs/2306.17709}{arXiv:\allowbreak 2306.17709}},
    eprint = {2306.17709},
    archivePrefix = {arXiv},
}

@misc{BDDHM25,
	author = {Noah Braeger and Arun Debray and Markus Dierigl and Jonathan J. Heckman and Miguel Montero},
	title = {Cobordism Utopia: {U}-Dualities, Bordisms, and the {S}wampland},
	year = {2025},
	note = {\href{https://arxiv.org/abs/2505.15885}{arXiv:\allowbreak 2505.15885}},
    eprint = {2505.15885},
    archivePrefix = {arXiv},
}

@misc{TY25,
	author = {Daniel Teixeira and Matthew Yu},
	title = {Mutual Influence of Symmetries and Topological Field Theories},
	year = {2025},
	note = {\href{https://arxiv.org/abs/2507.06304}{arXiv:\allowbreak 2507.06304}},
    eprint = {2507.06304},
    archivePrefix = {arXiv},
}

@misc{Fra11, 
        author = {John Francis},
        title = {Integrals on spin manifolds and the {$K$}-theory of {$\mathrm K(\mathbb Z, 4)$}},
        year = {2011},
        note = {\url{https://sites.math.northwestern.edu/~jnkf/writ/bspin2011.pdf}}
}

@article{Kne25,
	author="Christian Kneißl",
	title="Spin cobordism and the gauge group of type {I}/heterotic string theory",
	journal={J. High Energy Phys.},
	year="2025",
	volume="2025",
	number={1},
	note = {\href{https://arxiv.org/abs/2407.20333}{arXiv:\allowbreak 2407.20333}},
    doi = {10.1007/JHEP01(2025)181},
    pages = {181},
    fjournal = {Journal of High Energy Physics},
    eprint = {2407.20333},
    archivePrefix = {arXiv},
}

@book{BG10,
    AUTHOR = {Bruner, Robert R. and Greenlees, J. P. C.},
     TITLE = {Connective real {$K$}-theory of finite groups},
    SERIES = {Mathematical Surveys and Monographs},
    VOLUME = {169},
 PUBLISHER = {American Mathematical Society, Providence, RI},
      YEAR = {2010},
     PAGES = {vi+318},
      ISBN = {978-0-8218-5189-0},
   MRCLASS = {19L41 (13D45 19-02 19L47 19L64 55N15 55N91)},
  MRNUMBER = {2723113},
       DOI = {10.1090/surv/169},
       URL = {https://doi.org/10.1090/surv/169},
}

@article{Liu62,
    AUTHOR = {Liulevicius, Arunas L.},
     TITLE = {The cohomology of a subalgebra of the {S}teenrod algebra},
   JOURNAL = {Trans. Amer. Math. Soc.},
  FJOURNAL = {Transactions of the American Mathematical Society},
    VOLUME = {104},
      YEAR = {1962},
     PAGES = {443--449},
      ISSN = {0002-9947},
   MRCLASS = {55.34},
  MRNUMBER = {149476},
MRREVIEWER = {H. Suzuki},
       DOI = {10.2307/1993789},
       URL = {https://doi.org/10.2307/1993789},
}

@misc{Kur25,
	author = {Naoki Kuroda},
	title = {Computations of {Spin-Sp(4)}, {Spin-SU(8)}, and {Spin-Spin(16)} bordism groups
  in dimensions up to 7},
	year = {2025},
	note = {\href{https://arxiv.org/abs/2504.15014}{arXiv:\allowbreak 2504.15014}},
    eprint = {2504.15014},
    archivePrefix = {arXiv},
}

@misc{Bak20,
        author = {Andrew Baker},
        title = {Homotopy theory of modules over a commutative {$S$}-algebra: some tools
  and examples},
        year = {2020},
        note = {\href{https://arxiv.org/abs/2003.12003}{arXiv:\allowbreak 2003.12003}},
    eprint = {2003.12003},
    archivePrefix = {arXiv},
}

@article{Ada23,
 AUTHOR = {Adamyk, Katharine},
     TITLE = {Classifying and extending {$Q_0$} local {$\mathcal{A}(1)$}-modules},
   JOURNAL = {New York J. Math.},
  FJOURNAL = {New York Journal of Mathematics},
    VOLUME = {29},
      YEAR = {2023},
     PAGES = {939--980},
   MRCLASS = {55S10 (16T05 55T25)},
  MRNUMBER = {4627652},
        note = {\href{https://arxiv.org/abs/2107.02837}{arXiv:\allowbreak 2107.02837}},
    eprint = {2107.02837},
    archivePrefix = {arXiv},
}

@misc{BJ04,
	author = {Hans Joachim Baues and Mamuka Jibladze},
	title = {Computation of the {$E_3$}-term of the {A}dams spectral sequence},
	year = {2004},
	note = {\href{https://arxiv.org/abs/math/0407045}{arXiv:\allowbreak math/0407045}},
    eprint = {math/0407045},
    archivePrefix = {arXiv},
}

@article{SO12,
    AUTHOR = {Spitzweck, Markus and {\O}stvær, Paul Arne},
     TITLE = {Motivic twisted {$K$}-theory},
   JOURNAL = {Algebr. Geom. Topol.},
  FJOURNAL = {Algebraic \& Geometric Topology},
    VOLUME = {12},
      YEAR = {2012},
    NUMBER = {1},
     PAGES = {565--599},
      ISSN = {1472-2747},
   MRCLASS = {19L50 (14F42 55P43)},
  MRNUMBER = {2916287},
MRREVIEWER = {Niko Naumann},
       DOI = {10.2140/agt.2012.12.565},
	note = {\href{https://arxiv.org/abs/1008.4915}{arXiv:\allowbreak 1008.4915}},
    eprint = {1008.4915},
    archivePrefix = {arXiv},
}

@misc{vN25,
	author = {van Nigtevecht, Sven},
	title = {An introduction to filtered and synthetic spectra},
	year = {2025},
	note = {\href{https://arxiv.org/abs/2509.21127}{arXiv:\allowbreak 2509.21127}},
    eprint = {2509.21127},
    archivePrefix = {arXiv},
}

@misc{LWX24,
	author = {Weinan Lin and Guozhen Wang and Zhouli Xu},
	title = {On the Last {K}ervaire Invariant Problem},
	year = {2024},
	note = {\href{https://arxiv.org/abs/2412.10879}{arXiv:\allowbreak 2412.10879}},
    eprint = {2412.10879},
    archivePrefix = {arXiv},
}

@article{CD24,
	AUTHOR = {Carrick, Christian and Davies, Jack Morgan},
     TITLE = {A synthetic approach to detecting {$v_1$}-periodic families},
   JOURNAL = {Trans. Amer. Math. Soc.},
  FJOURNAL = {Transactions of the American Mathematical Society},
    VOLUME = {379},
      YEAR = {2026},
    NUMBER = {3},
     PAGES = {1985--2031},
      ISSN = {0002-9947,1088-6850},
   MRCLASS = {55T05 (55N20 55P42 55T15 55T25)},
  MRNUMBER = {5029817},
       DOI = {10.1090/tran/9526},
       URL = {https://doi.org/10.1090/tran/9526},
	note = {\href{https://arxiv.org/abs/2401.16508}{arXiv:\allowbreak 2401.16508}},
    eprint = {2401.16508},
    archivePrefix = {arXiv},
}

@book{Ada74,
    AUTHOR = {Adams, J. F.},
     TITLE = {Stable homotopy and generalised homology},
    SERIES = {Chicago Lectures in Mathematics},
 PUBLISHER = {University of Chicago Press, Chicago, Ill.-London},
      YEAR = {1974},
     PAGES = {x+373},
   MRCLASS = {55B20 (55E10)},
  MRNUMBER = {402720},
MRREVIEWER = {P. S. Landweber},
}

@article{BG03,
    AUTHOR = {Bruner, R. R. and Greenlees, J. P. C.},
     TITLE = {The connective {$K$}-theory of finite groups},
   JOURNAL = {Mem. Amer. Math. Soc.},
  FJOURNAL = {Memoirs of the American Mathematical Society},
    VOLUME = {165},
      YEAR = {2003},
    NUMBER = {785},
     PAGES = {viii+127},
      ISSN = {0065-9266},
   MRCLASS = {19L41 (19L47 20J05 55N91 55U25)},
  MRNUMBER = {1997161},
       DOI = {10.1090/memo/0785},
       URL = {https://doi-org.ezproxy.lib.utexas.edu/10.1090/memo/0785},
}

@incollection{Sul96,
    AUTHOR = {Sullivan, D. P.},
     TITLE = {Triangulating and smoothing homotopy equivalences and
              homeomorphisms. {G}eometric {T}opology {S}eminar {N}otes},
 BOOKTITLE = {The {H}auptvermutung book},
    SERIES = {$K$-Monogr. Math.},
    VOLUME = {1},
     PAGES = {69--103},
 PUBLISHER = {Kluwer Acad. Publ., Dordrecht},
      YEAR = {1996},
   MRCLASS = {57Q15 (57Q25 57R10)},
  MRNUMBER = {1434103},
MRREVIEWER = {Oliver Attie},
       DOI = {10.1007/978-94-017-3343-4_3},
       URL = {https://doi.org/10.1007/978-94-017-3343-4_3},
}

@article{Rub79,
    AUTHOR = {Rubinsztein, Ryszard L.},
     TITLE = {Some remarks on the cohomotopy of {${\bf R}P^{\infty }$}},
   JOURNAL = {Quart. J. Math. Oxford Ser. (2)},
  FJOURNAL = {The Quarterly Journal of Mathematics. Oxford. Second Series},
    VOLUME = {30},
      YEAR = {1979},
    NUMBER = {117},
     PAGES = {119--128},
      ISSN = {0033-5606},
   MRCLASS = {55P42 (55Q55)},
  MRNUMBER = {528896},
MRREVIEWER = {Stewart B. Priddy},
       DOI = {10.1093/qmath/30.1.119},
       URL = {https://doi.org/10.1093/qmath/30.1.119},
}

@article{Muk04,
    AUTHOR = {Mukerjee, Himadri Kumar},
     TITLE = {Classification of homotopy real {M}ilnor manifolds},
   JOURNAL = {Topology Appl.},
  FJOURNAL = {Topology and its Applications},
    VOLUME = {139},
      YEAR = {2004},
    NUMBER = {1--3},
     PAGES = {151--184},
      ISSN = {0166-8641},
   MRCLASS = {55P10 (55P25 57R19 57R67)},
  MRNUMBER = {2051103},
MRREVIEWER = {Alberto Cavicchioli},
       DOI = {10.1016/j.topol.2003.10.001},
       URL = {https://doi.org/10.1016/j.topol.2003.10.001},
}

@article{Muk06,
    AUTHOR = {Mukerjee, Himadri Kumar},
     TITLE = {Classification of homotopy {W}all's manifolds},
   JOURNAL = {Topology Appl.},
  FJOURNAL = {Topology and its Applications},
    VOLUME = {153},
      YEAR = {2006},
    NUMBER = {18},
     PAGES = {3467--3495},
      ISSN = {0166-8641},
   MRCLASS = {57R65 (55P10 55P25 57R19 57R67)},
  MRNUMBER = {2270599},
MRREVIEWER = {Laurence R. Taylor},
       DOI = {10.1016/j.topol.2006.03.002},
       URL = {https://doi.org/10.1016/j.topol.2006.03.002},
}

@article{Rud74,
    AUTHOR = {Rudyak, Yu. B.},
     TITLE = {The normal invariants of certain manifolds},
   JOURNAL = {Mat. Zametki},
  FJOURNAL = {Akademiya Nauk SSSR. Matematicheskie Zametki},
    VOLUME = {16},
      YEAR = {1974},
     PAGES = {763--769},
      ISSN = {0025-567X},
   MRCLASS = {57D12},
  MRNUMBER = {377916},
MRREVIEWER = {I. Pop},
}

@article{Wal68,
    AUTHOR = {Wall, C. T. C.},
     TITLE = {Free piecewise linear involutions on spheres},
   JOURNAL = {Bull. Amer. Math. Soc.},
  FJOURNAL = {Bulletin of the American Mathematical Society},
    VOLUME = {74},
      YEAR = {1968},
     PAGES = {554--558},
      ISSN = {0002-9904},
   MRCLASS = {57.05},
  MRNUMBER = {222905},
MRREVIEWER = {L. S. Husch},
       DOI = {10.1090/S0002-9904-1968-12006-3},
       URL = {https://doi.org/10.1090/S0002-9904-1968-12006-3},
}

@article{Muk03,
    AUTHOR = {Mukerjee, Himadri Kumar},
     TITLE = {Classification of homotopy {D}old manifolds},
   JOURNAL = {New York J. Math.},
  FJOURNAL = {New York Journal of Mathematics},
    VOLUME = {9},
      YEAR = {2003},
     PAGES = {271--293},
   MRCLASS = {57R67 (57R19)},
  MRNUMBER = {2016190},
MRREVIEWER = {R. E. Stong},
       URL = {http://nyjm.albany.edu:8000/j/2003/9_271.html},
}

@article{Wan22,
    AUTHOR = {Wang, Xueqi},
     TITLE = {On the classification of certain 1-connected 7-manifolds and
              related problems},
   JOURNAL = {Q. J. Math.},
  FJOURNAL = {The Quarterly Journal of Mathematics},
    VOLUME = {73},
      YEAR = {2022},
    NUMBER = {2},
     PAGES = {711--727},
      ISSN = {0033-5606},
   MRCLASS = {57K50 (55P15 57R15)},
  MRNUMBER = {4439806},
MRREVIEWER = {Yang Su},
       DOI = {10.1093/qmath/haab053},
       URL = {https://doi.org/10.1093/qmath/haab053},
	note = {\href{https://arxiv.org/abs/1810.08474}{arXiv:\allowbreak 1810.08474}},
    eprint = {1810.08474},
    archivePrefix = {arXiv},
}

@book{LdM71,
    AUTHOR = {L{\'{o}}pez de Medrano, S.},
     TITLE = {Involutions on manifolds},
    SERIES = {Ergebnisse der Mathematik und ihrer Grenzgebiete [Results in
              Mathematics and Related Areas], Band 59},
 PUBLISHER = {Springer-Verlag, New York-Heidelberg},
      YEAR = {1971},
     PAGES = {ix+103},
   MRCLASS = {57E30},
  MRNUMBER = {298698},
MRREVIEWER = {M. L. Curtis},
}

@book{MM79,
    AUTHOR = {Madsen, Ib and Milgram, R. James},
     TITLE = {The classifying spaces for surgery and cobordism of manifolds},
    SERIES = {Annals of Mathematics Studies, No. 92},
 PUBLISHER = {Princeton University Press, Princeton, NJ; University of Tokyo
              Press, Tokyo},
      YEAR = {1979},
     PAGES = {xii+279},
      ISBN = {0-691-08225-1},
   MRCLASS = {57N70 (55P47 55R35 57Q60)},
  MRNUMBER = {548575},
MRREVIEWER = {Yu. P. Solov\cprime ev},
}

@article{DY25,
	author = {Thibault D. Décoppet and Matthew Yu},
	 TITLE = {The classification of 3+1{D} symmetry enriched
              topological order},
   JOURNAL = {Comm. Math. Phys.},
  FJOURNAL = {Communications in Mathematical Physics},
    VOLUME = {407},
      YEAR = {2026},
    NUMBER = {11},
     PAGES = {Paper No. 238},
      ISSN = {0010-3616,1432-0916},
   MRCLASS = {18N10 (18M20)},
  MRNUMBER = {5135794},
       DOI = {10.1007/s00220-026-05753-8},
       URL = {https://doi.org/10.1007/s00220-026-05753-8},
	note = {\href{https://arxiv.org/abs/2509.10603}{arXiv:\allowbreak 2509.10603}},
    eprint = {2509.10603},
    archivePrefix = {arXiv},
}

@article{Sto94,
    AUTHOR = {Stolz, Stephan},
     TITLE = {Splitting certain {$M\, {\rm Spin}$}-module spectra},
   JOURNAL = {Topology},
  FJOURNAL = {Topology. An International Journal of Mathematics},
    VOLUME = {33},
      YEAR = {1994},
    NUMBER = {1},
     PAGES = {159--180},
      ISSN = {0040-9383},
   MRCLASS = {55N20 (53C21 55P42)},
  MRNUMBER = {1259520},
MRREVIEWER = {Donald M. Davis},
       DOI = {10.1016/0040-9383(94)90040-X},
       URL = {https://doi.org/10.1016/0040-9383(94)90040-X},
}

@incollection{JF25,
    AUTHOR = {Johnson-Freyd, Theo},
     TITLE = {{$(3+1)D$} topological orders with only a {$\mathbb{Z}_2$}-charged
              particle},
 BOOKTITLE = {Quantum symmetries: tensor categories, {TQFT}s, and vertex
              algebras},
    SERIES = {Contemp. Math.},
    VOLUME = {813},
     PAGES = {175--210},
 PUBLISHER = {Amer. Math. Soc., [Providence], RI},
      YEAR = {2025},
   MRCLASS = {18M20 (55S10 55T10 81T45)},
  MRNUMBER = {4886204},
       DOI = {10.1090/conm/813/16288},
       URL = {https://doi.org/10.1090/conm/813/16288},
	note = {\href{https://arxiv.org/abs/2011.11165}{arXiv:\allowbreak 2011.11165}},
    eprint = {2011.11165},
    archivePrefix = {arXiv},
}

@article{Mar74,
    AUTHOR = {Margolis, H. R.},
     TITLE = {Eilenberg-{M}ac {L}ane spectra},
   JOURNAL = {Proc. Amer. Math. Soc.},
  FJOURNAL = {Proceedings of the American Mathematical Society},
    VOLUME = {43},
      YEAR = {1974},
     PAGES = {409--415},
      ISSN = {0002-9939},
   MRCLASS = {55G10},
  MRNUMBER = {341488},
MRREVIEWER = {J. P. May},
       DOI = {10.2307/2038905},
}

@article{johnson2024minimal,
  title={Minimal nondegenerate extensions},
  author={Johnson-Freyd, Theo and Reutter, David},
  journal={J. Amer. Math. Soc.},
  fjournal={Journal of the American Mathematical Society},
  volume={37},
  number={1},
  pages={81--150},
  year={2024},
archivePrefix = {arXiv},
       eprint = {2105.15167},
 primaryClass = {math.QA},
       adsurl = {https://ui.adsabs.harvard.edu/abs/2021arXiv210515167J},
      note = {\href{https://arxiv.org/abs/2105.15167}{arXiv:\allowbreak 2105.15167}}
}

@ARTICLE{2025PhRvX..15b1062M,
       author = {{Ma}, Ruochen and {Zhang}, Jian-Hao and {Bi}, Zhen and {Cheng}, Meng and {Wang}, Chong},
        title = {Topological phases with average symmetries: the decohered, the disordered, and the intrinsic},
      journal = {Phys. Rev. X},
         year = 2025,
       volume = {15},
       number = {2},
          eid = {021062},
        pages = {021062},
          doi = {10.1103/PhysRevX.15.021062},
archivePrefix = {arXiv},
       note = {\href{https://arxiv.org/abs/2305.16399}{arXiv:\allowbreak 2305.16399}},
 primaryClass = {cond-mat.str-el},
       adsurl = {https://ui.adsabs.harvard.edu/abs/2025PhRvX..15b1062M},
    eprint = {2305.16399},
}

@misc{Debray:2025kfg,
    author = "Debray, Arun and Ye, Weicheng and Yu, Matthew",
    title = {How to build anomalous {(3+1)d} topological quantum field theories},
    eprint = "2510.24834",
    archivePrefix = "arXiv",
    primaryClass = "math-ph",
    year = "2025",
    note = {\href{https://arxiv.org/abs/2510.24834}{arXiv:\allowbreak 2510.24834}}
}

@misc{Car25,
	author = {Christian Carrick},
	title = {Slice spectral sequences through synthetic spectra},
	year = {2025},
	note = {\href{https://arxiv.org/abs/2510.19501}{arXiv:\allowbreak 2510.19501}},
    eprint = {2510.19501},
    archivePrefix = {arXiv},
}

@incollection{DG18,
    AUTHOR = {Debray, Arun and Gunningham, Sam},
     TITLE = {The {A}rf-{B}rown {TQFT} of pin{$^-$} surfaces},
 BOOKTITLE = {Topology and quantum theory in interaction},
    SERIES = {Contemp. Math.},
    VOLUME = {718},
     PAGES = {49--87},
 PUBLISHER = {Amer. Math. Soc., Providence, RI},
      YEAR = {2018},
   MRCLASS = {57R56 (57R15 81T25)},
  MRNUMBER = {3869641},
MRREVIEWER = {Carlos Segovia Gonz\'{a}lez},
       DOI = {10.1090/conm/718/14478},
        note = {\href{https://arxiv.org/abs/1803.11183}{arXiv:\allowbreak 1803.11183}},
    eprint = {1803.11183},
    archivePrefix = {arXiv},
}

@article{AM71,
    AUTHOR = {Adams, J. F. and Margolis, H. R.},
     TITLE = {Modules over the {S}teenrod algebra},
   JOURNAL = {Topology},
  FJOURNAL = {Topology. An International Journal of Mathematics},
    VOLUME = {10},
      YEAR = {1971},
     PAGES = {271--282},
      ISSN = {0040-9383},
   MRCLASS = {18G40 (55G10)},
  MRNUMBER = {294450},
MRREVIEWER = {A. K. Bousfield},
       DOI = {10.1016/0040-9383(71)90020-6},
       URL = {https://doi.org/10.1016/0040-9383(71)90020-6},
}

@misc{MO_Kn,
    TITLE = {Are the {AHSS} and {A}dams spectral sequence the same when computing connective {M}orava {K}-theory of a space?},
    AUTHOR = {Nicholas Kuhn and Tyler Lawson},
    NOTE = {\url{https://mathoverflow.net/q/359574}},
    year = {2020}
}

@incollection{Bru14,
    AUTHOR = {Bruner, Robert R.},
     TITLE = {Idempotents, localizations and {P}icard groups of {$
              A(1)$}-modules},
 BOOKTITLE = {An alpine expedition through algebraic topology},
    SERIES = {Contemp. Math.},
    VOLUME = {617},
     PAGES = {81--108},
 PUBLISHER = {Amer. Math. Soc., Providence, RI},
      YEAR = {2014},
   MRCLASS = {55S10 (19L41)},
  MRNUMBER = {3243394},
MRREVIEWER = {Constanze Roitzheim},
       DOI = {10.1090/conm/617/12346},
	note = {\href{https://arxiv.org/abs/1211.0213}{arXiv:\allowbreak 1211.0213}},
    eprint = {1211.0213},
    archivePrefix = {arXiv},
}

@article{Kit06,
    author = {Kitaev, Alexei},
    title = {Anyons in an exactly solved model and beyond},
    journal = {Ann. Physics},
    volume = {321},
    number = {1},
    pages = {2--111},
    year = {2006},
    doi = {10.1016/j.aop.2005.10.005},
    note = {\href{https://arxiv.org/abs/cond-mat/0506438}{arXiv:\allowbreak cond-mat/0506438}},
    eprint = {cond-mat/0506438},
    archivePrefix = {arXiv},
}

@inproceedings{Kit09,
    author = {Kitaev, Alexei},
    title = {Periodic table for topological insulators and superconductors},
    booktitle = {Advances in theoretical physics: Landau Memorial Conference},
    series = {AIP Conf. Proc.},
    volume = {1134},
    pages = {22--30},
    year = {2009},
    doi = {10.1063/1.3149495},
    note = {\href{https://arxiv.org/abs/0901.2686}{arXiv:\allowbreak 0901.2686}},
    eprint = {0901.2686},
    archivePrefix = {arXiv},
}

@phdthesis{Per78,
    AUTHOR = {Perloff, Michael},
     TITLE = {P{L}-Spin Cobordism},
    school = {Stanford University},
      YEAR = {1978},
}

@misc{2025arXiv251225069N,
       author = {{Ning}, Shang-Qiang and {Ren}, Xing-Yu and {Wang}, Qing-Rui and {Qi}, Yang and {Gu}, Zheng-Cheng},
        title = {Classification of interacting topological crystalline superconductors in three dimensions and beyond},
         year = 2025,
          eid = {arXiv:2512.25069},
archivePrefix = {arXiv},
       note = {\href{https://arxiv.org/abs/2512.25069}{arXiv:\allowbreak 2512.25069}},
 primaryClass = {cond-mat.str-el},
    eprint = {2512.25069},
}

@misc{Wan:2025ymd,
    author = "Wan, Zheyan and Wang, Juven",
    title = "{Anomalous (3+1)d Fermionic Topological Quantum Field Theories via Symmetry Extension}",
    eprint = "2512.25038",
    archivePrefix = "arXiv",
    primaryClass = "hep-th",
    month = "12",
    year = "2025",
    note = {\href{https://arxiv.org/abs/2512.25038}{arXiv:\allowbreak 2512.25038}}
}

@article{ZWYQG20,
  title = {Construction and classification of point-group symmetry-protected topological phases in two-dimensional interacting fermionic systems},
  author = {Zhang, Jian-Hao and Wang, Qing-Rui and Yang, Shuo and Qi, Yang and Gu, Zheng-Cheng},
  journal = {Phys. Rev. B},
  volume = {101},
  number = {10},
  pages = {100501},
  numpages = {6},
  year = {2020},
  publisher = {American Physical Society},
  doi = {10.1103/PhysRevB.101.100501},
  url = {https://link.aps.org/doi/10.1103/PhysRevB.101.100501}, 
    note = {\href{https://arxiv.org/abs/1909.05519}{arXiv:\allowbreak 1909.05519}},
    archivePrefix = "arXiv",
    primaryClass = "cond-mat.str-el",
    eprint = {1909.05519},
}

@article{JFY21,
    AUTHOR = {Johnson-Freyd, Theo and Yu, Matthew},
     TITLE = {Fusion 2-categories with no line operators are grouplike},
   JOURNAL = {Bull. Aust. Math. Soc.},
  FJOURNAL = {Bulletin of the Australian Mathematical Society},
    VOLUME = {104},
      YEAR = {2021},
    NUMBER = {3},
     PAGES = {434--442},
      ISSN = {0004-9727,1755-1633},
   MRCLASS = {18M20 (18N10)},
  MRNUMBER = {4338473},
       DOI = {10.1017/S0004972721000095},
       note = {\href{https://arxiv.org/abs/2010.07950}{arXiv:\allowbreak 2010.07950}},
    eprint = {2010.07950},
    archivePrefix = {arXiv},
}

@incollection{Ada69,
    AUTHOR = {Adams, J. F.},
     TITLE = {Lectures on generalised cohomology},
 BOOKTITLE = {Category {T}heory, {H}omology {T}heory and their
              {A}pplications, {III} ({B}attelle {I}nstitute {C}onference,
              {S}eattle, {W}ash., 1968, {V}ol. {T}hree)},
    SERIES = {Lecture Notes in Math., No. 99},
     PAGES = {1--138},
 PUBLISHER = {Springer, Berlin-New York},
      YEAR = {1969},
   MRCLASS = {55.30 (18.00)},
  MRNUMBER = {251716},
MRREVIEWER = {J.\ Siegel},
}

@article{JW73,
    AUTHOR = {Johnson, David Copeland and Wilson, W. Stephen},
     TITLE = {Projective dimension and {B}rown-{P}eterson homology},
   JOURNAL = {Topology},
  FJOURNAL = {Topology. An International Journal of Mathematics},
    VOLUME = {12},
      YEAR = {1973},
     PAGES = {327--353},
      ISSN = {0040-9383},
   MRCLASS = {57D90 (55B20)},
  MRNUMBER = {334257},
MRREVIEWER = {A.\ Liulevicius},
       DOI = {10.1016/0040-9383(73)90027-X},
       URL = {https://doi.org/10.1016/0040-9383(73)90027-X},
}

@article{Qui72,
    AUTHOR = {Quillen, Daniel},
     TITLE = {On the cohomology and {$K$}-theory of the general linear
              groups over a finite field},
   JOURNAL = {Ann. of Math. (2)},
  FJOURNAL = {Annals of Mathematics. Second Series},
    VOLUME = {96},
      YEAR = {1972},
     PAGES = {552--586},
      ISSN = {0003-486X},
   MRCLASS = {20J05 (18F25 55B15)},
  MRNUMBER = {315016},
MRREVIEWER = {D.\ W.\ Anderson},
       DOI = {10.2307/1970825},
       URL = {https://doi.org/10.2307/1970825},
}

@inproceedings{AH61,
    AUTHOR = {Atiyah, M. F. and Hirzebruch, F.},
     TITLE = {Vector bundles and homogeneous spaces},
 BOOKTITLE = {Proc. {S}ympos. {P}ure {M}ath., {V}ol. {III}},
     PAGES = {7--38},
 PUBLISHER = {Amer. Math. Soc., Providence, RI},
      YEAR = {1961},
   MRCLASS = {57.30 (14.52)},
  MRNUMBER = {139181},
MRREVIEWER = {R.\ Bott},
}

@article{Mau63,
    AUTHOR = {Maunder, C. R. F.},
     TITLE = {The spectral sequence of an extraordinary cohomology theory},
   JOURNAL = {Proc. Cambridge Philos. Soc.},
  FJOURNAL = {Proceedings of the Cambridge Philosophical Society},
    VOLUME = {59},
      YEAR = {1963},
     PAGES = {567--574},
      ISSN = {0008-1981},
   MRCLASS = {55.30 (55.52)},
  MRNUMBER = {150765},
MRREVIEWER = {Edgar\ H.\ Brown, Jr.},
       DOI = {10.1017/s0305004100037245},
       URL = {https://doi.org/10.1017/s0305004100037245},
}

@article{MNN17,
    AUTHOR = {Mathew, Akhil and Naumann, Niko and Noel, Justin},
     TITLE = {Nilpotence and descent in equivariant stable homotopy theory},
   JOURNAL = {Adv. Math.},
  FJOURNAL = {Advances in Mathematics},
    VOLUME = {305},
      YEAR = {2017},
     PAGES = {994--1084},
      ISSN = {0001-8708,1090-2082},
   MRCLASS = {55P91 (55P42)},
  MRNUMBER = {3570153},
MRREVIEWER = {Gregory\ Z.\ Arone},
       DOI = {10.1016/j.aim.2016.09.027},
       URL = {https://doi.org/10.1016/j.aim.2016.09.027},
        note = {\href{https://arxiv.org/abs/1507.06869}{arXiv:\allowbreak 1507.06869}},
    eprint = {1507.06869},
    archivePrefix = {arXiv},
}

@article{BK72,
    AUTHOR = {Bousfield, A. K. and Kan, D. M.},
     TITLE = {The homotopy spectral sequence of a space with coefficients in
              a ring},
   JOURNAL = {Topology},
  FJOURNAL = {Topology. An International Journal of Mathematics},
    VOLUME = {11},
      YEAR = {1972},
     PAGES = {79--106},
      ISSN = {0040-9383},
   MRCLASS = {55.52},
  MRNUMBER = {283801},
MRREVIEWER = {J.\ P.\ May},
       DOI = {10.1016/0040-9383(72)90024-9},
       URL = {https://doi.org/10.1016/0040-9383(72)90024-9},
}

@article{BK73a,
    AUTHOR = {Bousfield, A. K. and Kan, D. M.},
     TITLE = {A second quadrant homotopy spectral sequence},
   JOURNAL = {Trans. Amer. Math. Soc.},
  FJOURNAL = {Transactions of the American Mathematical Society},
    VOLUME = {177},
      YEAR = {1973},
     PAGES = {305--318},
      ISSN = {0002-9947,1088-6850},
   MRCLASS = {55H15},
  MRNUMBER = {372859},
MRREVIEWER = {Donald\ M.\ Davis},
       DOI = {10.2307/1996599},
       URL = {https://doi.org/10.2307/1996599},
}

@article{BK73b,
    AUTHOR = {Bousfield, A. K. and Kan, D. M.},
     TITLE = {Pairings and products in the homotopy spectral sequence},
   JOURNAL = {Trans. Amer. Math. Soc.},
  FJOURNAL = {Transactions of the American Mathematical Society},
    VOLUME = {177},
      YEAR = {1973},
     PAGES = {319--343},
      ISSN = {0002-9947,1088-6850},
   MRCLASS = {55H15},
  MRNUMBER = {372860},
MRREVIEWER = {Donald\ M.\ Davis},
       DOI = {10.2307/1996600},
       URL = {https://doi.org/10.2307/1996600},
}

@article{Bou03,
    AUTHOR = {Bousfield, A. K.},
     TITLE = {Cosimplicial resolutions and homotopy spectral sequences in
              model categories},
   JOURNAL = {Geom. Topol.},
  FJOURNAL = {Geometry \& Topology},
    VOLUME = {7},
      YEAR = {2003},
     PAGES = {1001--1053},
      ISSN = {1465-3060,1364-0380},
   MRCLASS = {55U35 (18G55 55P60 55T15)},
  MRNUMBER = {2026537},
MRREVIEWER = {Paul\ G.\ Goerss},
       DOI = {10.2140/gt.2003.7.1001},
       URL = {https://doi.org/10.2140/gt.2003.7.1001},
       note = {\href{https://arxiv.org/abs/math/0312531}{arXiv:\allowbreak math/0312531}},
    eprint = {math/0312531},
    archivePrefix = {arXiv},
}

@article{WG18,
  title = {Towards a Complete Classification of Symmetry-Protected Topological Phases for Interacting Fermions in Three Dimensions and a General Group Supercohomology Theory},
  author = {Wang, Qing-Rui and Gu, Zheng-Cheng},
  journal = {Phys. Rev. X},
  volume = {8},
  number = {1},
  pages = {011055},
  numpages = {29},
  year = {2018},
  publisher = {American Physical Society},
  doi = {10.1103/PhysRevX.8.011055},
  url = {https://link.aps.org/doi/10.1103/PhysRevX.8.011055},
    note = {\href{https://arxiv.org/abs/1703.10937}{arXiv:\allowbreak 1703.10937}},
    eprint = {1703.10937},
    archivePrefix = {arXiv},
}

@article{MathSmith,
    AUTHOR = {Debray, Arun and Devalapurkar, Sanath K. and Krulewski,
              Cameron and Liu, Yu Leon and Pacheco-Tallaj, Natalia and
              Thorngren, Ryan},
     TITLE = {The {S}mith fiber sequence and invertible field theories},
   JOURNAL = {Commun. Math. Phys.},
  FJOURNAL = {Communications in Mathematical Physics},
    VOLUME = {407},
      YEAR = {2026},
    NUMBER = {2},
     PAGES = {25},
      ISSN = {0010-3616,1432-0916},
   MRCLASS = {57R56 (55R12 57N70 81T45)},
  MRNUMBER = {5015202},
       DOI = {10.1007/s00220-025-05505-0},
       URL = {https://doi.org/10.1007/s00220-025-05505-0},
        note = {\href{https://arxiv.org/abs/2405.04649}{arXiv:\allowbreak 2405.04649}},
    eprint = {2405.04649},
    archivePrefix = {arXiv},
    pagetotal = {76},
}

@article{BR08,
    AUTHOR = {Baker, Andrew and Richter, Birgit},
     TITLE = {Uniqueness of {$E_\infty$} structures for connective covers},
   JOURNAL = {Proc. Amer. Math. Soc.},
  FJOURNAL = {Proceedings of the American Mathematical Society},
    VOLUME = {136},
      YEAR = {2008},
    NUMBER = {2},
     PAGES = {707--714},
      ISSN = {0002-9939,1088-6826},
   MRCLASS = {55P43 (19L41 55N15)},
  MRNUMBER = {2358512},
MRREVIEWER = {Jes\'{u}s\ Gonz\'{a}lez},
       DOI = {10.1090/S0002-9939-07-08984-8},
       URL = {https://doi.org/10.1090/S0002-9939-07-08984-8},
       note = {\href{https://arxiv.org/abs/math/0506422}{arXiv:\allowbreak math/0506422}},
    eprint = {math/0506422},
    archivePrefix = {arXiv},
}

@article{BLM23,
	  AUTHOR = {Beardsley, Jonathan and Luecke, Kiran and Morava, Jack},
     TITLE = {Brauer--{W}all groups and truncated {P}icard spectra of
              {K}-theory},
   JOURNAL = {Algebr. Geom. Topol.},
  FJOURNAL = {Algebraic \& Geometric Topology},
    VOLUME = {26},
      YEAR = {2026},
    NUMBER = {5},
     PAGES = {1749--1780},
      ISSN = {1472-2747,1472-2739},
   MRCLASS = {19L50 (55N15 55P42)},
  MRNUMBER = {5079112},
       DOI = {10.2140/agt.2026.26.1749},
       URL = {https://doi.org/10.2140/agt.2026.26.1749},
	note = {\href{https://arxiv.org/abs/2306.10112}{arXiv:\allowbreak 2306.10112}},
    eprint = {2306.10112},
    archivePrefix = {arXiv},
}

@misc{TY23a,
	author = {Yuji Tachikawa and Mayuko Yamashita},
	title = {Anderson duality of topological modular forms and its differential-geometric manifestations},
	year = {2023},
	note = {With an appendix by Sanath Devalapurkar. \href{https://arxiv.org/abs/2305.06196}{arXiv:\allowbreak 2305.06196}},
    eprint = {2305.06196},
    archivePrefix = {arXiv},
}

@article{TY23b,
    AUTHOR = {Tachikawa, Yuji and Yamashita, Mayuko},
     TITLE = {Topological modular forms and the absence of all heterotic
              global anomalies},
   JOURNAL = {Commun. Math. Phys.},
  FJOURNAL = {Communications in Mathematical Physics},
    VOLUME = {402},
      YEAR = {2023},
    NUMBER = {2},
     PAGES = {1585--1620},
      ISSN = {0010-3616,1432-0916},
   MRCLASS = {81T30 (81T45)},
  MRNUMBER = {4627327},
MRREVIEWER = {Farhang\ Loran},
       DOI = {10.1007/s00220-023-04761-2},
       URL = {https://doi.org/10.1007/s00220-023-04761-2},
       note = {\href{https://arxiv.org/abs/2108.13542}{arXiv:\allowbreak 2108.13542}},
    eprint = {2108.13542},
    archivePrefix = {arXiv},
}

@misc{TY25b,
	author = {Yuji Tachikawa and Kazuya Yonekura},
	title = {On invariants of two-dimensional minimally supersymmetric field theories},
	year = {2025},
	note = {\href{https://arxiv.org/abs/2508.04916}{arXiv:\allowbreak 2508.04916}},
    eprint = {2508.04916},
    archivePrefix = {arXiv},
}

@phdthesis{Sul66,
    AUTHOR = {Sullivan, Dennis Parnell},
     TITLE = {Triangulating homotopy equivalences},
      school = {Princeton University},
      YEAR = {1966},
}

@article{GJOS17,
    AUTHOR = {Gurski, Nick and Johnson, Niles and Osorno, Angélica M.
              and Stephan, Marc},
     TITLE = {Stable {P}ostnikov data of {P}icard 2-categories},
   JOURNAL = {Algebr. Geom. Topol.},
  FJOURNAL = {Algebraic \& Geometric Topology},
    VOLUME = {17},
      YEAR = {2017},
    NUMBER = {5},
     PAGES = {2763--2806},
      ISSN = {1472-2747,1472-2739},
   MRCLASS = {55S45 (18C20 18D05 19D23 55P42)},
  MRNUMBER = {3704242},
MRREVIEWER = {J\'{e}r\^{o}me\ Scherer},
       DOI = {10.2140/agt.2017.17.2763},
       URL = {https://doi.org/10.2140/agt.2017.17.2763},
        note = {\href{https://arxiv.org/abs/1606.07032}{arXiv:\allowbreak 1606.07032}},
    eprint = {1606.07032},
    archivePrefix = {arXiv},
}

@article{JO12,
    AUTHOR = {Johnson, Niles and Osorno, Angélica M.},
     TITLE = {Modeling stable one-types},
   JOURNAL = {Theory Appl. Categ.},
  FJOURNAL = {Theory and Applications of Categories},
    VOLUME = {26},
      YEAR = {2012},
     PAGES = {520--537},
      ISSN = {1201-561X},
   MRCLASS = {18B40 (55P42 55S45)},
  MRNUMBER = {2981952},
MRREVIEWER = {Geoffrey\ M. L. Powell},
    note = {\href{https://arxiv.org/abs/1201.2686}{arXiv:\allowbreak 1201.2686}},
    eprint = {1201.2686},
    archivePrefix = {arXiv},
    eid = {20},
}

@article{LKW17,
    AUTHOR = {Lan, Tian and Kong, Liang and Wen, Xiao-Gang},
     TITLE = {Modular extensions of unitary braided fusion categories and
              {$2+1{\rm D}$} topological/{SPT} orders with symmetries},
   JOURNAL = {Commun. Math. Phys.},
  FJOURNAL = {Communications in Mathematical Physics},
    VOLUME = {351},
      YEAR = {2017},
    NUMBER = {2},
     PAGES = {709--739},
      ISSN = {0010-3616,1432-0916},
   MRCLASS = {81R10 (81T08 81T40)},
  MRNUMBER = {3613518},
MRREVIEWER = {Andreas\ Arvanitoyeorgos},
       DOI = {10.1007/s00220-016-2748-y},
       URL = {https://doi.org/10.1007/s00220-016-2748-y},
       note = {\href{https://arxiv.org/abs/1602.05936}{arXiv:\allowbreak 1602.05936}},
    eprint = {1602.05936},
    archivePrefix = {arXiv},
}

@article{BGHNPRW17,
    AUTHOR = {Bruillard, Paul and Galindo, César and Hagge, Tobias and
              Ng, Siu-Hung and Plavnik, Julia Yael and Rowell, Eric C. and
              Wang, Zhenghan},
     TITLE = {Fermionic modular categories and the 16-fold way},
   JOURNAL = {J. Math. Phys.},
  FJOURNAL = {Journal of Mathematical Physics},
    VOLUME = {58},
      YEAR = {2017},
    NUMBER = {4},
     PAGES = {041704},
      ISSN = {0022-2488,1089-7658},
   MRCLASS = {81V35},
  MRNUMBER = {3641612},
MRREVIEWER = {Tim\ Raeymaekers},
       DOI = {10.1063/1.4982048},
       URL = {https://doi.org/10.1063/1.4982048},
       note = {\href{https://arxiv.org/abs/1603.09294}{arXiv:\allowbreak 1603.09294}},
    eprint = {1603.09294},
    archivePrefix = {arXiv},
    pagetotal = {31},
}

@article{JFS17,
    AUTHOR = {Johnson-Freyd, Theo and Scheimbauer, Claudia},
     TITLE = {({O}p)lax natural transformations, twisted quantum field
              theories, and ``even higher'' {M}orita categories},
   JOURNAL = {Adv. Math.},
  FJOURNAL = {Advances in Mathematics},
    VOLUME = {307},
      YEAR = {2017},
     PAGES = {147--223},
      ISSN = {0001-8708,1090-2082},
   MRCLASS = {18D10 (81R10)},
  MRNUMBER = {3590516},
MRREVIEWER = {Josep\ Elgueta},
       DOI = {10.1016/j.aim.2016.11.014},
       URL = {https://doi.org/10.1016/j.aim.2016.11.014},
       note = {\href{https://arxiv.org/abs/1502.06526}{arXiv:\allowbreak 1502.06526}},
    eprint = {1502.06526},
    archivePrefix = {arXiv},
}

@incollection{Lur09,
    AUTHOR = {Lurie, Jacob},
     TITLE = {On the classification of topological field theories},
 BOOKTITLE = {Current developments in mathematics, 2008},
     PAGES = {129--280},
 PUBLISHER = {Int. Press, Somerville, MA},
      YEAR = {2009},
   MRCLASS = {57R56 (18D10 18G30 57R15 57R75)},
  MRNUMBER = {2555928},
MRREVIEWER = {Julia Bergner},
        note = {\href{https://arxiv.org/abs/0905.0465}{arXiv:\allowbreak 0905.0465}},
    eprint = {0905.0465},
    archivePrefix = {arXiv},
}

@article{BD95,
    author = {Baez, John C. and Dolan, James},
    title = {Higher-dimensional algebra and topological quantum field theory},
    journal = {J. Math. Phys.},
    fjournal = {Journal of Mathematical Physics},
    volume = {36},
    number = {11},
    year = {1995},
    pages = {6073--6105},
    issn = {0022-2488},
    doi = {10.1063/1.531236},
    note = {\href{https://arxiv.org/abs/q-alg/9503002}{arXiv:\allowbreak q-alg/9503002}},
    eprint = {q-alg/9503002},
    archivePrefix = {arXiv},
}

@article{Gun16,
        title = "Spin {H}urwitz numbers and topological quantum field theory",
        journal = "Geom. Topol.",
        fjournal = "Geometry \& Topology",
        volume = "20",
        number = "4",
        pages = "1859--1907",
        year = "2016",
        author = "Sam Gunningham",
        note = {\href{https://arxiv.org/abs/1201.1273}{arXiv:\allowbreak 1201.1273}},
    eprint = {1201.1273},
    archivePrefix = {arXiv},
}

@misc{Ant24,
	author = {Benjamin Antieau},
	title = {Spectral sequences, décalage, and the {B}eilinson {$t$}-structure},
	year = {2024},
	note = {\href{https://arxiv.org/abs/2411.09115}{arXiv:\allowbreak 2411.09115}},
    eprint = {2411.09115},
    archivePrefix = {arXiv},
}

@article{BCC93,
        title={Cohomology with coefficients in symmetric cat-groups. {A}n extension of {E}ilenberg-{MacLane}'s classification theorem},
        volume={114},
        DOI={10.1017/S0305004100071498},
        number={1},
        journal={Math. Proc. Cambridge Philos. Soc.},
        fjournal={Mathematical Proceedings of the Cambridge Philosophical Society},
        publisher={Cambridge University Press},
        author={Bullejos, M. and Carrasco, P. and Cegarra, A. M.},
        year={1993},
        pages={163--189}
}

@incollection{Dri06,
        author="Drinfeld, Vladimir",
        editor="Etingof, Pavel
        and Retakh, Vladimir
        and Singer, I. M.",
        title="Infinite-Dimensional Vector Bundles in Algebraic Geometry",
        bookTitle="The Unity of Mathematics: In Honor of the Ninetieth Birthday of I.M. Gelfand",
        year="2006",
        publisher="Birkhäuser Boston",
        address="Boston, MA",
        pages="263--304",
        isbn="978-0-8176-4467-3",
        doi="10.1007/0-8176-4467-9_7",
        note = {\href{https://arxiv.org/abs/math/0309155}{arXiv:\allowbreak math/0309155}},
    eprint = {math/0309155},
    archivePrefix = {arXiv},
}

@article{GK14,
        author="Ganter, Nora and Kapranov, Mikhail",
        title="SYMMETRIC AND EXTERIOR POWERS OF CATEGORIES",
        journal="Transform. Groups",
        fjournal="Transformation Groups",
        year="2014",
        volume="19",
        number="1",
        pages="57--103",
        issn="1531-586X",
        doi="10.1007/s00031-014-9255-z",
        note = {\href{https://arxiv.org/abs/1110.4753}{arXiv:\allowbreak 1110.4753}},
    eprint = {1110.4753},
    archivePrefix = {arXiv},
}

@article{Pat12,
        author="Patel, Deepam",
        title="De {Rham} {$\mathcal E$}-factors",
        journal="Invent. Math.",
        fjournal="Inventiones Mathematicae",
        year="2012",
        volume="190",
        number="2",
        pages="299--355",
        issn="1432-1297",
        doi="10.1007/s00222-012-0381-8",
        note = "\url{http://www.math.purdue.edu/~patel471/deRhamep.pdf}"
}

@article{HS05,
        author = "Hopkins, M.J. and Singer, I.M.",
        doi = "10.4310/jdg/1143642908",
        fjournal = "Journal of Differential Geometry",
        journal = "J. Differential Geom.",
        number = "3",
        pages = "329--452",
        publisher = "Lehigh University",
        title = "Quadratic functions in geometry, topology, and {M}-theory",
        url = "https://doi.org/10.4310/jdg/1143642908",
        volume = "70",
        year = "2005",
        note = {\href{https://arxiv.org/abs/math/0211216}{arXiv:\allowbreak math/0211216}},
    eprint = {math/0211216},
    archivePrefix = {arXiv},
}

@article{GJO19,
    AUTHOR = {Gurski, Nick and Johnson, Niles and Osorno, Angélica M.},
     TITLE = {The 2-dimensional stable homotopy hypothesis},
   JOURNAL = {J. Pure Appl. Algebra},
  FJOURNAL = {Journal of Pure and Applied Algebra},
    VOLUME = {223},
      YEAR = {2019},
    NUMBER = {10},
     PAGES = {4348--4383},
      ISSN = {0022-4049},
   MRCLASS = {55P42 (18D05 18D50 19D23 55P15 55P48)},
  MRNUMBER = {3958095},
MRREVIEWER = {Constanze Roitzheim},
       DOI = {10.1016/j.jpaa.2019.01.012},
        note = {\href{https://arxiv.org/abs/1712.07218}{arXiv:\allowbreak 1712.07218}},
    eprint = {1712.07218},
    archivePrefix = {arXiv},
}

@article{MOPSV22,
    AUTHOR = {Moser, Lyne and Ozornova, Viktoriya and Paoli, Simona and
              Sarazola, Maru and Verdugo, Paula},
     TITLE = {Stable homotopy hypothesis in the {T}amsamani model},
   JOURNAL = {Topology Appl.},
  FJOURNAL = {Topology and its Applications},
    VOLUME = {316},
      YEAR = {2022},
     PAGES = {108106},
      ISSN = {0166-8641,1879-3207},
   MRCLASS = {18N10 (18B40 18M05 18N20 19D23 55P10 55P42 55U35)},
  MRNUMBER = {4438951},
MRREVIEWER = {Tilman\ Bauer},
       DOI = {10.1016/j.topol.2022.108106},
       URL = {https://doi.org/10.1016/j.topol.2022.108106},
       note = {\href{https://arxiv.org/abs/2001.05577}{arXiv:\allowbreak 2001.05577}},
    eprint = {2001.05577},
    archivePrefix = {arXiv},
    pagetotal = {40},
}

@article{JFY22,
    AUTHOR = {Johnson-Freyd, Theo and Yu, Matthew},
     TITLE = {Topological orders in {$(4+1)$}-dimensions},
   JOURNAL = {SciPost Phys.},
  FJOURNAL = {SciPost Physics},
    VOLUME = {13},
      YEAR = {2022},
    NUMBER = {3},
     PAGES = {068},
      ISSN = {2542-4653},
   MRCLASS = {81T45 (18M20 18N10)},
  MRNUMBER = {4492330},
       DOI = {10.21468/scipostphys.13.3.068},
       URL = {https://doi.org/10.21468/scipostphys.13.3.068},
       note = {\href{https://arxiv.org/abs/2104.04534}{arXiv:\allowbreak 2104.04534}},
    eprint = {2104.04534},
    archivePrefix = {arXiv},
    pagetotal = {20},
}

@article{Sha23,
    AUTHOR = {Sharma, Amit},
     TITLE = {Picard groupoids and {$\Gamma$}-categories},
   JOURNAL = {Cah. Topol. G\'{e}om. Diff\'{e}r. Cat\'{e}g.},
  FJOURNAL = {Cahiers de Topologie et G\'{e}om\'{e}trie Diff\'{e}rentielle
              Cat\'{e}goriques},
    VOLUME = {64},
      YEAR = {2023},
    NUMBER = {3},
     PAGES = {287--324},
      ISSN = {1245-530X,2681-2363},
   MRCLASS = {18M05 (18N40)},
  MRNUMBER = {4625820},
    note = {\href{https://arxiv.org/abs/2002.05811}{arXiv:\allowbreak 2002.05811}},
    eprint = {2002.05811},
    archivePrefix = {arXiv},
}

@article{BP25,
    AUTHOR = {Burklund, Robert and Pstrągowski, Piotr},
     TITLE = {Quivers and the {A}dams spectral sequence},
   JOURNAL = {Adv. Math.},
  FJOURNAL = {Advances in Mathematics},
    VOLUME = {471},
      YEAR = {2025},
     PAGES = {110270},
      ISSN = {0001-8708,1090-2082},
   MRCLASS = {55T15 (16G20 18N99 19L20)},
  MRNUMBER = {4890707},
MRREVIEWER = {Yangyang\ Ruan},
       DOI = {10.1016/j.aim.2025.110270},
       URL = {https://doi.org/10.1016/j.aim.2025.110270},
    note = {\href{https://arxiv.org/abs/2305.08231}{arXiv:\allowbreak 2305.08231}},
    eprint = {2305.08231},
    archivePrefix = {arXiv},
    pagetotal = {76},
}

@article{Hys23,
	AUTHOR = {Hyslop, Logan},
     TITLE = {Whitehead filtrations for computations in topological
              {H}ochschild homology},
   JOURNAL = {Algebr. Geom. Topol.},
  FJOURNAL = {Algebraic \& Geometric Topology},
    VOLUME = {26},
      YEAR = {2026},
    NUMBER = {7},
     PAGES = {2443--2474},
      ISSN = {1472-2747,1472-2739},
   MRCLASS = {55R20 (19D55 55T25)},
  MRNUMBER = {5125657},
       DOI = {10.2140/agt.2026.26.2443},
       URL = {https://doi.org/10.2140/agt.2026.26.2443},
	note = {\href{https://arxiv.org/abs/2311.06717}{arXiv:\allowbreak 2311.06717}},
}

@misc{Wu26,
	author = {Yuchen Wu},
	title = {Burklund-{L}in-{W}ang-{X}u Methods in the Cofiber-of-Tau Formalism and Applications to Equivariant Slice Differentials},
	year = {2026},
	note = {\href{https://arxiv.org/abs/2606.02779}{arXiv:\allowbreak 2606.02779}},
    eprint = {2606.02779},
    archivePrefix = {arXiv},
}

@phdthesis{Hed20,
    AUTHOR = {Alice Hedenlund},
     TITLE = {Multiplicative {T}ate spectral sequences},
      school = {University of Oslo},
      YEAR = {2020},
      note = {\url{https://www.mn.uio.no/math/personer/vit/rognes/theses/hedenlund-thesis.pdf}}
}

@misc{KPY26,
	author = {Ryohei Kobayashi and Abhinav Prem and Matthew Yu},
	title = {Sixteen-Fold Way for Fermionic Topological Orders},
	year = {2026},
	note = {\href{https://arxiv.org/abs/2606.28682}{arXiv:\allowbreak 2606.28682}},
    eprint = {2606.28682},
    archivePrefix = {arXiv},
}

@article{May88,
    AUTHOR = {Ellen {Maycock Parker}},
     TITLE = {The {B}rauer group of graded continuous trace
              {$C^*$}-algebras},
   JOURNAL = {Trans. Amer. Math. Soc.},
  FJOURNAL = {Transactions of the American Mathematical Society},
    VOLUME = {308},
      YEAR = {1988},
    NUMBER = {1},
     PAGES = {115--132},
      ISSN = {0002-9947},
   MRCLASS = {46L05 (16A16 22D25 55R10)},
  MRNUMBER = {946434},
MRREVIEWER = {J. A. Mingo},
       DOI = {10.2307/2000953},
       URL = {https://doi.org/10.2307/2000953},
}

@article{FHM03,
    AUTHOR = {Fausk, H. and Hu, P. and May, J. P.},
     TITLE = {Isomorphisms between left and right adjoints},
   JOURNAL = {Theory Appl. Categ.},
  FJOURNAL = {Theory and Applications of Categories},
    VOLUME = {11},
      YEAR = {2003},
     PAGES = {107--131},
      ISSN = {1201-561X},
   MRCLASS = {18D10 (14A10 55U35)},
  MRNUMBER = {1988072},
MRREVIEWER = {Gabriele\ Vezzosi},
    note = {\href{https://arxiv.org/abs/math/0206079}{arXiv:\allowbreak math/0206079}},
    eprint = {math/0206079},
    archivePrefix = {arXiv},
    eid = {4},
}

@article{Wal62,
    AUTHOR = {Wall, C. T. C.},
     TITLE = {A characterization of simple modules over the {S}teenrod
              algebra {${\rm mod}\ 2$}},
   JOURNAL = {Topology},
  FJOURNAL = {Topology. An International Journal of Mathematics},
    VOLUME = {1},
      YEAR = {1962},
     PAGES = {249--254},
      ISSN = {0040-9383},
   MRCLASS = {55.34},
  MRNUMBER = {148061},
MRREVIEWER = {R.\ Brown},
       DOI = {10.1016/0040-9383(62)90108-8},
       URL = {https://doi.org/10.1016/0040-9383(62)90108-8},
}

@article{Kan58a,
    AUTHOR = {Kan, Daniel M.},
     TITLE = {A combinatorial definition of homotopy groups},
   JOURNAL = {Ann. of Math. (2)},
  FJOURNAL = {Annals of Mathematics. Second Series},
    VOLUME = {67},
      YEAR = {1958},
     PAGES = {282--312},
      ISSN = {0003-486X},
   MRCLASS = {55.00},
  MRNUMBER = {111032},
       DOI = {10.2307/1970006},
       URL = {https://doi.org/10.2307/1970006},
}

@article{Kan58b,
    AUTHOR = {Kan, Daniel M.},
     TITLE = {On homotopy theory and c.s.s. groups},
   JOURNAL = {Ann. of Math. (2)},
  FJOURNAL = {Annals of Mathematics. Second Series},
    VOLUME = {68},
      YEAR = {1958},
     PAGES = {38--53},
      ISSN = {0003-486X},
   MRCLASS = {55.00},
  MRNUMBER = {111033},
MRREVIEWER = {J. A. Zilber},
       DOI = {10.2307/1970042},
       URL = {https://doi.org/10.2307/1970042},
}

@article{PP26,
    AUTHOR = {Patchkoria, Irakli and Pstrągowski, Piotr},
     TITLE = {Adams {S}pectral {S}equences and {F}ranke's {A}lgebraicity
              {C}onjecture},
   JOURNAL = {Mem. Amer. Math. Soc.},
  FJOURNAL = {Memoirs of the American Mathematical Society},
    VOLUME = {321},
      YEAR = {2026},
    NUMBER = {1636},
     PAGES = {v+140},
      ISSN = {0065-9266,1947-6221},
      ISBN = {978-1-4704-8178-0; 978-1-4704-8708-9},
   MRCLASS = {Prelim},
  MRNUMBER = {5101334},
       DOI = {10.1090/memo/1636},
       URL = {https://doi.org/10.1090/memo/1636},
       note = {\href{https://arxiv.org/abs/2110.03669}{arXiv:\allowbreak 2110.03669}},
    eprint = {2110.03669},
    archivePrefix = {arXiv},
}

@article{Pst23,
    AUTHOR = {Pstrągowski, Piotr},
     TITLE = {Synthetic spectra and the cellular motivic category},
   JOURNAL = {Invent. Math.},
  FJOURNAL = {Inventiones Mathematicae},
    VOLUME = {232},
      YEAR = {2023},
    NUMBER = {2},
     PAGES = {553--681},
      ISSN = {0020-9910,1432-1297},
   MRCLASS = {55P43 (14F42 18N60 55T15)},
  MRNUMBER = {4574661},
MRREVIEWER = {Geoffrey\ M. L. Powell},
       DOI = {10.1007/s00222-022-01173-2},
       URL = {https://doi.org/10.1007/s00222-022-01173-2},
       note = {\href{https://arxiv.org/abs/1803.01804}{arXiv:\allowbreak 1803.01804}},
    eprint = {1803.01804},
    archivePrefix = {arXiv},
}

@article{LW24,
    AUTHOR = {Livernet, Muriel and Whitehouse, Sarah},
     TITLE = {Homotopy theory of spectral sequences},
   JOURNAL = {Homology Homotopy Appl.},
  FJOURNAL = {Homology, Homotopy and Applications},
    VOLUME = {26},
      YEAR = {2024},
    NUMBER = {1},
     PAGES = {69--86},
      ISSN = {1532-0073,1532-0081},
   MRCLASS = {18G40 (18N40)},
  MRNUMBER = {4709809},
       DOI = {10.4310/hha.2024.v26.n1.a5},
       URL = {https://doi.org/10.4310/hha.2024.v26.n1.a5},
       note = {\href{https://arxiv.org/abs/2207.02156}{arXiv:\allowbreak 2207.02156}},
    eprint = {2207.02156},
    archivePrefix = {arXiv},
}

@article{CESLW20,
    AUTHOR = {Cirici, Joana and Egas Santander, Daniela and Livernet, Muriel
              and Whitehouse, Sarah},
     TITLE = {Model category structures and spectral sequences},
   JOURNAL = {Proc. Roy. Soc. Edinburgh Sect. A},
  FJOURNAL = {Proceedings of the Royal Society of Edinburgh. Section A.
              Mathematics},
    VOLUME = {150},
      YEAR = {2020},
    NUMBER = {6},
     PAGES = {2815--2848},
      ISSN = {0308-2105,1473-7124},
   MRCLASS = {18N40 (55U35)},
  MRNUMBER = {4190091},
MRREVIEWER = {Scott\ Balchin},
       DOI = {10.1017/prm.2019.45},
       URL = {https://doi.org/10.1017/prm.2019.45},
       note ={\href{https://arxiv.org/abs/1805.00374}{arXiv:\allowbreak 1805.00374}},
    eprint = {1805.00374},
    archivePrefix = {arXiv},
}

@article{CLW25,
    AUTHOR = {Cirici, Joana and Livernet, Muriel and Whitehouse, Sarah},
     TITLE = {Model category structures on truncated multicomplexes for
              complex geometry},
   JOURNAL = {Bull. Lond. Math. Soc.},
  FJOURNAL = {Bulletin of the London Mathematical Society},
    VOLUME = {57},
      YEAR = {2025},
    NUMBER = {12},
     PAGES = {4163--4177},
      ISSN = {0024-6093,1469-2120},
   MRCLASS = {18N40 (18G40 32C99 32Q55)},
  MRNUMBER = {5007401},
       DOI = {10.1112/blms.70239},
       URL = {https://doi.org/10.1112/blms.70239},
       note = {\href{https://arxiv.org/abs/2501.13509}{arXiv:\allowbreak 2501.13509}},
    eprint = {2501.13509},
    archivePrefix = {arXiv},
}

@misc{LW24a,
	author = {Muriel Livernet and Sarah Whitehouse},
	title = {Spectral sequences via linear presheaves},
	year = {2024},
	note = {\href{https://arxiv.org/abs/2406.02777}{arXiv:\allowbreak 2406.02777}},
    eprint = {2406.02777},
    archivePrefix = {arXiv},
}

@article{Bro25,
    AUTHOR = {Brotherston, James A.},
     TITLE = {A distributive lattice of model structures relating to
              spectral sequences},
   JOURNAL = {Homology Homotopy Appl.},
  FJOURNAL = {Homology, Homotopy and Applications},
    VOLUME = {27},
      YEAR = {2025},
    NUMBER = {2},
     PAGES = {91--122},
      ISSN = {1532-0073,1532-0081},
   MRCLASS = {18N40 (06D99 18G40 18N55)},
  MRNUMBER = {4971731},
MRREVIEWER = {Philippe\ Gaucher},
       DOI = {10.4310/hha.2025.v27.n2.a5},
       URL = {https://doi.org/10.4310/hha.2025.v27.n2.a5},
       note = {\href{https://arxiv.org/abs/2402.09893}{arXiv:\allowbreak 2402.09893}},
    eprint = {2402.09893},
    archivePrefix = {arXiv},
}

@misc{Bro24,
	author = {James A. Brotherston},
	title = {Monoidal model structures on filtered chain complexes relating to spectral sequences},
	year = {2024},
	note = {\href{https://arxiv.org/abs/2402.09207}{arXiv:\allowbreak 2402.09207}},
    eprint = {2402.09207},
    archivePrefix = {arXiv},
}

@article{MR19,
    AUTHOR = {Muro, Fernando and Roitzheim, Constanze},
     TITLE = {Homotopy theory of bicomplexes},
   JOURNAL = {J. Pure Appl. Algebra},
  FJOURNAL = {Journal of Pure and Applied Algebra},
    VOLUME = {223},
      YEAR = {2019},
    NUMBER = {5},
     PAGES = {1913--1939},
      ISSN = {0022-4049,1873-1376},
   MRCLASS = {55U35 (18G35)},
  MRNUMBER = {3906533},
MRREVIEWER = {Jan\ Spali\'{n}ski},
       DOI = {10.1016/j.jpaa.2018.08.007},
       URL = {https://doi.org/10.1016/j.jpaa.2018.08.007},
       note = {\href{https://arxiv.org/abs/1802.07610}{arXiv:\allowbreak 1802.07610}},
    eprint = {1802.07610},
    archivePrefix = {arXiv},
}

@article{FGLW22,
    AUTHOR = {Fu, Xin and Guan, Ai and Livernet, Muriel and Whitehouse,
              Sarah},
     TITLE = {Model category structures on multicomplexes},
   JOURNAL = {Topology Appl.},
  FJOURNAL = {Topology and its Applications},
    VOLUME = {316},
      YEAR = {2022},
     PAGES = {108104},
      ISSN = {0166-8641,1879-3207},
   MRCLASS = {18N40 (18G40)},
  MRNUMBER = {4438949},
MRREVIEWER = {Emma\ Lepri},
       DOI = {10.1016/j.topol.2022.108104},
       URL = {https://doi.org/10.1016/j.topol.2022.108104},
       note = {\href{https://arxiv.org/abs/2001.10873}{arXiv:\allowbreak 2001.10873}},
    eprint = {2001.10873},
    archivePrefix = {arXiv},
    pagetotal = {26},
}

@phdthesis{Ait24,
	author = {Colin Aitken},
	title = {Model Structures on Infinity-Categories of Filtrations},
	year = {2024},
    school = {The University of Chicago},
	note = {\href{https://arxiv.org/abs/2402.17921}{arXiv:\allowbreak 2402.17921}},
    eprint = {2402.17921},
    archivePrefix = {arXiv},
}

@article{FN23,
    AUTHOR = {Freed, Daniel S. and Neitzke, Andrew},
     TITLE = {The dilogarithm and abelian {C}hern-{S}imons},
   JOURNAL = {J. Differential Geom.},
  FJOURNAL = {Journal of Differential Geometry},
    VOLUME = {123},
      YEAR = {2023},
    NUMBER = {2},
     PAGES = {241--266},
      ISSN = {0022-040X,1945-743X},
   MRCLASS = {53C07},
  MRNUMBER = {4571802},
MRREVIEWER = {Paul\ F.\ Bracken},
       DOI = {10.4310/jdg/1680883577},
       URL = {https://doi.org/10.4310/jdg/1680883577},
       note = {\href{https://arxiv.org/abs/2006.12565}{arXiv:\allowbreak 2006.12565}},
    eprint = {2006.12565},
    archivePrefix = {arXiv},
}

@misc{LPT25,
	author = {Guchuan Li and Sarah Petersen and Elizabeth Tatum},
	title = {A spectrum-level splitting of the {$ku_\mathbb{R}$}-cooperations algebra},
	year = {2025},
	note = {\href{https://arxiv.org/abs/2503.17149}{arXiv:\allowbreak 2503.17149}},
    eprint = {2503.17149},
    archivePrefix = {arXiv},
}

@misc{Mor25,
	author = {Jackson Morris},
	title = {Rings of cooperations for hermitian {K}-theory over finite fields},
	year = {2025},
	note = {\href{https://arxiv.org/abs/2509.02786}{arXiv:\allowbreak 2509.02786}},
    eprint = {2509.02786},
    archivePrefix = {arXiv},
}

@misc{MPT25,
	author = {Jackson Morris and Sarah Petersen and Elizabeth Tatum},
	title = {Splittings of truncated motivic {Brown--Peterson} cooperations algebras},
	year = {2025},
	note = {\href{https://arxiv.org/abs/2509.19542}{arXiv:\allowbreak 2509.19542}},
    eprint = {2509.19542},
    archivePrefix = {arXiv},
}

@article{Mor25a,
	AUTHOR = {Morris, Jackson},
     TITLE = {On the ring of cooperations for real hermitian {K}-theory},
   JOURNAL = {Ann. K-Theory},
  FJOURNAL = {Annals of K-Theory},
    VOLUME = {11},
      YEAR = {2026},
    NUMBER = {4},
     PAGES = {633--708},
      ISSN = {2379-1683,2379-1691},
   MRCLASS = {14F42 (19G38 55Q10 55T15)},
  MRNUMBER = {5148298},
       DOI = {10.2140/akt.2026.11.633},
       URL = {https://doi.org/10.2140/akt.2026.11.633},
	note = {\href{https://arxiv.org/abs/2506.16672}{arXiv:\allowbreak 2506.16672}},
    eprint = {2506.16672},
    archivePrefix = {arXiv},
}

@article{CDDH26,
    AUTHOR = {Chakrabhavi, Vivek and Debray, Arun and Dierigl, Markus and
              Heckman, Jonathan J.},
     TITLE = {Reflection branes, bordisms, and {$U$} dualities},
   JOURNAL = {Phys. Rev. D},
  FJOURNAL = {Physical Review D},
    VOLUME = {113},
      YEAR = {2026},
    NUMBER = {6},
     PAGES = {066015},
      ISSN = {2470-0010,2470-0029},
   MRCLASS = {83E50},
  MRNUMBER = {5066220},
       DOI = {10.1103/99d8-hyhs},
       URL = {https://doi.org/10.1103/99d8-hyhs},
       note = {\href{https://arxiv.org/abs/2509.03573}{arXiv:\allowbreak 2509.03573}},
    eprint = {2509.03573},
    archivePrefix = {arXiv},
    pagetotal = {30},
}

@phdthesis{Guo18,
    author = "Meng Guo",
    title = "Some Calculations of Cobordism Groups and Their Applications in Physics",
    school = "Harvard University",
    year = "2018",
    note = "\url{https://dash.harvard.edu/bitstream/handle/1/40050106/GUO-DISSERTATION-2018.pdf?sequence=4}"
}

@article{Wil73,
    AUTHOR = {Wilson, W. Stephen},
     TITLE = {A new relation on the {S}tiefel-{W}hitney classes of spin
              manifolds},
   JOURNAL = {Illinois J. Math.},
  FJOURNAL = {Illinois Journal of Mathematics},
    VOLUME = {17},
      YEAR = {1973},
     PAGES = {115--127},
      ISSN = {0019-2082},
   MRCLASS = {55F40},
  MRNUMBER = {310883},
MRREVIEWER = {R.\ De Sapio},
       URL = {http://projecteuclid.org/euclid.ijm/1256052043},
}

@phdthesis{Sie13,
    author = "Christian Siegemeyer",
    title = "On the {G}romov-{L}awson-{R}osenberg Conjecture for Finite Abelian 2-Groups of Rank 2",
    school = "Westfälische Wilhelms-Universität Münster",
    year = "2013",
    note = "\url{https://repositorium.uni-muenster.de/document/miami/4882c034-0925-432c-b743-edf299b5192b/diss_siegemeyer.pdf}"
}

@article{She25,
    AUTHOR = {Sheckler, Aiden},
     TITLE = {Anomaly matching in 6d {$\mathcal N=(2,0)$} {SCFT}s from {M}5
              cobordism},
   JOURNAL = {J. High Energy Phys.},
  FJOURNAL = {Journal of High Energy Physics},
      YEAR = {2025},
    VOLUME = {2025},
     PAGES = {116},
      ISSN = {1126-6708,1029-8479},
   MRCLASS = {81T40 (81T30 81T50)},
  MRNUMBER = {5004602},
       DOI = {10.1007/jhep12(2025)116},
       URL = {https://doi.org/10.1007/jhep12(2025)116},
       note = {\href{https://arxiv.org/abs/2501.04785}{arXiv:\allowbreak 2501.04785}},
    number = {12},
    pagetotal = {28},
    eprint = {2501.04785},
    archivePrefix = {arXiv},
}

@misc{GLWZ26,
	author = {Zhenbang Gu and Ran Luo and Yi-Nan Wang and Yi Zhang},
	title = {Understanding Non-Split 2-Group Symmetry: {(3+1)D} {SymTFT}, Anomaly and Bordism},
	year = {2026},
	note = {\href{https://arxiv.org/abs/2608.08766}{arXiv:\allowbreak 2608.08766}},
    eprint = {2608.08766},
    archivePrefix = {arXiv},
}

@article{GIHR20,
    AUTHOR = {Guillou, Bertrand J. and Hill, Michael A. and Isaksen, Daniel
              C. and Ravenel, Douglas Conner},
     TITLE = {The cohomology of {$C_2$}-equivariant {$\mathcal A(1)$} and the
              homotopy of {${\rm ko}_{C_2}$}},
   JOURNAL = {Tunis. J. Math.},
  FJOURNAL = {Tunisian Journal of Mathematics},
    VOLUME = {2},
      YEAR = {2020},
    NUMBER = {3},
     PAGES = {567--632},
      ISSN = {2576-7658,2576-7666},
   MRCLASS = {14F42 (55Q91 55T15)},
  MRNUMBER = {4041284},
MRREVIEWER = {Bj\o rn\ Ian\ Dundas},
       DOI = {10.2140/tunis.2020.2.567},
       URL = {https://doi.org/10.2140/tunis.2020.2.567},
       note = {\href{https://arxiv.org/abs/1708.09568}{arXiv:\allowbreak 1708.09568}},
    eprint = {1708.09568},
    archivePrefix = {arXiv},
}

@article{CQ21,
    AUTHOR = {Culver, Dominic Leon and Quigley, J. D.},
     TITLE = {{$kq$}-resolutions {I}},
   JOURNAL = {Trans. Amer. Math. Soc.},
  FJOURNAL = {Transactions of the American Mathematical Society},
    VOLUME = {374},
      YEAR = {2021},
    NUMBER = {7},
     PAGES = {4655--4710},
      ISSN = {0002-9947,1088-6850},
   MRCLASS = {14F42 (55Q10 55Q45 55Q51 55T15)},
  MRNUMBER = {4273174},
MRREVIEWER = {Gabriel\ Angelini-Knoll},
       DOI = {10.1090/tran/8367},
       URL = {https://doi.org/10.1090/tran/8367},
       note = {\href{https://arxiv.org/abs/1905.11952}{arXiv:\allowbreak 1905.11952}},
    eprint = {1905.11952},
    archivePrefix = {arXiv},
}

@misc{DFNQY25,
	author = {Dhankhar, Prerna and Field, Rebecca and Nigam, Arjun and Quigley, J. D. and Yang, Albert Jinghui},
	title = {Very effective algebraic and hermitian {K}-theory of the cyclic group of order two
},
	year = {2025},
	note = {\href{https://arxiv.org/abs/2509.25006}{arXiv:\allowbreak 2509.25006}},
    eprint = {2509.25006},
    archivePrefix = {arXiv},
}

@article{CKQ21,
    AUTHOR = {Culver, Dominic Leon and Kong, Hana Jia and Quigley, J. D.},
     TITLE = {Algebraic slice spectral sequences},
   JOURNAL = {Doc. Math.},
  FJOURNAL = {Documenta Mathematica},
    VOLUME = {26},
      YEAR = {2021},
     PAGES = {1085--1119},
      ISSN = {1431-0635,1431-0643},
   MRCLASS = {14F42 (55P42 55P91 55T05 55T15)},
  MRNUMBER = {4324461},
MRREVIEWER = {Matthias\ Wendt},
    note = {\href{https://arxiv.org/abs/2007.08682}{arXiv:\allowbreak 2007.08682}},
    eprint = {2007.08682},
    archivePrefix = {arXiv},
}

@article{Hil11,
    AUTHOR = {Hill, Michael A.},
     TITLE = {Ext and the motivic {S}teenrod algebra over {$\mathbb{R}$}},
   JOURNAL = {J. Pure Appl. Algebra},
  FJOURNAL = {Journal of Pure and Applied Algebra},
    VOLUME = {215},
      YEAR = {2011},
    NUMBER = {5},
     PAGES = {715--727},
      ISSN = {0022-4049,1873-1376},
   MRCLASS = {55T15 (14F42 19D45 55P43 55S10)},
  MRNUMBER = {2747214},
MRREVIEWER = {Kyle\ M.\ Ormsby},
       DOI = {10.1016/j.jpaa.2010.06.017},
       URL = {https://doi.org/10.1016/j.jpaa.2010.06.017},
       note = {\href{https://arxiv.org/abs/0904.1998}{arXiv:\allowbreak 0904.1998}},
    eprint = {0904.1998},
    archivePrefix = {arXiv},
}

@article{Orm11,
    AUTHOR = {Ormsby, Kyle M.},
     TITLE = {Motivic invariants of {$p$}-adic fields},
   JOURNAL = {J. K-Theory},
  FJOURNAL = {Journal of K-Theory. K-Theory and its Applications in Algebra,
              Geometry, Analysis \& Topology},
    VOLUME = {7},
      YEAR = {2011},
    NUMBER = {3},
     PAGES = {597--618},
      ISSN = {1865-2433,1865-5394},
   MRCLASS = {14F42 (19D50)},
  MRNUMBER = {2811717},
MRREVIEWER = {Matthias\ Wendt},
       DOI = {10.1017/is011004017jkt153},
       URL = {https://doi.org/10.1017/is011004017jkt153},
       note = {\href{https://arxiv.org/abs/1002.5007}{arXiv:\allowbreak 1002.5007}},
    eprint = {1002.5007},
    archivePrefix = {arXiv},
}

@article{IS11,
    AUTHOR = {Isaksen, Daniel C. and Shkembi, Armira},
     TITLE = {Motivic connective {$K$}-theories and the cohomology of
              {A}(1)},
   JOURNAL = {J. K-Theory},
  FJOURNAL = {Journal of K-Theory. K-Theory and its Applications in Algebra,
              Geometry, Analysis \& Topology},
    VOLUME = {7},
      YEAR = {2011},
    NUMBER = {3},
     PAGES = {619--661},
      ISSN = {1865-2433,1865-5394},
   MRCLASS = {19L41 (14F42 55N15 55T15)},
  MRNUMBER = {2811718},
MRREVIEWER = {Michael\ J.\ Paluch},
       DOI = {10.1017/is011004009jkt154},
       URL = {https://doi.org/10.1017/is011004009jkt154},
       note = {\href{https://arxiv.org/abs/1002.2368}{arXiv:\allowbreak 1002.2368}},
    eprint = {1002.2368},
    archivePrefix = {arXiv},
}

@article{Qui19,
    AUTHOR = {Quigley, J. D.},
     TITLE = {The motivic {M}ahowald invariant},
   JOURNAL = {Algebr. Geom. Topol.},
  FJOURNAL = {Algebraic \& Geometric Topology},
    VOLUME = {19},
      YEAR = {2019},
    NUMBER = {5},
     PAGES = {2485--2534},
      ISSN = {1472-2747,1472-2739},
   MRCLASS = {55P42 (14F42 55T15)},
  MRNUMBER = {4023321},
MRREVIEWER = {Rui\ Miguel\ Saramago},
       DOI = {10.2140/agt.2019.19.2485},
       URL = {https://doi.org/10.2140/agt.2019.19.2485},
       note = {\href{https://arxiv.org/abs/1801.06035}{arXiv:\allowbreak 1801.06035}},
    eprint = {1801.06035},
    archivePrefix = {arXiv},
}

@article{Ghe18,
    AUTHOR = {Gheorghe, Bogdan},
     TITLE = {The motivic cofiber of {$\tau$}},
   JOURNAL = {Doc. Math.},
  FJOURNAL = {Documenta Mathematica},
    VOLUME = {23},
      YEAR = {2018},
     PAGES = {1077--1127},
      ISSN = {1431-0635,1431-0643},
   MRCLASS = {14F42 (55P43 55S10 55S35)},
  MRNUMBER = {3874951},
MRREVIEWER = {Drew\ Heard},
       DOI = {10.4171/DM/642},
       URL = {https://doi.org/10.4171/DM/642},
       note = {\href{https://arxiv.org/abs/1701.04877}{arXiv:\allowbreak 1701.04877}},
    eprint = {1701.04877},
    archivePrefix = {arXiv},
}

@article{OO13,
    AUTHOR = {Ormsby, Kyle M. and {\O}stvær, Paul Arne},
     TITLE = {Motivic {B}rown-{P}eterson invariants of the rationals},
   JOURNAL = {Geom. Topol.},
  FJOURNAL = {Geometry \& Topology},
    VOLUME = {17},
      YEAR = {2013},
    NUMBER = {3},
     PAGES = {1671--1706},
      ISSN = {1465-3060,1364-0380},
   MRCLASS = {14F42 (19D50 19E15 55P43)},
  MRNUMBER = {3073932},
MRREVIEWER = {Andrei\ D.\ Halanay},
       DOI = {10.2140/gt.2013.17.1671},
       URL = {https://doi.org/10.2140/gt.2013.17.1671},
       note = {\href{https://arxiv.org/abs/1208.5007}{arXiv:\allowbreak 1208.5007}},
    eprint = {1208.5007},
    archivePrefix = {arXiv},
}

@article{Bac17,
    AUTHOR = {Bachmann, Tom},
     TITLE = {The generalized slices of {H}ermitian {$K$}-theory},
   JOURNAL = {J. Topol.},
  FJOURNAL = {Journal of Topology},
    VOLUME = {10},
      YEAR = {2017},
    NUMBER = {4},
     PAGES = {1124--1144},
      ISSN = {1753-8416,1753-8424},
   MRCLASS = {14F42 (14C25 19E08 55P42)},
  MRNUMBER = {3743071},
MRREVIEWER = {Guillermo\ Corti\~{n}as},
       DOI = {10.1112/topo.12032},
       URL = {https://doi.org/10.1112/topo.12032},
       note = {\href{https://arxiv.org/abs/1610.01346}{arXiv:\allowbreak 1610.01346}},
    eprint = {1610.01346},
    archivePrefix = {arXiv},
}

@article{ARO20,
    AUTHOR = {Ananyevskiy, Alexey and R\"{o}ndigs, Oliver and {\O}stv{\ae}r, Paul Arne},
     TITLE = {On very effective hermitian {$K$}-theory},
   JOURNAL = {Math. Z.},
  FJOURNAL = {Mathematische Zeitschrift},
    VOLUME = {294},
      YEAR = {2020},
    NUMBER = {3--4},
     PAGES = {1021--1034},
      ISSN = {0025-5874,1432-1823},
   MRCLASS = {14F42 (19G38)},
  MRNUMBER = {4074031},
MRREVIEWER = {Heng\ Xie},
       DOI = {10.1007/s00209-019-02302-z},
       URL = {https://doi.org/10.1007/s00209-019-02302-z},
       note = {\href{https://arxiv.org/abs/1712.01349}{arXiv:\allowbreak 1712.01349}},
    eprint = {1712.01349},
    archivePrefix = {arXiv},
}

@article{Ati66,
    AUTHOR = {Atiyah, M. F.},
     TITLE = {{$K$}-theory and reality},
   JOURNAL = {Quart. J. Math. Oxford Ser. (2)},
  FJOURNAL = {The Quarterly Journal of Mathematics. Oxford. Second Series},
    VOLUME = {17},
      YEAR = {1966},
     PAGES = {367--386},
      ISSN = {0033-5606,1464-3847},
   MRCLASS = {55.30 (57.30)},
  MRNUMBER = {206940},
MRREVIEWER = {J.\ F.\ Adams},
       DOI = {10.1093/qmath/17.1.367},
       URL = {https://doi.org/10.1093/qmath/17.1.367},
}

@article{Lew95,
    AUTHOR = {Lewis, Jr., L. Gaunce},
     TITLE = {Change of universe functors in equivariant stable homotopy
              theory},
   JOURNAL = {Fund. Math.},
  FJOURNAL = {Fundamenta Mathematicae},
    VOLUME = {148},
      YEAR = {1995},
    NUMBER = {2},
     PAGES = {117--158},
      ISSN = {0016-2736,1730-6329},
   MRCLASS = {55P91 (55P20 55Q91)},
  MRNUMBER = {1360142},
MRREVIEWER = {J.\ M.\ Boardman},
       DOI = {10.4064/fm-148-2-117-158},
       URL = {https://doi.org/10.4064/fm-148-2-117-158},
}

@incollection{Gre18,
    AUTHOR = {Greenlees, J. P. C.},
     TITLE = {Four approaches to cohomology theories with reality},
 BOOKTITLE = {An alpine bouquet of algebraic topology},
    SERIES = {Contemp. Math.},
    VOLUME = {708},
     PAGES = {139--156},
 PUBLISHER = {Amer. Math. Soc., [Providence], RI},
      YEAR = {2018},
      ISBN = {978-1-4704-2911-9},
   MRCLASS = {55P91 (18G60 19L47 55N91 55Q91)},
  MRNUMBER = {3807754},
MRREVIEWER = {Jose\ Cantarero},
       DOI = {10.1090/conm/708/14261},
       URL = {https://doi.org/10.1090/conm/708/14261},
       note = {\href{https://arxiv.org/abs/1705.09365}{arXiv:\allowbreak 1705.09365}},
    eprint = {1705.09365},
    archivePrefix = {arXiv},
}

@article{Dug05,
    AUTHOR = {Dugger, Daniel},
     TITLE = {An {A}tiyah-{H}irzebruch spectral sequence for {$KR$}-theory},
   JOURNAL = {$K$-Theory},
  FJOURNAL = {$K$-Theory. An Interdisciplinary Journal for the Development,
              Application, and Influence of $K$-Theory in the Mathematical
              Sciences},
    VOLUME = {35},
      YEAR = {2005},
    NUMBER = {3--4},
     PAGES = {213--256},
      ISSN = {0920-3036,1573-0514},
   MRCLASS = {19L64 (55S45 55T25)},
  MRNUMBER = {2240234},
MRREVIEWER = {G\'{e}rald\ Gaudens},
       DOI = {10.1007/s10977-005-1552-9},
       URL = {https://doi.org/10.1007/s10977-005-1552-9},
       note = {Published in 2006. \href{https://arxiv.org/abs/math/0304099}{arXiv:\allowbreak math/0304099}},
    eprint = {math/0304099},
    archivePrefix = {arXiv},
}

@incollection{Voe02,
    AUTHOR = {Voevodsky, Vladimir},
     TITLE = {Open problems in the motivic stable homotopy theory. {I}},
 BOOKTITLE = {Motives, polylogarithms and {H}odge theory, {P}art {I}
              ({I}rvine, {CA}, 1998)},
    SERIES = {Int. Press Lect. Ser.},
    VOLUME = {3},
     PAGES = {3--34},
 PUBLISHER = {Int. Press, Somerville, MA},
      YEAR = {2002},
      ISBN = {1-57146-090-X},
   MRCLASS = {14F35 (55P42)},
  MRNUMBER = {1977582},
MRREVIEWER = {Vladimir\ I.\ Guletski\u{\i}},
note = {\url{https://www.math.ias.edu/Voevodsky/files/files-annotated/Dropbox/Published_papers/Motives/Open/Stage3/newopen.pdf}}
}

@inproceedings{Voe98,
    AUTHOR = {Voevodsky, Vladimir},
     TITLE = {{$\mathbf A^1$}-homotopy theory},
 BOOKTITLE = {Proceedings of the {I}nternational {C}ongress of
              {M}athematicians, {V}ol. {I} ({B}erlin, 1998)},
   JOURNAL = {Doc. Math.},
  FJOURNAL = {Documenta Mathematica},
      YEAR = {1998},
    SERIES = {Extra Vol. I},
     PAGES = {579--604},
      ISSN = {1431-0635,1431-0643},
   MRCLASS = {14F35 (14A15 55U35)},
  MRNUMBER = {1648048},
MRREVIEWER = {Mark\ Hovey},
note = {\url{https://www.math.ias.edu/vladimir/sites/math.ias.edu.vladimir/files/A1_homotopy_ICM_1998_Berlin_published.pdf}}
}

@article{DW25,
    AUTHOR = {Davis, Donald M. and Wilson, W. Stephen},
     TITLE = {The connective {$KO$}-theory of the {E}ilenberg-{M}ac{L}ane
              space {$K(\mathbb{Z}_2,2)$}, {I}: the {$E_2$} page},
   JOURNAL = {J. Homotopy Relat. Struct.},
  FJOURNAL = {Journal of Homotopy and Related Structures},
    VOLUME = {20},
      YEAR = {2025},
    NUMBER = {4},
     PAGES = {511--521},
      ISSN = {2193-8407,1512-2891},
   MRCLASS = {55S10 (55N15 55N20 55T15)},
  MRNUMBER = {4984184},
MRREVIEWER = {Geoffrey\ M. L. Powell},
       DOI = {10.1007/s40062-025-00379-4},
       URL = {https://doi.org/10.1007/s40062-025-00379-4},
       note = {\href{https://arxiv.org/abs/2410.23094}{arXiv:\allowbreak 2410.23094}},
    eprint = {2410.23094},
    archivePrefix = {arXiv},
}

@book{Dav26,
    AUTHOR = {Davis, Donald M.},
     TITLE = {The connective {$KO$} theory of the {E}ilenberg-{M}ac{L}ane
              space {$K(\mathbb Z/2,2)$}},
    SERIES = {New York Journal of Mathematics. NYJM Monographs},
    VOLUME = {7},
 PUBLISHER = {State University of New York, University at Albany, Albany,
              NY},
      YEAR = {2026},
     PAGES = {88},
   MRCLASS = {55T15 (55N15 55N20 57R20)},
  MRNUMBER = {5062716},
       URL = {https://nyjm.albany.edu/m/2026/7.htm},
       note = {\href{https://arxiv.org/abs/2502.14982}{arXiv:\allowbreak 2502.14982}},
    eprint = {2502.14982},
    archivePrefix = {arXiv},
}

@article{BCQ25,
    AUTHOR = {Balderrama, William and Culver, Dominic Leon and Quigley, J.
              D.},
     TITLE = {The motivic lambda algebra and motivic {H}opf invariant one
              problem},
   JOURNAL = {Geom. Topol.},
  FJOURNAL = {Geometry \& Topology},
    VOLUME = {29},
      YEAR = {2025},
    NUMBER = {3},
     PAGES = {1489--1570},
      ISSN = {1465-3060,1364-0380},
   MRCLASS = {14F42 (55Q25 55Q45 55S10 55T15)},
  MRNUMBER = {4918111},
MRREVIEWER = {Yangyang\ Ruan},
       DOI = {10.2140/gt.2025.29.1489},
       URL = {https://doi.org/10.2140/gt.2025.29.1489},
       note = {\href{https://arxiv.org/abs/2502.14982}{arXiv:\allowbreak 2112.07479}}
}

@misc{Bak26,
	author = {Andrew Baker},
	title = {Poincaré duality spaces related to the {J}oker},
	year = {2026},
	note = {\href{https://arxiv.org/abs/2603.29425}{arXiv:\allowbreak 2603.29425}},
    eprint = {2603.29425},
    archivePrefix = {arXiv},
}

@misc{Sal26,
	author = {Andrew Salch},
	title = {Ext-groups of graded comodules},
	year = {2026},
	note = {\href{https://arxiv.org/abs/2609.26465}{arXiv:\allowbreak 2609.26465}},
    eprint = {2609.26465},
    archivePrefix = {arXiv},
}

@misc{San26,
	author = {Sean Sanford},
	title = {On the {M}ext groups of {sVect$_\R$ and sVect$_{\mathbb H}$}},
	year = {2026},
	note = {\href{https://arxiv.org/abs/2609.20153}{arXiv:\allowbreak 2609.20153}},
    eprint = {2609.20153},
    archivePrefix = {arXiv},
}

@misc{Fre12,
    title = {Lectures on twisted {$K$}-theory and orientifolds},
    author = {Daniel S. Freed},
    year = {2012},
    note = {\url{https://people.math.harvard.edu/~dafr/vienna.pdf}}
}

@article{HNP25,
    AUTHOR = {Harpaz, Yonatan and Nuiten, Joost and Prasma, Matan},
     TITLE = {On {$k$}-invariants for {$(\infty,n)$}-categories},
   JOURNAL = {Algebr. Geom. Topol.},
  FJOURNAL = {Algebraic \& Geometric Topology},
    VOLUME = {25},
      YEAR = {2025},
    NUMBER = {2},
     PAGES = {721--790},
      ISSN = {1472-2747,1472-2739},
   MRCLASS = {18N65 (55S45)},
  MRNUMBER = {4930552},
MRREVIEWER = {Philippe\ Gaucher},
       DOI = {10.2140/agt.2025.25.721},
       URL = {https://doi.org/10.2140/agt.2025.25.721},
       note = {\href{https://arxiv.org/abs/2011.12723}{arXiv:\allowbreak 2011.12723}},
    eprint = {2011.12723},
    archivePrefix = {arXiv},
}

@article{Mah82,
    AUTHOR = {Mahowald, Mark},
     TITLE = {The image of {$J$} in the {$EHP$} sequence},
   JOURNAL = {Ann. of Math. (2)},
  FJOURNAL = {Annals of Mathematics. Second Series},
    VOLUME = {116},
      YEAR = {1982},
    NUMBER = {1},
     PAGES = {65--112},
      ISSN = {0003-486X},
   MRCLASS = {55Q40},
  MRNUMBER = {662118},
       DOI = {10.2307/2007048},
        note = {Correction in \textit{Annals of Mathematics}, 120:399--400, 1984.}
}

@book{Koc96,
    AUTHOR = {Kochman, S. O.},
     TITLE = {Bordism, stable homotopy and {A}dams spectral sequences},
    SERIES = {Fields Institute Monographs},
    VOLUME = {7},
 PUBLISHER = {American Mathematical Society, Providence, RI},
      YEAR = {1996},
     PAGES = {xiv+272},
      ISBN = {0-8218-0600-9},
   MRCLASS = {55P42 (55N20 55N22 55Q45 55R40 55T15)},
  MRNUMBER = {1407034},
MRREVIEWER = {J. P. C. Greenlees},
       DOI = {10.1090/fim/007},
       URL = {https://doi.org/10.1090/fim/007},
}

@article{BMSS15,
    AUTHOR = {Robert Bruner and Khairia Mira and Laura Stanley and Victor Snaith},
     TITLE = {Ossa's Theorem via the {K}unneth Formula},
   JOURNAL = {Mathematics and Statistics},
  FJOURNAL = {Mathematics and Statistics},
    VOLUME = {3},
      YEAR = {2015},
    NUMBER = {3},
     PAGES = {58--64},
     note = {\href{https://arxiv.org/abs/1008.0166}{arXiv:\allowbreak 1008.0166}},
    eprint = {1008.0166},
    archivePrefix = {arXiv},
}

@misc{KPY26a,
	author = {Ryohei Kobayashi and Abhinav Prem and Matthew Yu},
	title = {Realizing Braided Fusion 2-Categories as {(3+1)D} Mixed-State Topological Orders},
	year = {2026},
	note = {\href{https://arxiv.org/abs/2610.03864}{arXiv:\allowbreak 2610.03864}},
    eprint = {2610.03864},
    archivePrefix = {arXiv},
}

@misc{Dha26,
	author = {Prerna Dhankhar},
	title = {On the Motivic {$K$}-Theories of Classifying Spaces of Cyclic 2-Groups},
	year = {2026},
	note = {\href{https://arxiv.org/abs/2610.27016}{arXiv:\allowbreak 2610.27016}},
    eprint = {2610.27016},
    archivePrefix = {arXiv},
}

@article{BHS26,
    AUTHOR = {Burklund, Robert and Hahn, Jeremy and Senger, Andrew},
     TITLE = {Galois reconstruction of {A}rtin-{T}ate {$\R$}-motivic
              spectra},
   JOURNAL = {Geom. Topol.},
  FJOURNAL = {Geometry \& Topology},
    VOLUME = {30},
      YEAR = {2026},
    NUMBER = {5},
     PAGES = {1625--1717},
      ISSN = {1465-3060,1364-0380},
   MRCLASS = {14F42},
  MRNUMBER = {5106758},
       DOI = {10.2140/gt.2026.30.1625},
       URL = {https://doi.org/10.2140/gt.2026.30.1625},
       note = {\href{https://arxiv.org/abs/2010.10325}{arXiv:\allowbreak 2010.10325}},
}

@misc{AP25,
	author = {Keita Allen and Lucas Piessevaux},
	title = {Synthetic equivariant spectra for finite abelian groups and motivic homotopy theory},
	year = {2025},
	note = {\href{https://arxiv.org/abs/2510.20197}{arXiv:\allowbreak 2510.20197}},
    eprint = {2510.20197},
    archivePrefix = {arXiv},
}

@misc{Emm26,
	author = {Konstantin Emming},
	title = {Filtrations in {$\C$}-motivic stable homotopy theory},
	year = {2026},
	note = {\href{https://arxiv.org/abs/2608.04877}{arXiv:\allowbreak 2608.04877}},
    eprint = {2608.04877},
    archivePrefix = {arXiv},
}

@misc{HRW22,
	author = {Jeremy Hahn and Arpon Raksit and Dylan Wilson},
	title = {A motivic filtration on the topological cyclic homology of commutative ring spectra},
	year = {2022},
	note = {\href{https://arxiv.org/abs/2206.11208}{arXiv:\allowbreak 2206.11208}},
    eprint = {2206.11208},
    archivePrefix = {arXiv},
}

@article{GIKR22,
    AUTHOR = {Gheorghe, Bogdan and Isaksen, Daniel C. and Krause, Achim and
              Ricka, Nicolas},
     TITLE = {{$\C$}-motivic modular forms},
   JOURNAL = {J. Eur. Math. Soc. (JEMS)},
  FJOURNAL = {Journal of the European Mathematical Society (JEMS)},
    VOLUME = {24},
      YEAR = {2022},
    NUMBER = {10},
     PAGES = {3597--3628},
      ISSN = {1435-9855,1435-9863},
   MRCLASS = {14F42 (55N34 55Q45 55S10 55T15)},
  MRNUMBER = {4432907},
MRREVIEWER = {Zhouli\ Xu},
       DOI = {10.4171/jems/1171},
       URL = {https://doi.org/10.4171/jems/1171},
       note = {\href{https://arxiv.org/abs/1810.11050}{arXiv:\allowbreak 1810.11050}},
    eprint = {1810.11050},
    archivePrefix = {arXiv},
}

@article{BMS19,
    AUTHOR = {Bhatt, Bhargav and Morrow, Matthew and Scholze, Peter},
     TITLE = {Topological {H}ochschild homology and integral {$p$}-adic
              {H}odge theory},
   JOURNAL = {Publ. Math. Inst. Hautes \'{E}tudes Sci.},
  FJOURNAL = {Publications Math\'{e}matiques. Institut de Hautes \'{E}tudes
              Scientifiques},
    VOLUME = {129},
      YEAR = {2019},
     PAGES = {199--310},
      ISSN = {0073-8301,1618-1913},
   MRCLASS = {14F30 (13A35)},
  MRNUMBER = {3949030},
MRREVIEWER = {Lance\ Edward\ Miller},
       DOI = {10.1007/s10240-019-00106-9},
       URL = {https://doi.org/10.1007/s10240-019-00106-9},
        note = {\href{https://arxiv.org/abs/1802.03261}{arXiv:\allowbreak 1802.03261}},
    eprint = {1802.03261},
    archivePrefix = {arXiv},
}

@misc{CDN24,
	author = {Christian Carrick and Jack Morgan Davies and van Nigtevecht, Sven},
	title = {The descent spectral sequence for topological modular forms},
	year = {2024},
	note = {\href{https://arxiv.org/abs/2412.01640}{arXiv:\allowbreak 2412.01640}},
    eprint = {2412.01640},
    archivePrefix = {arXiv},
}

@misc{Ari21,
	author = {Stefano Ariotta},
	title = {Coherent cochain complexes and {B}eilinson t-structures},
	year = {2021},
	note = {With an appendix by Achim Krause. \href{https://arxiv.org/abs/2109.01017}{arXiv:\allowbreak 2109.01017}},
    eprint = {2109.01017},
    archivePrefix = {arXiv},
}

@mastersthesis{Cre25,
  author	= {Julie Creemers},
  title = {Spectral sequences obtained from towers in the {$\infty$}-category of spectra: décalage, exact couples and {M}assey products},
  school	= {Utrecht University},
  year   = {2025},
  note = {\url{https://studenttheses.uu.nl/items/6a22ab2e-906a-4e19-b079-caf0dc3fbd8b}}
}

@misc{AP25a,
	author = {Toni Annala and Piotr Pstrągowski},
	title = {A note on weight filtrations at the characteristic},
	year = {2025},
	note = {\href{https://arxiv.org/abs/2502.19626}{arXiv:\allowbreak 2502.19626}},
    eprint = {2502.19626},
    archivePrefix = {arXiv},
}

@misc{BJM24,
	author = {J. Francis Baer and Maxwell Johnson and Peter Marek},
	title = {Stable comodule deformations and the synthetic {A}dams-{N}ovikov spectral sequence},
	year = {2024},
	note = {\href{https://arxiv.org/abs/2402.14274}{arXiv:\allowbreak 2402.14274}},
    eprint = {2402.14274},
    archivePrefix = {arXiv},
}

@article{CDvN25,
    AUTHOR = {Carrick, Christian and Davies, Jack Morgan and van Nigtevecht,
              Sven},
     TITLE = {Descent spectral sequences through synthetic spectra},
   JOURNAL = {Int. Math. Res. Not. IMRN},
  FJOURNAL = {International Mathematics Research Notices. IMRN},
      YEAR = {2025},
    VOLUME = {15},
     PAGES = {Paper No. rnaf216, 26},
      ISSN = {1073-7928,1687-0247},
   MRCLASS = {14D23 (19L99 55N34)},
  MRNUMBER = {4937581},
       DOI = {10.1093/imrn/rnaf216},
       URL = {https://doi.org/10.1093/imrn/rnaf216},
       note = {\href{https://arxiv.org/abs/2407.01507}{arXiv:\allowbreak 2407.01507}},
    eprint = {2407.01507},
    archivePrefix = {arXiv},
}

@phdthesis{Pie26,
    author = "Lucas Piessevaux",
    title = "Synthetic Equivariant Spectra and Complex Bordism",
    school = "University of Bonn",
    year = "2026",
    note = "\url{https://bonndoc.ulb.uni-bonn.de/xmlui/handle/20.500.11811/14463}"
}

@misc{Pst23a,
	author = {Piotr Pstrągowski},
	title = {Perfect even modules and the even filtration},
	year = {2023},
	note = {\href{https://arxiv.org/abs/2304.04685}{arXiv:\allowbreak 2304.04685}},
    eprint = {2304.04685},
    archivePrefix = {arXiv},
}

@misc{Ant25,
	author = {Benjamin Antieau},
	title = {Filtrations and cohomology {III}: cohomology of {$E_\infty$} rings},
	year = {2025},
	note = {\href{https://arxiv.org/abs/2512.15509}{arXiv:\allowbreak 2512.15509}},
    eprint = {2512.15509},
    archivePrefix = {arXiv},
}

@misc{AR24,
	author = {Benjamin Antieau and Noah Riggenbach},
	title = {Cyclotomic synthetic spectra},
	year = {2024},
	note = {\href{https://arxiv.org/abs/2411.19929}{arXiv:\allowbreak 2411.19929}},
    eprint = {2411.19929},
    archivePrefix = {arXiv},
}

@article{HM25,
    AUTHOR = {Hedenlund, Alice and Moulinos, Tasos},
     TITLE = {The synthetic {H}ilbert additive group scheme},
   JOURNAL = {Selecta Math. (N.S.)},
  FJOURNAL = {Selecta Mathematica. New Series},
    VOLUME = {31},
      YEAR = {2025},
    NUMBER = {5},
     PAGES = {Paper No. 89, 64},
      ISSN = {1022-1824,1420-9020},
   MRCLASS = {14F42 (13F20)},
  MRNUMBER = {4963699},
       DOI = {10.1007/s00029-025-01085-y},
       URL = {https://doi.org/10.1007/s00029-025-01085-y},
       note = {\href{https://arxiv.org/abs/2411.17441}{arXiv:\allowbreak 2411.17441}},
    eprint = {2411.17441},
    archivePrefix = {arXiv},
}

@misc{Ant25a,
	author = {Benjamin Antieau},
	title = {Filtrations and cohomology {II}: the {G}auss-{M}anin connection},
	year = {2025},
	note = {\href{https://arxiv.org/abs/2511.14731}{arXiv:\allowbreak 2511.14731}},
    eprint = {2511.14731},
    archivePrefix = {arXiv},
}

@misc{QZ25,
	author = {Ryan Quinn and Qi Zhu},
	title = {Multiplicative Equivariant {T}hom Spectra \& Structured {R}eal Orientations},
	year = {2025},
	note = {\href{https://arxiv.org/abs/2512.15573}{arXiv:\allowbreak 2512.15573}},
    eprint = {2512.155731},
    archivePrefix = {arXiv},
}

@incollection{AS06,
    AUTHOR = {Atiyah, Michael and Segal, Graeme},
     TITLE = {Twisted {$K$}-theory and cohomology},
 BOOKTITLE = {Inspired by {S}. {S}. {C}hern},
    SERIES = {Nankai Tracts Math.},
    VOLUME = {11},
     PAGES = {5--43},
 PUBLISHER = {World Sci. Publ., Hackensack, NJ},
      YEAR = {2006},
      ISBN = {978-981-270-062-9; 981-270-062-5},
   MRCLASS = {55N15 (19L99 55R40 55S25 55S30 55T25)},
  MRNUMBER = {2307274},
MRREVIEWER = {Francis\ Clarke},
       DOI = {10.1142/9789812772688\{_}0002}

@article{Dec24a,
    AUTHOR = {D{\'{e}}coppet, Thibault D.},
     TITLE = {The 2-{D}eligne tensor product},
   JOURNAL = {Kyoto J. Math.},
  FJOURNAL = {Kyoto Journal of Mathematics},
    VOLUME = {64},
      YEAR = {2024},
    NUMBER = {1},
     PAGES = {1--29},
      ISSN = {2156-2261,2154-3321},
   MRCLASS = {18M20 (18M05 18N10 18N25)},
  MRNUMBER = {4677193},
MRREVIEWER = {R.\ H.\ Street},
       DOI = {10.1215/21562261-2023-0005},
       URL = {https://doi.org/10.1215/21562261-2023-0005},
note = {\href{https://arxiv.org/abs/abs2103.16880}{arXiv:\allowbreak abs/2103.16880}},
    eprint = {2103.16880},
    archivePrefix = {arXiv},
}

@misc{KKSW26,
	author = {Masashi Kawahira and Shotaro Kawanago and Shota Saito and Hiroki Wada},
	title = {Anomalies and discrete torsion in {T}ype {II} and {T}ype {I} worldsheet theories},
	year = {2026},
	note = {\href{https://arxiv.org/abs/2608.30774}{arXiv:\allowbreak 2608.30774}},
    eprint = {2608.30774},
    archivePrefix = {arXiv},
}

@article{BC76,
 author = {Edgar H. Brown and Michael Comenetz},
 journal = {American Journal of Mathematics},
 number = {1},
 pages = {1--27},
 publisher = {Johns Hopkins University Press},
 title = {Pontrjagin Duality for Generalized Homology and Cohomology Theories},
 volume = {98},
 year = {1976}
}

@article{FMS07,
        author="Freed, Daniel S.
        and Moore, Gregory W.
        and Segal, Graeme",
        title="The Uncertainty of Fluxes",
        journal="Communications in Mathematical Physics",
        year="2007",
        month="Apr",
        day="01",
        volume="271",
        number="1",
        pages="247--274",
        issn="1432-0916",
        doi="10.1007/s00220-006-0181-3",
        url="https://doi.org/10.1007/s00220-006-0181-3",
        note = "\url{https://arxiv.org/abs/hep-th/0605198}"
}

@misc{DS26,
	author = {Amartya Shekhar Dubey and Luuk Stehouwer},
	title = {The {$J$}-Homomorphism From Localising The Cobordism Hypothesis},
	year = {2026},
	note = {\href{https://arxiv.org/abs/2611.19615}{arXiv:\allowbreak 2611.19615}},
    eprint = {2611.19615},
    archivePrefix = {arXiv},
}

\end{document}